%% file: dgg-v3.tex
\documentclass[reqno]{amsart}
\usepackage[all]{xy}
\usepackage{amsmath}
\usepackage{amssymb,latexsym}
\usepackage{eucal}
\usepackage{graphics}
\newtheorem{theorem}{Theorem}[section]
\newtheorem{proposition}[theorem]{Proposition}
\newtheorem{lemma}[theorem]{Lemma}

\newtheorem{definition}[theorem]{Definition}
\newtheorem{remark}[theorem]{Remark}

\def\re{{\it resp.}}
\def\bee{\begin{equation}}
\def\ga{\gamma}
\def\ka{\kappa}
\def\Ga{\Gamma}
\def\UGa{\underline{\Gamma}}
\def\uga{\underline{\gamma}}
\def\De{\Delta}
\def\UDe{\underline{\Delta}}
\def\ude{\underline{\delta}}
\def\de{\delta}
\def\uv{\underline{v}}
\def\uu{\underline{u}}
\def\om{\omega}
\def\Om{\Omega}
\def\UOm{\underline{\Omega}}

\def\cc{{\mathcal{C}}}
\def\ddc{\mathcal{D}}
\def\kc{\mathcal{K}}

\def\pc{{\mathcal{P}}}
\def\qc{{\mathcal{Q}}}
\def\rc{{\mathcal{R}}}
\def\gc{{\mathcal{G}}}
\def\pa{{P}^{\mathrm{ad}}}
\def\pac{{\mathcal{P}^{\mathrm{ad}}}}
\def\epc{\mathcal{E}q(\pc, \pc')}
\def\ep2c{\mathcal{E}q(\pc', \pc'')}

\def\oxc{{\mathcal{O}_X}}

\def\h2c{{{\Psi^2_{X/S}}}}

\def\J2{{J_{02}^{(2)}}}
\def\la{\longrightarrow}

\def\to{\mathrm{Tors}}

\def\ot{\otimes}

\def\sp{\mathrm{Spec}}

\def\pra{{p_{0}^{\ast}}}
\def\prb{{p_{1}^{\ast}}}
\def\prc{{p_{2}^{\ast}}}
\def\prd{{p_{3}^{\ast}}}

\def\prab{{p_{01}^{\ast}}}
\def\prbc{{p_{12}^{\ast}}}

\def\prac{{p_{02}^{\ast}}}

\def\dea{{\Delta^{1}_{X/S}}}
\def\deb{{\Delta^{2}_{X/S}}}
\def\dec{{\Delta^{3}_{X/S}}}
\def\ded{{\Delta^{4}_{X/S}}}
\def\li{\mathrm{lien}}

\def\AL{\mathrm{Ar\; Lie}}
\def\AI{\mathrm{Ar_I\; Lie}}
\def\OL{\mathrm{Ob\; Lie}}
\def\ap{\mathrm{AP}}
\def\op{\mathrm{OP}}
\def\ae{\mathrm{AE}}
\def\oe{\mathrm{OE}}

\renewcommand{\include}{\input}
\begin{document}

\title{Differential geometry of gerbes}
\author{Lawrence Breen \and William Messing}
\date{}

\address{LB: UMR CNRS 7539 \newline
\indent Institut Galil\'ee \newline
\indent Universit\'e Paris 13 \newline
\indent  F-93430  Villetaneuse, France}
  \email{breen@math.univ-paris13.fr}
\address{ WM: School of Mathematics \newline
  \indent University of Minnesota \newline
  \indent 127 Vincent Hall \newline
  \indent 206 Church Street S.E. \newline
  \indent Minneapolis, MN 55455, USA}
  \email{messing@math.umn.edu}

\maketitle

\tableofcontents

\setcounter{section}{-1}

\include{dgg-ch0}

\include{dgg-ch1}

\include{dgg-ch2}

\include{dgg-ch3}

\include{dgg-ch4}

\include{dgg-ch5}

\include{dgg-ch6}

\include{dgg-ch7}

\include{dgg-appendixa}
\include{dgg-appendixb}

  \end{document}

%% file: dgg-ch0.tex
\section*{Introduction}
\label{intro}

%\addtocounter{subsec}{1}%
\setcounter{equation}{0}%
%{\bf \noindent \thesection.\thesubsec}
A foundational concept
in differential geometry is that of a 
connection on
a $G$-principal bundle, which  
embodies the notion of parallel transport between infinitesimally 
close fibers of the bundle. The  curvature    of the connection
 measures the compatibilities of these parallel transports around 
 infinitesimal closed loops in the base space $X$.
It 
is a global 2-form on $X$ with values in the Lie algebra of the gauge 
group 
of $P$, and 
satisfies the Bianchi identity. When the structure group 
$G$  is abelian, the adjoint action of $P$ is trivial, and $\om$
 is simply a closed $\mathrm{Lie}\;G$-valued form on $X$.

\bigskip

Just as principal bundles (also known as torsors) are  geometric 
realizations 
of $G$-valued \v{C}ech 1-cocycles, a gerbe on $X$ is the geometric 
embodiment of a 2-cocycle on $X$ with values in a non-abelian group. 
We will show here how to define, in a very global manner, the 
corresponding  concepts of connection and curvature 
for a gerbe $\pc$ on $X$. This generalizes the construction 
given  by 
J.-L. Brylinski \cite{bry}
in the case of gerbes with coefficients 
in the  abelian group $\mathbb{C}^\ast$. The point of view presented 
here, though compatible with that of \cite{bry} is somewhat more global.
A  connection on $\pc$ once more describes the parallel transport between 
infinitesimally close fibers.
The fibers of $\pc$ are now categories,   and 
the 
additional degree of freedom provided by natural transformations 
between functors allows one to choose, for a fixed
connection on $\pc$, an additional curving. 
While the connection is now a functor  above a 
 line segment between a pair of infinitesimally close points, 
we view  the curving  as a higher mode of parallel 
 transport around the boundary of an infinitesimal 2-simplex
 between induced first level  parallel transports. 
The 3-curvature of the curving  data provided by
 the connection and the
 curving is defined as the
 obstruction to coherence between the induced curvings
 above the faces of a 3-simplex.
 
 \bigskip

 The analogue of the Bianchi identity which the 3-curvature satisfies 
 essentially 
states that the pullback of these obstructions are themselves coherent
above any infinitesimal 
4-simplex. A complication arises, however, when one passes from 
 principal bundles to gerbes. While a connection on a 
 principal bundle is defined by an 
 isomorphism between neighboring fibers, a connection on a gerbe 
 merely  provides one, in categorical style,  with an equivalence 
 between neighboring fiber categories. In particular, the 
 inverse functor of a connection is no longer 
 uniquely defined. While the complications which 
ensue can  be dispensed with in the case of 
 abelian gerbes, this is no longer the case in the more general 
 situation considered here, and it is necessary to introduce
an additional 1-arrow which we call the fake curvature.
 There results an extra term in the equation \eqref{relnufi} satisfied by 
 the 3-curvature, so that this equation differs from the   naive 
 generalization from 2- to 3-forms of the Bianchi identity. 

\bigskip

In this increasingly complicated situation, a very efficient tool is 
provided by A. Kock's combinatorial approach to differential calculus
\cite{AK:Synth}-\cite{AK:Bianchi}. While 
this was carried out by him  within the context of synthetic geometry, we 
have shown in \cite{cdf}
that it remains valid, up to 
some minor differences, in an  algebro-geometric setting. We show 
here that the notions of connections, 
curving data, and associated curvature forms, and the
formulas describing the relations between them, become quite transparent 
when formulated in this combinatorial language.  Once local 
trivializations of the gerbe have been chosen, these can  be 
translated into  more familiar terms as equations satisfied by
traditional Lie algebra-valued 
differential forms. As could have been expected from the cocyclic 
description of a gerbe in \cite{fest}, \cite{lb:2-gerbes}, one can no 
longer consider in isolation the local 
 3-curvature forms attached to the connective structure of the 
 gerbe $\pc$.
These are now accompanied by auxiliary lower degree 
forms. The local description of 
the 
connection and curving, and the associated 3-curvature is 
therefore more complicated than in the principal bundle case, and we refer
to 
theorem \ref{th:omi}  for a precise statement.

\bigskip

The idea of  geometric quantization has provided a strong link between 
differential forms and principal bundles. 
This can be extended from principal bundles to gerbes.
As explained in \cite{bry}, to which we refer for further 
information on this subject,
P. A. M. Dirac's classic study of the magnetic field 
of a monopole may be interpreted as the construction of 
a gerbe on $S^{3}$ with 
connective structure. Recently, the concept of gerbes and $n$-gerbes 
with connections and holonomy
has  surfaced in a  number of areas in the mathematics and
physics literature  (see for example 
\cite{freed}, 
\cite{murray}, 
\cite{kalkkinen}, \cite{mackaay}, 
\cite{attal}), mostly  in the context of gerbes  
with abelian 
coefficients. The advent of string theory and  more recently the 
progress in understanding its solitonic excitations such as 
D(irichlet) and Neveu-Schwarz branes,
has brought to the fore certain differential forms of  degree 
2 and higher, for which it  
is of interest to provide a geometrical interpretation. We refer in 
particular  to \S 6 of \cite{fw} for a very down to earth description
of the abelian theory 
in terms of \v{C}ech 2- and 3-cocycles. In that 
language, our  local 2-forms $B_{i}$ \eqref{defbigki},
 which describe  the curving 
data, are known as Neveu-Schwarz $B$-fields. Here they take their
values in (the Lie algebras of) non-abelian groups,
and the 
 geometrically defined
 global 3-curvature form $\Om$ \eqref{defomeg} is the associated 
 field strength $H$.  It has been suggested that such  
  forms with values in the non-abelian group $U(k)$ 
 arise when  $k$  Neveu-Schwarz branes coincide (see \cite{dijk}).

\bigskip

We now describe in more detail the content of the present text. While 
we have placed ourselves firmly within the context of algebraic 
geometry, in which the concepts of gerbes and stacks have been to date most 
fully developed, and for which we have worked out the more general
combinatorial differential calculus in \cite{cdf}, we wish to keep in mind 
the applications mentioned above. We have therefore 
adopted a somewhat hybrid style of exposition.  The results  concerning 
torsors and gerbes  are valid  when working with any 
site. Nevertheless we use the language of algebraic geometry and refer 
to sheaves for an unspecified topology. In order to make the theory 
easily comparable with the $C^{\infty}$-theory, we speak of open sets 
and their intersections, when we of course mean fiber products.
Starting with section 
 \ref{sec:curv-Bian}, we work exclusively in the \'etale 
topology. More precisely, when we consider sheaves of groups, torsors, 
bitorsors, gerbes and $gr$-stacks we will work on the big \'etale site, while 
in our Lie-theoretic considerations of differential forms we will 
restrict ourselves to the small \'etale site. In relating the two we 
will frequently use the big site on an \'etale open set. We will always state 
our 
results for the relative situation $X/S$, even though the 
consideration of families of  topological spaces with a parameter 
space $S$ is not as standard in topology as in algebraic geometry.

\bigskip

In the 
first   section, we present a very explicit, cocyclic, discussion of 
the theory of principal bundles and connections. 
There are many equivalent definitions of a connection, and we choose 
here
the one which is most standard in algebraic geometry. 
Once $P$ has 
been locally trivialized, the connection is described by a connection 
1-form, yielding by covariant differentiation a local value for the 
curvature. This description of the curvature  in terms of the 
connection form of the curvature is classical, in the $C^{\infty}$ context 
of differential geometry.
While the corresponding assertion is well known in 
algebraic geometry when the structure group $G$ is the 
multiplicative group (as in \cite{mm}, \cite{m1}), it is harder to come by
when the group $G$ is  arbitrary. We can refer for this to \cite{onishchik}
\S  5, where non-abelian cocycles are introduced,   but otherwise 
the preferred method of exposition
paraphrases the classical $C^{\infty}$ theory (see  for example 
\cite{gm} \S 5).
We believe that our very short proof 
of this assertion based on \cite{cdf}, though it is some sense close to the 
discussion carried out in \cite{berth} in a more additive context,
is a good illustration
of the power of the techniques presented here. For future reference, we
also discuss in this section the notion of a connection on group. In the 
 situation of a relative scheme $X/S$, the standard notion of a 
 connection on a $G$-principal bundle requires that $G$ be the 
 pullback to $X$ of an $S$-group scheme. We  introduce here the 
 corresponding notion in the more general case of an $X$-group 
 scheme $G$ endowed with a connection, and examine its basic 
 properties. In order to do so, it is necessary to extend from $S$-group
schemes to $X$-group schemes some of the results of \cite{cdf}. This is
 carried out in Appendix A.

\bigskip

 We then give  a 
 short description of the 2-category of gerbes on a space. 
 In particular, a very explicit description of the cocycle 
and coboundary relations for a fully trivialized gerbe is   given.
We have found it very convenient to display this
diagrammatically, as in 
\cite{breen-labesse} in the more limited context of Galois cohomology,
 rather than  simply as a set of 
impenetrable equations. In   section \ref{sec:morita0}
 we review 
Giraud's analogue \cite{Gir} IV proposition 5.2.5 in the 
principal bundle context 
of the Morita description of natural 
transformations between categories of modules (as in \cite{bass} 
chapter II), and the corresponding bitorsor cocycle description of 
a gerbe with a chosen family of local objects $x_i$ as a family of bitorsors
$P_{ij}$, together with coherent isomorphisms $\Psi_{ijk}$ \eqref{1bcoc}.
We  
introduce  the stack $\pac = \mathcal{E}q(\pc,\, \pc)$ 
of self-equivalences of a gerbe $\pc$. This group-like stack 
deserves to be called 
the  gauge stack of the gerbe $\pc$ since it is the analogue of the 
gauge group $P^\mathrm{ad}:= \mathrm{Aut}_{G}(P)$ of base space 
preserving 
self-transformations of a $G$-principal 
bundle $P$. We give in proposition \ref{pcad} a description of 
$\pac$
which closely parallels the well-known construction \eqref{0coch1} 
of the gauge group
$\pa$ of $P$  as a twisted form of the group $G$, for the adjoint action of 
$P$ on $G$. 
Our result is easiest to state
 in the special case of $G$-gerbes, in other words gerbes attached 
to a globally defined coefficient group $G$. In that case $\pac$ 
is a twisted form of the group-like stack $\mathrm{Bitors}(G)$
of bitorsors on $G$, with
the
inner twisting   provided, just as in the torsor case,
by the cocycle data arising from  the gerbe $\pc$. 
 Once more, these cocycles are systematically displayed
in a ``2-dimensional 
algebra'' diagrammatic setting in which their significance and 
various verifications are best 
understood.

\bigskip

We now come to our  main object of study. We  define in section 
 \ref{sec:curv-Bian} the 
concepts of connection and curving data on a gerbe $\pc$. Our 
definitions of these notions 
are somewhat more intrinsic than those of Brylinski, since they do 
not depend on the choice of local trivializing objects. They are also
more general, even for gerbes with
abelian coefficient groups. Since Brylinski  only considered  the traditional 
 cohomology with values in an abelian group, 
he in effect worked in the 2-category of
abelian gerbes on a space $X$, as defined  in 
\cite{lb:2-gerbes} definition 2.9, rather than in the full subcategory 
of the category of gerbes consisting of all those gerbes whose 
coefficient groups are abelian.
Once the curving data associated to the curving data 
 has been defined, one shows that 
it defines a 3-curvature 3-form $\Om$ \eqref{defomeg} on $X$,
with values in the arrows of the  gauge stack 
$\pac$ introduced above. The form  $\Om$  
satisfies a higher Bianchi-type identity, which is expressed by the 
commutativity of the cubic diagram  of 2-arrows \eqref{bianchicube}. 
 As we have already stated, the result
\eqref{omcocy1}, \eqref{bianchi3} is somewhat 
more complicated than could have been expected, since the term on the 
right-hand side of \eqref{bianchi3} is in general non-trivial. The 
proof of the identity is once more geometric, and follows from the 
analysis of the four-dimensional hypercube \eqref{hypercube}. We also 
describe here in geometric  terms the notion of a morphism between a 
pair of connections on a gerbe  \eqref{def:xa},
 \eqref{def:a}, \eqref{prism}
and the corresponding notion of a natural transformation \eqref{def:r}
between two such morphisms. 
 Generalizing this, there is a natural notion 
of a morphism between gerbes endowed with connection (really connection 
triples) and of morphisms between such morphisms.
Finally, we give a de Rham type description of these concepts of
connection,  curving and induced 3-curvature, as well as of the 
action  of a transformation triple. In this description,  the 
complication regarding the non-triviality  of the right-hand term
in  \eqref{bianchi3} becomes quite natural.
The full  de Rham description displayed in diagram
\eqref{dr4}  is 
of necessity somewhat elaborate, since we are in a context in which
the connection on $\pc$  is not integrable. It rests upon 
the notion of a differential form with values in a $gr$-stack, which
we introduce in appendix \ref{sec:conger}.  This extends the concept of
a group-valued differential form, as in \cite{cdf}.

\bigskip

The next two sections provide us with  successively more explicit descriptions
of the concepts introduced in section \ref{sec:curv-Bian}. We have seen  that
 the choice of local objects of the gerbe $\pc$ provides us with a description
of $\pc$ in terms of bitorsors. We now show that the curving pair 
$(\epsilon,\, K)$  and the
 induced 3-curvature $\Omega$ may be described in a similar manner. 
The result is stated in proposition \ref{prop:omi}, and reduces in the
abelian case, when in addition the fake curvature $\kappa$ is trivial
and the group $G$  is the multiplicative group $G_m$,
to Hitchin's description  of a gerbe with connection 
\cite{hitch}. In section \ref{sec:fulldec}, we consider  additional
 trivializations  involving the choice of certain arrows in the
gerbe. In that case, a fully cocyclic description of the concepts 
introduced in section \ref{sec:curv-Bian} is obtained. The result for
cocycles is  stated in theorem \ref{th:omi}. The corresponding
coboundary relations, and the relations between these, are respectively
described in paragraphs  \ref {disc:cobound} and \ref{disc:cob-cob}.

\bigskip

In the final section, we look at two  special cases
of the theory. We first assume that  the given gerbe $\pc$ is trivial. This 
 hypothesis is not without interest, as it is the full extension
 from torsors to gerbes of the well-known 
 assertion that a connection on a trivial $G$-torsor is defined by a 
 global $G$-valued 1-form on $X$, and  its curvature  by 
 the induced Maurer-Cartan 2-form \eqref{glmc}.
 The cocyclic data and the 3-curvature are now given by global $G$- and 
 $\mathrm{Aut}\, G$-valued forms, satisfying   four 
 quite reasonable equations (displayed in \S \ref{special2}). Some of the 
  coboundary  relations which we find in this trivial gerbe case  also occur,
 in a physics context,
 in the recent 
 preprint \cite{IC} of I. Chepelev.
As a second  example of our general theory,
we  assume that the $G$-gerbe is abelian, with a 
 connection which respects this abelian gerbe structure, and 
 suppose 
 that the fake curvature 
 1-arrow $\kappa$ is trivial. In that case the $G$-valued cocycles
associated to the trivial gerbe 
  are seen  to be the expected ones, with values in the ordinary $G$-valued
\v{C}ech-de Rham complex of $X$.
 
\bigskip

As will be apparent from this summary, we   give  four distinct
 descriptions  of the theory of gerbes with a  connective structure, and their
 associated  curvature forms. The first, which is global and geometric,
does not  require any auxilliary choices.  While it may be considered somewhat
abstract, it is the one which displays in the clearest form the various
phenomena which are encountered.  The second description
 is of a homological nature, and 
introduces a \v{C}ech-de Rham theory of $gr$-stacks with connections. This
 is  somewhat complicated since we do not restrict ourselves to the
 case of  abelian gerbes, but it encodes very efficiently the notions 
of connection, curving, fake curvature and 3-curvature which we are dealing 
with. The third description  is 
semi-local and enables us to extend Hitchin's approach in \cite{hitch} to
 our non-abelian context. It is a compromise between the global theory and a  
purely local one, since it expresses the theory of gerbes in terms of the 
simpler concept of bitorsors. However, certain aspects of 
 the differential geometry are  however somewhat obscured 
 within this framework. The fourth description is purely local.
 Here, all geometric objects which occur have been locally trivialized, so 
that we can display in full detail the cocycle and coboundary conditions for
 gerbes with a  connective structure. These
 are a priori described in the language of combinatorial forms, but we have 
then translated each of these equations into the traditional language of
 differential calculus.  It is our contention that each of these four
 descriptions has its advantages, and sheds some light upon the other three.

\bigskip

We would like to thank Victor Ginzburg for bringing the preprint 
\cite{hitch} to our attention. We are grateful to Christiaan 
Hofman for pointing out a missing term in one of the formulas in the 
first version of this text.
The first author also thanks   Patrick Polo
for a useful comment regarding Lie theory, and  Bernard Julia and  Jussi 
Kalkkinen for their continued interest and observations regarding
this project.

%% file: dgg-ch1.tex
\addtocounter{section}{1}
 \section{Connections on groups and torsors}
\label{torsors}

\subsection{}
\setcounter{equation}{0}%
\label{sub:1}
Let $X$ be a  scheme, and $G$ a sheaf of groups on $X$ (we will also 
say that $G$ is an $X$-group). We will
in general  suppose that the topology on $X$ is the \'etale topology and 
speak of open sets of $X$. The main features of our theory remain 
valid for a Grothendieck topology,
as in \cite{sga3} expos\'e IV, \cite{milne} chapter II. This is in 
particular the case for all of \S \ref{sub:1}.

\bigskip

Let $P$ be a $G$-torsor (or 
locally trivial
principal $G$- bundle)  on $X$.
For any open cover 
$\mathcal{U}=(U_{i})_{i \in I}$ of $X$ 
 and any family of
local sections $s_{i} \in \Gamma(U_{i}, P)$, the torsor $P$
is described by the 
$G$-valued 1-cochain $g_{ij}$ defined by 
\[  s_{j} = s_{i}g_{ij}
   \]
    and which satisfies the cocycle condition
    \begin{equation}
        \label{1coc}
        g_{ij}\,g_{jk} \:=\: g_{ik}\:.
         \end{equation}
        For another choice of local sections $s'_{i}$  of $P$
         on the same open cover 
         $\mathcal{U}$ of $X$, and associated 
         1-cocycles $g'_{ij}$, 
          we may set 
         \[ s_{i}= s'_{i}\ga_{i}\] 
         for a $G$-valued 0-cochain $\ga_{i}$. The 
         coboundary relation
         \[ g'_{ij} = \ga_{i}\, g_{ij} \, \ga_{j}^{-1}\]
         is then satisfied.

         \bigskip
         
Consider 
the sheaf $\mathrm{Aut}\,P$ of $G$-equivariant self-morphisms of $P$ which
leave 
the base $X$ fixed. This is a sheaf of groups $\pa$ on $X$, 
 for  the group law given by composition of morphisms,
whose global sections  constitute the gauge
group\footnote{In the physics literature, this is generally called the 
group of gauge transformations and $G$ itself is referred to as the 
gauge group of $P$.} of 
$P$. For any  section   $u \in \pa$ we may compare the 
        sections $u(s_{i})$ and $s_{i}$ of $P$. Setting
    \bee
    \label{1u}
    u(s_{i}) = s_{i}\,\gamma_{i}\:,
    \end{equation}
   one finds that
   \bee
    \label{0coch}
        g_{ij} = \ga_{i}\,g_{ij}\,\ga_{j}^{-1}\:.
        \end{equation}
        Conversely, a  $G$-valued 0-cochain $\ga_{i} \in \Ga(U_{i}, 
        \, G)$
        satisfying (\ref{0coch})  determines  by \eqref{1u} 
         a global section $u$ of $\pa$.
         
        \bigskip
        
We may  rewrite (\ref{0coch}) as 
   \bee
   \label{0coch1}  
   {}^{g_{ij\,}}\!\ga_{j} = \ga_{i}
   \end{equation} 
   where $ {}^{g_{ij\,}}\!\ga_{j}:= g_{ij}\,\ga_{j}\,g_{ij}^{-1}$. 
The restrictions of the sheaves of groups $\pa$ and  $G$ above the 
open sets $U_{i}$ are isomorphic 
{\it via} the maps
\[
\begin{array}{ccc}
    G_{\mid U_{i}} & \simeq  & \pa_{\mid U_{i}} \\
    \ga_{i} & \mapsto & (g \mapsto \ga_{i} g)\:,
    \end{array}\]
  and $\pa$ can be recovered by gluing 
  together local copies of the sheaf $G$ according to the rule   
  (\ref{0coch1}), in other words {\it via} the adjoint action of the 
  group $G$ upon itself. This assertion may also be  expressed
  as the isomorphism of sheaves of groups
   \bee
   \label{defpa}
   \xymatrix@R=3pt{ P \wedge^{G}G
    \ar[r]^(.55){\sim}&\pa\\(p,g)
    \ar@{|->}[r]& (p \mapsto pg) 
   }\end{equation}
   of \cite{Gir} III proposition 
  2.3.7, 
  and identifies $\pa$ with the  gauge group 
  $\mathrm{Ad}\, P$ of $P$ defined in \cite{AB} \S 2.
  We  can also write \eqref{defpa} as an isomorphism 
  \bee
  \label{padj}
  \pa  \simeq {}^{P}\!G
\end{equation}
in order to emphasize that the basic object being  considered here is $G$ 
rather than $P$.
  
  \bigskip
  
 More generally, consider  a pair of $G$-torsors $P$ and 
  $P'$ on $X$, and chosen families of local sections $s_{i},\, s'_{i}$ 
  of $P$ and  $P'$ with respect to the same open cover $\mathcal{U}$ 
  of $X$, with  
corresponding $G$-valued  1-cocycles $g_{ij},\,g'_{ij}$.
A global section $u$ of 
  the sheaf $\mathrm{Isom}_{G}(P,\, P')$, in other words 
  a morphism of $G$-torsors $u: P \la P'$ on $X$,  is  described  by 
  the $G$-valued  0-cochains
  $g_{i}$ defined by
  \bee
  \label{gi}
  u(s_{i}) = s'_{i}\,g_{i}  
\end{equation} 
and which satisfy the 
  equations
  \bee
  \label{0coch2}
        g'_{ij} = g_{i}\,g_{ij}\,g_{j}^{-1}\:.
  \end{equation}
  Rewritten in the style of (\ref{0coch1}), this becomes 
  \bee 
  \label{0coch3}
  g'_{ij}\,g_{j}\,g_{ij}^{-1} = g_{i}\:.
  \end{equation}
  When $P$ and $P'$ are  a pair of trivial torsors, with global 
  sections $s$ and $s'$, the sheaf of sets
  $\mathrm{Isom}_{G}(P,\, P')$ is isomorphic
  to the underlying sheaf of  pointed sets of $G$, under  the rule
  \[\begin{array}{ccc}
      G & \simeq & \mathrm{Isom}_{G}(P, P')\\
      g & \mapsto & (s \mapsto s'g)\:.
      \end{array}\]
    Since the  $G$-torsors $P$ and $P'$ are both locally isomorphic
    to the underlying sheaf 
   of sets of $G$, equation (\ref{0coch3}) expresses the fact that 
   $\mathrm{Isom}_{G}(P,\, P')$ may be constructed  by gluing
   local copies $G_{\mid U_{i}}$ of the sheaf of set $G$
   above $U_{ij}$ according to the rule
  \bee
  \label{gi1}
   \begin{array}{ccc}
       G_{\mid U_{ij}} & \simeq &  G_{\mid U_{ij}}\\
       g_{j} & \mapsto & g'_{ij}\,g_{j}\, g_{ij}^{-1}\:.
       \end{array}\end{equation}
By composition, the sheaf  $\mathrm{Isom}_{G}(P,\, P')$
is both a left   
${P'}^{\mathrm{ad}}$ and a right  $\pa$-torsor  on $X$. These left and right 
actions commute, so that  $\mathrm{Isom}_{G}(P,\, P')$ is in fact a 
       $({P'}^{\mathrm{ad}}, \pa)$-bitorsor on $X$ 
       (\cite{Gir}  III\S 1.5,\cite{fest}).
       The right 
       $\pa$-torsor structure is described locally by the 
       rule 
       \[(g_{i} )\: (\ga_{i}) := (g_{i}\,\ga_{i})\]
    and is easily verified to be compatible with the gluing data. 
     So is the rule
       \[ (\ga'_{i})\:(g_{i}) := (\ga'_{i}\,g_{i}) \]
     for the  corresponding left $P{'}^{\mathrm{ad}}$-torsor structure.
       
       \begin{remark}
           \label{rem:twisom}
            {\rm
         For a given 
       isomorphism  $\mu: G \la G'$ of sheaves of groups on $X$, 
       the previous 
       discussion extends to a description of  the bitorsor 
       $\mathrm{Isom}_{\mu}(P,\, P')$ of $\mu$-equivariant 
       isomorphisms from 
       $P$ to $P'$. A section $u:P \la P'$ of this sheaf 
       corresponds to an isomorphism of  $G'$-torsors
       \[\begin{array}{ccc}
       P \wedge^{G}G' & \simeq & P'\\
       (p, g') & \mapsto & u(p)g'
       \end{array}
       \]
       By  (\ref{0coch2}), such an $u$ is described  by a 
       family of  0-cochains $g'_{i} \in \Gamma (U_{i}, G')$ such 
       that
     \[         g'_{ij} = g'_{i}\,\mu(g_{ij})\,(g'_{j})^{-1}\:.\]
  To this corresponds, in the  (\ref{0coch3}) style, the 
        equivalent equation
 \bee
 \label{eq:twisom}
 g'_{ij}\,g'_{j}\,\mu(g_{ij})^{-1} = g'_{i}\end{equation}
 which tells us that the sheaf $\mathrm{Isom}_{\mu}(P, P')$ may 
 be constructed from local copies of the $G'_{\mid U_{i}}$,
 glued  above the open sets $U_{ij}$ according to the rule
 \[g' \mapsto g'_{ij}\, g' \,\mu(g_{ij})^{-1}\:.\] The right 
 $\pa$-torsor  ({\it resp.}  the left ${P'}^{\mathrm{ad}}$-torsor) structure 
 on $\mathrm{Isom}_{\mu}(P, P')$ is 
now described in local terms by 
   \[(g'_{i}) \: (\ga_{i}) := (g'_{i}\,\mu(\ga_{i}))\]
    ({\it resp.} 
    \[(\ga'_{i}) \: (g'_{i}) := (\ga'_{i}\,g'_{i})\quad )\:.\]
   }\end{remark}

\bigskip

 \subsection{}
\setcounter{equation}{0}%
\label{subpar61}
We work in the \'etale topology. 
Let $S$ be a scheme, $X$ be an $S$-scheme, and $G$ be a
sheaf of groups on $X$.
As in \cite{cdf} (1.10.1), we denote by $\Delta^{n}_{X/S}$
the $S$-scheme which parametrizes  
$(n+1)$-tuples of first order infinitesimally close points of $X$. 
For $n >1$, the notion of infinitesimal proximity which we have in mind 
here is somewhat refined, but  the naive notion is adequate whenever we may
ignore 2-torsion phenomena, for example whenever 2 is invertible on 
the base scheme $S$. In that case the infinitesimal neighborhood 
$\Delta^{n}_{X/S}$  of  the diagonal embedding of $X$ 
in its $n+1$-fold product 
may be replaced by the somewhat coarser, but more intuitive,
neighborhood $\Delta^{(n)}_{X/S}$ of {\it loc. cit.} (1.4.8).
For $n=1$, both notions  coincide, and $\Delta^{1}_{X/S}$ is simply 
the scheme of first order principal 
parts on $X$, endowed with the projections $p_{0}$ and 
$p_{1}$ onto the first and second factor  $X$. 
When $\De^{n}_{X/S}$ is regarded as an 
$X$-scheme, it will always be {\it via} the projection $p_{0}$. 
Let 
\bee
\label{defdel}
\Delta: X \la \Delta^{1}_{X/S}
\end{equation}
be the diagonal immersion.
If $p_{\alpha}: \Delta^n_{X/S} \la \Delta^m_{X/S}$
is a projection, we will denote 
by $E_{\alpha}$ the inverse image
 $p_{\alpha}^{\ast}E$ of an object $E$ on $\Delta^m_{X/S}$, 
and similarly for morphisms. 

\bigskip

For $T$ an $S$-scheme, we refer to $(x_0, \ldots,x_n) \in
\Delta^n_{X/S}(T)$ as a $T$-valued infinitesimal $n$-simplex, and
say that it is degenerate if $x_{i} = x_j$ for some $i \neq
j$. We denote the subscheme of degenerate infinitesimal
simplices  by $s\Delta^n_{X/S}$. In \cite{cdf} 1.12 this was
denoted by $\partial\Delta^n_{X/S}$.

\bigskip

\begin{definition}
    \label{def:gcon}
    A 
    connection on the $X$-group $G$ is a group isomorphism
    \bee
    \label{grcon}
    \mu: p^{\ast}_{1}G \la p^{\ast}_{0}G\end{equation}
    satisfying $\Delta^{\ast}(\mu) = 1_{G}$. 
     \end{definition}
     
     \bigskip
     
     \noindent{\bf Examples} \hspace{3pt}{\it i}) If $G$ is the base 
     change of 
     an $S$-group, then $G$ has a canonical connection.
     
     \medskip
     
     \qquad \qquad  {\it ii}) If $G$ is a reductive $X$-group 
     scheme and $X$ is affine, then there are connections on $G$, 
     since  $\mathrm{Isom}(\prb G,\, \pra G)$ is a smooth
     $\dea$-scheme (\cite{sga3} XXIV corollary 1.8).

\bigskip
     
   Let $\mathrm{Aut}(G)$ denote the sheaf of group automorphisms of an $X$-group
   $G$.
A connection $\mu$ on $G$ induces by transport of structure a connection 
$i_{\mu}$ on $\mathrm{Aut}(G)$, which we will also denote by 
$\mu^{\mathrm{ad}}$. The image $\mu^{\mathrm{ad}}(u)$ of an element $u 
\in \prb\mathrm{Aut}(G)$ by this 
connection is therefore defined by the commutativity of the following 
square:
\begin{equation}
    \label{diag:muad}
\xymatrix{
\prb G \ar[r]^{u} \ar[d]_{\mu}  & \prb G \ar[d]^{\mu} \\
\pra G \ar[r]_{\mu^{\mathrm{ad}}(u)}& \pra G
}\end{equation}

Varying the open set $U $, we obtain the sheaf
of connections  $\mathrm{Conn}(G)$. Composition on the left with pointed
automorphisms of $G$ equips  $\mathrm{Conn}(G)$ with the structure of a left 
torsor
under the $X$-group
$\mathrm{Lie}(\mathrm{Aut}(G),\, \Om^1_{X/S})$.

\begin{definition}
Let $\mu$ be a connection on the $X$-group $G$. 
  The curvature $\kappa_{\mu}$ of  $\mu$ is the section
of $\pra \mathrm{Aut}(G)$ on $\deb$
 defined     
 by 
\begin{equation}
    \label{defkmu}
    \kappa_{\mu}:= 
(p_{01}^{\ast}\mu)\,(p_{12}^{\ast}\mu)\,(p_{02}^{\ast}\mu)^{-1}\,.
\end{equation}
\end{definition}
\noindent By definition, we have a commuting  diagram 
    \bee
     \label{def:kmu}
\xymatrix{
p^{\ast}_{2}G \ar[r]^{p_{12}^{\ast}\mu}
\ar[d]_{p_{02}^{\ast}\mu} & p^{\ast}_{1}G 
\ar[d]^{p_{01}^{\ast}\mu}\\
 p^{\ast}_{0}G \ar[r]_{\kappa_{\mu}}
&  p^{\ast}_{0}G \:\:.
}\end{equation}
 The restriction of $\kappa_{\mu}$ to
 $s\Delta^{2}_{X/S} $ is the identity so that
$\kappa_{\mu}$ is a section of $\mathrm{Lie} (\mathrm{Aut}(G),\, \Om^2_{X/S})$.
  
  \bigskip

\subsection{}
\setcounter{equation}{0}%
Let $(G,\, \mu)$ be an $X$-group  with connection and $P$ a $G$-torsor on $X$.
\begin{definition} 
    \label{def:con-0}
  {\it i}) A connection on $P$ is 
a $\mu$-equivariant isomorphism on $\dea$:
\bee
\label{def:con}
\epsilon : p_{1}^{\ast}P \la 
p_{0}^{\ast}P\:,\end{equation}
 satisfying $\Delta^{\ast}\epsilon = \mathrm{id}_{P}$.

\medskip

\hspace{2cm} {\it ii}) The curvature  of $\epsilon$ is the
$\kappa_{\mu}$-equivariant automophism of  $\pra (P)$ over $\deb$ 
\bee
 \kappa_{\epsilon}:= (\prab\epsilon)\, (\prbc\epsilon) \, (\prac \epsilon)^{-1}
\end{equation} 
\end{definition}

As the restriction of $\kappa_{\epsilon}$ to $s\deb$ is $\mathrm{id}_{\pra (P)}$,
we can view  $\kappa_{\epsilon}$  as a section of
$\mathrm{Lie}(\mathrm{Aut}_{\kappa_{\mu}}(P),\,\Om^2_{X/S})$ (here we abuse
notation
as $\mathrm{Aut}_{\kappa_{\mu}}(P)$ is defined over $\deb$ rather than
over $X$ and we have supressed
$(p_{0})_{\ast}$).
When $G$ comes from an 
$S$-group and $\mu$ is its canonical connection, this definition of a connection 
on the  $G$-torsor $P$ coincides with Ehresmann's notion of  a connection on a
principal bundle {\it via} parallel transport, and
in this case $\kappa_{{\epsilon}}$ 
is a $\mathrm{Lie}(\pa)$-valued 2-form on $X$.
The following square, whose commutativity defines $\kappa_{\epsilon}$,
is the analogue of \eqref{def:kmu}
\begin{equation}
        \label{def:kappa1}
        \xymatrix@=11pt{
   \prc P \ar[rr]^{\epsilon_{12}} \ar[dd]_{\epsilon_{02}}
    && \prb P \ar[dd]^{\epsilon_{01}}\\&&\\\pra 
    P\ar[rr]_{\kappa_{\epsilon}}&& \pra 
    P
    }\end{equation} 
    The connection $\epsilon$ induces a connection
 \begin{equation}
     \label{mudef}
 \epsilon^{\mathrm{ad}}: \prb \pa \la \pra \pa
 \end{equation}
on the group $\pa$,  determined by the commutativity of the squares
\bee
\label{def:ieps}
\xymatrix{
p^{\ast}_{1}P \ar[r]^{u} \ar[d]_{\epsilon} & p^{\ast}_{1} P
\ar[d]^{\epsilon}\\ p^{\ast}_{0}P \ar[r]_{\epsilon^{\mathrm{ad}}(u)} & 
p^{\ast}_{0}P}
\end{equation}
 for all sections $u$ of $\prb(\pa)$. 
 
 \bigskip
 
\subsection{}
\setcounter{equation}{0}%
\label{sec:contors}
 In this section, we do not limit  ourselves to the \'etale topology. 
 We refer to  \cite{Gir}, \cite{fest}  for a general discussion 
     of bitorsors and simply recall here  that for a pair of $X$-groups $G$ 
     and $H$, a $(G,\, H)$-bitorsor  on $X$ is a sheaf $P$ on $X$,  
     together with a left action of $G$ on $P$ and a right action of 
     $H$ on $P$, which commute with each other and respectively 
     define on $P$ a 
      left $G$- and a right $H$-torsor structure. When $G$ and $H$ 
      are equal, we will simply say that  $P$ is  a $G$-bitorsor.
      
      \bigskip

      We choose in this paragraph to view a 
      $(G,\, H)$-bitorsor as a left $G$-bitorsor with additional 
      structure. A right multiplication by $H$ on the  left $G$-torsor
      $P$ on $X$ is 
      determined, for each section $p \in P$, by the morphism
      $u_{p} \in \mathrm{Isom}(H,\, G)$ defined by
      \begin{equation}
          \label{defu}
          ph = u_{p}(h) p\end{equation}
      for all $h \in H$.  For any other section $p'$ of $P$, we may set
      \begin{equation}
          \label{defmor}
          p' = \ga p
      \end{equation}
      for some $\ga \in G$. It is readily verified that 
      \[u_{p'} = i_{\ga}\, u_{p}\,.\]
      It follows that the giving of a $(G,\, H)$-bitorsor structure 
      on a left $G$-torsor $P$ is equivalent to that of 
      a morphism of sheaves on  $X$ 
      \[ u: P \la \mathrm{Isom}(H,\, G)\]
      which is equivariant with respect to the inner conjugation map 
      \[ i: G \la \mathrm{Aut}(G)\:.\]
      \noindent In particular, if a bitorsor $P$ is  trivial as a 
      left  $G$-torsor, with 
      chosen global section $p$, then the bitorsor structure of $P$ is 
      entirely described by the corresponding element $u_{p} \in 
      \mathrm{Isom}(H,\, G)$.
      If $K$ is another $X$-group, and  $Q$ 
      an $(H,\, K)$-bitorsor
      on $X$,  whose right multiplication by $K$ is described by the 
      corresponding morphism 
      \[ v: Q \la \mathrm{Isom}(K, \, H)\,.\]
      The contracted product $P \wedge^{H} Q$ of $P$
      and $Q$  is a $(G,\, K)$-bitorsor. It may be 
      considered  
      as   a left $G$-torsor on $X$  for which the 
      right multiplication by $K$ is described by the morphism 
      induced by the composite map
      \[ \xymatrix{P \times Q  \ar[r] & \mathrm{Isom}(H,\, G) \times  
      \mathrm{Isom}(K,\, H) \ar[r] & \mathrm{Isom}(K,\, G)
      }\]
      where the second map is the composition of isomorphisms.
      If $(P,\, u)$ and $(P',\, u')$ are a pair of $(G,\, H)$-bitorsors,
      and 
       $f: P \la P'$ is a morphism of left $G$-torsors, then $f$  is 
       in fact a bitorsor morphism if and only if 
       \begin{equation}
           \label{bitormor}
           u = u' \, f\,.\end{equation}
      \noindent    Let 
      \[ f: (P,p) \la (P', \, p') \]
     be a morphism between a pair of $(G, \, H)$-bitorsors with chosen global 
     sections $p$ and $p'$. As  a morphism of left-torsors, $f$ is 
      is determined by the section $g \in \Gamma (X, \, G)$
     defined by the equation
     \begin{equation} 
         \label{def:f}
         f(p) = g p'\,.
         \end{equation}
       The equation \eqref{bitormor}, which ensures that  $f$ 
       is compatible with the right action of $H$, is now equivalent to
       the condition
     \begin{equation}
         \label{defuprime}
         u_p = i_{g}\, u'_{p'}\,.
         \end{equation}

    \bigskip
    
     Much of the previous discussion regarding torsors can  be 
     extended to bitorsors. For example, to a $(G,H)$-bitorsor $P$ on $X$ 
     is associated the $X$-group $\pa$ of bitorsor isomorphisms from 
     $P$ to itself. Similarly, given a  pair of $X$-groups with connnection
     $(G,\, \mu)$ and $(H,\,\nu)$, a connection on a $(G,H)$-bitorsor
     $P$ on $X$ is a $(\mu,\, \nu)$-equivariant morphism $\epsilon$
     \eqref{def:con}
     satisfying the condition
     $\Delta^{\ast}\epsilon = 1_{P}$. In particular, the underlying 
     left $G$- and right $H$-torsors of $P$ are then both endowed with a 
     connection.
 
\bigskip

\subsection{}
\setcounter{equation}{0}%
 We now give the local expressions for a connection and its curvature
on a right 
$G$-torsor $P$. We will use here some of the definitions and results of
 \cite{cdf} extended from $S$-groups to $X$-groups with connections, 
 as justified in
 Appendix A (to which we also refer for additional notation).

\bigskip

Let us  choose local sections $s_{i}$ of $P$ on an open cover 
$\mathcal{U}= (U_{i})_{i \in I}$ of $X$, which determine  a family of 
1-cocycles $g_{ij}$ (\ref{1coc}).
The connection $\epsilon$
is described in 
local terms by the family of sections 
$\om_{i} \in \Gamma(\Delta^{1}_{U_{i}/S},\, G)$, defined by
\bee
\label{1con}
\epsilon(p^{\ast}_{1}s_{i}) = 
(p^{\ast}_{0}s_{i})\,\om_{i}\:.\end{equation} 
Reasoning as in \eqref{eq:twisom}, we find that
      \begin{align}
 \om_{j} 
&= (p_{0}^{\ast}g_{ij})^{-1}\,\, \om_{i}\, 
\,\mu(p_{1}^{\ast}g_{ij}) \label{omcoc1}\\
&=\om_{i}^{\:{\ast}_{\mu}\,g_{ij}}\,,\label{omcoc-alt}
\end{align}
 so that by \eqref{astmu} we have 
in additive notation
\bee
\label{omcoc11}
\om_{j}= \om_{i}^{g_{ij}} + \delta^{0}_{\mu}(g_{ij})\:.\end{equation}
If we write
\[\kappa_{\epsilon\,\mid U_i}(\pra s_i) = 
(\pra s_i) \kappa_i\,,\]
we find that the local curvatures $\kappa_{i} \in \Gamma(U_{i},\, 
\mathrm{Lie}(\mathrm{Aut}(G),\, \Omega^{2}_{X/S})$
are expressed in  terms of the local connection 1-forms $\om_{i}$
 by the rule
\bee
\label{loccurv}
\kappa_{i}= \delta^{1}_{\mu}\,\om_{i} \,.  \end{equation}
They  satisfy 
\begin{equation}
\label{k-twist}
\kappa_{j}= (p^{\ast}_{0}g_{ij})^{-1} \,\kappa_{i}\,\, 
\kappa_{\mu}(p^{\ast}_{0}g_{ij})
\end{equation}
 or in additive notation, 
\bee
\label{k-twist-1}\kappa_{j} =
 \kappa_{i}^{\,g_{ij}}\,+\,  [g_{ij}, \kappa_{\mu}]\:,
\end{equation}
where the bracket pairing is defined by \eqref{lie3}.
When $G$ comes from $S$, and is endowed with its canonical connection, 
the local curvature equations simplify:
\begin{eqnarray}
\label{kd1}
\kappa_{i} & = &
\delta^1 (\om_{i}) \\
\label{kd11}
&= &d\om_{i} + [\om_{i}]^{(2)}\:.
\end{eqnarray}
The second equality is the structural equation of E. Cartan, as is in 
  \cite{cdf} theorem 3.3, with the expression $ [\om_{i}]^{(2)}$ more 
commonly  written as $1/2[\om_i,\, \om_i]$, or as $\om_i \wedge \om_i\,$.
The local connection and curvature forms $\kappa_{i}$ are now related 
to each other by 
the following simpler versions of \eqref{k-twist} and 
\eqref{k-twist-1}:
 \begin{align}
\label{gluecurv}
\om_{j}&= \om_{i}^{g_{ij}} + \delta^{0}(g_{ij})\\\notag
\kappa_{j} & = (\pra g_{ij})^{-1}\,\kappa_{i}\, (\pra g_{ij})\\
\label{gluecurv2}
&= 
\kappa_{i}^{\,g_{ij}}\:.\end{align}
The last  equation  confirms that the local curvature 2-forms $\kappa_{i}$ glue
to a globally defined 2-form on $X$, with values in  (the Lie algebra 
of) the
gauge bundle of $P$.
For $P$ representable, the previous discussion can
be carried out in a Grothendieck 
topology context, for the single covering open set $P \la X$ of 
$X$. This is also the point of view  taken in the 
differential geometry texts, where 
the curvature $\kappa$ is generally  described 
as a  relative $\mathrm{Lie} \,G$-valued 2-form $\kappa_{P}$ on $P/S$. 
The gluing condition (\ref{gluecurv2}) for the
local curvature forms $\kappa_{i}$
then  correspond  to  descent conditions from $P$ to $X$. These may be 
expressed
 as
 the equivariance and horizontality conditions for the curvature 
form $\kappa_{P}$
(see for example  \cite{gm} 5.2). 
\begin{remark}{\rm We still assume that $G$ comes from $S$.
    
    \medskip
    
    \qquad \qquad \quad {\it i}) 
     One can study with the same combinatorial techniques 
    the effect
of a change  of local sections of the torsor $P$ on 
 the connection and curvature forms, and recover the corresponding 
 classical formulas.  

 \medskip
 
 \qquad \qquad \quad {\it ii}) When the group  $G$ is abelian,  
    the  discussion 
     reduces to that in \cite{mm} \S 3.1.
      The adjoint action of $g_{ij}$ on $\omega_{i}$ is trivial,  
 so that
the first equation in  \eqref{gluecurv}
reduces to 
 \[ \om_{j} = \om_{i} + \delta^{0}(g_{ij}) \:.\]
The gauge group $\pa$ of $P$ is now canonically isomorphic to the group  
$G$ itself, so that 
 the pair $(\om_{i}, g_{ij})$ is now a 1-cocycle on $X$ with values 
 in the  truncated  version 
 $ (G \stackrel{\delta^{0}}{\la}
 \mathrm{Lie} \,(G) \ot \Om^{1}_{X/S})$ of the $G$-valued de Rham 
 complex on $X$. Similarly, the
formula \eqref{kd11} for
 the local  curvature
form $\kappa_i$
reduces to the
       simpler formula 
 \[\kappa_{i} = d \om_{i}\:.\] }

    \end{remark}
    
\subsection{}
\setcounter{equation}{0}%
The  first part of each of the 
following two  propositions is a   global versions of the classical Bianchi
identity. We   prove them by a diagrammatic argument which  is
 a cubic variant of the  tetrahedral proof in  \cite{AK:Bianchi}, 
theorem 9.1. 

\begin{proposition}
\label{bianchigroup}
 {\it i})
Let $G$ be an group scheme defined over $X$, and 
endowed with a connection $\mu$ \eqref{grcon}
and associated curvature $\kappa_{\mu}$ \eqref{defkmu}. Then
  $\kappa_{\mu} \in \mathrm{Lie}(\mathrm{Aut}(G),\, \Omega^2_{X/S})$ satisfies
the equation 
 \begin{equation}
        \label{defkapmu0}
        \delta^{2}_{\mu^{\mathrm{ad}}} \kappa_{\mu} = 0
        \end{equation}
where $\mu^{\mathrm{ad}}$ is 
 the connection on  $\mathrm{Aut}(G)$ induced as in \eqref{diag:muad} by $\mu$. 

\medskip

\hspace{2.3cm}  {\it ii)} Let $\mu'$ be another connection on $G$, so that
 $\mu'= \alpha \, \mu$ for some 1-form  $\alpha \in 
 \mathrm{Lie}(\mathrm{Aut}(G), \, \Omega^1_{X/S})$,
 and let $\kappa_{\mu'}$ be its curvature. Then the
relation 
\[\kappa_{\mu'} = ( \delta^{1}_{\mu^{\mathrm{ad}}} \alpha) \:\kappa_{\mu}\]
is satisfied in $ \mathrm{Lie}(\mathrm{Aut}(G), \, \Omega^2_{X/S})$. 
\end{proposition}

\bigskip

\begin{proposition}
\label{bianchitors}
  {\it i})
  Let $G$ be an $S$-group scheme,   
  $P$ be a  $G$-torsor on an $S$-scheme $X$, and 
    $\prb P \stackrel{\epsilon}{\la} \pra P$ a connection on $P$ with 
   associated  curvature the $\mathrm{Lie}\:{\pa}$-valued 2-form 
   $\kappa_{\epsilon}$. 
    Then $\kappa_{\epsilon}$ satisfies the equation 
    \begin{equation}
        \label{defkap0}
        \delta^{2}_{\epsilon^{\mathrm{ad}}} \kappa_{\epsilon} = 0\,.
        \end{equation}

\medskip

\hspace{2.3cm}  {\it ii}) Let  $\epsilon'$ be another connection on the
 $G$-torsor $P$, so
 that   $\epsilon' =  h \, \epsilon$ for some element 
 $h \in \mathrm{Lie}(\pa, \,
 \Omega^1_{X/S})$, and let $\kappa_{\epsilon'}\in \mathrm{Lie}(\pa, \,
 \Omega^2_{X/S})$  be the  
 curvature of $\epsilon'$. Then the relation 
\begin{equation}
\label{cobcap}
 \kappa_{\epsilon'} = (\delta^1_{\epsilon^{\mathrm{ad}}}\, h) \, 
 \kappa_{\epsilon}\end{equation}
is satisfied.
 \end{proposition}

\bigskip

We  now  prove the first part of   proposition \ref{bianchigroup}.
Consider the following diagram
of groups  above $\dec$:

\begin{equation}
      \label{bianchigrdiag}
    \includegraphics{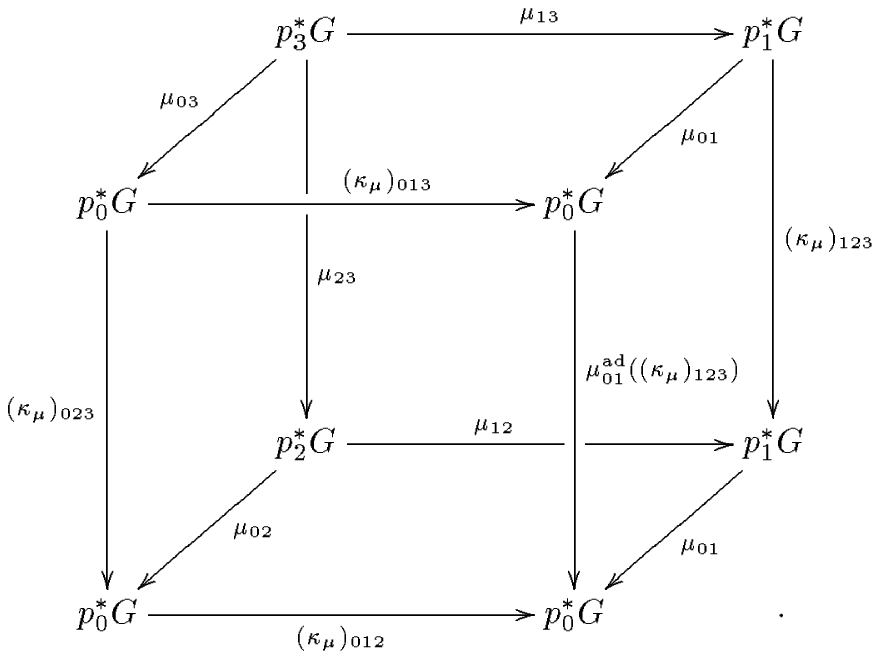} 
\end{equation}  

    Denoting by $S_{ijk}$ pullback of the commutative square
    \eqref{def:kmu}
    by the projection $p_{ijk}: \dec \la \deb$, we observe that 
    the top and
   bottom faces  of \eqref{bianchigrdiag}  are
   respectively the commutative squares $S_{013}$ and $S_{012}$, while the
 left  vertical face is  $S_{023}$, and the 
    back vertical one is $S_{123}$. Each of these four faces commute,
    and  so does
    the right hand vertical one, with $\mu^{\mathrm{ad}}:\prb 
    \mathrm{Aut}(G) \la \pra \mathrm{Aut}(G) $
    defined 
   as in 
    \eqref{diag:muad}.  Since all arrows in 
    this diagram are invertible,  it follows that the  remaining 
    front vertical square in diagram \eqref{bianchigrdiag} is also 
    commutative. 
   Its commutativity expresses 
    combinatorially the Bianchi identity \eqref{defkapmu0}.

\bigskip

The second part of proposition \ref{bianchigroup} can also be proved 
diagrammatically, by considering the diagram

\begin{equation}
\label{compcon1}
\xymatrix{
\prc G \ar[r]^{\mu_{12}}
\ar[d]_{\mu_{02}}
 & \prb G  \ar[r]^{\alpha_{12}} \ar[d]_{\mu_{01}}
 &
\pra G  \ar[d]_{\mu_{01}} \ar[dr]^{\mu'_{01}}
  &\\
\pra G   \ar[r]^{\kappa_{\mu}}  \ar[d]_{\alpha_{02}} &\pra G 
\ar[r]_{\mu^{\mathrm{ad}}_{01}(\alpha_{12})}
 \ar[d]_{\alpha_{02}} & \pra G \ar[r]_{\alpha_{01}} & \pra G\\
\pra G  \ar[r]_{\kappa_{\mu}} &\pra G & & 
}
\end{equation}
 The two upper squares in this diagram commute
 by definition,  and so does 
 the lower one by lemma 2.8 of \cite{cdf}, applied to the 
 $\mathrm{Aut}(G)$-valued forms
$\alpha$ and  $\kappa_{\mu}$. The assertion now follows from
the commutativity of the
 outer diagram and the definition of $\kappa_{\mu}'\,.$ 

    \begin{flushright}
        $\Box$
        \end{flushright}
        
  The proof of proposition \ref{bianchitors} is similar. We consider 
instead of \eqref{bianchigrdiag}  the 
 diagram of principal bundles 
  \begin{equation}
      \label{bianchidiag}
    \includegraphics{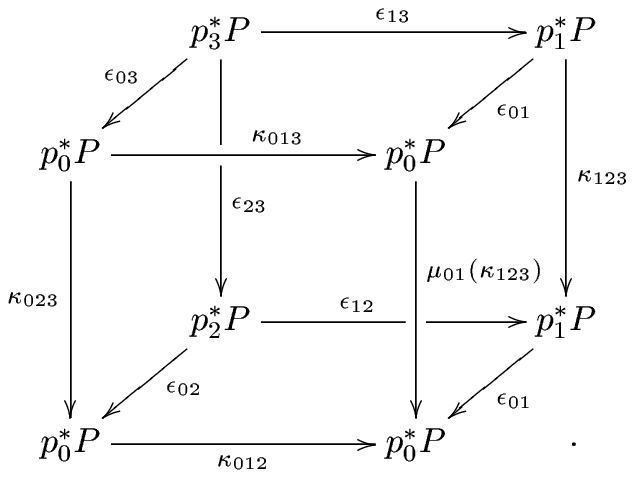} 
\end{equation}
Once more,  four of its  faces   commute by definition, since they are of the
pullbacks of  \eqref{def:kappa1}, and so does the right-hand one, 
which is of the form \eqref{def:ieps} Therefore, so does the front one,
whose commutativity
expresses  the equation \eqref{defkap0}. 
 We refer to this front square 
   \begin{equation}
       \label{bsq}
    \xymatrix{\pra P \ar[r]^{\kappa_{013}} \ar[d]_{\kappa_{023}} & \pra  
    P \ar[d]^{\epsilon^{\mathrm{ad}}_{01}(\kappa_{123})}\\
    \pra P \ar[r]_{\kappa_{012}}&\pra P
   }\end{equation}
   as the Bianchi square. 
   The second part of the proposition
 is proved by considering the diagram analogous to \eqref{compcon1} in which
all the groups $G$ have been replaced by the corresponding torsors $P$.
\begin{flushright}
        $\Box$
        \end{flushright}

\bigskip

\begin{remark}
    \label{abcase} {\rm  Let us assume once more that $G$ comes from an 
    $S$-group scheme, and is endowed wih its canonical connection.
    
    \medskip
    
  \qquad \qquad \quad    {\it i})  When
the group $\mathrm{Aut}(G)$ is representable, proposition \ref{bianchitors} 
     reduces locally, by  \cite{cdf} theorem 3.7,
     to the assertion that the local curvature 
     2-forms $\kappa_{i}$ \eqref{loccurv} associated to $\kappa$
     satisfy the classical Bianchi
     identity
     \bee
\label{bianchi:cl}\delta_{\omega_{i}}^{2}(\kappa_{i}) = d\kappa_{i}  + 
[\om_{i},\, \kappa_{i}] = 0\:,\end{equation}
with
$\om_{i}$  the corresponding 
local connection 1-form  associated as in \eqref{1con} to the 
connection $\epsilon$. This is in particular the case when $P = T_{G}$
is a trivial $G$-bundle,  with connection described
by a global $\mathrm{Lie}\, G$-valued 
    1-form $\om$ and associated curvature 2-form 
    \begin{equation}
        \label{glmc}
    \kappa:= d\om + 
    [\om]^{(2)}\:.\end{equation}
    
    \medskip

 \qquad \qquad \quad {\it ii}) Let  $G$  be the gauge group $\pa$ of a
$G$-torsor  with connection  $(P,\,\epsilon)$,  and
$\epsilon^{\mathrm{ad}}$ the induced   connection on $G$ induced by 
$\epsilon$. In that case 
proposition \ref{bianchigroup} for the pair $(G,\, \epsilon^{\mathrm{ad}})$ 
follows from
 proposition \ref{bianchitors}, by applying the construction 
 \eqref{mudef}
 to the  Bianchi square \eqref{bsq}.

\medskip

 \qquad \qquad \quad {\it iii}) Part {\it i}) of
 proposition \ref{bianchigroup} for $G= \pa$
is a formal consequence of part {\it i}) of proposition \ref{bianchitors}
since a necessary condition for  diagram \eqref{bsq} to commute is 
that both paths have the same equivariance property with respect 
to the left $\pa$ actions. This is precisely the content of \eqref{defkapmu0}

\bigskip

  }\end{remark}

\bigskip

Let $(G,\, \mu)$  be an $X$-group equipped with a connection.
 We consider the crossed module $[G \stackrel{i}{\la} \mathrm{Aut}(G)]$
 (with $i$ the 
inner conjugation homomorphism)
as a complex with $G$ placed in degree $-1$, and attach to it 
the commutative
diagram 
\begin{equation}
\label{interp-G} \xymatrix{
   \mathrm{Lie} \,G  \otimes \Omega^1_{X/S} 
   \ar[r]^{\delta^1_{\mu}}\ar@{|->}[d]
   & \mathrm{Lie} \,G  \otimes \Omega^2_{X/S} \ar[d]^{i}
\ar[r]^{\delta^2_{\mu}} &  \mathrm{Lie}\, G \otimes \Omega^3_{X/S}\ar[d]^{i} 
\ar[r]^(.7){\delta^3_{\mu}}&\: \cdots   \\ \{\mu\}
 \ar[r]^(.4){\delta^1} 
 & \mathrm{Lie}\,\mathrm{Aut}( G)  \otimes \Omega^2_{X/S}
\ar[r]^{\delta^2_{\mu^{\mathrm{ad}}}} &  \mathrm{Lie}\, 
\mathrm{Aut}(G) \otimes \Omega^3_{X/S}
\ar[r]^(.7){\delta^3_{\mu^{\mathrm{ad}}}} &\: \cdots
}\end{equation}
where the the term $ \{\mu\}$ is a section of  
$\mathrm{Conn(G)}$.  and the morphism $\delta^{1}$ assigns to a connection 
its curvature. The first vertical arrow represents the  action on 
$\mathrm{Conn}\,(G)$ of the sheaf
$\mathrm{Lie}\,(\mathrm{Aut}( G)  \otimes
\Omega^1_{X/S})$, so that an element $\alpha$ and sends $\mu$ to 
another section $\mu'$ of $\mathrm{Conn(G)}$.  
The commutativity of the first square expresses
that map $\delta^1$ is
$\delta^1_{\mu}$-equivariant. Since $\mu \in \mathrm{Conn}\,(G)$
functions  both as a term  in this diagram and  
as a parameter which determines the arrows,  
formally the diagram displayed 
should be viewed as the fiber over $\mu$ of a universal 
diagram defined over the sheaf $\mathrm{Conn}\,(G)$.

\bigskip

Let $P$ be a torsor on $X$ under an $S$-group and suppose that $P$ is 
endowed with a connection $\epsilon$. To this we associate a  similar
diagram
\begin{equation}
\label{interp} \xymatrix{
  \mathrm{Lie}\, \pa  \otimes \Omega^1_{X/S} \ar@{|->}[d]
  \ar[r]^{\delta^1_{\mu_{\epsilon}}}  &
  \mathrm{Lie}\, \pa  \otimes \Omega^2_{X/S} \ar[d]^{i}
\ar[r]^{\delta^2_{\mu_{\epsilon}}} & 
\mathrm{Lie}\, \pa \otimes \Omega^3_{X/S}\ar[d]^{i} 
\ar[r]^(.7){\delta^3_{\mu_{\epsilon}}}&\: \cdots   \\
\{\mu_{\epsilon}\}
 \ar[r]^(.4){\delta^1} 
 & \mathrm{Lie}\,\mathrm{Aut}( \pa)  \otimes \Omega^2_{X/S}
\ar[r]^{\delta^2_{(\mu_{\epsilon})^{\mathrm{ad}}}} &  \mathrm{Lie}\, 
\mathrm{Aut}(\pa) \otimes \Omega^3_{X/S}
\ar[r]^(.7){\delta^3_{(\mu_{\epsilon})^{\mathrm{ad}}}} &\: \cdots
}\end{equation}
where  we have set
\[\mu_{\epsilon}:= \epsilon^{\mathrm{ad}}\,.\]
The map  which associates to a connection $\mu_{\epsilon}$ on
 $\pa$ its curvature, 
 is  denoted once more  by $\delta^1$.
The  pair $(\mu_{\epsilon},\, \kappa_{\epsilon})$ attached to the 
conection $\epsilon$ may  be viewed as a 1-cocycle 
for the total complex associated to  diagram 
\eqref{interp}, since the  cocycle condition for this 1-cochain
is given by the pair of equations 
\[\delta_{\mu_{{\epsilon}}}^2 \, \kappa_{\epsilon} = 0\]
\[\delta^1\mu_{\epsilon}  = i(\kappa_{\epsilon})\,.\]
The former is the Bianchi identity, and the latter follows from the
functoriality of $i$ with respect to group isomorphisms. 
When $\epsilon$ is replaced by a new connection $\epsilon' = h \, 
\epsilon$, the  cocycle pair $(\mu_{\epsilon},\, \kappa_{\epsilon})$
is transformed  into a pair 
$(\mu_{\epsilon'},\,\kappa_{\epsilon'})$,
with   $\kappa_{\epsilon'}$ defined by \eqref{cobcap} and
\[ \mu_{\epsilon'} = i_{h}\,\mu_{\epsilon}\,.\] 
This action  of the element $h$ in 
the group $\mathrm{Lie}(\pa,\, \Omega^{1}_{X/S})$ on the cocycle 
pairs may be thought of as a coboundary relation. Once more, it 
transform the 2-form $\kappa_{\epsilon}$ into a 2-form $\kappa_{\epsilon'}$
which lives in the corresponding diagram indexed by $\mu'$, and 
therefore satisfies the corresponding 
Bianchi identity $\delta^{2}_{\mu_{\epsilon'}}\kappa_{\epsilon'} = 0$. 
In terms of the original pair $(\mu_{\epsilon},\, \kappa_{\epsilon})$ 
this may be restated (in additive notation) as
\[ \delta^{2}_{\mu + i_{h}}(\kappa_{\mu} + \delta^{1}h) = 0\,.\]
When $G$ comes from $S$ and $\mu$ is the canonical connection, this 
formula reduces to the formula
\[\delta^{2}_{i_{h}}(\delta^{1}h) = 0\]
of \cite{cdf} lemma 3.5.

%% file: dgg-ch2.tex
  \section{Gerbes and their gauge  stacks}
    \label{gerbes}
\subsection{}
\setcounter{equation}{0}%
\label{gerbes1}
In the present section, we do not assume that we are working in the 
\'etale topology. The original reference for the concept of a gerbe is 
\cite{Gir}. For the reader's convenience, we will now briefly review the 
description of gerbes 
in terms of  explicit non-abelian 2-cocycles,
given in  \cite{lb:2-gerbes}. We refer to that text for a  discussion 
of  fibered categories and stacks,  an alternative  source being \cite{lmb}. 

\bigskip

 Let $\pc$ be a gerbe on a scheme $X$. We
   choose a family of  objects $x_{i}$ in the fiber categories  
   $\pc_{U_{i}}$, for 
    some open cover $\mathcal{U}=(U_{i})_{i \in I}$ of $X$. These 
    objects determine  a  family of coefficient sheaves of groups
    $G_{i} = \mathrm{Aut}_{\pc}(x_{i})$,
    defined above the open sets $U_{i}$.
   Assume that we can choose, for each pair 
    $i, j \in I$, an arrow
     \bee
    \label{path}
    \phi_{ij}: (x_{j})_{\mid U_{ij}} \la (x_{i})_{\mid U_{ij}}
    \end{equation}
    in $ \pc_{U_{ij}}$.  The gerbe axioms do not 
    guarantee that there are such  paths, and their existence will mean 
    that the cocycle classifying $\pc$ lives in \v{C}ech cohomology, 
     rather than in the more general sheaf cohomology corresponding 
     to hypercovers of $X$. 
    The gerbe axioms 
    do however ensure that such paths exist locally. As explained in
    \cite{lb:2-gerbes},  this introduces 
      new families of upper indices for the cocycle pairs, and 
      coboundary formulas and  describes in a very 
    concrete manner  the cohomology with respect 
    to a hypercover of $X$ (rather than with respect to an  ordinary
    cover $\mathcal{U}$ of $X$ as in \v{C}ech cohomology). 
   When  $X$ is a scheme which is quasi-projective over a 
    ring and if we work in the \'etale topology, then by M. Artin's 
    theorem, \cite{martin} theorem 4.1, \v{C}ech covers are cofinal 
    among hypercovers and there is no need for upper  indices. 
    Similarly if we work in the topological context with  a 
    paracompact space, then the proof of \cite{spanier} chapter 6, 
    section 9,
    lemma 12 shows that upper indices are again 
    unnecessary. We will from now on always make the assumption that 
    the paths $\phi_{ij}$ \eqref{path} exist so that 
    this is indeed the case.
     
     \bigskip
     
   Under this assumption, we  may associate to  the gerbe $\pc$, 
    with chosen trivializing data 
    $(x_{i}, \phi_{ij})$, the cochains 
    \begin{equation}
        \label{def:lij0}
        \lambda_{ij}: G_{j} \la G_{i}\end{equation}
    in $\mathrm{Isom}((G_{j})_{\mid U_{ij}}, (G_{i})_{\mid 
    U_{ij}})$
    and $g_{ijk} \in \Ga(U_{ijk},\, G_{i})$ defined by\footnote{The 
    terms $\phi_{ij}$ in the definition of  $\lambda_{ij}$ in
    formula (2.4.1) of  \cite{lb:2-gerbes} are incorrectly ordered, 
    but the rest of the discussion there is correct.}
   \begin{alignat}{2}
       \label{def:lij}
        \lambda_{ij}(g) &:= (\phi_{ij})_{\ast}(g)
        =  \phi_{ij}\,g\, \phi_{ij}^{-1}\\
\label{defgijk} g_{ijk} &:= \phi_{ij}\, \phi_{jk }\,(\phi_{ik})^{-1}\:.
    \end{alignat}   
    The  pairs $(\lambda_{ij}, g_{ijk})$ satisfy the following 
    cocycle conditions,
    respectively defined above the open sets $U_{ijk}$ and 
    $U_{ijkl}$:
     \begin{alignat}{2}
        \lambda_{ij} \: \lambda_{jk} &= i_{g_{ijk}} 
        \lambda_{ik}\label{coclam}\\
        \lambda_{ij}(g_{jkl})\,g_{ijl}  &= g_{ijk} \, 
        g_{ikl}\:.\label{cocg}
        \end{alignat}
   Here $ i_{g_{ijk}}$ is the inner conjugation automorphism of 
   $G_{i}$, defined by 
    \bee
    \label{inconj}
    i_{g}(\gamma) := g\,\ga \, g^{-1}
    \end{equation}
    for any section $\ga$ of 
    $G_{i}$. 
    
    \bigskip

    Suppose now that we consider other choices $x'_{i}$ and $\phi_{ij}'$ 
    for the local objects and arrows in $\pc$ associated to the same 
    open cover $\mathcal{U}$ of $X$. These yield a different family 
    of local groups $G'_{i}\,$, and of 
    corresponding cocycle pairs $(\lambda'_{ij},\, g'_{ijk})$ with 
    $\lambda'_{ij}:G'_{j}\la G'_{i}$ and $g'_{ijk} \in \Gamma(U_{i},\, 
    G'_{i})$.
   With the same blanket assumption as above,  we can choose an arrow 
\bee
\label{chidef}
\chi_{i}: x_{i} \la x'_{i}\end{equation}
in 
$\pc_{U_{i}}$.   Conjugation by $\chi_{i}$ then  defines a group isomorphism
    \bee
    \label{def:imu}
\begin{array}{rcl}
    m_{i}: G_{i} &\la &G'_{i}\\
    g & \mapsto &\chi_{i} \, g \,\chi_{i}^{-1}\:.
    \end{array}
    \end{equation} 
The arrow $\chi_{i}$ also defines  a section
 $\delta_{ij} \in \Ga(U_{ij},\, G'_{i}) $,  by 
the equation
    \bee
    \label{def:delij}
    \phi_{ij}'\chi_{j} \: = \: \delta_{ij}\: \chi_{i} \, \phi_{ij}\:,
    \end{equation}
     so that  $\delta_{ij} $ 
     measures the defect in commutativity of the diagram 
    \[  \xymatrix{
    x_{j} \ar[r]^{\phi_{ij}} \ar[d]_{\chi_{j}}& x_{i} 
    \ar[d]^{\chi_{i}}\\
   x'_{j} \ar[r]_{\phi'_{ij}}   & x'_{i}
    }\]
     in 
    $\pc_{U_{ij}}$.  The pair $(m_{i}, \delta_{ij})$ satisfies 
    the following coboundary conditions (\cite{lb:2-gerbes}  (2.4.16) and 
    (2.4.17))\footnote{In the present context, \cite{lb:2-gerbes}
    (2.4.17) would 
    literally read
    \[  g_{ijk}' \, \delta_{ik} = \delta_{ij} 
    \,\Lambda_{ij}(\delta_{jk})\, m_{i}(g_{ijk})\]
   in the notation of \cite{lb:2-gerbes}. We make use here of
   (\ref{cocd1}) in order to transform this expression 
    into the more  pleasant (\ref{cocd2}).}: 
    \begin{alignat}{2}
        \lambda_{ij}'\,m_{j}&=  i_{\delta_{ij}}\, m_{i}\, \lambda_{ij}
        \label{cocd1}\\ g_{ijk}'\, \delta_{ik}&=  
        {\lambda}_{ij}'(\delta_{jk})\,\delta_{ij}\, 
        m_{i}(g_{ijk})\:.\label{cocd2}
        \end{alignat}

    \begin{remark}
        {\rm
         A difficulty , in dealing  with general gerbes, is 
that the cocycles correspond to some sort of cohomology with local 
coefficients, since they take their values in the family of sheaves 
of groups  $G_{i}$. The reader may  prefer 
 to assume at first  that the gerbes being considered are 
 $G$-gerbes, for a given $X$-group $G$. In that 
 case, the groups $G_{i}$ are the restrictions of  $G$ to the open sets 
 $U_{i}$. We still have a cocycle pair 
  $(\lambda_{ij}, \, g_{ijk})$ satisfying equations \eqref{coclam} 
  and \eqref{cocg}, with now $\lambda_{ij} \in \Gamma (U_{ij}, \,
  \mathrm{Aut}(G))$, and $g_{ijk} 
  \in \Ga(U_{ijk},\, G)$. Let alternate trivializing data $(x'_{i},\, 
  \phi'_{ij})$ for $\pc$ determine another  cocycle pair 
     $(\lambda'_{ij}\,,g'_{ijk})$ . The choice of a family of  arrows
     $\chi_{i}$ \eqref{chidef} now determines  sections $m_{i} \in 
     \Gamma(U_{i}, \, \mathrm{Aut}(G))$ \eqref{def:imu}
     and sections $\delta_{ij} \in 
     \Gamma(U_{ij},\, G)$ \eqref{def:delij}
  determining the coboundary relations  \eqref{cocd1} and \eqref{cocd2}.
   This shows that $G$-gerbes are 
  classified by elements in the \v{C}ech  cohomology set 
  $H^{1}(\mathcal{U},\, G \la \mathrm{Aut}(G))$  with 
  values,  as in \cite{fest}, in the coefficient crossed module
  $( G \la \mathrm{Aut}(G))$ determined by the inner conjugation 
  homomorphism \eqref{inconj}.
 The definition of such a non-abelian cohomology set associated to a 
 sheaf of group $G$ is due to Dedecker \cite{dede}, with a shift in degree.
  }
        \end{remark}
    
 \subsection{}
\setcounter{equation}{0}%
\label{adif}
Instead of examining  the change in cocycle pairs which a 
new  choice of the trivializing data  produces, we may consider an 
    equivalence (\cite{Gir} I Corollaire 1.5.2)  between a pair
    of   gerbes on $X$. This is  a morphism  of gerbes $u:\pc \la 
    \pc'$, together with  a morphism $u^{-1}: \pc' \la \pc$, 
    and a natural equivalence
\begin{equation}
\label{def:Theta}
\Theta: u^{-1}u \Longrightarrow 
    1_{\pc}\:.
\end{equation}
 We refer to $u^{-1}$ as a quasi-inverse of $u$. It is by 
    no means unique, nor is, for a fixed $u^{-1}$, the associated 2-arrow
    $\Theta$. The 2-arrow $\Theta$ determines an adjoint 2-arrow 
(\cite{maclane} IV \S 4, \cite{marcum} 1.5) 
\begin{equation}
\label{def:Theta-tilde}
\widetilde{\Theta}: u u^{-1} \Longrightarrow 1_{\pc'}\:,
\end{equation}
 such that the composite 2-arrow in each of 
the following  diagrams is the identity 2-arrow:
    \begin{equation}
        \label{diag:adj}
    \begin{array}{cccc}
        \xymatrix{
    \ar@{}[r]^(.99){\,}="6"&\pc \ar[dr] \ar@{}@<-.5ex>[dr]^{u} \ar[rr]^{1_{\pc}} 
    \ar@{}[rr]_(.45){\,}="1"&& \pc
    \\
     \pc'\ar[ur]
\ar@{}@<-1.5ex>[ur]^{u^{-1}} \ar[rr]_{1_{\pc'}}\ar@{}[rr]^{\,}="5" &&\pc'
\ar@{}[u]^(.01){\,}="2" \ar[ur]_{u^{-1}}&
    \ar@{}"5";"6"_(.2){\,}="7"
\ar@{}"5";"6"_(.6){\,}="8"
\ar@{=>}"7";"8"^{\widetilde{\Theta}^{-1}}
    \ar@{}"2";"1"_(.4){\,}="3"
\ar@{}"2";"1"_(.8){\,}="4"
\ar@{=>}"3";"4"_{\Theta}
} &\qquad&\qquad& \xymatrix{
\pc \ar[dr]_{u} \ar[rr]^{1_{\pc}}  \ar@{}[rr]_(.45){\,}="1"&& \pc 
\ar@{}[d]_(.01){\,}="6"\ar[dr]^{u} &\\
& \pc'\ar@{}[u]^(.01){\,}="2" \ar[ur]
\ar@{}@<-1.5ex>[ur]^{u^{-1}} \ar[rr]_{1_{\pc'}}
\ar@{}[rr]^(.45){\,}="5"&&\pc'\:.
\ar@{}"2";"1"_(.4){\,}="3"
\ar@{}"2";"1"_(.8){\,}="4"
\ar@{=>}"3";"4"^{\Theta}
\ar@{}"5";"6"_(.2){\,}="7"
\ar@{}"5";"6"_(.6){\,}="8"
\ar@{=>}"7";"8"_{\widetilde{\Theta}^{-1}}
}
\end{array}
\end{equation}

\bigskip

We have only been dealing so far  with 
    morphisms between the gerbes $\pc$ and $\pc'$, but 
    this discussion carries over to their restrictions
 $\pc_{\mid U}$ and $\pc'_{\mid U}$ over 
    varying open sets $U$ of $X$. We denote the 
    resulting fibered category of equivalences between $\pc$ and $\pc'$
 by $\epc$. By \cite{Gir} II corollary 2.1.5,
    $\epc$ is   a stack,
    just as  the presheaf $\mathrm{Hom}(F,\, F')$ 
    associated to sheaves  of groups $F$ and $F'$ on $X$ is automatically 
    a sheaf. It follows that the given description  of a global object $u 
    \in  \epc$ is correct as it stands, and does
 not call for  an additional level of sheafification.
 
 \bigskip

 \bigskip

There is a horizontal composition of 2-arrows 
\bee
\label{comp1}
\begin{array}{ccc}
\ep2c   \, \times \,\epc &\la &\mathcal{E}q(\pc, \pc'')\\
(u^{2},u^{1})& \mapsto & u^{2}\,u^{1}\\
(t^{2}, t^{1})& \mapsto & t^{2}\ast t^{1} 
\end{array}
\end{equation}
which associates to a composable pair of natural equivalences $t_{1}$ 
and $t_{2}$ 
\bee
\label{comp2}
 \xymatrix{
\pc
\ar@/^1pc/[rr]^{u^{1}}_{}="1" 
\ar@/_1pc/[rr]_{v^{1}}_{}="2" &
& \pc' \ar@/^1pc/[rr]^{u^{2}}_{}="3" 
\ar@/_1pc/[rr]_{v^{2}}_{}="4"&&
\pc''
\ar@{=>}"1";"2"^{t^{1}}
\ar@{=>}"3";"4"^{t^{2}}}
\end{equation}
the horizontally composed diagram
\bee
\label{comp3}
 \xymatrix{
\pc
\ar@/^1pc/[rrr]^{u^{2}u^{1}}_{}="1" 
\ar@/_1pc/[rrr]_{v^{2}v^{1}}_{}="2" &&
& \pc''\:. 
\ar@{=>}"1";"2"^{t^{2}\ast \,t^{1}}
}
\end{equation} 
The operation of horizontal composition on the left or on the right 
with an identity arrow $\mathrm{id}_{u}$ will be referred to as a 
whiskering.

\bigskip

The stack  of self-equivalences 
$\mathcal{E}q(\pc,\, \pc)$ of the gerbe $\pc$ will be called the gauge 
stack of $\pc$, and denoted $\pac$.
Setting $\pc = \pc' = \pc''$, the rule \eqref{comp1}  for horizontal 
composition completely  describes  the monoidal 
structure on the gauge stack $\pac$ of a gerbe $\pc$. While the associativity
 and  unit constraints are strict, this is not the case for the inverse law
in $\pac$.
 As we have seen, for every object $u  \in \pac$, there exists a (quasi-)
 inverse $u^{-1}$, together with a pair of compatible arrows \eqref{def:Theta} 
and \eqref{def:Theta-tilde}. Thus  $\pac$ is a group-like monoidal stack
 in groupoids. A 
stack in groupoids with  a group-like monoidal structure
is  sometimes called a  $gr$-stack (\cite{schreier} D\'efinition  3.1.1).
When we merely set $\pc'' = \pc'$ ({\it resp.}  $\pc = \pc'$)
 in \eqref{comp2}, 
we obtain  the left 
${\mathcal{P}'}^{\mathrm{ad}}$-torsor structure on the stack $\epc$ 
({\it resp.}
the right ${\mathcal{P}'}^{\mathrm{ad}}$-torsor structure  on the
stack $\ep2c$). We refer to \cite{fest} definition 6.1 for the 
precise definition of a  torsor under a $gr$-stack.

\bigskip

To  an object $(u,\, u^{-1},\, \Theta)$  of $\mathcal{E}q(\pc,\, \pc')$
corresponds,  functorially in 
$u$,
an
 equivalence of {\it gr}-stacks 
\begin{equation}
\label{def:ad}
\begin{array}{ccc}
\pac &\stackrel{u^{ad}}{\la} & \mathcal{P'}^{\mathrm{ad}}\\
w & \mapsto & uwu^{-1}
\end{array} \end{equation}
We will exhibit part of this data by  giving, for any object $w \in \pac$,
  a conjugation 2-arrow 
\begin{equation}
\label{msquare}
\xymatrix{
\pc \ar[rr]^w_(.99){\,}="2"  \ar[d]_u &&  \mathcal{P}
\ar[d]^u
\\
\mathcal{P'} \ar@{}[rr]_(.01){\,}="1" \ar[rr]_{u^{\mathrm{ad}}(w)}
 &&
\mathcal{P'}
\ar@{}"1";"2"_(.4){\,}="3"
\ar@{}"1";"2"_(.65){\,}="4"
\ar@{=>}"3";"4"^{M_u(w)\:}
}
\end{equation}
whose construction is functorial in $w$ and 
compatible with the {\it gr}-structure of $\pac$
and with the composition of 1-arrows $u$.
For any arrow
$t: u \Longrightarrow v$ in $\mathcal{E}q(\pc,\, \pc')$, 
we are also provided  with
a 2-arrow
\begin{equation}
    \label{def:tad}
    t^{\mathrm{ad}}(w): u^{\mathrm{ad}}(w) \Longrightarrow 
    v^{\mathrm{ad}}(w)
     \end{equation}
for which the following  degenerate prism, 
with left and right face $t$  and  bottom face $ t^{\mathrm{ad}}(w)$,
is commutative:
\begin{equation}
    \label{t-prism}
\xymatrix{
& \pc  \ar[rrrr]^{w}
\ar@{}[rrrr]_(.2){\,}="1"
 \ar[ldd]_(.2){u}
  \ar@{}[ldd]^(.7){\,}="5"
 \ar[d]^{v}
 \ar@{}[d]^{\,}="6"
\ar@{}[d]^(.3){\,}="2"
 &&&& \pc 
\ar[ldd]_(.2){u}  \ar@{}[ldd]_(.8){\,}="3"
\ar[d]^{v} 
\\
& \pc'  \ar[rrrr]^{v^{\mathrm{ad}}(w)} &&&& \pc'
\\
\pc' \ar@{=}[ur] 
 \ar[rrrr]^{u^{\mathrm{ad}}(w)} 
 \ar@{}[rrrr]^(.8){\,}="4"
 \ar@{=}[ur] &&&& \pc'
 \ar@{=}[ur] & \:.
\ar@{=>}"2";"1"_{M_{v} (w)}
\ar@{=>}"4";"3"^{M_{u} (w)}
\ar@{}"5";"6"
}
\end{equation} The  construction
\eqref{t-prism}
 is  compatible
with  both  vertical and horizontal composition of such 2-arrows
 $t:u \Longrightarrow  v$.
 Suppose now that $u \in \pac$, in other words that $\pc = \pc'$. In 
 that case the construction \eqref{def:ad} defines a monoiddal  ``inner 
 conjugation'' morphism
 \begin{equation}
     \label{def:jmap}
     \begin{array}{ccc}
 \pac  & \stackrel{j}{\la} & \mathcal{E}q(\pac)\\
 u & \mapsto & u^{\mathrm{ad}}\end{array}\end{equation}
 which we will prefer to call $j$ rather than the customary $i$ in 
 order to avoid confusion with the inner conjugation map $i$ for groups.  
 We will sometimes write $j(u)$ as $j_{u}$ when we wish to emphasize 
 the similarity between $j$ and $i$. 

\bigskip

Since the construction of $u^{\mathrm{ad}}$ is functorial in $u$, there
corresponds to any diagram in the 2-category of gerbes $\pc$ on $X$ (with 
equivalences of gerbes $u$ as 1-arrows and natural equivalences $t$ between
these as 2-arrows) a similar diagram in the  2-category of 
{\it gr}-categories, whose objects, 1- and 2-arrows are respectively 
replaced by $\pac, \, u^{\mathrm{ad}}$ and $t^{\mathrm{ad}}$. 
We will refer to this process as the {\it adification} of the original diagram.
A related construction is that given by whiskering. A typical instance 
in which we will make use of whiskering is the following. Consider a diagram 
\begin{equation}
\label{def:alph}
\xymatrix{
 \qc
 \ar@/^1pc/[rr]^u
 \ar@/_1pc/[rr]_v
 \ar@{}@/^1pc/[rr]_{}="1" 
 \ar@{}@/_1pc/[rr]^{}="2" 
 && \pc 
\ar@{}"1";"2"^(.1){\,}="3" 
\ar@{}"1";"2"^(.9){\,}="4" 
\ar@<-1ex>@{=>}"3";"4"^{\:\alpha}
}
\end{equation}
where $u$ is an equivalence.
We can  whisker it on the left by the given (quasi-)inverse
   $u^{-1}$ of $u$ and  compose  with the 2-arrow 
 $\widetilde{\Theta}$ \eqref{def:Theta-tilde}:
\begin{equation}
\label{def:alph-t}
\xymatrix{
 \pc \ar[rr]_{u^{-1}}
\ar@{}[rr]^(.999){\,}="6"
\ar@/^2pc/[rrrr]^(.49){1}
\ar@{}@/^2pc/[rrrr]_(.49){\,}="5"
&&
 \qc
 \ar@/^1pc/[rr]^(.4){u}
 \ar@/_1pc/[rr]_v
 \ar@{}@/^1pc/[rr]_{}="1" 
 \ar@{}@/_1pc/[rr]^{}="2" 
 &&  \pc  \:, 
\ar@{}"1";"2"^(.2){\,}="3" 
\ar@{}"1";"2"^(.9){\,}="4" 
\ar@<-1ex>@{=>}"3";"4"^{\:\alpha}
\ar@{=>}"5";"6"_(.35){\widetilde{\Theta}\,}
}
\end{equation}
thereby defining    a 2-arrow
 $\tilde{\alpha}:
 I \Longrightarrow v\,u^{-1}$ sourced at  $I= \mathrm{id}_{\pc}$.
 Conversely, the 2-arrow
$\alpha$  can be recovered from  $\tilde{\alpha}$ by
 whiskering $\tilde{\alpha}$ by $u$  and then composing with the 2-arrow $\Theta$
 \eqref{def:Theta} adjoint to $\widetilde{\Theta}$. 

 \bigskip
 
 Such a construction can  be applied in a number 
 of related situations, as we illustrate by the following example.
 Suppose   that we start from a 
 2-arrow
  \begin{equation}
\label{def:alph-2}
\xymatrix{
 \rc
 \ar@/^1pc/[rr]^{v \,u}
 \ar@/_1pc/[rr]_{w \,u}
 \ar@{}@/^1pc/[rr]_{}="1" 
 \ar@{}@/_1pc/[rr]^{}="2" 
 && \pc 
\ar@{}"1";"2"^(.1){\,}="3" 
\ar@{}"1";"2"^(.9){\,}="4" 
\ar@<-1ex>@{=>}"3";"4"^{\:\beta}
}
\end{equation}
with $u: \rc \la \qc$ once more an equivalence,  and $v$ and $w$ a pair 
of 
arrows from $\qc$ to $\pc$. Whiskering on the left 
by $u^{-1}$ and 
applying $\widetilde{\Theta}$ and its inverse, this induces a diagram   
\begin{equation}
\label{def:alph-3}
\xymatrix{
 \qc
 \ar@/^1pc/[rr]^v
 \ar@/_1pc/[rr]_w
 \ar@{}@/^1pc/[rr]_{}="1" 
 \ar@{}@/_1pc/[rr]^{}="2" 
 && \pc 
\ar@{}"1";"2"^(.1){\,}="3" 
\ar@{}"1";"2"^(.9){\,}="4" 
\ar@<-1ex>@{=>}"3";"4"^{\:\widetilde{\beta}}\,.
}
\end{equation}
It follows by adjointness of $\Theta$ and $\widetilde{\Theta}$ that 
$\widetilde{\beta}$ is the unique 2-arrow inducing $\beta$ when 
whiskered on the left by $u$. The giving of $\beta$ and  $\widetilde{\beta}$
are therefore equivalent. We will in the sequel refer to 
any such  composite construction   as a whiskering.

\bigskip

\subsection{}
\setcounter{equation}{0}%
\label{algar}
Let us now choose 
    trivializations of both gerbes $\pc$ and  $\pc'$ on the same 
    family of 
    open sets $\mathcal{U}$  by objects 
    $x_{i},\:x_{i}'$ and corresponding 
    arrows 
     \[
\begin{array}{ccc}
    \phi_{ij}: x_{j} \la x_{i}& \qquad \qquad &\phi_{ij}': x_{j}' \la 
    x_{i}'\:.
   \end{array}
    \] 
Once more, the trivializing pairs $(x_{i},\, \phi_{ij})$
and $(x'_{i},\, \phi'_{ij})$ determine corresponding local groups 
$G_{i}$ and $ G'_{i}$, together with
associated cocycle pairs $(\lambda_{i}, \, 
g_{ijk}), \, (\lambda'_{i}, \, 
g'_{ijk})$
The discussion which led up to  the coboundary terms $(m_{i}, \, 
\de_{ij})$
satisfying equations \eqref{cocd1} and \eqref{cocd2}
 may be reproduced in the present context, and now leads to a 
cocyclic description of the equivalence $u$.
Let us begin by   choosing a family of arrows
    \begin{equation}
        \label{upath}
        \ga_{i}:u(x_{i})  \la  x_{i}'\end{equation}
    in $\mathrm{Ar}(\pc')_{U_{i}}$, as we may always do by 
    backtracking and refining 
  the original open cover $\mathcal{U}$ of $X$  if necessary.
 To any such a family of
  arrows 
  $\ga_{i}$ is
    associated a family of pairs $(m_{i}, \delta_{ij})$, 
    with $m_{i} \in 
    \mathrm{Isom}(G_{i},\, G_{i}')$ defined by 
   \[m_{i}(g) := \ga_{i}\,u(g)\, \ga_{i}^{-1} \] for 
    any $g: x_{i }\la x_{i}$ in $G_{i}$, 
    and $\delta_{ij} \in \Gamma(U_{ij}, G_{i}')$ determined by the 
    equation 
    \bee
    \label{def:delij1}
    \phi_{ij}'\ga_{j} \: = \: \delta_{ij}\: \ga_{i} \, u(\phi_{ij})
    \end{equation}
    in $\pc'_{U_{ij}}$ analogous to \eqref{def:delij}.  The pairs
    $(m_{i}, \delta_{ij})$
    completely determine the arrow $u$, and satisfy anew the 
    coboundary conditions \eqref{cocd1} and \eqref{cocd2}.
        We will use the shorthand notation
        \bee
        \label{displ1}
        \xymatrix{
(  \lambda_{ij}\,,\,  g_{ijk} )
\ar[rr]^{(m_{i}\,,\, \delta_{ij})}
 && (  \lambda_{ij}'\,,\, g_{ijk}' )}\end{equation}
to display   this cocyclic description
        of the  equivalence  of gerbes $u: \pc \la \pc'$.

 \bigskip
    
    A natural 
    transformation between a pair of equivalences $u,\,  v$ from $\pc$ 
    to
    $\pc'$
    can be described in similar cocyclic
    terms. This is discussed  in 
    \cite{lb:tc} \S  5 when $X$ is the spectrum of a field endowed 
    with the \'etale topology, 
    and can be restated as follows for an open cover 
    $\mathcal{U}$ of $X$.
    Given such an arrow $t$
 \bee  \label{diagP} \xymatrix{
 \pc \ar@/^1pc/[rr]^{u}_{}="1" 
\ar@/_1pc/[rr]_{v}_{}="2" && \pc'
\ar@{=>}"1";"2"^{t}}\end{equation}
    in $\mathcal{E}q(\pc, \pc')$,
    and assume that the morphisms $\ga^{u}_{i}$ (\ref{upath}) 
    associated to $u$ and the corresponding morphisms
    \bee \ga^{v}_{i}: v(x_{i}) \la x'_{i}    \end{equation}
    associated to $v$ have been chosen, yielding corresponding 
    families of the  pairs 
    $({m}_{i}^{u}, {\delta}^{u}_{ij})$ and
    $(m^{v}_{i}, \delta^{v}_{ij})$.
    The natural 
    transformation $t: u \Longrightarrow v$ determines, for each 
    object $x_{i}$, an arrow
    \[\xymatrix{u(x_{i}) \ar[r]^{t(x_{i})} & v(x_{i})}
    \] in $\pc'$.
The  transformation    $t$ is  entirely 
    described by the associated  0-cochain $\theta_{i} \in \Gamma (U_{i}, 
    G_{i}')$ defined by 
     \begin{equation}
        \label{defthet}
        \theta_{i}:= \ga^{v}_{i}\, t(x_{i}) \,(\ga_{i}^{u})^{-1}\:,
        \end{equation}
        and  the diagram
        \[\xymatrix{
        u(x_{i}) \ar[r]^{t(x_{i})} \ar[d]_{\ga^{u}_{i}} & v(x_{i})
        \ar[d]^{\ga^{v}_{i}}\\
        x'_{i}\ar[r]_{\theta_{i}}& x'_{i}
                }\] therefore  commutes. 
Since $t$ is a natural transformation,   so does the  diagram
\[
\xymatrix{
& G_{i} \ar[ld]_{u} \ar[rd]^{v}&\\
\mathrm{Aut}_{\pc'}(u(x_{i})) 
\ar[rr]_{t(x_{i})_{\ast}}&&\mathrm{Aut}_{\pc'}(v(x_{i}))\,,
}\]
and (\ref{defthet}) therefore yields by conjugation
a  condition
\begin{equation}
    \label{cocthet1}
    m^{v}_{i} \:=\: i_{\theta_{i}} m^{u}_{i}\end{equation} 
which the $\theta_{i}$ must satisfy. 
Comparing the terms $\delta^{u}_{ij}$ and $\delta^{v}_{ij}$
respectively associated by \eqref{def:delij1} to the equivalences $u$ 
and $v$ yields the further  condition
\begin{equation}
    \label{cocthet2}
    \lambda_{ij}'(\theta_{j}) = \delta^{v}_{ij}\, \theta_{i}\, 
(\delta_{ij}^{u})^{-1}\:.\end{equation}
 We  will find it convenient to 
 display symbolically the diagram \eqref{diagP} as
\bee
\label{defthi}
 \xymatrix{
(  \lambda_{ij}\,,  g_{ijk} )
\ar@/^1pc/[rr]^{(m^{u}_{i}\,,\, \delta_{ij}^{u})}_{}="1" 
\ar@/_1pc/[rr]_{(m^{v}_{i}, \delta^{v}_{ij})}_{}="2" &
& (  \lambda_{ij}'\,,\, g_{ijk}' )\:.
\ar@{=>}"1";"2"^{\theta_{i}}}
\end{equation}
In writing such a diagram it is understood that the local objects
 $x_i$ in $\pc$ and $\pc'$ 
have been fixed, as well as the arrows $\phi_{ij}$ \eqref{path},
 $\gamma_{i}^{u}$, $\gamma_{i}^{v}$ \eqref{upath}.
The expressions
$(  \lambda_{ij}\,,  g_{ijk} )$,$\:(  \lambda_{ij}'\,,\, g_{ijk}' )$
then satisfy (\ref{coclam}), (\ref{cocg}),
$(m^{u}_{i}\,,\, \delta_{ij}^{u})$,  $(m^{v}_{i}, \delta^{v}_{ij})$ 
satisfy  (\ref{cocd1}), (\ref{cocd2}) and  $\theta_{i}$ satisfies 
\eqref{cocthet1}, \eqref{cocthet2}. 
The vertical composition of the 2-arrows
 in diagram 
 \[
\xymatrix{
(  \lambda_{ij}\,,\,  g_{ijk} )
\ar@/^2pc/[rrr]^{(m^{u}_{i}\,,\, \delta^{u}_{ij})}="1" 
\ar@{}[rrr]_{}="2"
\ar[rrr]_>>>>>>{(m^{v}_{i}, \delta^{v}_{ij})}
\ar@/_2pc/[rrr]_{(m^{w}_{i}, \delta^{w}_{ij})}="3"
 &&& (  \lambda_{ij}'\,,\, g_{ijk}' )
\ar@{}"1";"2"^(0.3){}="4"
\ar@{}"1";"2"^(0.9){}="5"
\ar@{=>}"4";"5"_{\theta_{i}\,}
\ar@{}"2";"3"^(0.1){}="6"
\ar@{}"2";"3"^(0.7){}="7"
\ar@{=>}"6";"7"_{\tilde{\theta}_{i}\,}}\]
{\it i.e.} 
composition of arrows in the category $\epc$, is given by  the rule
      \bee
      \label{compvert}
 \xymatrix{
(  \lambda_{ij}\, , \,  g_{ijk} )
\ar@/^1pc/[rr]^{(m_{i}^{u}\,,\, \delta_{ij}^{u})}_{}="1" 
\ar@/_1pc/[rr]_{(m^{w}_{i}, \delta^{w}_{ij})}_{}="2" &
& (  \lambda_{ij}'\,,\, g_{ijk}' )\:.
\ar@{=>}"1";"2"^{\tilde{\theta}_{i}\theta_{i}}}
\end{equation}

 When the horizontally composable 
diagram \eqref{comp2} is displayed as
\bee
\label{comp5}
 \xymatrix{
(  \lambda_{ij},\,  g_{ijk} )
\ar@/^1pc/[rr]^{(m_{i}^{u^1}\,,\, \delta^{u^1}_{ij})}_{}="1" 
\ar@/_1pc/[rr]_{(m^{v^1}_{i}, \delta^{v^1}_{ij})}_{}="2" &
& (  {\lambda'}_{ij}\,,\, {g'}_{ijk} )
\ar@/^1pc/[rr]^{(m^{u^2}_{i}\,,\, 
\delta^{u^2}_{ij})}_{}="3" 
\ar@/_1pc/[rr]_{(m^{v^2}_{i}, \delta^{v^2}_{ij})}_{}="4"&&
(   \lambda_{ij}''\,,\, g_{ijk}'' ) 
\ar@{=>}"1";"2"^{\theta^{1}_{i}}
\ar@{=>}"3";"4"^{\theta^{2}_{i}}}\:,
\end{equation}  
the horizontally composed diagram \eqref{comp3} becomes
 \bee
\label{comp6}
 \xymatrix{
(  \lambda_{ij}\,,\,  g_{ijk} )
\ar@/^1pc/[rrrr]^{(m_{i}^{u^2}m_{i}^{u^1}\,,\, 
\delta^{u^2}_{ij}m_{i}^{u^2}(\delta^{u^1}_{ij}))}_{}="1" 
\ar@/_1pc/[rrrr]_{(m^{v^2}_{i}m^{v^1}_{i}\,,\,
\delta^{v^2}_{ij}
m_{i}^{v^2}(\delta^{v^1}_{ij}))}_{}="2" &
&&& (  \lambda_{ij}''\,,\, g_{ijk}'' )
\ar@{=>}"1";"2"^{\,\theta^{2}_{i}m_{i}^{u^2}(\theta^{1}_{i})}}\:.
\end{equation}
In particular,  setting $\theta_{i}^{2}$ or $\theta_{i}^{1}$ equal to 
1, we see that the diagram 

\bee
    \label{comp5a}
     \xymatrix{
(  \lambda_{ij},\,  g_{ijk} )
\ar@/^1pc/[rr]^{(m_{i}^{u}\,,\, \delta^{u}_{ij})}_{}="1" 
\ar@/_1pc/[rr]_{(m^{v}_{i}, \delta^{v}_{ij})}_{}="2" &
& (  {\lambda'}_{ij}\,,\, {g'}_{ijk} ) \ar[rr]^{(m^{u'}_{i}\,,\, 
\delta^{u'}_{ij})}_{}="3" 
&&
(  \lambda_{ij}''\,,\, g_{ijk}'' ) 
\ar@{=>}"1";"2"^{\theta}}
\end{equation}
composes to 
 \bee
\label{comp6a}
 \xymatrix{
(  \lambda_{ij}\,,\,  g_{ijk} )
\ar@/^1pc/[rrrr]^{(m_{i}^{u'}m_{i}^{u}\,,\, 
\delta^{u'}_{ij}m_{i}^{u'}(\delta^{u}_{ij}))}_{}="1" 
\ar@/_1pc/[rrrr]_{(m^{u'}_{i}m^{v}_{i}\,,\,
\delta^{u'}_{ij}
m_{i}^{u'}(\delta^{v}_{ij}))}_{}="2" &
&&& (  \lambda_{ij}''\,,\, g_{ijk}'' )
\ar@{=>}"1";"2"^{\,m_{i}^{u'}(\theta_{i})}}\:,
\end{equation}
and the diagram 
\bee
\label{comp5b}
 \xymatrix{
(  \lambda_{ij},\,  g_{ijk} )
\ar[rr]^{(m_{i}^{u}\,,\, \delta^{u}_{ij})}_{}="1" 
&& (  {\lambda'}_{ij}\,,\, {g'}_{ijk} ) \ar@/^1pc/[rr]^{(m^{u'}_{i}\,,\, 
\delta^{u'}_{ij})}_{}="3" 
\ar@/_1pc/[rr]_{(m^{v'}_{i}, \delta^{v'}_{ij})}_{}="4"&&
(   \lambda_{ij}''\,,\, g_{ijk}'' ) 
\ar@{=>}"3";"4"^{\theta'_{i}}}
\end{equation}  
to 
 \bee
\label{comp6b}
 \xymatrix{
(  \lambda_{ij}\,,\,  g_{ijk} )
\ar@/^1pc/[rrrr]^{(m_{i}^{u'}m_{i}^{u}\,,\, 
\delta^{u'}_{ij}m_{i}^{u'}(\delta^{u}_{ij}))}_{}="1" 
\ar@/_1pc/[rrrr]_{(m^{v'}_{i}m^{u}_{i}\,,\,
\delta^{v'}_{ij}
m_{i}^{v'}(\delta^{u}_{ij}))}_{}="2" &
&&& (  \lambda_{ij}''\,,\, g_{ijk}'' )
\ar@{=>}"1";"2"^{\,\theta'_{i}}}\:.
\end{equation}

\bigskip

We now specialize from $\epc$ to the gauge stack  $\pac$ the cocyclic 
 description of objects and arrows  which we  obtained in \eqref{defthi}.
In the present situation, we have $x_{i} = 
x_{i}'$ and $\phi_{ij} = \phi_{ij}'$, so that $g_{ijk} = g_{ijk}'$ 
(\re $\lambda_{ij} = \lambda_{ij}'$). The elements $m_{i}^{u}$ and 
$m_{i}^{v}$  both live in 
$\mathrm{Aut}(G_{i})$,  $\theta_{i}$ is a section of $G_{i}$ above
$U_{i}$,  and 
$\delta_{ij}^{u},\, \delta_{ij}^{v} \in \Gamma (U_{ij}, G_{i})$.  
 Condition \eqref{cocthet1} remains unchanged. The conditions
 (\ref{cocd1}), (\ref{cocd2})  associated 
to the diagram 
\bee
\label{comp6c}
 \xymatrix{
(  \lambda_{ij}\,,\,  g_{ijk} )
\ar@/^1pc/[rr]^{(m_{i}^{u}\,,\, \delta_{ij}^{u})}_{}="1" 
\ar@/_1pc/[rr]_{(m^{v}_{i}, \delta^{v}_{ij})}_{}="2" &
& (  \lambda_{ij}\,,\, g_{ijk} )
\ar@{=>}"1";"2"^{\theta_{i}}}
\end{equation}
of objects and arrow in $\pac$ now read as 
 \begin{alignat}{2}
        \lambda_{ij}\,m_{j}^{u} &= i_{\delta_{ij}}\,
        m_{i}^{u}\,   \lambda_{ij}
        \label{coc3}\\
        g_{ijk}\, \delta^{u}_{ik} &=  
        {\lambda}_{ij}(\delta^{u}_{jk})\,\delta^{u}_{ij}\, 
        m_{i}^{u}(g_{ijk})\label{coc4}
      \intertext{and similar conditions are satisfied by the pair
        $(m_{i}^{v},\,\delta_{ij}^{v})$. 
        The   counterpart of   \eqref{cocthet2} is the equation}
        \lambda_{ij}(\theta_{j})  &=  {\delta}^{v}_{ij}\, \theta_{i}\, 
(\delta_{ij}^{u})^{-1}   \:.\label{cocthet3}
        \end{alignat}  
        which relates 
$\theta_{i}$ to $\theta_{j}$.

        \bigskip

Vertical composition, {\it i.e.} composition of arrows in $\pac$ is 
given by the rule \eqref{compvert}, in other words by ordinary 
multiplication. Horizontal composition, in other words the monoidal 
structure on $\pac$,  is given by specializing the composition formulas 
for 1- and 2-arrows in \eqref{comp6}. Once more, as in \eqref{comp6b},
whiskering on the 
left has no effect on the cochain  attached to a 2-arrow, whereas  
whiskering on the right does, as in  
\eqref{comp6a}.

%% file: dgg-ch3.tex
\section{Morita theory for locally trivialized gerbes}
\label{sec:morita0}
 \subsection{}
\setcounter{equation}{0}
\label{sec:morita}
Before we reinterpret the gauge stack, we will first review Giraud's Morita 
theorem \cite{Gir} IV proposition  5.2.5, paying particular attention to  questions
of variance.  Recall that in its most basic form, it asserts that any 
equivalence 
\bee
\label{mor1}
\mathrm{Tors} (X, \, G) \stackrel{\Phi}\la \mathrm{Tors} (X,\,H)
\end{equation}
between a pair of trivial gerbes on a scheme 
$X$ is described by the $(H,G)$-bitorsor
$P_{\Phi}$ defined on $X$  by 
\bee 
\label{mor2}P_{\Phi} := \mathrm{Isom}(\Phi(T_{G}), T_{H})\:,\end{equation}
where $T_{G}$ is the trivial right $G$-torsor on $X$.
The actions of $H$ and $G$ 
 on $P_{\Phi}$ come from the composition of isomorphisms, taking 
into account the isomorphism
{\renewcommand{\arraystretch}{1.3}
\begin{equation}\label{mor2.5}
    \begin{array}{ccc}
G &\simeq &\mathrm{Isom}_{G}(T_{G},\, 
T_{G})\\
g& \mapsto & (1 \mapsto g)\end{array}\end{equation}}
\noindent  and the corresponding isomorphism for $H$.
 A natural 
 transformation \bee
 \label{phi1}\psi: \Phi_{1} \Longrightarrow 
 \Phi_{2}\end{equation} between  a pair of such functors
 determines a morphism of $H$-torsors
 \[\Phi_{1}(T_{G}) \la \Phi_{2}(T_{G})\] 
 and therefore by composition a morphism of bitorsors 
 \bee
 \label{phi2}
 P_{\psi}: P_{\Phi_{2}} \la P_{\Phi_{1}}\:,\end{equation}
 so that the functor $\Phi \mapsto P_{\Phi}$ is contravariant in the 
 equivalence $\Phi$.
 
 \bigskip

A quasi-inverse functor  
associates to a bitorsor  $P$  the morphism $\Phi_{P}$ 
\eqref{mor1} defined by   
\bee
\label{mor3}\Phi_{P} (Q) = Q  \wedge P^{0}\:.\end{equation}
Here $P^{0}$ is the $(G,H)$-bitorsor opposite to $P$, whose
 underlying sheaf of sets coincides with that of $P$, but with the right
 $G$-action 
on $P$ transferred in the usual manner to a left action on $P^{0}$
by the  rule $g \ast p := pg^{-1}$, 
and  the  
  action of $H$ similarly transferred from left to right.
  We will sometimes write $P^{-1}$ for $P^{0}$.  The 
  ``opposite'' map $P\mapsto P^{0}$ is actually a contravariant functor
{\renewcommand{\arraystretch}{1.3}
  \[\begin{array}{ccc}
  \mathrm{Bitors}(H,\, G) &\la &  \mathrm{Bitors}(G,\, H)\\
  P & \mapsto & P^{0}\end{array}
  \] }
 which associates to the $(H,\,G)$-bitorsor morphism $f: P\la P'$ the 
 bitorsor morphism 
 \bee
 \label{def:po}
 f^0: (P')^{0} \la  P^{0}\end{equation}
 defined by
 \[f^{0}(p') = p\qquad \Longleftrightarrow  \qquad f(p)= 
 p'\:.\]
 This functor is compatible with the contracted product in 
 the following sense: given an $(H,G)$-bitorsor $P_{1}$ and an 
 $(K,H)$-bitorsor $P_{2}$, there is a canonical 
 isomorphism $(G,K)$-bitorsors 
{\renewcommand{\arraystretch}{1.3}
\bee
 \label{opos}
\begin{array}{ccc}
   (P_{2} \wedge^H  P_{1})^{0} & \simeq & P_{1}^{0} \wedge^{H} 
   P_{2}^0\\
   (p_{2},p_{1}) & \mapsto  & (p_{1},p_{2})\:,
   \end{array}
   \end{equation}}
and this isomorphism is compatible 
with the associativity isomorphisms \[P_{3}\wedge^{K} (P_{2}
\wedge^{H} P_{1}) 
\simeq (P_{3}\wedge^{K}  P_{2}) \wedge^{H} P_{1}\:.\]
When $P$ has a global section, and is therefore described
by an isomorphism $u: H \la G$, the associated functor
$\Phi_P$ is the familiar ``extension of the structural group'' functor
$v_\ast$ associated to $v:=u^{-1}$.
 
\bigskip

The  Morita theorem can now  be  summarized as 
follows (we henceforth denote the stack of $G$-torsors on $X$ 
by $\mathrm{Tors}(G)$ when the context is clear):
\begin{proposition}
    \label{morpro}The map $\Phi 
\mapsto P_{\Phi}$ \eqref{mor2} defines an anti-equivalence 
between the stack  of torsor equivalences
%\linebreak
 $\mathcal{E}q(\mathrm{Tors}(G), \, 
 \mathrm{Tors}(H))$  
 and the stack  of $(H,G)$-bitorsors on $X$, with a  quasi-inverse functor   
 defined by $P \mapsto \Phi_{P}$ \eqref{mor3}. This anti-equivalence 
 is coherently compatible with the composition of equivalences 
 between stacks of torsors.\end{proposition}
 In particular,
 the natural 
transformation $\Phi_{f}$ \eqref{phi1} associated to  a bitorsor isomorphism
$f: P_{2} \la P_{1}$ 
is defined by 
   \begin{equation}
\label{def:phi-f}
\xymatrix@R=6pt{
 \Phi_{P_{1}}(Q) \ar[rr]^{\Phi_{f}(Q)}&& \Phi_{P_{2}}(Q)\\
 Q \wedge (P_{1})^0 \ar@{|->}[rr]^{1_{Q}\wedge f^{0}}
&& Q \wedge (P_{2})^0}
\end{equation}
for any $G$-torsor $Q$.

\bigskip

 The correspondence of proposition \ref{morpro}
is multiplicative, 
in the sense that for any pair of composable
equivalences
\[\xymatrix{\to(G) \ar[r]^{\Phi_{1}} & \to(H) \ar[r]^{\Phi_{2}}& \to 
(K)}\]
there is a natural transformation 
{\renewcommand{\arraystretch}{1.5}
\[\begin{array}{ccc}  P_{\Phi_{2}}\wedge P_{\Phi_{1}} & \simeq &
P_{\Phi_{2}\circ \Phi_{1}}\\
u_{2} \wedge u_{1} & \mapsto & u_{2}\:\Phi_{2}(u_{1})\end{array}\]
}
and this  is compatible with the associativity isomorphism for the 
wedge product of bitorsors. Conversely, for any pair of multipliable 
bitorsors $P_{2}$ and 
$P_{1}$, 
the canonical  bitorsor isomorphism \eqref{opos} determines a natural 
transformation 
  \begin{equation}
\label{comp-phi}
 \Phi_{P_{2}} \,\circ\,\Phi_{P_{1}} \simeq \Phi_{P_{2}\wedge P_{1}}
\end{equation}
compatible with the associativity isomorphisms
\[P_{3} \wedge (P_{2}\wedge P_{1})
\simeq (P_{3} \wedge P_{2}) \wedge P_{1}\]
between contracted products  of bitorsors.

\bigskip

\subsection{}
\setcounter{equation}{0}
 Instead of starting  here from the  global 
 definition of a gerbe $\pc$, or from the purely local description of 
 \S 2, in which  both a family of local
objects $x_{i} \in \mathrm{ob}\,\pc_{U_{i}}$ and a family of arrows 
$\phi_{ij}$ (\ref{path}) have been chosen, we adopt here 
a semi-local approach, in which we merely 
choose the locally trivializing objects $x_{i}$, but no arrows (\ref{path}). 
This point of view is presented in \cite{fest} (see also
\cite{lb:tc} \S 2.4-2.5), as well as in the work of K.-H. Ulbrich
 \cite{ul}, \cite{ul1}). It  also occurs in \cite{K1} \S 3.13-3.16
in an additive
 category context framework in which the 
 stacks $\mathcal{C}_{U_{i}}$ are stacks  of modules rather 
 than of torsors. More recently, it has 
 been advocated  for abelian gerbes
in related contexts  by N. Hitchin  
\cite{hitch}, and by M.K. Murray  under the terminology of bundle 
gerbes \cite{murray}.

\bigskip

Let $\pc$ be a gerbe on $X$.
The choice of  objects  $x_{i} \in \mathrm{ob}\,\pc_{U_{i}}$
determines once more for each $i \in I$ an $U_{i}$-group  $G_{i}: = 
\mathrm{Aut}_{\pc}(x_{i})$, and 
an equivalence
of gerbes 
\bee
\label{triv}
\begin{array}{ccc}
\pc_{\mid U_{i}}&\stackrel{ \phi_{i}}{\la}&\mathrm{Tors}(U_{i},\, G_{i})\\
y& \mapsto & \mathrm{Isom}_{\pc_{\mid U_{i}}}(x_{i}, y)
\end{array}\end{equation}
above 
$U_{i}$, which sends $x_{i}$ itself to the trivial $G_{i}$-torsor.
A gerbe $\pc$ on $X$ may be 
constructed by gluing together trivial
 gerbes 
$\mathrm{Tors}(G_{i})$. The gluing 
data consists here in  a family of morphism of gerbes above $U_{ij}$
\begin{equation}
    \label{phij}
    \phi_{ij}: \mathrm{Tors}(G_{j})_{\mid U_{ij}} \la
\mathrm{Tors}(G_{i})_{\mid U_{ij}}
\end{equation}
together with natural
transformations  
\bee
\label{psijk}
\psi_{ijk}: \phi_{ij}\,\phi_{jk} \Longrightarrow \phi_{ik}
\end{equation}
defined above $U_{ijk}$, and which satisfies the tetrahedral coherence
condition
identifying to each other the pair of natural transformations 
{\renewcommand{\arraystretch}{1.5}
\begin{equation}
       \label{tetra}
     \begin{array}{ccc}
                 (\phi_{ij}\, \phi_{jk})\, \phi_{kl} &
                 \Longrightarrow \phi_{ik}\, 
         \phi_{kl} & \Longrightarrow \phi_{il}\\
          \phi_{ij}\,( \phi_{jk}\, \phi_{kl}) & \Longrightarrow \phi_{ij}\, 
         \phi_{jl} & \Longrightarrow \phi_{il} 
         \end{array}\end{equation}}
\noindent  defined above the quadruple intersection
$U_{ijkl}$. A global  object in $\pc$  is determined by a family of local 
objects $P_{i} \in \mathrm{Tors}(U_{i}, \,G_{i})$, together with 
$G_{i}$-torsor isomorphisms
\bee
\label{obP}
f_{ij}: \phi_{ij}(P_{j}) \la P_{i}\end{equation}
such that the diagram
 \bee\label{ob2P}\xymatrix{   
  \phi_{ij}\phi_{jk}(P_{k}) \ar[rr]^{\phi_{ij}(f_{jk})}
  \ar[d]_{\psi_{ijk}(P_{k})}  && \phi_{ij}(P_{j})\ar[d]^{f_{ij}}\\
  \phi_{ik}(P_{k}) \ar[rr]_{f_{ik}}&& P_{i} }\end{equation}
    commutes. Similarly, an arrow
   $a: (P_{i}, f_{ij}) \la ({P}_{ij}^{'}, {f}_{ij}^{'})$ 
    between a pair of such objects 
    in the glued 
    gerbe  is determined by a family of  isomorphism
    \[a_{i}: P_{i}\la P'_{i}\] in $\mathrm{Tors}(G_{i})$  such that 
    the square
     \bee
     \label{arP}
     \xymatrix{   
  \phi_{ij}(P_{j}) \ar[d]_{\phi_{ij}(a_{j})} \ar[r]^{f_{ij}}& P_{i} 
  \ar[d]^{a_{i}}\\
  \phi_{ij}(P'_{j}) \ar[r]_{f'_{ij}}& P'_{i}
  }\end{equation}
 of $G_{i}$-torsors on $U_{ij}$  commutes.
  
  \bigskip
  
 An object $u$ in the stack of  equivalences between $\pc$ 
  and $\pc'$ may be  described in similar terms by a family of local
equivalences $u_{i}: \to (G_{i}) \la \to (G'_{i})$
  together with  natural equivalences    $\Gamma_{ij}$ 
 above $U_{ij}$:
 \begin{equation}
\label{defui1} 
\xymatrix{
\to(G_{j}) \ar[r]^{u_{j}}\ar[d]_{ \phi_{ij}} & 
\to(G_{j}') \ar[d]^{\phi_{ij}'}
\ar@{}[d]_(.3){\,}="1"
\\
\to(G_{i}) \ar[r]_{u_{i}}
\ar[r]^(.3){\,}="2"& \to(G_{i}')\:.
\ar@{}"1";"2"^(.2){\,}="3"
\ar@{}"1";"2"^(.6){\,}="4"
\ar@{=>}"3";"4"_{\Ga_{ij}}
}\end{equation} 
such  that 
the composite 2-arrow in the pasting diagram
\begin{equation}
    \label{compat}
\xymatrix{
\to(G_{k}) \ar[rrr]^{u_{k}}="1" \ar[dr]
\ar[dd] \ar@{}[dd]^(.7){\,}="11"&&& \to(G_{k}') \ar[dl] \ar[dd]
\ar@{}[dd]_(.7){\,}="10"\\
&\to(G_{j}) \ar[r]^{u_{j}}="2"
 \ar@{}[r]_{\,}="9"\ar[dl] \ar@{}[dl]_(.2){\,}="12"& \to(G_{j}')\ar[dr]
 \ar@{}[dr]^(.4){\,}="7"&\\
\to(G_{i}) \ar[rrr]_{u_{i}}
\ar@{}[rrr]^{\,}="8"&&&\to(G_{i}')
\ar@{}"1";"2"^(.3){\,}="3"
\ar@{}"1";"2"^(.7){\,}="4"
\ar@{=>}"3";"4"^{\,\Ga_{jk}}
\ar@{}"9";"8"^(.2){\,}="5"
\ar@{}"9";"8"^(.7){\,}="6"
\ar@{=>}"5";"6"^{\,\Ga_{ij}}
\ar@{}"7";"10"^(.2){\,}="15"
\ar@{}"7";"10"^(.8){\,}="16"
\ar@{=>}"16";"15"_{({\psi'}_{ijk})^{-1}}
\ar@{}"11";"12"^(.3){\,}="13"
\ar@{}"11";"12"^(.7){\,}="14"
\ar@{=>}"14";"13"_{\psi_{ijk}}
}\end{equation}
above $U_{ijk}$ is equal to the 2-arrow $\Ga_{ik}$.

 \bigskip

Similarly, to a given morphism $t: u \Longrightarrow v$ 
\eqref{diagP} corresponds a family of morphisms $t_{i}$
\bee
\label{defti}\xymatrix{
\to(G_{i})\ar@/^1pc/[r]^{u_{i}}="3"
     \ar@/_1pc/[r]_{v_{i}}="4"
     & \to(G_{i}')
\ar@{}"3";"4"^(.2){\,}="7"
     \ar@{}"3";"4"^(.7){\,}="8"
\ar@{=>}"7";"8"_{t_{i}}}
\end{equation}
and the compatibility of the $t_{i}$'s with the 2-arrows $\Ga_{ij}$  
\eqref{defui1} may be exhibited as the commutativity of the following 
cyclindrical
diagram of 2-arrows
 \bee
 \label{cyl}
 \xymatrix{
    \to(G_{j})\ar@/^1pc/[r]^{u_{j}}="1"
    \ar@/_1pc/[r]_{v_{j}}="2"
    \ar[dd]_{\phi_{ij}}&
    \to(G'_{j}) \ar[dd]^{\phi'_{ij}}\\
    &&\\
    \to(G_{i})\ar@/^1pc/[r]^{u_{i}}="3"
     \ar@/_1pc/[r]_{v_{i}}="4"
     & \to(G_{i}')
     \ar@{}"1";"2"^(.2){\,}="5"
     \ar@{}"1";"2"^(.7){\,}="6"
     \ar@{}"3";"4"^(.2){\,}="7"
     \ar@{}"3";"4"^(.7){\,}="8"
     \ar@{=>}"5";"6"_{t_{j}}
      \ar@{=>}"7";"8"_{t_{i}}
      }\end{equation}
 in which the  the back and  the  front  2-arrow are the 
 2-arrows $\Ga_{ij}^u$ and $\Ga_{ij}^v$ respectively attached, as in    
 \eqref{defui1}, to $u$ and  to $v$.

 \subsection{}
\setcounter{equation}{0}
 The Morita theorem gives us an alternate description of the gluing data 
 for the 
 glued stack $\pc$.
 By proposition \ref{morpro}, the morphism 
(\ref{phij}) corresponds to the $(G_{i}, G_{j})$-bitorsor 
\bee
\label{morita}
P_{ij} := \mathrm{Isom}(\phi_{ij}(T_{j})_{\mid U_{ij}},\, (T_{i})_{\mid 
U_{ij}})
\end{equation}
above $U_{ij}$, where $T_{i}$ is the trivial $G_{i}$-torsor $T_{G_{i}}$
on $U_{i}$. 
Note that the canonical $(G_{j},\, G_{i})$ bitorsor isomorphism
\[ \phi_{ij}(T_{j}) \stackrel{\sim}{\la} \phi_{i}(x_{j}) = 
\mathrm{Isom}_{\pc}(x_{i}, x_{j})\:,\]
determines a $(G_{i}, \, G_{j})$-bitorsor isomorphism 
\bee
\label{defpij}
\mathrm{Isom}_{\pc}(x_{j}, x_{i})   \stackrel{\sim}{\la}  P_{ij}
\end{equation}
above $U_{ij}$.
The natural transformation $\psi_{ijk}$ is described by the 
$(G_{i},G_{k})$-bitorsor isomorphism 
 \bee
 \label{1bcoc}
 \begin{array}{ccc}
P_{ik} & \stackrel{\Psi_{ijk}}{\la}& P_{ij} \wedge^{G_{j}}P_{jk}
\end{array}
\end{equation}
above $U_{ijk}$ induced by the isomorphism
\[\psi_{ijk}(T_{k}):\phi_{ij}\phi_{jk}(T_{k}) \la \phi_{ik}(T_{k})\:. \]
The 
coherence condition \eqref{tetra} now asserts that the 
the diagram of bitorsors
       \[\xymatrix{
    P_{il} \ar[rr]^(.45){\Psi_{ijl}} \ar[d]_{\Psi_{ikl}} && P_{ij} 
    \wedge^{G_{j}} P_{jl}\ar[d]^{1\wedge \Psi_{jkl}}\\
   P_{ik}\wedge^{G_{k}} P_{kl} \ar[rr]_(.45){ \Psi_{ikl} \wedge 1}&&
   P_{ij} \wedge^{G_{j}} P_{jk} \wedge^{G_{k}} P_{kl}
    }\]
 above $U_{ijkl}$ is commutative. We will call such a family of 
 pairs $(P_{ij},\, \Psi_{ijk})$, or of triples
 $(G_{i},\,P_{ij},\, \Psi_{ijk})$, a bitorsor cocycle on $X$, even 
 though this terminology is generally used for
the equivalent data provided by the 
pair $(P_{ij},\, \Psi_{ijk}^{-1})$. 

\bigskip

By \eqref{ob2P}, an object in the glued gerbe $\pc$ may now  be described,
in terms of 
the bitorsors $P_{ij}$ attached to the gluing data, 
 by a family of torsors $P_{i} 
\in \mathrm{Tors}(U_{i},\,G_{i})$, together with the
isomorphisms
\bee
\label{fij1}
f_{ij}: P_{j} \wedge^{G_{j}} P_{ij}^{0} \la P_{i}\end{equation}
in $\mathrm{Tors}(U_{ij}, \:G_{i})$ describing the maps \eqref{obP}. 
The condition corresponding to the commutativity of diagram \eqref{ob2P} 
which the maps $f_{ij}$ must satisfy are most pleasantly described in terms
of the associated morphisms of $G_{j}$-torsors
 \bee
    \label{fij2}
    \tilde{f}_{ij}: P_{j} \la P_{i} \wedge^{G_{i}}P_{ij}\end{equation}
    defined by 
   \[ \tilde{f}_{ij}(p_{j})= (f_{ij}(p_{j},\,p_{ij}),\, p_{ij}),\]
   a definition which makes sense since the right-hand term is 
   independent of the choice of a section $p_{ij}$ of $P_{ij}$.
The condition states that each diagram
 \bee
\label{fij3}
\xymatrix{
P_{k} \ar[rr]^{\tilde{f}_{jk}} \ar[d]_{\tilde{f}_{ik}} && P_{j}\wedge P_{jk} 
\ar[d]^{\tilde{f}_{ij}\wedge P_{jk}}\\
P_{i}\wedge P_{ik} \ar[rr]_{P_{i}\wedge \Psi_{ijk}}
&& P_{i} \wedge P_{ij}\wedge 
P_{jk}}\end{equation}
  commutes. Similarly, a morphism $a:(P_{i}, f_{ij}) \la (P_{i}', 
  f_{ij}')$  between two such objects of $\pc$
   is given, by \eqref{arP}, by a family of $G_{i}$-torsor 
  isomorphisms $a_{i}: P_{i} \la P'_{i}$ for which each of the 
  following diagrams commute
  \bee
\label{fij4}
\xymatrix{
P_{j}\ar[r]^(.4){\tilde{f}_{ij}}\ar[d]_{a_{j}} &
 P_{i}\wedge P_{ij}\ar[d]^{a_{i}\wedge  P_{ij}}
\\
P_{j}' \ar[r]_(.4){\tilde{f}_{ij}'}& P_{i}'\wedge P_{ij}\:.
}\end{equation}

 \begin{remark}{\rm 
   When $\pc$ is an abelian $G_{m}$-gerbe in the sense of 
   \cite{lb:2-gerbes},
the bitorsor structure on the 
corresponding $G_{m}$-torsors is the obvious one which doesn't 
distinguish between the left and the right action, so that the 
bitorsor cocycle structure is defined by pairs $(P_{ij},\, 
\Psi_{ijk}^{-1})$ with $P_{ij}$ simply a $G_{m}$-torsor  on $U_{ij}$, 
or its associated line bundle. 
What is referred to in \cite{hitch} as a gerbe is now seen to be 
one possible description of 
an abelian $G_{m}$-gerbe $\pc$ on $X$, together with a chosen family of local 
trivializations $(x_{i})_{i \in I}$ above the open sets $U_{i}$.
}
\end{remark}
     
The  equivalences $u_{i}$ attached to an object $u$ in $\epc$ are described
by the 
$(G_{i}', G_{i})$-bitorsors 
\bee
\label{defgai}\Gamma_{i} := \mathrm{Isom}(u_{i}(T_{i}), 
T_{i}') \stackrel{\sim}{\la} 
\mathrm{Isom}_{\pc'}(u(x_{i}),\, x'_{i})\end{equation}
defined on the open sets $U_{i}$. The natural transformation 
$\Ga_{ij}$ \eqref{defui1} corresponds to a $(G_{i}',\, G_{j})$-bitorsor 
isomorphism
\bee
\label{ardij0}
\xymatrix{
\Gamma_{i}\wedge P_{ij} \ar[r]^{\tilde{g}_{ij}}& P_{ij}' \wedge 
\Gamma_{j}}\end{equation}
on $U_{ij}$ and the latter can also be written  as a $(G_{i}',\, 
G_{i})$-isomorphism\footnote{There should be no confusion between the 
present $g_{ij}$, which is an arrow, and the $g_{ij}$ occuring in 
\eqref{1coc}, which is a $G$-valued 1-cocycle.}
\bee
\label{ardij}
\xymatrix{
\Gamma_{i} \ar[r]^(.3){g_{ij}}& P_{ij}' \wedge 
\Gamma_{j}\wedge P_{ij}^{0}}\:.\end{equation}
In such a  form it is the  analogue of \eqref{0coch3}.
In terms of the morphisms $\tilde{g}_{ij}$, the compatibility 
condition \eqref{compat} is now expressed by the commutativity of the 
diagram
  \[\xymatrix{
 \Ga_{i}\wedge P_{ik}\ar[dd]_{\tilde{g}_{ik}}
  \ar[rr]^<<<<<<<<<<{\Gamma_{i} \wedge 
\Psi_{ijk}}&&
\Ga_{i}\wedge P_{ij}\wedge P_{jk} \ar[d]^{\tilde{g}_{ij}\wedge 
P_{jk}}\\
&&P_{ij}'\wedge\Ga_{j}\wedge P_{jk}
\ar[d]^{P_{ij}' \wedge \tilde{g}_{jk}}\\ 
P_{ik}'\wedge\Ga_{k}\ar[rr]_<<<<<<<<<<{\Psi_{ijk}'\wedge \Ga_{k}}&&P_{ij}' 
\wedge P_{jk}' \wedge \Ga_{k}
}\]
This property is better stated in terms 
of the corresponding morphisms $g_{ij}$ \eqref{ardij}. Neglecting the 
canonical isomorphism 
\[(P_{ij} \wedge P_{jk})^{0} \simeq  P_{jk}^{0} \wedge P_{ij}^{0}\:,\]
it then asserts that the 
diagram
\bee
\label{compdij}
\xymatrix{
\Ga_{i} \ar[rrrr]^{g_{ij}}\ar[d]_{g_{ik}}&&&&
P_{ij}' \wedge \Ga_{j}\wedge P_{ij}^{0}
\ar[d]^{P_{ij}'\wedge g_{jk} \wedge P_{ij}^{0}}
\\
P_{ik}' \wedge \Ga_{k} \wedge P_{ik}^{0} 
\ar[rrrr]_<<<<<<<<<<<<<<<<<<<{{\Psi'}_{ijk} \wedge 1_{\Ga_{k}}
\wedge (\Psi^{-1}_{ijk})^{0}}
& &&& 
(P_{ij}' \wedge P_{jk}') \wedge \Ga_{k} \wedge (P_{ij} \wedge 
P_{jk})^{0}
}\end{equation}
commutes.

\bigskip

We now set $\Ga_{i}^{u} := \Ga_{i}$ and denote by $\Ga^{v}_{i}$ the
corresponding family of $(G_{i}', G_{i})$-bitorsors 
\[\Ga^{v}_{i}:= \mathrm{Isom}(v_{i}(T_{i}),T_{i}') \simeq 
\mathrm{Isom}(v_{i}(x_{i}),\, x'_{i})\]
 attached to $v: \pc 
\la \pc'$. The morphisms $t_{i}$ \eqref{defti} then correspond to the 
bitorsor isomorphisms 
\[
\label{def:thetai1}
\begin{array}{ccc}
    \Ga_{i}^{v} &\stackrel{ \Theta_{i}}{\la} &\Ga_{i}^{u}\\
    \ga^{v}_{i} & \mapsto & \ga^{v}_{i} \, t(x_{i})\,.
    \end{array}\]
%so that, in the notation of \eqref{defthet},
%\bee
%\label{thetatheta}
%\Theta_{i}(\ga^{v}_{i}) = \theta_{i}\ga^{u}_{i}
%\end{equation}
The commutativity of  \eqref{cyl}, when expressed in terms of the arrows 
$g^{u}_{ij}$ and $g^{v}_{ij}$ associated as in  \eqref{ardij} 
to the pair of morphisms $u$ and $v$, asserts that the diagram 
     \bee
     \label{comthet}
     \xymatrix{
\Ga_{i}^{v} \ar[r]^(.3){g^v_{ij}} \ar[d]_{\Theta_{i}}
& P_{ij}' \wedge \Gamma_{j}^{v} \wedge 
P_{ij}^{0} \ar[d]^{ P_{ij}' \wedge \Theta_{j} \wedge P_{ij}^{\mathrm{0}}}\\
\Ga_{i}^{u}
\ar[r]_(.3){g^u_{ij}}& P_{ij}' \wedge \Gamma_{j}^{u} \wedge 
P_{ij}^{0} }\end{equation}
commutes.

\bigskip

This discussion may be summarized as follows.
\begin{proposition}
    \label{pcad}
   Consider a pair of  locally trivialized gerbes $(\pc,\, x_{i})$ 
  (\re\:
    $(\pc',\, x_{i}')$) with associated locally defined groups 
    $G_{i}$ (\re $ \:G_{i}'$),
 bitorsor cocycles $P_{ij}$ (\re $P_{ij}'$) and bitorsor isomorphisms 
$\:\Psi_{ijk}$  (\re $\: \Psi_{ijk}'$) \eqref{1bcoc}. The stack $\epc$ is 
 anti-equivalent to the stack on $X$ obtained by gluing the bitorsor 
 stacks 
\bee
\label{defdi}
\mathcal{D}_{i} := \mathrm{Bitors}_{U_{i}}(G_{i}',,\, 
G_{i})\end{equation}
 by the $(P_{ij}', P_{ij}^{0})$-gluing
 data  $(r_{ij},\, s_{ijk})$ 
 defined by 
  \bee
   \label{defrij}
   \xymatrix@R=3pt{
(\ddc_{j})_{\mid U_{ij}} \ar[r]^{r_{ij}}& (\ddc_{i})_{\mid U_{ij}}\\
\Gamma_{j}\ar@{|->}[r]& P_{ij}'\wedge \Ga_{j}\wedge P_{ij}^{0}
}\end{equation}
and with  $s_{ijk}$ the natural transformation
\[ \Psi'_{ijk}\wedge (-) \wedge ({\Psi^{0}}_{ijk})^{-1}\]
acting as in the lower line of  \eqref{compdij}.
\end{proposition}
    \begin{flushright}
        $\Box$
        \end{flushright}
    
        The most interesting case is that of the gauge stack $\pac$, 
        in 
        which $\pc'$ is equal to $\pc$. The contracted product of bitorsors 
        then defines on  
        $\ddc_{i}= \mathrm{Bitors}(G_{i})$  a 
        $gr$-stack structure, 
        and the gluing data $r_{ij}$ is defined  by the conjugation action 
        \bee
        \label{def:twpac}
        \Ga_{j} \mapsto P_{ij} \wedge \Ga_{j} \wedge (P_{ij})^{0}\:.
        \end{equation}
        Since this adjoint action preserves the 
        monoidal structure, as does any conjugation, the glued stack $\ddc$ 
        on $X$ which it defines 
 has an induced $gr$-stack structure, for which the anti-equivalence 
        $\ddc \simeq \pac$ is a monoidal functor. In a more compact 
        form and  in analogy with \eqref{padj}, this can be expressed as
        as
        \bee 
        \label{gadj}\pac \simeq 
        {}^{(P_{ij},\Psi_{ijk})\,}\!(\mathrm{Bitors}(G_{i}))_{i\in 
        I}\:.\end{equation}
 When $\pc$ is a $G$-gerbe,  the groupoid-type data 
$(\mathrm{Bitors}(G_{i}))_{i\in I}$  reduces to the 
$gr$-stack $\mathcal{G}:= \mathrm{Bitors}(G)$ of $(G,\,G)$-bitorsors on $X$, 
so that  the anti-equivalence
\eqref{gadj} simply  becomes 
\bee
\label{gadj0}
\pac \simeq     {}^{(P_{ij},\Psi_{ijk})\,}\!\mathcal{G}\:.\end{equation}

\bigskip

Returning to the general gerbe (rather than $G$-gerbe) case, let us
 observe that the right $\pac$ and left   ${\pc'}^{\mathrm{ad}}$-torsor
 structures 
on the stack $\epc$ can both  be recovered from their local description. The 
right $\pac$-torsor structure, for example, is given above $U_{i}$ by 
the contracted product pairing
\[
 \xymatrix@R=3pt{\mathrm{Bitors}(G'_{i},\, G_{i}) \times 
\mathrm{Bitors}(G_{i}) \ar[r]&
\mathrm{Bitors}(G'_{i},\, G_{i})\\
(\Ga_{i}^{1},\, \Ga_{i}^{2}) \ar@{|->}[r] & \Ga_{i}^{1} \,\Ga_{i}^{2}}\]
where we have dropped the $\wedge$ symbol from the contracted 
product, and this pairing is compatible with the gluing morphisms 
$r_{ij}$ \eqref{defrij} {\it via} the coherent family of isomorphisms
\[ P_{ij}'\,(\Ga_{j}^{1}\,\Ga_{j}^{2})\, P_{ij}^{0} \simeq 
(P_{ij}'\,\Ga_{j}^{1} \,P_{ij}^{0} )\:(P_{ij}\,\Ga_{j}^{2} \,
P_{ij}^{0})\:.\]
More generally, the  horizontal pairing functor 
\eqref{comp1} is now anti-equivalent to the one constructed by 
 twisting appropriately the contracted product pairings
\[ \mathrm{Bitors}(G''_{i},\, G'_{i}) \times
\mathrm{Bitors}(G'_{i},\, G_{i}) \la
\mathrm{Bitors}(G''_{i},\, G_{i})\:.\]

 \begin{remark}{\rm
It is tempting to express the right-hand term of \eqref{gadj0} even  more 
compactly, by analogy with \eqref{padj}. For this we set 
 \[ \gc := \mathcal{E}q (\mathrm{Tors}(G), \, \mathrm{Tors}(G))\:.\]
As in \cite{fest} 
proposition 7.3, there exists a canonical equivalence between the 
2-stack of $G$-gerbes  and the 2-stack of right $\gc$-torsors  on 
$X$, which to any $G$-gerbe $\pc$ associates the $\gc$-torsor
\[ \tilde{\pc} := \mathcal{E}q( \mathrm{Tors}(G),\,\pc)\:.\]
This provides an equivalence of the gauge stack $\pac$ of $\pc$ with  
the $gr$-stack 
\[\tilde{\pc}^{\mathrm{ad}}:=
\mathcal{E}q_{\gc}(\tilde{\pc},\, \tilde{\pc} ) \]
of $\gc$-equivariant self-equivalences of 
$\tilde{\pc}$. The anti-equivalence \eqref{gadj0} can now be stated 
as  
\[\pac \simeq  \tilde{\pc}^{\mathrm{ad}} 
\simeq \tilde{\pc} \wedge^{\gc} \gc\,, \]
where $\gc$ acts on itself by conjugation, which 
expresses the   $gr$-stack $\pac$ as  a twisted
inner form 
of the $gr$-stack   $\mathcal{G}$. 
A morphism 
   $\tilde{\pc}  
\wedge^{\mathcal{G}}\mathcal{G} \la \tilde{\pc}^{\mathrm{ad}}$
can be constructed directly,
by factoring the morphism of stacks
\[\xymatrix@R=3pt{ \tilde{\pc} \times \mathcal{G} \ar[r] & 
\tilde{\pc}^{\mathrm{ad}}  \\
(p,g) \ar@{|->}[r] &   (p \mapsto pg) }\]
 through the contracted product, as defined in  \cite{fest} 
\S 6.7. More generally, the stack $\epc$,
(with $\pc$ a $G$-gerbe and $\pc'$ a 
$G'$-gerbe on $X$) is described by an equivalence of stacks
\bee
\label{bitortw}
\epc \simeq \tilde{\pc}'\, \wedge^{\mathcal{G}'} \,\mathrm{Bitors}(G', 
\, G)\, 
\wedge^{\mathcal{G}}\, \tilde{\pc}^{0}\end{equation} which expresses 
the stack  $\epc$ as 
a doubly twisted form  of the stack $\mathrm{Bitors}(G', \, G) $.
}\end{remark}

%% file: dgg-ch4.tex
\section{Connections, curving data and the higher Bianchi identity}
\label{sec:curv-Bian}
\subsection{}
\label{subsec:curv-Bian1}
\setcounter{equation}{0}%
Let   $\pc$ be a gerbe on an $S$-scheme  $X$.
\begin{definition}
     \label{defin: equicon}{\it i}) A
     connection on the gerbe $\pc$ relative to $X/S$ is an equivalence
     of gerbes $(\epsilon,\, \epsilon^{-1}, \,
\Theta)$  on $\Delta^{1}_{X/S}$ \eqref{def:Theta}:
     \bee
     \label{gercon}
\begin{array}{cccc}   \epsilon: p_{1}^{\ast}\pc \la
p_{0}^{\ast}\pc && \epsilon^{-1}:p_{0}^{\ast}\pc \la  p_{1}^{\ast}\pc\,,
\end{array}
\end{equation}
\begin{equation}
\label{gercon2}
\Theta: \epsilon^{-1}\, \epsilon  \Longrightarrow 1_{\pra \pc}\,,
\end{equation}
together with  a  natural equivalence  $\eta$:
  \bee
  \label{2ar:eta}
  \xymatrix{
  \pc  \ar[rr]^{1_{\pc}}="1"\ar[d]_(.45){\wr} && \pc\ar[d]^(.45){\wr}\\
  \Delta^{\ast}(p_{1}^{\ast}\pc)
  \ar[rr]^{\Delta^{\ast}\epsilon}="2" && \Delta^{\ast}(p_{0}^{\ast}\pc)
\,,
  \ar@{}"1";"2"^(.30){\,}="3"
    \ar@{}"1";"2"^(.70){\,}="4"
\ar@{=>}"3";"4"^{\,\eta}\,}
\end{equation}
where the  vertical arrows in \eqref{2ar:eta} are the canonical equivalences
  induced, in
the notation of \S \ref{subpar61}, by the simplicial identities
$p_{1}\,\Delta = p_{0}\,\Delta = 1_{X}$ on $\De^{\ast}_{X/S}$.

\medskip

\qquad \qquad  \quad {\it ii}) A morphism
  between two such  connections
$(\epsilon,\, \eta)$ and $(\epsilon',\, \eta')$ on $\pc$ is a natural
equivalence
\bee
\label{def:zeta}
\zeta: \epsilon \Longrightarrow \epsilon'\:,\end{equation}
on $X$ such that the composite 2-arrow
  \bee
  \label{2ar:eta1}
  \xymatrix{
  \pc  \ar[rr]^{1_{\pc}}="1"\ar[d]_(.45){\wr} && \pc\ar[d]^(.45){\wr}\\
  \Delta^{\ast}(p_{1}^{\ast}\pc)
  \ar[rr]^{\Delta^{\ast}\epsilon}="2"
  \ar@{}[rr]_{\,}="5"
  \ar@/_2pc/[rr]_{\Delta^{\ast}\epsilon'}="6"&&
  \Delta^{\ast}(p_{0}^{\ast}\pc)
\:,
\ar@{}"1";"2"^(.30){\,}="3"
    \ar@{}"1";"2"^(.70){\,}="4"
\ar@{=>}"3";"4"^{\,\eta}
\ar@{}"5";"6"^(.12){\,}="7"
\ar@{}"5";"6"^(.62){\,}="8"
\ar@{=>}"7";"8"^{\,\Delta^{\ast}\zeta}}
\end{equation}
coincides with $\eta'$.

\end{definition}

With these definitions the functor given by forgetting $\eta$
is a faithful functor from the stack of
connections on $\pc$ to the stack  $\mathcal{E}q(p_{1}^{\ast}\pc,\,
p_{0}^{\ast}\pc)$.
We  denote by $\tilde{\Theta}$
  the  adjoint equivalence
  \eqref{def:Theta-tilde} of $\Theta$ \eqref{gercon2},
  and by $\tilde{\eta}: 1_{\pc} \Longrightarrow
  \Delta^{\ast}\epsilon^{-1}$ the unique 2-arrow such that  $\eta$ and
  $\tilde{\eta}$
are compatible with $\Delta^{\ast}\Theta$. We will often simply refer
to a gerbe with a connection $(\epsilon,\, \epsilon^{-1},\,\Theta)$
simply as $(\pc,\, \epsilon)$, without making explicit the
  quasi-inverse $\epsilon^{-1}$ of $\epsilon$, or the 2-arrow
$\Theta$ \eqref{gercon2}. We now introduce additional curving data.
For this, we will denote by $\epsilon_{ij}$ the pullback
$p^{\ast}_{ij}   \epsilon$ of $\epsilon$ by the projection $p_{ij}$
from $\deb$ to $\dea$.

\begin{definition}
{\it i)}  Let $(\pc, \, \epsilon)$  be a gerbe with connection
  on $X$. An arrow 
\begin{equation}
     \label{defK0}
  \widetilde{K}:\hspace{.3cm}  \ast\,  \Longrightarrow
  \epsilon_{12}
\,  \epsilon_{01} \, \epsilon^{-1}_{02}
\end{equation}
in  the stack $ \mathrm{Lie}(\pac,\, \Om^2_{X/S})$ 
(see definition \ref{def:b2}) with target
$ \prab \epsilon\,  \prbc \epsilon \,\prac \epsilon^{-1}$
is called a curving for $(\pc, \, \epsilon)$.

\medskip

\qquad \qquad \quad  \quad   {\it  ii)} The source
$\kappa \in  \mathrm{Lie}(\pac,\, \Om^2_{X/S})$
of the arrow $K$ is called the fake curvature associated to the
  pair $(\epsilon,\,\widetilde{K})$.
\end{definition}

\bigskip

The  curving $\widetilde{K}$ with source $\kappa$
\begin{equation}
\label{defK1}
\xymatrix{
\pra \pc  \ar@/^1pc/[rr]^{\kappa}
\ar@/_1pc/[rr]_{ \,\epsilon_{01} \,\,
  \epsilon_{12}\,\, \epsilon^{-1}_{02}}
\ar@{}@/^1pc/[rr]_{\,}="1"
\ar@{}@/_1pc/[rr]^{\,}="2"
  &&\pra \pc
\ar@<-.5ex>@{=>}"1";"2"^(.3){\,\widetilde{K}}
}
\end{equation}
determines by whiskering  with $\epsilon_{02}$
  a 2-arrow
\begin{equation}
\label{defkK}
\xymatrix{
&\pra \pc\ar@/^/[dr]^{\kappa}\ar@{}[d]^(.25){\,}="1"
\ar@{}[d]^(.75){\,}="2"
  &\\
\prc \pc \ar@/^/[ur]^{\epsilon_{02}}  \ar[r]_{\epsilon_{12}}& \prb \pc
\ar[r]_{\epsilon_{01}}&\pra \pc\:.
\ar@<-1ex>@{=>}"1";"2"^(.4){\,K}
}\end{equation}
  and the 2-arrow $\widetilde{K}$ may be recovered from $K$ by
  whiskering by  $\epsilon_{02}^{-1}$, as explained in
\S  \ref{adif}. The giving of a curving $\eqref{defK1}$ with source $\kappa$
is thus equivalent to that of a 2-arrow \eqref{defkK}. We will henceforth
refer either to $\widetilde{K}$ or to $K$ as a curving, and to
$(\epsilon,\, \widetilde{K})$ or $(\epsilon,\, K)$ as a curving pair. It
is often convenient to include in the description of $K$ the name $\kappa$
  of  its
source 1-arrow. We will call $(\epsilon,\, \widetilde{K},\,
\kappa)$ or $(\epsilon,\, K,\,
\kappa)$ a connection triple on the gerbe $\pc$ even though this is
  somewhat redundant, since strictly
speaking the giving of the 2-arrow $\widetilde{K}$ \eqref{defK0}
determines its source  1-arrow $\kappa$.

\bigskip

The 1-arrow $\kappa$ is an object
in the fiber category  of the gauge stack $\pra\pac$ of
$\pc$ above  $\deb$,
with given quasi-inverse
$\kappa^{-1} $ and associated  2-arrow
\bee
\label{defEta}
\kappa^{-1}\kappa \stackrel{H}{\Longrightarrow}1.
\end{equation}
  The restriction  $t^{\ast}\kappa$ of $\kappa$
to the degenerate  subsimplex $s\deb$  of $\deb$ is  endowed with a 2-arrow
   \bee
\label{katriv}
    \xymatrix{
    t^{\ast}\kappa \,\,\ar@{=>}[r]^{\tau} & 1}
     \end{equation}
    so that $\kappa$ may actually be considered
    as an object in the fiber on $X$ of the
    stack $\mathrm{Lie}(\pac,\, \Om^{2}_{X/S})$ of $\pac$-valued
    2-forms on $X$.
The 2-arrow $\tau$ is uniquely determined by its restrictions
along each  the two degeneracy maps
\begin{equation}
\label{degmap}
s_0, s_1:\dea \la \deb \,.
\end{equation}
\bigskip
and these restrictions
  coincide above
    $X  \hookrightarrow \dea$.

    \bigskip

\begin{remark}\label{process}
{\rm
   \hspace{2cm}  {\it i}) \hspace{.25cm}
    The  diagram
  \[ \xymatrix{
&&\pra \pc \\
\prc \pc \ar@/^1pc/[urr]^{\epsilon_{02}}
\ar[r]_{\epsilon_{12}}&\prb \pc \ar[r]_{\epsilon_{01}}
&\pra \pc
}\]
  is analogous to a
horn,  for which the 2-arrow $K$ would be a Kan
filler

\bigskip

\bee
      \label{kanfil} \xymatrix{
&&\pra \pc \ar@{}[d]_{\,}="3"
\ar@{-->}[d]^{\kappa}
\ar@{}[dl]^(.3){\,}="7"
\ar@{}[dl]^(.7){\,}="8"
\\
\prc \pc \ar@/^1pc/[urr]^{\prac \epsilon}
\ar@{}@/^1pc/[urr]_{\,}="1"
\ar[r]_{\prbc \epsilon}
\ar@{}[r]^(.99){\,}="2"
& \prb \pc \ar[r]_{\prab \epsilon}
&\pra \pc
\ar@{}"1";"2"^{\,}="4"
\ar@{}"4";"3"^(.01){\,}="5"
\ar@{}"4";"3"^(.7){\,}="6"
\ar@{=>}"7";"8"_{K}
\:.
}\end{equation}
In view of the chosen  orientation for
the  1-arrows, such a 2-arrow
does not have an immediate interpretation, and its rigorous
significance is given by  diagram \eqref{defK1}.
To borrow a metaphor from \cite{B-D} \S1, we may think of $K$ as a
`process' which allows one to subtract  the upper arrow of \eqref{kanfil}
from the lower composite one, and of $\kappa$ as the result of this 
subtraction  in
the 2-stack of gerbes on $\Delta^{2}_{X/S}$.

\bigskip

\hspace {2cm}{\it ii})  Consider the lien functor, which associates
to a gerbe $\pc$ on
$X$ its lien, $\mathrm{lien}(\pc)$.
By \cite{Gir} IV 2.1.5.2, the existence of a natural
transformation  $t$ \eqref{diagP} between a pair of arrows $u$ and $v$ implies
that the induced morphisms $\mathrm{lien}(u)$ and $\mathrm{lien}(v)$
between $\mathrm{lien}(\pc)$
and $\mathrm{lien}(\pc')$ are equal. It follows that the lien functor
associates to a connection $\epsilon$ \eqref{gercon} on $\pc$ a
connnection
\bee
\label{lien1}
\prb \, \li (\pc) \stackrel{\lambda}{\la} \pra \, \li
(\pc)\end{equation}
on  $\mathrm{lien}(\pc)$, and
transforms  diagram \eqref{defkK} into a commutative diagram of liens
\[ \xymatrix{\prc \,\li (\pc) \ar[r]^{\prac \lambda}
\ar[d]_{\prbc \lambda} & \pra\, \li (\pc)\ar[d]^{\li (\kappa)} \\
\prb \, \li (\pc) \ar[r]_{\prab \lambda}& \pra \, \li (\pc)\:.}\]
The fake curvature $\kappa$ can thus be interpreted as a
lifting to the 2-stack of gerbes of the curvature of the lien
connection $\lambda$.
}\end{remark}

     \vspace{.75cm}

The connection $(\epsilon, \epsilon^{-1} \Theta)$ \eqref{gercon}
   determines a morphism of $gr$-stacks \eqref{def:ad}
\bee
\label{defmuad}
\begin{array}{cccc}
    \mu := \epsilon^{\mathrm{ad}}:& \prb \pc^\mathrm{ad} &\la &\pra
\pc^\mathrm{ad}\\
&u& \mapsto & \epsilon \, u \, \epsilon^{-1}
\end{array}\end{equation}
which is a (group structure preserving)
connection on the  gauge stack $\pac$. This is a monoidal functor
whose pullback to $X$ is provided with a monoidal 2-arrow \cite{maclane}
XI \S 2.
$\eta_{\mu}:
1_{\pac} \Longrightarrow \Delta^{\ast}(\mu)$.

\bigskip

Similarly, the ``inner conjugation'' $\kappa^{\mathrm{ad}}:u \mapsto
\kappa\, u\, \kappa^{-1}$ is an object of the
stack $\mathcal{E}q(\pra \pac, \,\pra \pac)$. We will denote it by
$i_{\kappa}$, though it is not strictly speaking an
inner conjugation by $\kappa$, since $\kappa^{-1}$ is not the strict
inverse of $\kappa$.  By the discussion in \S \ref{adif}, diagram
\eqref{defkK}, induces by adification a diagram of $gr$-stacks \bee
\label{defkKad}
\xymatrix{
&\pra \pac \ar@/^/[dr]^{i_{\kappa}}
\ar@{}[d]^(.25){\,}="1"
\ar@{}[d]^(.75){\,}="2"
\ar@<-1ex>@{=>}"1";"2"^{\kc} &\\
\prc \pac \ar@/^/[ur]^{\mu_{02}}  \ar[r]_{\mu_{12}}& \prb
\pac
\ar[r]_{\mu_{01}}&\pra \pac
}\end{equation}
above $\deb$ , where  $\mathcal{K}$ is the 2-arrow $K^{\mathrm{ad}}$
induced by $K$.
In particular, the natural
transformation $\kc$  determines  an arrow
\bee
\label{curly0}
\kc (x):\, i_{\kappa}\,\,\,\mu_{02}(x) \la \mu_{01}\,
\mu_{12}(x)\end{equation}
in $\pra\, \pac$ for every object $x \in \mathrm{ob}(\prc\, \pac)$.
We will denote  by $\mathrm{Conn}(\pac)$ and $\mathrm{Curv}(\pac)$
the stacks $\mathcal{E}q(\prb \pac,\, \pra \pac)$ and
 $\mathcal{E}q(\prc \pac, \, \pra \pac)$ of which $\mu$ and $\kc$ 
are respectively an object and an arrow.

\bigskip

With the same notations as in  \eqref{def:kappa1},
we will now display diagram \eqref{defkK} as the square

      \begin{equation}
         \label{def:kappa2}
         \xymatrix@=13pt{
    \prc \pc \ar[rr]^{\epsilon_{12}} \ar[dd]_{\epsilon_{02}}
     && \prb \pc \ar[dd]^{\epsilon_{01}}
     \ar@{}[dd]_(.55){\, }="1"\\&&\\\pra \pc \ar[rr]_{\kappa}
     \ar@{}[rr]^(.55){\,}="2"&& \pra  \pc \:.
      \ar@{}"1";"2"^(.2){\,}="3"
      \ar@{}"1";"2"^(.8){\,}="4"
      \ar@{=>}"4";"3"^{K}
     }\end{equation}

\bigskip

Just as the square \eqref{def:kappa1}
has now been replaced, when passing from torsors to gerbes,
by the 2-arrow $K$ \eqref{def:kappa2}, we may
now enrich the commutative diagram \eqref{bianchidiag}
to the following diagram of 2-arrows in the 2-category  of gerbes 
above $\dec$:

  \bigskip

  \begin{equation}
          \label{cube}
\includegraphics{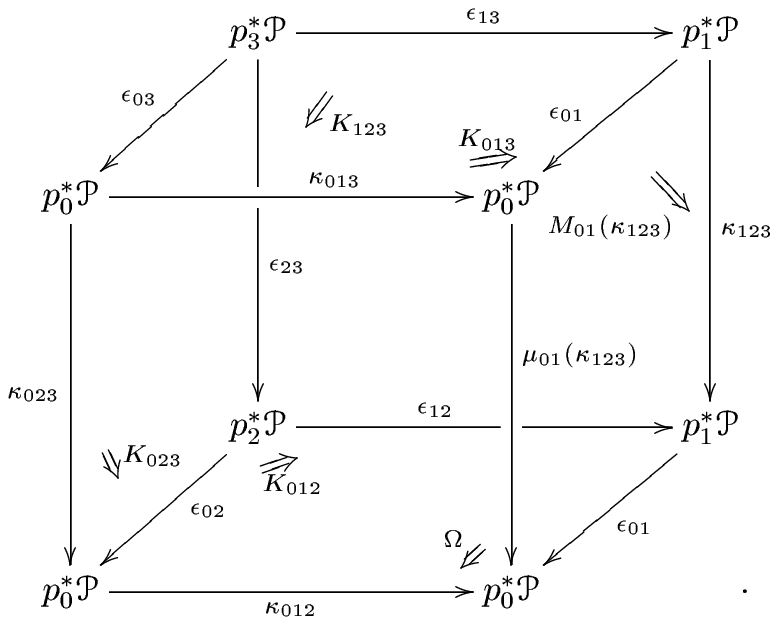}
\end{equation}

\bigskip

Each of the four displayed 2-arrows $K_{ijk}$ is the  pullback of the
2-arrow $K$ \eqref{defkK}
by the corresponding  projection from
$\Delta^{3}_{X/S}$ to $\deb$, and  the right hand vertical face is
$M_{01}(\kappa_{123})$. Here $M_{01}$ is the pullback from $\dea$
to $\dec$ of the conjugation 2-arrow $M$ associated, as in
\eqref{msquare},
to $\epsilon_{01}$.
  Finally
we can now define the front  2-arrow
\bee
\label{defomeg}
\xymatrix{\pra \pc \ar[r]^{\kappa_{013}} \ar[d]_{\kappa_{023}}
&\pra \pc \ar[d]^{\mu_{01}(\kappa_{123})}\ar@{}[d]_(.3){\,}="1"\\
\pra \pc  \ar[r]_{\kappa_{012}}\ar[r]^(.3){\,}="2" &\pra \pc
\ar@{}"1";"2"^(.2){\,}="3"
\ar@{}"1";"2"^(.7){\,}="4"
\ar@{=>}"3";"4"_{\Om}
}
\end{equation}
as the unique 2-arrow in the fiber category of $\pac$
above $\dec$
for which the cubic diagram of 2-arrows \eqref{cube} commutes. It is
explicitly  constructed  by pasting
together  the five other faces (or in some instances
their inverses), as in:
\begin{equation}
\label{}
\xymatrix{
\prd \pc  \ar@/^1pc/[rrrr]^{\mu_{01}(\kappa_{123}) \, \kappa_{013}\,
  \epsilon_{03}}
\ar@/_1pc/[rrrr]_{\kappa_{012}\, \kappa_{023} \,  \epsilon_{03} }
\ar@{}@/^1pc/[rrrr]_{\,}="1"
\ar@{}@/_1pc/[rrrr]^{\,}="2" &&&&\pra \pc\,
\ar@{=>}"1";"2"
}
\end{equation}
and whiskering on the left with $\epsilon_{03}^{-1}$.

\bigskip

The restriction
  $t^{\ast} \Om$  of $\Om$ to the
  degenerate subsimplex $t: \,s\dec\hookrightarrow \dec$
  is identified {\it via} \eqref{katriv}
   with the 2-arrow $1_{\sigma}$, where $\sigma$ is the arrow in the
   fiber of $\pra \pc$ above $s\dec$ whose restriction along each of
   the degeneracies $s_{j}$ is $\kappa$. By whiskering, $\Omega$
corresponds to a 2-arrow
  \[1 \Longrightarrow \kappa_{012}\, \kappa_{023}
  \, (\kappa_{013})^{-1}   \, (\prab\mu)(\kappa_{123})^{-1} \] above
  $\dec$  whose restriction to $s\dec$
may be identified with
the identity 2-arrow  on $1_{\pra(\pc)}\,$,
   and which  therefore  defines  an arrow in
  fiber  of the stack  $\mathrm{Lie}(\pac , \Om^3_{X/S})$ of
   $\mathrm{Lie}(\pac)$-valued 3-forms on $X$.
   We denote it  by
   $\widehat{\Omega}$,  and will call it  the
   3-curvature of the gerbe $\pc$ with curving pair  $(\epsilon, \, K)$.
In order not to overburden the notation, we at times simply write it as
$\Om$, since the giving of $\Om$ and $\widehat{\Omega}$ are equivalent.

   \bigskip

\bigskip

   By the adification process described in \S \ref{adif},
   diagram \eqref{cube} induces a diagram of {\it
   gr}-stacks

   \begin{equation}
       \label{cubead}
     \includegraphics{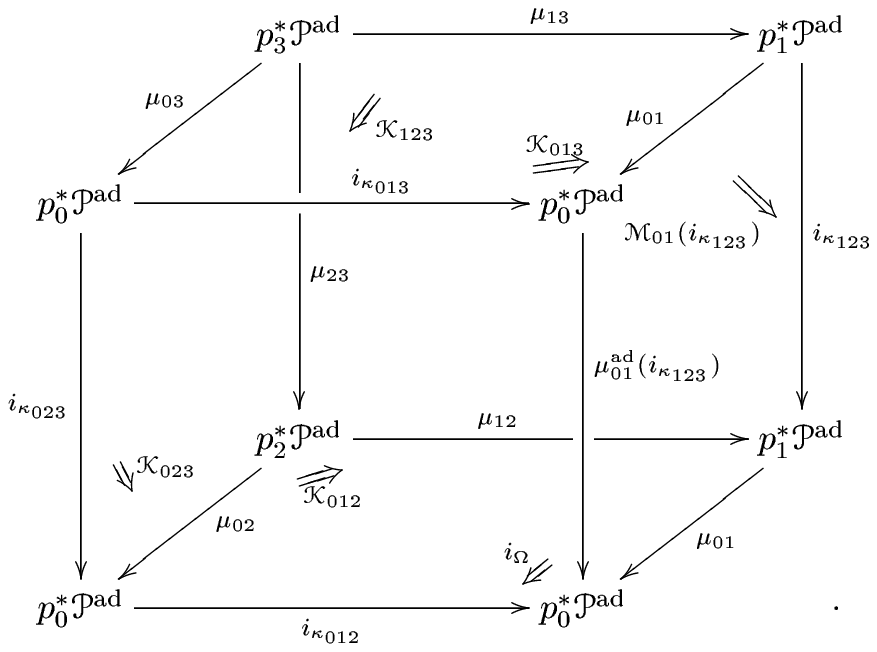}
       \end{equation}
      This  describes the relation between
the various pullbacks $\kc_{ijk}$ of  $\kc$ to $\dec$
and   the adifications
  $i_{\Omega}$ and
       $\mathcal{M}$   of the 2-arrows $\Omega$ and  $M$.

\bigskip

    In addition to this relation,
  the 3-curvature form $\Omega$
   satisfies   a  higher Bianchi identity, which we will now discuss.
   This  expresses a compatibility between the   pullbacks  $\Om_{ijkl}$
  of $\Om$ by the various
projections $p_{ijkl}:\Delta^{4}_{X/S}\la \Delta^{3}_{X/S}$.

  \bigskip

  \begin{theorem}
       \label{th:gl}
       Let $\pc$ be a gerbe on $X$.
A  curving pair $(\epsilon,\, K)$
  determines a  quadruple
       $(\mu, \, \kc,\, \kappa,\, \Omega)$, with $\mu$  \eqref{defmuad}
       a (group-structure
       preserving) connection
    on the gauge $gr$-stack $\pac$ of $\pc$, $\kc$
    an arrow \eqref{defkKad} in the fiber of
    $\mathcal{E}q(\prc \pac,\, \pra
    \pac)$ above $\deb$, with a trivialization on the degenerate
    subsimplex $s\deb$ of $\deb$,  $\kappa$
    an object  \eqref{defkK} in the fiber on $X$ of the stack
    $\mathrm{Lie} (\pac, \, \Om^{2}_{X/S})$,
    and  $\Om$  an arrow \eqref{defomeg} in the fiber on $X$ of the
    stack
    $\mathrm{Lie}(\pac,\, \Om^{3}_{X/S})$. The 3-curvature arrow
  $\Om$  satisfies
    the two  conditions respectively expressed by the commutativity of the
  diagram of 2-arrows
\eqref{cubead} and of  the following
    diagram of 2-arrows:
   \end{theorem}

    \begin{equation}
          \label{bianchicube}
\includegraphics{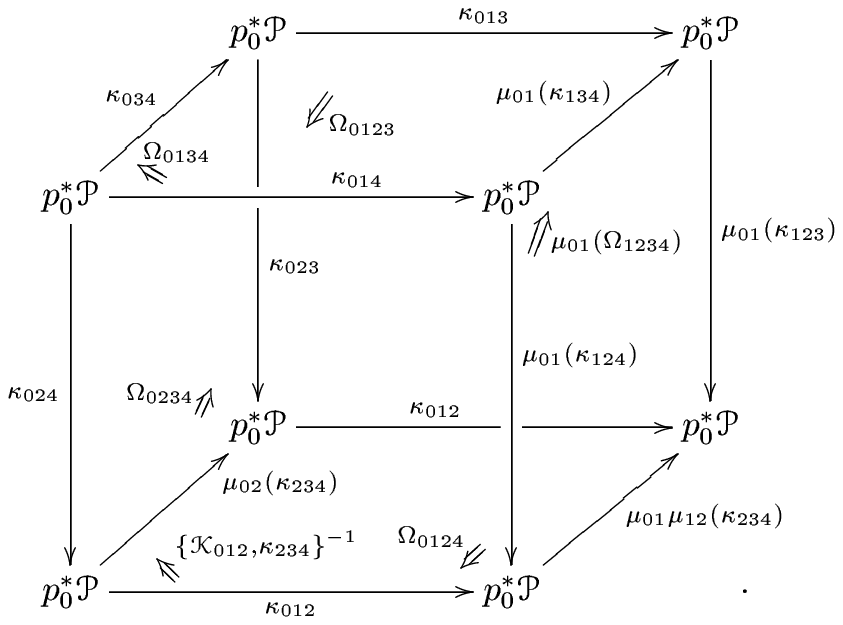}
         \end{equation}

         We call  diagram \eqref{bianchicube} in the fiber of $\pc$
 above $\ded$ the (higher) Bianchi cube.
         As indicated,
         five of its faces are pullbacks $\Omega_{ijkl}$ from $\dec$
         to $\ded$ of the 2-arrow
         $\Omega$. The bottom  face, labelled
         $\{\mathcal{K}_{012},\,\kappa_{234}\}^{-1}$
         is the  inverse of  the arrow
         \bee
\label{curly}
\{\kc_{012},\, \kappa_{234}\}:\, \kappa_{012}\, \mu_{02}(\kappa_{234}) \la
  (\mu_{01}\,  \mu_{12})(\kappa_{234})\, \kappa_{012}
\end{equation}
obtained by pasting the  conjugation 2-arrow
$M_{\kappa_{012}}(\mu_{02}(\kappa_{234}))^{-1}$ \eqref{msquare}
and the
2-arrow $\kc_{012}(\kappa_{234})$
\eqref{curly0}
         associated  to the natural transformation
         $\mathcal{K}_{012}$ \eqref{defkKad}
         and the object $\kappa_{234}$
         in  $\prc \pac\:$:
\begin{equation}
             \label{def:kc;k}
     \xymatrix@=60pt{\pra \pc
\ar@{}[rrd]^(.3){\,}="9"
\ar@{}[rrd]^(.45){\,}="10"
 \ar[rr]^{\kappa_{012}}&& \pra \pc \\
\pra \pc \ar[u]^{\mu_{02}(\kappa_{234})}
\ar@{}[u]_(.48){\,}="1"
     \ar[rr]_{\kappa_{012}}  &&\pra \pc\:.
     \ar@/^1.7pc/[u]^(.4){i_{\kappa_{012}}\mu_{02}(\kappa_{234})}
     \ar@{}@/^1.7pc/[u]_{\,}="5"
     \ar@{}@/^1.7pc/[u]^(.48){\,}="2"
     \ar@/_1.7pc/[u]_{\mu_{01}\mu_{12}(\kappa_{234})}
     \ar@{}@/_1.7pc/[u]^{\,}="6"
     \ar@{}"1";"2"^(.3){\,}="3"
      \ar@{}"1";"2"^(.4){\,}="4"
      \ar@{}"5";"6"^(.1){\,}="7"
  \ar@{}"5";"6"^(.9){\,}="8"
     \ar@{=>}"7";"8"^{\kc (\kappa_{234})}  
\ar@{=>}"9";"10"^{M^{-1}_{\kappa_{012}}} }
 \end{equation}
By whiskering, the giving of the arrow \eqref{curly} in $\pac$ is
 equivalent to that of an arrow 
\begin{equation}
\label{curly1}
\mathrm{id}_{\pra \pac} \la  (\mu_{01}\,  \mu_{12})(\kappa_{234})\,\,
 \kappa_{012}\,
  \mu_{02}(\kappa_{234}^{-1})\,\kappa_{012}^{-1}
\end{equation}
sourced at the identity object of $\pra \pac$,
 and which we will denote $[\kc_{012},\, \kappa_{234}]$
or simply $[\kc,\, \kappa]$.

\bigskip

The only part of theorem \ref{th:gl} which remains to be  proved
is the commutativity of
   the Bianchi cube
\eqref{bianchicube}. This is proved by inserting it  as one of
  the constituent cubes in a
  4-dimensional hypercube above $\ded$.
  Such a 4-cube
  may be viewed, when projected
  onto 3-space, as constructed from an ``inner'' 3-cube, with six cubes
  attached to it along its six  faces, and a final ``outer'' cube.  The
  4-cube which we will consider is the following one, in which for
  greater legibility, the 2-arrows along the faces have been omitted:

  \begin{equation}
         \label{hypercube}
\includegraphics{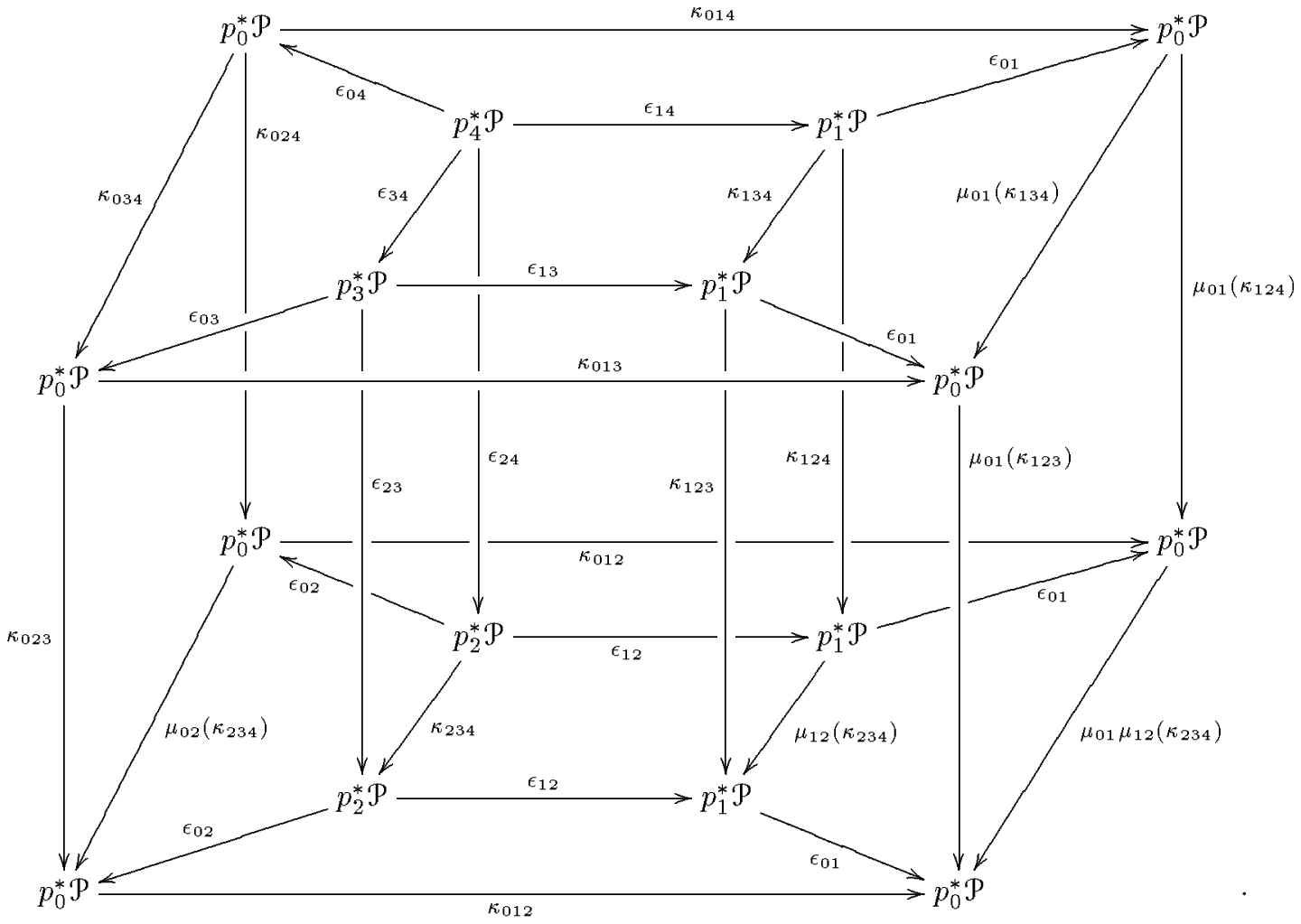}
\end{equation}

    Denoting by $C$ the cube \eqref{cube} by which we defined the
    3-curvature $\Omega$, and by $C_{ijkl}$ the corresponding pullback
    to $\ded$,  the construction of this hypercube is carried out by
    choosing $C_{1234}$ as inner cube, and by attaching to it the
    cubes $C_{ijkl}$ and the Bianchi cube according to the following
    table

.

     \begin{table}[ht]
      \begin{center}
         \renewcommand{\arraystretch}{1.3}
         \begin{tabular}{|c|c|c|c|c|c|c|c|}\hline
           inner  & left & right &top & bottom & front & back& outer\\
             \hline
             $C_{1234}$ & $C_{0234}$ &$M$-cube & $C_{0134}$
             & $\{\:,\:\}$ & $C_{0123}$ & $C_{0124}$ & Bianchi  \\ \hline
              &$K_{234}$ &$\Omega_{1234}$
               &$K_{134}$&$M_{12}(\kappa_{234})$& $K_{123}$&$K_{124}$&
             \\ \hline
       \end{tabular}
\caption{}
\end{center}
\end{table}
  The bottom line in the table
    describes the face along which each of these cubes
    is attached to the inner cube.
  We refer to the right-hand  cube as the $M$-cube,
    since  each  of its
faces involve a 2-arrow $M$ \eqref{msquare}. Its commutativity
follows from the compatibility of the squares \eqref{msquare} with
  the 1-arrows  $u$. It could also have been denoted $M_{01}(\Om_{1234})$
since it describes the effect of conjugation by $\epsilon_{01}$
  on the 2-arrow  $\Om_{1234}$. The bottom cube  in
    \eqref{hypercube}, denoted $\{\,,\,\}$ in the table, is the cube
    \begin{equation}
        \label{conjcube}
    \includegraphics{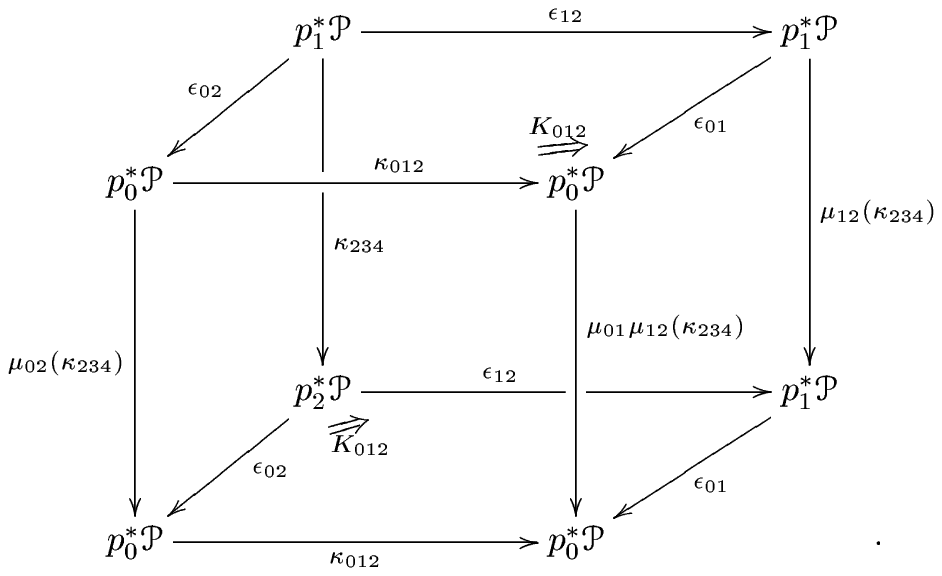}
\end{equation}
    which we will refer to  as $\mathrm{Conj}_{K_{012}}(\kappa_{234})$.
   Its left, back and right faces are
respectively
  the conjugation faces
  $M_{02}(\kappa_{234})$,  $M_{12}(\kappa_{234})$
and  $M_{01}(\mu_{12}(\kappa_{234}))$\eqref{msquare}.
It  may be viewed as an instance of the  prism \eqref{t-prism},
where we have set $u:= \kappa_{012}\, \epsilon_{02}$, $v:=
\epsilon_{01}\,\epsilon_{12}$      and $t:= K_{012}$.
Its front face $\{\kc _{012},\, \kappa_{234}\}^{-1}$
is the inverse of the  pasting   \eqref{def:kc;k}  of
  $M_{\kappa_{012}}(\mu_{02}(\kappa_{234}))$ and the 2-arrow
$t^{\mathrm{ad}}(\kappa_{234}):= \kc_{012}(\kappa_{234})$.

\bigskip

  The proof
         of the higher Bianchi identity now
  parallels that of Bianchi identity
         \eqref{bianchi:cl}: since seven of the faces of the hypercube
         \eqref{hypercube} are
         commutative, and all 2-arrows occuring are invertible, so is the
         outer (Bianchi) cube \eqref{bianchicube}.
         \begin{flushright}
             $\Box$
             \end{flushright}

The commutativity of the Bianchi cube
may be stated  algebraically as the following  non-abelian 3-cocycle
condition,
in which a whiskering on the right of a 2-arrow $\Omega$ by a 1-arrow
$\kappa$ is denoted ${}^{\ka\,}\!\Om$:
\begin{align}
\label{omcocy}
\Om_{0123} \,\, {}^{\mu_{01}(\ka_{123})\,}\!\Om_{0134}\, \,
\mu_{01}(\Om_{1234}) & =
{}^{\ka_{012}\,}\!\Om_{0234}\,\,
\{\kc_{012},\,\ka_{234}\}^{-1}\,\,
{}^{\mu_{01}\mu_{12}(\ka_{234})\,}\!\Om_{0124}\,.
\end{align}

\bigskip

  The  twisted 3-cocycle condition \eqref{omcocy}
  is not quite as formidable as it appears at first
     glance. We henceforth suppose that the gerbe $\pc$ is locally trivialized by
     a family  of
     objects $(x_{i})_{i \in I}$ above  an open
     cover $\mathcal{U}= (U_{i})_{i \in I}$ of $X$,
     and that  the corresponding sheaves
     $G_{i}:= \mathrm{Aut}_{\pc}(x_{i})$  are representable
     by flat $U_{i}$-group schemes. This  assumption
    will almost always be satisfied in practice, for example whenever $\pc$ is
     a $G$-gerbe for some flat $S$-group scheme $G$.  Since the
     restriction of the $G_{i}$-gerbe $\pc$ to $U_{i}$ is equivalent to
     $\mathrm{Tors}(G_{i})$, the restriction of $\pac$
     to $U_{i}$
     is  equivalent by \cite{Gir} IV proposition 5.2.5
     to the
     $gr$-stack $\mathrm{Bitors}(G_{i})$ of $G_{i}$-bitorsors. The
     latter is the  stack associated to the crossed module  $G_{i}
     \la \mathrm{Aut}(G_{i})$,  whose components are infinitesimally
     pushout reversing by \cite{cdf} proposition 2.3. For such
     sheaves, lemma 2.8 of {\it loc. cit.} (and its extension from
     commutator pairings to the pairing introduced there in (2.8.2))
     imply that
     there are canonical isomorphisms between those constituent
     terms of \eqref{omcocy} involving a  left action of a pullback
     of $\kappa$
     and the corresponding terms in which this action is omitted. These
     also involve  canonical isomorphisms permuting certain 1-arrows
     in $\pac$. Such canonical isomorphisms can additionally be used to
     permute terms in \eqref{omcocy}, so that
     the  condition \eqref{omcocy} may finally be replaced
     by the  simpler identity
     \bee
      \label{omcocy1}
     \Om _{0123}\, \, \Om_{0134}\,\, \mu_{01}(\Om_{1234})
   \:=\:  \Om_{0234}\,\, \Om_{0124}\, \,\{\kc_{012},\,  \kappa_{234}\}^{-1}.
      \end{equation}
        This can be expressed even more compactly as
  \bee
  \label{bianchi3}
  \delta^{3}_{\mu}\Om = \{\kc_{012},\,
    \kappa_{234}\}^{-1}\:,\end{equation}
or in additive notation as
\bee
\label{bianchi4}
\delta^{3}_{\mu}\Om +  \{\kc_{012},\, \kappa_{234}\} = 0\,.
\end{equation}
    Here $\delta^{3}_{\mu}$  is the third de Rham differential for the
    $gr$-stack $\pac$ with connection $\mu$, defined on objects as
    in  (\ref{def:d3mu}).
     For   $\kappa$
    trivial,
     this equation reduces
    to  the  more familiar  $\mu$-twisted higher Bianchi identity
    \bee
    \label{bianord}\delta^{3}_{\mu}\Om = 0
\end{equation}
    for $\Om$.
The
diagram \eqref{cubead} may be expressed algebraically as
\begin{equation}
   \label{relKkad}
   \kc_{012} \,\,\, \kc_{023}\,\,\,\,j_{\Omega} =
   \mu_{01}(\kc_{123})\, \,\,\,
  {}^{\mu_{01}(\kappa_{123})\,}\kc_{013}
\end{equation}
with $j$ as in \eqref{def:jmap}.
   Since the factors commute, this is equivalent to the equation
\begin{equation}
\label{relKkad1}
j_{\,\Omega} = \delta^2_{\mu^{\mathrm{ad}}}\,\kc
\end{equation}
in $\AL (\mathcal{E}q(\pac,\, \pac), \,\, \Omega^3_{X/S})$.

\bigskip

\subsection{} 
\label{hgen}
\setcounter{equation}{0}%
Consider two
curving pairs $(\epsilon',\, K')$ and $(\epsilon, \, K)$ on $\pc$.
 There always exists an object $h \in \pac$ and a
2-arrow $x$:
\begin{equation}
         \label{def:xa}
    \xymatrix{\prb \pc  \ar[rr]^{\epsilon'} \ar@{}[rr]_{\,}="1"
    \ar@/_/[rd]_{\epsilon}& & \pra \pc\\
    &\pra \pc   \ar@{}[r]^(.01){\,}="2" \ar@/_/[ru]_{h} &
      \ar@{}"1";"2"^(.2){\,}="3"
    \ar@{}"1";"2"^(.8){\,}="4"
    \ar@{=>}"3";"4"_{x} }
    \end{equation} 
We can  always assume that 
   there is  a 2-arrow   $\chi$
in the diagram
     \bee
  \label{def:chi1}
  \xymatrix{
  \pc  \ar[rr]^{1_{\pc}}="1"\ar[d]_(.45){\wr} && \pc\ar[d]^(.45){\wr}\\
  \Delta^{\ast}(p_{0}^{\ast}\pc)
  \ar[rr]^{\Delta^{\ast}h}="2" && \Delta^{\ast}(p_{0}^{\ast}\pc)
\:,
  \ar@{}"1";"2"^(.30){\,}="3"
    \ar@{}"1";"2"^(.70){\,}="4"
\ar@{=>}"3";"4"^{\,\chi}}
\end{equation}
so that  $(h, \chi)$ is an object  in  the stack of
  $\mathrm{Lie}(\pac)$-valued 1-forms on $X$. Similarly,
 $x$ is now an arrow  in the
fiber of  $\mathcal{E}q(\prb \pc,\, \pra \pc)$  above $\dea$, whose restriction above $X \hookrightarrow \dea$ is the identity. 
We also have an   arrow $a$   in the fibre of $\pac$
    above $\deb$, which lives in a diagram
    \begin{equation}
         \label{def:a}
   \xymatrix@R=15pt@C+4pt{\pra \pc & \pra \pc \ar[l]_{\kappa'}
   \ar@{}[ldd]^(.4){\,}="1"
    \ar@{}[ldd]^(.6){\,}="2"
   \ar@{=>}"1";"2"^{a}\\
   \pra \pc \ar[u]^{h_{01}}&\\
   \pra \pc \ar[u]^{\mu_{01}(h_{12})}
   & \pra \pc \ar[l]^{\kappa} \ar[uu]_{h_{02}}}
\end{equation}
    above  $\deb$. It is  characterized by  the commutativity of
the  following diagram of 2-arrows
\begin{equation}
\label{prism}
\includegraphics{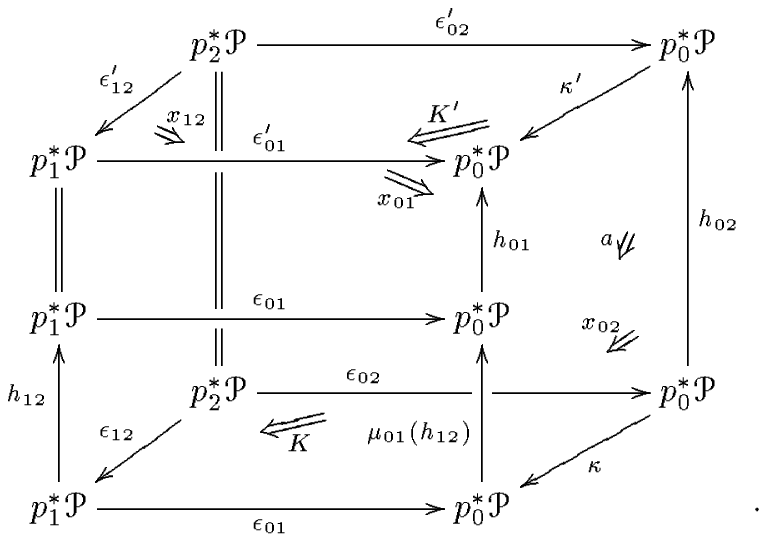}
\end{equation}
so that in particular the restriction of the $a$ above the degenerate
subsimplex $s\deb$ of $\deb$ is the identity.
    The unlabelled lower front square is the  2-arrow
    square which defines the arrow
   $\mu_{01}(h_{12})$ as in  the definition of $u^{\mathrm{ad}}$ 
\eqref{msquare}.
     We will call  any triple
$(x,h,a)$ as in \eqref{def:xa},\eqref{def:a}, for which diagram
\eqref{prism} is commutative and the degeneracy conditions  satisfied,
 a transformation triple.
By adification,  diagram \eqref{def:xa} induces  a diagram
        of monoidal stacks
         \begin{equation}
         \label{def:xi}
    \xymatrix{\prb \pac  \ar[rr]^{\mu'} \ar@{}[rr]_{\,}="1"
    \ar@/_/[rd]_{\mu}& & \pra \pac\\
    &\pra \pac   \ar@{}[r]^(.01){\,}="2" \ar@/_/[ru]_{i_{h}} &\qquad .
      \ar@{}"1";"2"^(.2){\,}="3"
    \ar@{}"1";"2"^(.8){\,}="4"
    \ar@{=>}"3";"4"_{\xi}}
    \end{equation}
 Diagram \eqref{prism}
induces  by adification
a diagram
\begin{equation}
\label{prism-1}
\includegraphics{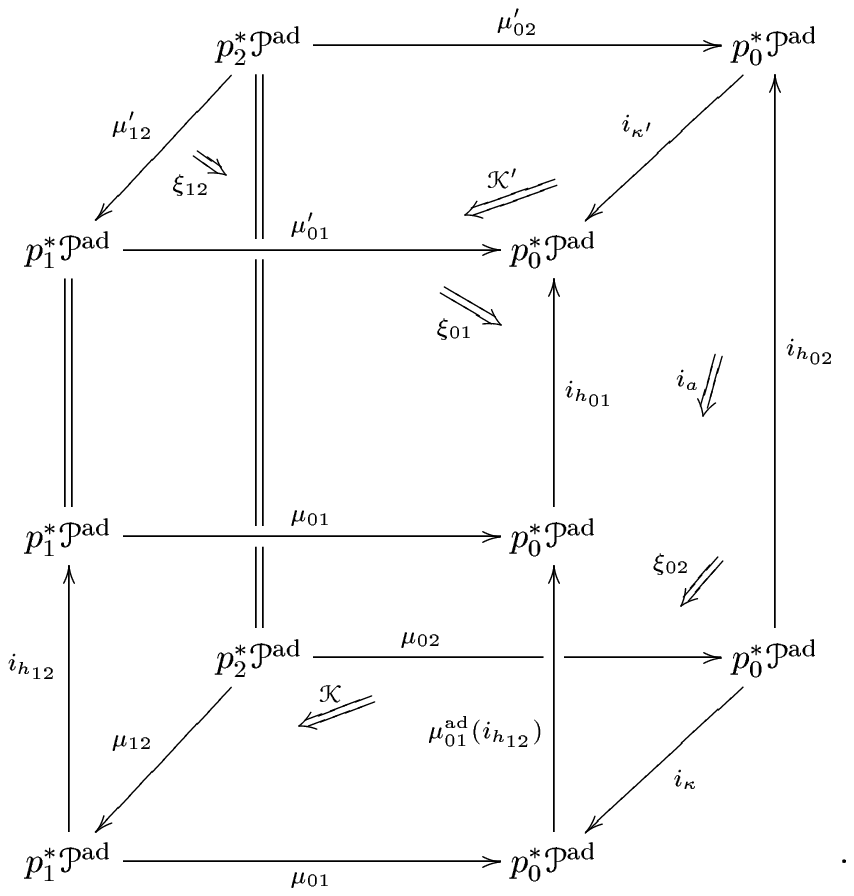}
\end{equation}
 and the   of 3-curvature forms
  $\Omega'$  and $\Omega$ respectively  determined by the curving
pairs $(\epsilon',\, K')$ and $(\epsilon,\, K)$ are compared to each
other by
 the diagram
  \begin{equation}
          \label{cob-prism}
  \includegraphics{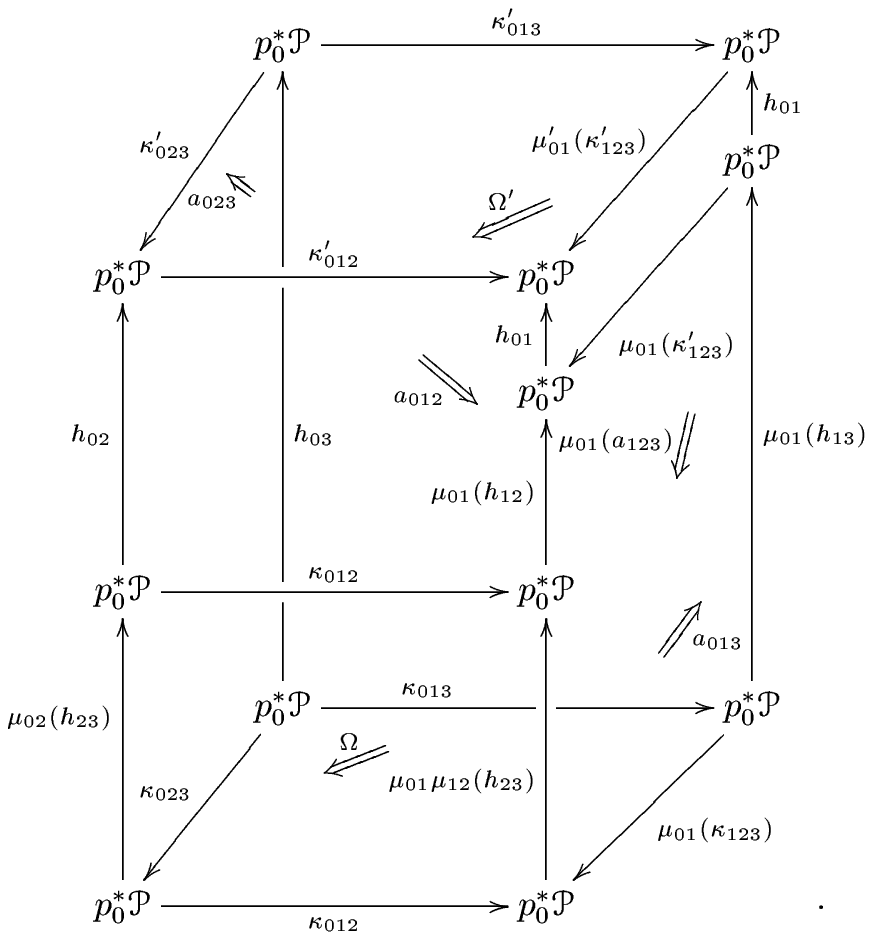}
        \end{equation}
       All 2-arrows in this diagram have already been
       introduced, except for the  lower front 2-arrow, and
       the upper right-hand unlabelled one. The former is the pullback
       from $\deb$ to  $\dec$ of the
       pasting 2-arrow  $\{\kc_{012},\, h_{23}\}$
       constructed as in \eqref{def:kc;k}.
   Similarly, for any object $u \in \prb \pac$,
    we denote by $\{\xi,\,\, u\}$  the pasting 2-arrow
  \begin{equation}
             \label{def:xiu}
     \xymatrix@=60pt{
     \pra \pc \ar[rr]^{h_{01}}
     && \pra \pc \\
\pra \pc \ar[u]^{\mu_{01}(u_{123})}
\ar@{}[u]_(.48){\,}="1"
     \ar[rr]_{h_{01}}  &&\pra \pc
     \ar@/^1.7pc/[u]^(.4){i_{h_{01}}\mu_{01}(u_{123})}
     \ar@{}@/^1.7pc/[u]_{\,}="5"
     \ar@{}@/^1.7pc/[u]^(.48){\,}="2"
     \ar@/_1.7pc/[u]_{\mu'_{01}(u_{123})}
     \ar@{}@/_1.7pc/[u]^{\,}="6"
     \ar@{}"1";"2"^(.3){\,}="3"
      \ar@{}"1";"2"^(.5){\,}="4"
      \ar@<-.5ex>@{}"4";"3"_{M_{h_{01}}}
     \ar@{=>}"4";"3"
     \ar@{}"6";"5"^(.1){\,}="7"
     \ar@{}"6";"5"^(.9){\,}="8"
  \ar@<-.5ex>@{}"7";"8"_(.4){\xi(\,u_{123})}
     \ar@{=>}"7";"8"}
\end{equation}
  which the 2-arrow $\xi$
\eqref{def:xi} induces  in $\pra \pac$.
    With this terminology, the 2-arrow pertaining to the unlabelled
    square in diagram \eqref{cob-prism} is $\{\xi,\,\, \kappa'_{123}\}$.

\bigskip

The commutativity of the diagram of 2-arrows \eqref{cob-prism}
is proved
 by considering the following pasting  diagram,
from which the 2-arrows have been omitted:
 \begin{equation}
         \label{hypercube3}
\includegraphics{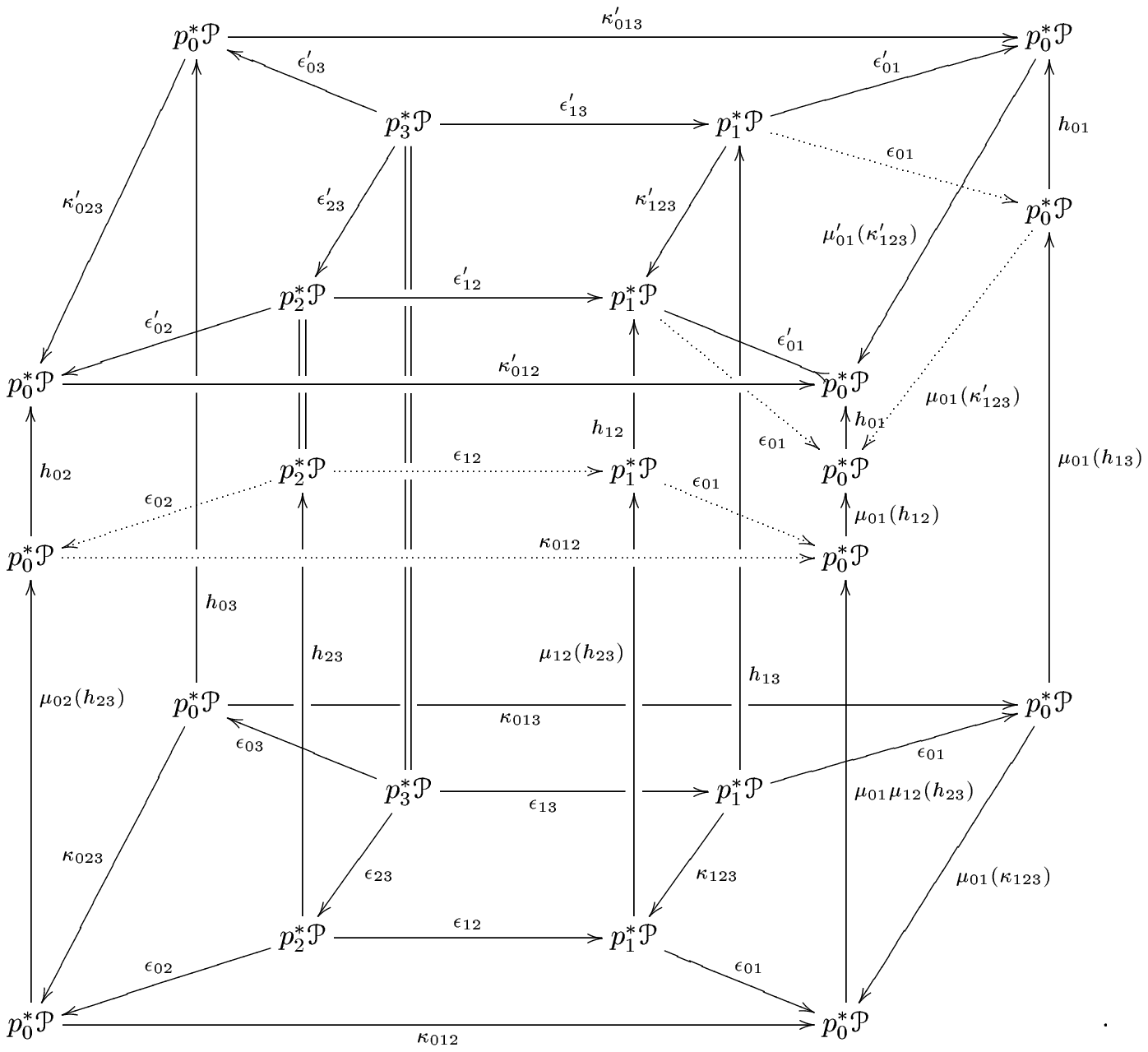}
\end{equation}
   \bigskip

   The eight cubes occuring here may be listed as follows:

      \begin{table}[ht]
      \begin{center}
         \renewcommand{\arraystretch}{1.3}
         \begin{tabular}{|c|c|c|c|c|c|c|c|}\hline
           inner  & left & right &top & bottom & front & back& outer\\
             \hline
             $P_{123}$ & $P_{023}$ &$\mathrm{Conj}_{x_{01}}(\kappa'_{123})$
             & $C'$
             &$ C$ &$P_{012}$ & $P_{013}$ & \eqref{cob-prism} \\ \hline
              &  &$\mathrm{Triv}_{\epsilon_{01}}(a_{123})$
              & && $\mathrm{Conj}_{K_{012}}(h_{23})$ & &
             \\ \hline
       \end{tabular}
\caption{}
\end{center}
\end{table}
\noindent An expression $P_{ijk}$  stands here for the corresponding 
 pullback from
  $\dec$ to $\ded$ of the diagram \eqref{prism}, and the cubes
  $C'$ and $C$ are those defining, as in \eqref{cube}, the 3-curvature
terms  $\Om'$ and $\Om$. Finally the  commutative prism
  $\mathrm{Conj}_{x_{01}}(\kappa'_{123})$  at the top right-hand side is
the degenerate prism \eqref{t-prism} associated to the 2-arrow $x: \epsilon'
  \Longrightarrow h\, \epsilon$ \eqref{def:xa} and the object
  $\kappa'_{123}$ in $\prb \pac$.
Since the other pasting diagrams listed here commute, so does the outer one
 \eqref{cob-prism}.
Diagram \eqref{cob-prism} may now be interpreted, together with diagrams
\eqref{def:a}, \eqref{def:xi} and \eqref{prism-1}, as the geometric
version of the rule by which the triple $(\xi,\,h,\,a)$  induced by
the transformation  triple $(x,\, h,\,a)$
transforms a cocycle quadruple $(\mu',\, \kc',\, \kappa',\,\Omega')$
into an equivalent quadruple $(\mu,\, \kc,\, \kappa,\, \Omega)$.

\bigskip

  In order to complete this global discussion of
   transformation  triples  for a gerbe $\pc$ it remains to describe
  the sense in which any two such   triples
  $(\xi,\, h,\,a)$  and  $(\xi',\, h',\,a')$
  are equivalent.
For any pair of  triples $(x,\,h,\,a)$ and $(x',\,h',\,a')$
there exists  a unique
  2-arrow
    \begin{equation}
        \label{def:r}
       \xymatrix{
  \pra \pc \ar@/^1pc/[rr]^{h}_{}="1"
\ar@/_1pc/[rr]_{h'}_{}="2" &&  \pra \pc
\ar@{=>}"1";"2"^{r}}
\end{equation}
above $\deb$
 for which the composed diagram of two arrows
  \begin{equation}
         \label{cond1:r}
    \xymatrix{\prb \pc  \ar[rr]^{\epsilon'} \ar@{}[rr]_{\,}="1"
    \ar@/_/[rd]_{\epsilon}& & \pra \pc\\
    &\pra \pc   \ar@{}[l]^(.01){\,}="2" \ar[ru]^{h}_{\,}="5"
    \ar@/_2pc/[ru]_{h'}^{\,}="6"  &
      \ar@{}"1";"2"^(.1){\,}="3"
    \ar@{}"1";"2"^(.7){\,}="4"
    \ar@{=>}"3";"4"_{x}
     \ar@{=>}"5";"6"_{r}}
    \end{equation}
     coincides with $x'$.
It then follows that the pullback
  $\Delta^{\ast}r$ of $r$  by the diagonal embedding is
compatible with  2-arrows $\chi$ and $\chi'$ \eqref{def:chi1}
associated to $h$ and $h'$, and
   also with  the obvious
    compatibility   between the 2-arrows $a$ and $a'$
    induced above $\deb$ by the 2-arrows
    $r_{01},\,\mu_{01}( r_{12})$ and $r_{02}$. We can now view the
    2-arrow $r$ as defining a transformation between the triple
    $(\xi',\, h',\,a')$ and the triple $(\xi,\, h,\,a)$.

\begin{remark}
{\rm
When the element $h \in \pac$ is trivial, diagram \eqref{def:xa} reduces
to a diagram 
 \begin{equation}
        \label{def:x}
       \xymatrix{
  \prb \pc \ar@/^1pc/[rr]^{\epsilon'}_{}="1"
\ar@/_1pc/[rr]_{\epsilon}_{}="2" &&  \pra \pc
\ar@{}"1";"2"^(.1){\,}="3"
\ar@{}"1";"2"^(.9){\,}="4"
\ar@<-1ex>@{=>}"3";"4"^{x}
}
\end{equation} 
The entire previous discussion may then be carried out in this restricted 
context, by collapsing  those
  edges  in  the previous
 diagrams which were defined
in terms of   $h$. In particular, the 2-arrow  \eqref{def:a} simply becomes
a 2-arrow
 \begin{equation}
        \label{def:a1}
       \xymatrix{
  \pra \pc \ar@/^1pc/[rr]^{\kappa'}_{}="1"
\ar@/_1pc/[rr]_{\kappa}_{}="2" &&  \pra \pc
\ar@{}"1";"2"^(.1){\,}="3"
\ar@{}"1";"2"^(.9){\,}="4"
\ar@<-1ex>@{=>}"3";"4"^{a}\,,
}
\end{equation}
While it is natural to study the effect on the 3-curvature of morphisms
 \eqref{def:x}, since these are  the    arrows in the stack of connections
on $\pc$, it should be noted that for  an arbitrary pair of connections
 $\epsilon'$ and $\epsilon$, there does not necessarily exist a 2-arrow $x$
 \eqref{def:x} of which they are source and target.
}\end{remark}

\bigskip

\subsection{}
\setcounter{equation}{0}%
\label{alg}
The
previous geometric 
 discussion will now be interpreted in a  more algebraic manner.
Let $\gc$ be a $gr$-stack. We will now define a connection on $\gc$
as we did in definition \ref{def:gcon} for groups.

\begin{definition}
\label{def:grstcon}
A connection on a $gr$-stack $\gc$ on  $X$ is a monoidal equivalence
\[\mu: \prb \gc \la \pra \gc\]
of $gr$-stacks above $\dea$, together with a monoidal
2-morphism \[1_{\gc} \Longrightarrow  \De^{\ast} \mu\] 
\end{definition}
For $n>0$, we  define  $\mu$-twisted differential functors 
\begin{equation}
   \label{defdiff}
\delta^n_\mu:\mathrm{Lie}(\gc,\, \Omega^n_{X/S}) \la
  \mathrm{Lie}(\gc,\, \Omega^{n+1}_{X/S})
\end{equation}
between     the stacks
  of  $\mathrm{Lie}(\gc)$-valued     forms on $X$
by the same formulas as in \eqref{def:dmu},
but now applied to both objects and
    arrows in the Picard stacks $
\mathrm{Lie}(\gc,\, \Omega^{n }_{X/S})$ (definition \ref{def:b2}).
When $\gc$ is the associated $gr$-stack of a pushout reversing crossed module,
the functors $\de^i_{\mu}$ are morphisms of Picard stacks for $n > 1$. This
is in particular  the case when $\gc$ is the stack of $G$-bitorsors,
 with $G$ a flat  $X$-group.

\bigskip

The diagram of pointed   stacks
\begin{equation}
     \label{dr0p}
\xymatrix{\ar[r]^(.2){\delta^1_{\mu}} &
\mathrm{Lie}(\gc,\, \Omega^{2}_{X/S}) \ar[r]^{\delta^2_{\mu}}
&  \mathrm{Lie}(\gc,\, \Omega^{3}_{X/S})\ar[r]^{\delta^3_{\mu}}&
  \mathrm{Lie}(\gc,\, \Omega^{4}_{X/S}) \ar[r]^(.75){\delta^4_{\mu}}&{}
}\end{equation}
can be displayed as follows by separating  the objects and arrows:
     \begin{equation}
\label{drp}
\xymatrix{\ar[r]^(.2){\delta^1_{\mu}} &
\AL(\gc,\, \Omega^{2}_{X/S})  \ar@<1ex>[d]^t \ar@<-1ex>[d]_s
  \ar[r]^{\delta^2_{\mu}}
&  \AL(\gc,\, \Omega^{3}_{X/S})
   \ar@<1ex>[d]^t \ar@<-1ex>[d]_s
\ar[r]^{\delta^3_{\mu}}&
  \AL(\gc,\, \Omega^{4}_{X/S})
\ar@<1ex>[d]^t \ar@<-1ex>[d]_s
  \ar[r]^(.8){\delta^4_{\mu}} &
\\
\ar[r]^(.2){\delta^1_{\mu}} &
\OL(\gc,\, \Omega^{2}_{X/S}) \ar[r]^{\delta^2_{\mu}}
&  \OL(\gc,\, \Omega^{3}_{X/S})\ar[r]^{\delta^3_{\mu}}&
  \OL(\gc,\, \Omega^{4}_{X/S})  \ar[r]^(.8){\delta^4_{\mu}} &
}\end{equation}
with $s$ and $t$ the source and target maps.
Since any 2-arrow is equivalent by whiskering on the left
as in  \eqref{def:alph-t} to one 
whose source is
the identity 1-arrow,
we can restrict ourselves without loosing any information
from the sets  $\AL(\gc,\, \Omega^{i}_{X/S})$ to the subsets
  $\AI(\gc,\, \Omega^{i}_{X/S})$ consisting of those arrows which are
  sourced at the
identity object $I$ of $\gc$.
Diagram \eqref{drp} can be  be   replaced  by the simpler
  diagram 
with commutative squares:
  \begin{equation}
\label{dr1p}
\xymatrix{\ar[r]^(.2){\delta^1_{\mu}} &
\AI(\gc,\, \Omega^{2}_{X/S})  \ar[d]^{t}
  \ar[r]^{\delta^2_{\mu}}
&  \AI(\gc,\, \Omega^{3}_{X/S})
   \ar[d]^{t}
\ar[r]^{\delta^3_{\mu}}&
  \AI(\gc,\, \Omega^{4}_{X/S})
\ar[d]^{t}
  \ar[r]^(.8){\delta^4_{\mu}} &
\\
\ar[r]^(.2){\delta^1_{\mu}} &
\OL(\gc,\, \Omega^{2}_{X/S}) \ar[r]^{\delta^2_{\mu}}
&  \OL(\gc,\, \Omega^{3}_{X/S})\ar[r]^{\delta^3_{\mu}}&
  \OL(\gc,\, \Omega^{4}_{X/S})  \ar[r]^(.8){\delta^4_{\mu}} &\, .
}\end{equation}
We now set $\gc:= \pac$ and $\mu:= \epsilon^{\mathrm{ad}}$
\eqref{defmuad}, 
and  enrich diagram \eqref{dr1p}
to a diagram
\begin{equation}
\label{dr2p}
\xymatrix{\ar[r]^(.2){\delta^1} &
\AI(\pac,\, \Omega^{2}_{X/S})     \ar[d]^{t}
  \ar[r]^{\delta^2_{\mu}}
&  \AI(\pac,\, \Omega^{3}_{X/S})
   \ar[d]^(.3){t}
\ar[r]^{\delta^3_{\mu}}&
  \AI(\pac,\, \Omega^{4}_{X/S})
\ar[d]^{t}
  \ar[r]^(.8){\delta^4_{\mu}} &
\\
\ar[r]_(.2){\delta^1}
  &
\OL(\pac,\, \Omega^{2}_{X/S})\ar@{-->}[rru]^(.35){[\kc,\,-]}
  \ar[r]_{\delta^2_{\mu}}
&  \OL(\pac,\, \Omega^{3}_{X/S})\ar[r]_{\delta^3_{\mu}}&
  \OL(\pac,\, \Omega^{4}_{X/S})  \ar[r]_(.8){\delta^4_{\mu}} &\,.
}\end{equation}
The dotted arrow sends  $u \in \OL(\pac,\, \Omega^{2}_{X/S})$
to the element $[\kc_{012},\,u_{234}]$ in  $ \AI(\pac,\, \Omega^{4}_{X/S})$
defined by the same proceedure as in \eqref{curly1},  but now applied to an 
arbitrary  $u \in  \OL(\pac,\, \Omega^{2}_{X/S})$ rather than simply 
to the fake curvature $\kappa$. More generally, this construction
determines an arrow
\bee
\label{brack2}
 [\kc,\, -]: \OL(\pac,\, \Omega^{i}_{X/S}) \la  \AI(\pac,\, \Omega^{i+2}_{X/S})
\end{equation}
and which associates to $u \in  \OL(\pac,\, \Omega^{i}_{X/S})$
an
element $[\kc_{012},\, u_{2,\cdots , i+2}]$  which we will also simply denote
 by 
$[\kc,\, u]$.

\bigskip

By examining the 1-skeleton of the  cube \eqref{bianchicube}, we see
that the equation
\begin{equation}
\label{reldk}
\delta^3_{\mu}\, \delta^{2}_{\mu}(\kappa) = t[\kc, \kappa]
\end{equation}
  is satisfied whenever $(\mu,\, \kc,\, \kappa)$ are as in
\eqref{defkKad}, for example when they are as in theorem \ref{th:gl}
  the first three
terms of a quadruple $(\mu,\, \kc,\, \kappa,\, \Omega)$.
This makes up to some extent
  for the fact that the composite morphism
  $\delta^3_{\mu} \circ \delta^2_{\mu}$ is not
  trivial on objects, since we are here   in the presence of a
non-integrable connection  $\mu$ on $\pac$.
   Part of theorem \ref{th:gl} may now
  be restated
as the assertion that the pair  $(\kappa,\, \Omega)$
in  a quadruple  $(\mu,\, \kc,\, \kappa,\, \Omega)$  determines  a cocycle
in  the diagram \eqref{dr2p}.
  By very definition
of the 3-curvature
  $\Omega$  \eqref{defomeg}, the corresponding element
of $ \AI(\pac,\, \Omega^{3}_{X/S})$, which we  denote, as in 
\S \ref{subsec:curv-Bian1} by  $\widehat{\Omega}$,
satisfies the equation
\begin{equation}
     \label{add3}
   t\, \widehat{\Omega} +   \delta^{2}_{\mu}(\kappa) =  0\,.\end{equation}
The second  equation which
the pair $(\kappa,\, \Omega)$ must satisfy in order  to be a
cocycle in this de Rham complex is, in additive notation,    the equation
\eqref{bianchi4} which we restate here as
\begin{equation}
\label{neomegeq}
\delta^3_{\mu}\widehat{\Omega} + [\kc,\, \kappa] = 0\,.
\end{equation}

\medskip 

\begin{remark}
{\rm
Returning   to the diagram \eqref{dr2p}, let us    observe
that a  cocycle pair  $(\kappa,\, \Omega)$ in a quadruple $(\mu, \,
\kc,\, \kappa,\, \Omega)$ is endowed with  an arrow
\[\xymatrix@C=8pt{I \ar[rrr]^(.4){[\kc,\, \kappa]} &&&
  \delta^3_{\mu}\circ  \delta^2_{\mu}\, \kappa }\]
in the category $\mathrm{Lie} (\pac, \, \Omega^4_{X/S})$,
and the 3-curvature $\widehat{\Omega}$ is
an arrow \[\xymatrix@C=8pt{I \ar[rr]^{\widehat{\Omega}} && 
\delta^2_{\mu}}(\kappa^{-1})\]
in $\mathrm{Lie}\,(\pac, \Omega^3_{X/S})$.
The higher Bianchi identity
  \eqref{neomegeq}
   asserts      that the composite arrow
\begin{equation}
\label{bianchicat}
    \xymatrix@C=8pt{I \ar[rr] &&\delta^3_{\mu}(I)
\ar[rrr]^(.4){\delta^3_{\mu}\,\widehat{\Omega}}&& &\delta^3_{\mu}
  \delta^2_{\mu}(\kappa^{-1}) \ar[rrrr]^(.6){[\kc,\, \kappa^{-1}]^{-1}}& &&& I}
  \end{equation}
is the identity arrow, a very natural assertion  in the present categorical
homological algebra context (for a similar condition, see for example
\cite{ul2} (0.1))} \end{remark}

\vspace{.75cm}

In addition to \eqref{add3} and \eqref{bianchi4},
the quadruple $(\mu,\, \kc,\, \kappa,\, \Omega)$ introduced
in theorem \ref{th:gl} satisfies two additional conditions. The first
of these describes diagram \eqref{defkKad}. In the present context, this
asserts that
\bee 
\label{add01}
 t \kc = \delta^{1}\mu - i_{\kappa}\:,\end{equation}
an expression whose right-hand side has   meaning in
$\OL (\pac,\, \Omega^{3}_{X/S})$  even though
$\mathrm{Ob} \, \mathcal{E}q_{.}(\prc \pc, \,
\pra \pc)$ is not canonically identified with the stack
  $\OL (\pac,\, \Omega^{2}_{X/S})$. The inner conjugation map $\gc \la
  \mathcal{E}q(\gc)$ associated to a monoidal category $\gc$ will
  as in \eqref{def:jmap} now  be denoted $j$ rather than $i$.
 It should not be confused  with
  the target map $t$, which in the category $\gc$ of $G$-bitorsors
  corresponds to the inner conjugation  map for the group $G$. Equation
\eqref{add01}  is therefore rewritten  as
  \bee
  \label{add1}
   t \kc = \delta^{1}\mu -j(\kappa)\,.
   \end{equation}
  The second condition
     expresses the commutativity of the diagram \eqref{cubead}
and is simply  equation \eqref{relKkad1},
an expression whose right-hand side again has
  meaning,  even though
$\mathrm{Ob} \, \mathcal{E}q_{.}(\prc \pc, \,
\pra \pc)$ is not canonically identified with the stack
  $\OL (\pac,\, \Omega^{2}_{X/S})$.

\bigskip

  Instead of  viewing the pair $(\kappa,\, \widehat{\Omega})$ simply
as a cocycle in the diagram \eqref{dr2p}, we now view
the full quadruple $(\mu,\, \kc,\, \kappa,\, \widehat{\Omega})$ as
a cocycle in the extended  diagram

\begin{equation}
     \label{dr7}
     \xymatrix@R=48pt{
&&&&&&\ap_{4} \ar[r]^{t}
\ar[d]^(.25){j} |!{[dl];[r]}\hole
& \op_{4}\ar[d]^{j}\\
&&&&*+[F]{\ap_{3}}\ar[r]\ar[urr]^{\delta^{3}_{\mu}} \ar[d]
|!{[dl];[r]}\hole |!{[dl];[urr]}\hole
\ar@{}[d]^(.5){}="1"
\ar@{}[d]^(.928){}="2"
& \op_{3}
\ar[urr]_(.7){\delta^{3}_{\mu}}
 \ar[d]^(.25){j} & \ae_{4}\ar[r]^{t} &
\oe_{4}\\
&&&
*+[F]{\op_{2}}  
\ar[urr]_(.7){\delta^{3}_{\mu}}
\ar[d]_(.35){j}
\ar@{.>}[rrruu]_(.75){[\widehat{\kc},\,]} 
&
\ae_{3}
\ar[urr]_(.75){{\delta^{2}_{\mu^{\mathrm{ad}}}}}  |!{[ur];[r]}\hole
\ar[r]_(.42){t}
&\oe_{3}\ar[urr]_{\delta^{3}_{\mu^{\mathrm{ad}}}}
&&\\
&&*+[F]{\widehat{\kc}}
\ar[urr]_(.75){{\delta^{2}_{\mu^{\mathrm{ad}}}}}  |!{[ur];[r]}\hole
\ar[r]_(.42){t}& 
\mathrm{OCu}\ar[urr]_{{\delta^{2}_{\mu^{\mathrm{ad}}}}}  &&&&\\
  & *+[F]{\mu}\ar[urr]_{\delta^{1}}
&&&&&&
\ar@{->}"1";"2"^(.4){j}
}
\end{equation}

\bigskip

\noindent whose vertices have been abbreviated according to the following table

\bigskip

\begin{table}[ht]
     \label{abb1}
      \begin{center}
         \renewcommand{\arraystretch}{1.3}
         \begin{tabular}{|c|c|}\hline
           abbreviation & full description \\
             \hline
            $AP_{i}$&$\mathrm{Ar}_{I}\ \mathrm{Lie}(\pac,\,\Om^{i}_{X/S})$
            \\ \hline
            $OP_{i}$ & $\mathrm{Ob}\ \mathrm{Lie}(\pac,\,\Om^{i}_{X/S})$
            \\ \hline
               $AE_{i}$ &
$ \mathrm{Ar}_{I}\ \mathrm{Lie}(\mathcal{E}q(\pac),\,\Om^{i}_{X/S})$
          \\ \hline  $OE_{i}$ &
$  \mathrm{Ob}\ \mathrm{Lie}(\mathcal{E}q(\pac),\,\Om^{i}_{X/S})$
             \\ \hline  $\mathrm{OCu}$&$\mathrm{Ob}\ \mathrm{Curv}(\pac)$
             \\ \hline
       \end{tabular}
\caption{}
\end{center}
\end{table}
\noindent The terms $\mu$ and $\kc$ respectively live in the
sheaves
$\mathrm{Ob}\ \mathrm{Conn}$ and  $\mathrm{Ar}\ \mathrm{Curv}$.
The two upper lines in diagram \eqref{dr7} are those of the de Rham
diagram \eqref{dr2p} and the two lower ones constitute a portion of
the de Rham diagram  associated in a similar manner to the $gr$-stack
  $\mathcal{E}q(\pac)$ with connection $\mu^{\mathrm{ad}}$. The
  vertical maps $j$ are  the maps induced at the Lie level by the
  inner conjugation functors
  .
  The first one of these however (with source $\mathrm{Ob}\ 
\mathrm{Lie}(\pac, \,
  \Om^{2}_{X/S})$),
  is not a monoidal functor since its target is not a monoidal
  category. Instead, it is
  an action of its source on the sheaf $\mathrm{Ar}\ \mathrm{Curv}$ of curvings.

  \bigskip

  The situation is now analogous to that described in the torsor case by
diagram \eqref{interp}.
The elements $(\mu, \kc)$  in
  the  two lower lines are
  parameters   whose values  determine the two sorts of 
arrows occuring  in the upper lines of the diagram.
     The  additional conditions \eqref{add1} and \eqref{relKkad1}
imply that the  quadruple $(\mu,\, \kc,\, \kappa,\,
\widehat{\Omega})$, whose terms live in each of the four framed
locations,
is indeed a cocycle for  the full diagram
  \eqref{dr7}.

  \bigskip

  As in the  torsor
case    \eqref{interp}, we may now interpret the transformation
  conditions for quadruples embodied in the discussion of \S \ref{hgen} as a
  coboundary condition in this diagram. Consider a transformation triple
  $(x, \, h, \, a)$,
   which according to diagrams
  \eqref{def:xa}-\eqref{prism} and \eqref{def:xi}-\eqref{cob-prism}
  transforms a quadruple  $(\mu',\, \kc',\, \kappa',\, \Omega')$ into a new
  quadruple  $(\mu,\, \kc,\, \kappa,\, \Omega)$. The effect of each of
  these transformations on the terms in diagram \eqref{dr7} can be read
 off from
  the  geometric discussion in \S \ref{hgen},
  and we will simply record the result  here.

\bigskip

  The  coboundary transformations  in question are
\begin{align}
\label{trans1}
\mu &=\mu' + t\xi - j_h \\
\kappa & = \kappa' + ta - \delta^{1}_{\mu}h\\
\kc &= \kc' +\delta^{1}_{\mu'}(\xi) - j_a\\
\Om &= \Om' -\delta^{2}_{\mu}(a) - [\xi,\, \kappa'] - [\kc,\, h]
-[h,a]
\end{align}
These relations are those
implicit for coboundaries in the extended diagram
\begin{equation}
     \label{dr4}
\xymatrix@R=48pt@C=36pt{
&&&&&&\ap_{4} \ar[r]^{t}
\ar[d]_(.3){j} |!{[dl];[r]}\hole
& \op_{4}\ar[d]^{j}\\
&&&&*+[F]{\ap_{3}}\ar[r]\ar[urr]^{\delta^{3}_{\mu'}}
\ar[d]|!{[dl];[r]}\hole |!{[dl];[urr]}\hole
\ar@{}[d]^(.1){}="1"
\ar@{}[d]^(.3){}="2"
& \op_{3}\ar[urr]_(.3){\delta^{3}_{\mu'}} \ar[d]_(.35){j} & \ae_{4}\ar[r]^{t} &
\oe_{4}\\
&&*+[F--]{\ap_{2}}\ar[urr]^{\delta^{2}_{\mu'}}\ar@{|->}[d]^(.65){j}
\ar[r]_(.6){t}& *+[F]{\op_{2}} \ar@{-->}[ru]
\ar[urr]_(.3){\delta^{2}_{\mu'}}\ar[d]_(.35){j}
\ar@{.>}[rrruu]_(.75){[\widehat{\kc}',\,\,]} &
\ae_{3}
\ar[urr]_(.75){{\delta^{2}_{{\mu'}^{\mathrm{ad}}}}}  |!{[ur];[r]}\hole
\ar[r]_{t}
&\oe_{3}\ar[urr]_{\delta^{3}_{{\mu'}^{\mathrm{ad}}}}
&&\\
&*+[F--]{\op_{1}}\ar@{.>}[rrruu]^(.2){[\widehat{\kc}',\,\,\,]}
\ar[urr]\ar@{|->}[d]_(.35){j}
&*+[F]{\widehat{\kc}'}
\ar[urr]_(.75){{\delta^{2}_{{\mu'}^{\mathrm{ad}}}}}  |!{[ur];[r]}\hole
\ar[r]_(.42){t}& \mathrm{OCu}\ar[urr]_{{\delta^{2}_{{\mu'}^{\mathrm{ad}}}}}
&&&&\\
*+[F--]{\widehat{\xi}} \ar[r]_{t}\ar[urr]^(.3){\delta^{1}}|!{[ur];[r]}\hole &
*+[F]{\mu'}\ar[urr]_{\delta^{1}}
&&&&&&
\ar@{-}"1";"2"
}
\end{equation}
obtained from \eqref{dr7} by adding the  extra vertices $\mathrm{OP}_{1}$
and $\ap_{2}$ and an extra parameter
$\hat{\xi}$ in $\mathrm{Ar}_I\ \mathrm{Conn}(\pac)$ together with the
corresponding edges. The arrow $j$ with source $OP_1$ 
represents an action of its source on its target, the  objects
in the 
 stack  of connections on the $gr$-stack $\pac$.
An additional  enrichment has been 
 provided in this  diagram. This is
the unlabelled dotted arrow $[\xi,\,\,]$ from $\op_{2}$ to
$\ap_{3}$ which the parameter $\xi$ generates, just as the parameter
$\kc'$ generated 
 the arrow $[\kc',\,\, ]$ in \eqref{dr2p}.
The positions in this
diagram which the
three coboundary  terms $(\xi,\,h,\,a)$ occupy are highlighted by the
three dotted arrow frames and do indeed represent the three possible
slots for coboundary terms.

\begin{remark}{\rm The term $r$  \eqref{def:r}
which generates a transformation
between a pair of   triples $(\xi, \, h, \, a)$ and
$(\xi', \, h', \, a')$ can be represented in similar terms by adding
to diagram \eqref{dr4} a final vertex $\ap_{1}:=\mathrm{Ar}\
\mathrm{Lie}(\pac,\, \Om^{1}_{X/S})$  and the corresponding edges,
in order to complete the composite
parallelepiped. The transformations  which the section $r \in \ap_{1}$
determines  can then either  be read off  from diagram \eqref{def:xi}
and  the diagram expressing the additional compatibilities for $r$
mentioned at the end of \S \ref{hgen},
or simply from the completed parallellepiped. In particular, it is
seen that $r$  acts respectively on each of the terms of the triple
through $j$, through $t$ and through $\delta^{1}_{\mu'}$.
}  \end{remark}

%% file: dgg-ch5.tex
\section{The partial decomposition of gerbes with curving pairs}
\label{partly}
\subsection{}
\setcounter{equation}{0}%
The  global description of the quadruple $(\mu,\, \kc,\, \kappa,\, \Om)$
 associated by theorem \ref{th:gl}   to a  curving pair $(\epsilon,\, K)$ on a
gerbe $\pc$ above  an $S$-scheme $X$  can be made more explicit if we choose a
 family of local sections $x_i \in \mathrm{ob}(\pc_{U_i})$,
 since we can then apply Morita theory.
 We have already  seen that the 
    gerbe $\pc$ is then  described by the 
    associated family of bitorsor cocycles $(G_{i}, \, P_{ij}, 
    \,\Psi_{ijk})$ (\ref{morita}, \ref{1bcoc}). 
 By (\ref{defgai}),
a connection $\epsilon$ on $\pc$ 
is defined locally  by  the 
$(p_{0}^{\ast}G_{i},\, p_{1}^{\ast}G_{i})$-bitorsor
\bee
\label{def:gai}\Gamma_{i}:=
\mathrm{Isom}(\epsilon_{\mid U_{i}}(p_{1}^{\ast}T_{i}),\, 
p_{0}^{\ast}T_{i}) 
\simeq \mathrm{Isom}_{p_{0}^{\ast}\pc}(\epsilon(p^{\ast}_{1}x_{i}),\, 
p_{0}^{\ast}x_{i})\end{equation}
on $\Delta^{1}_{U_{i}/S}$.
The 2-arrow $\eta$ \eqref{2ar:eta} determines an isomorphism
$G_{i}$-bitorsors
\bee
\label{def:triv}
\xi_{i}: T_{i} \stackrel{\sim}{\la} \Delta^{\ast}\Gamma_{i}\:. \end{equation}
with  source the trivial $G_{i}$-bitorsor $T_{i}$ above $U_i$,  which 
identifies the pullback of $\Ga_{i}$ by the diagonal 
embedding $\De: U_{i} \hookrightarrow \De^{1}_{U_{i}/S}$
with the trivial $G_{i}$-bitorsor 
$\mathrm{Isom}_{\pc}(x_{i},\, x_{i})$.
For varying indices
 $i$ and $j$, the pullbacks
 to $\Delta^{1}_{U_{ij}/S}$ are 
related to each other by  isomorphisms of 
$(p_{0}^{\ast}G_{i},\,p_{1}^{\ast}G_{i})$-bitorsors\footnote{We are really dealing with the restrictions
of $\Ga_{i}$ and $\Ga_{j}$
to $\De^{1}_{U_{ij}}$, but we will omit  this from the notation here
and in similar situations.}
(\ref{ardij}) 
\begin{equation}
    \label{defdij1}
\xymatrix{
\Gamma_{i} \ar[r]^(.4){g_{ij}}& {}^{ P_{ij}\ast\,\,}\!\Gamma_{j}
}
\end{equation}
 compatible with the pointings of the source and target determined 
 by $\xi_{i}$ and $\xi_{j}$,  
where the target of $g_{ij}$ is the $(\pra G_{i},\, \prb G_{i})$-bitorsor  
above $\De^{1}_{U_{ij}}$ defined by 
\[{}^{ P_{ij}\ast\,\,}\!\Gamma_{j} :=
p_{0}^{\ast}P_{ij}  \,\wedge^{\pra G_{j}} \,\Gamma_{j}\wedge^{\prb 
G_{j}}
(p_{1}^{\ast}P_{ij})^{0}\:.\]
 The upper ${\ast}$ sign has been 
inserted, as in \eqref{omcoc-alt}, in order to remind us that
this twisting  by $P_{ij}$ of the bitorsor $\Ga_{j}$ 
is  a  \emph{twisted} adjoint action.
 Finally, for all 
indices $(i, j, k)$, the diagram
\bee
\label{compdij3}
\xymatrix{
\Ga_{i} \ar[rrr]^{g_{ij}}\ar[d]_{g_{ik}}&&&
{}^{P_{ij}\ast\,}\!\Ga_{j}
\ar[d]^{{}^{P_{ij}\ast\,\,}\! g_{jk} }
\\
{}^{P_{ik}\ast\,\,}\! \Ga_{k} 
\ar[rrr]_<<<<<<<<<<<<<<<<<<<{{}^{\Psi_{ijk}\ast\,}\!(1_{\Ga_{k}})}
& &&
{}^{(P_{ij} \wedge P_{jk})\ast\,\,}\!
\Ga_{k}
}\end{equation}
above $\De^{1}_{U_{ijk}/S}$ commutes.

\bigskip

   Similarly, an equivalence $\zeta$ 
(\ref{def:zeta}) between a pair of connections  $\epsilon,\, \epsilon'$ 
on a locally trivialized gerbe $(\pc,\, (x_{i})_{i \in I})$  is described
by  a family of 
 $ (p_{0}^{\ast}G_{i},\,p_{1}^{\ast}G_{i})$
-bitorsor isomorphisms 
\begin{equation}
    \label{equivzi}
    \zeta_{i}: \Ga_{i}^{\epsilon'} \la 
\Gamma_{i}^{\epsilon}
\end{equation} 
between the bitorsors above $U_{i}$ associated to $\epsilon'$ and 
$\epsilon$, and for which for which the induced diagrams 
\bee \label{comzeta}
     \xymatrix{
\Ga_{i}^{\epsilon'} 
\ar@<-2pt>[r]^(.4){g_{ij}^{\epsilon'}} 
\ar[d]_{\zeta_{i}}
& {}^{ P_{ij}\ast\,\,}\!\Gamma_{\,j}^{\epsilon'}
\ar[d]^{{}^{ P_{ij}\ast\,\,}\!\zeta_{\,j}}\\
\Ga_{i}^{\epsilon}
\ar@<-2pt>[r]_(.4){g_{ij}^{\epsilon}} 
&{}^{ P_{ij}\ast\,\,}\!\Gamma_{\,j}^{\epsilon}
}
\end{equation}
commute.
    
\bigskip

Assuming $G_{i}$ is a $U_{i}$-group scheme, or more generally is 
pushout reversing, we consider, as in \cite{cdf} proposition 2.2,  
the  exact sequence of sheaves of groups with abelian kernel 
on $\Delta^{1}_{U_{i}/S}$
\bee
\label{abker}
0 \la \Delta_{\ast}(\mathrm{Lie}(G_{i}, \Om^{1}_{U_{i}/S})) 
\la p^{\ast}_{0}G_{i} 
 \la \Delta_{\ast}G_{i}\:.\end{equation}
 The $\xi_{i}$-pointed sections of the left
$p_{0}^{\ast}G_{i}$-torsor $\Gamma_{i}$ constitute a sheaf  on 
$\De^{1}_{U_{i}/S}$, which  is the direct image under the diagonal 
embedding $\Delta$ of 
a $\mathrm{Lie}(G_{i}, 
\Om^{1}_{U_{i}/S})$-torsor $\UGa_{i}$ above $U_{i}$.
The  right action of $p^{\ast}_{1}G_{i}$
on $\Gamma_i$ induces an additional    right action of 
$\mathrm{Lie}(G_{i},\, \Om^{1}_{U_{i}/S} )$ on $\underline{\Ga}_{i}$, 
so that $\underline{\Ga}_{i}$ is  in fact a  bitorsor under 
the sheaf of abelian groups
$\mathrm{Lie}(G_{i},\,\Om^{1}_{U_{i}/S}).$
In more 
concrete terms, 
$\underline{\Gamma}_{i}$ 
 is isomorphic to the sheaf 
 \bee
 \label{gamunder}
 \mathrm{Isom}_{ p_{0}^{\ast}\pc,\,\mathrm{pt}}
 (\epsilon p_{1}^{\ast}\,x_{i}, \,p_{0}^{\ast}\,x_{i})\:.\end{equation}
 of $1_{x_{i}}$-pointed 
 isomorphisms from $\epsilon \prb x_{i}$ to $\pra x_{i}$.
 The left and 
 right actions of the abelian sheaf 
 $\mathrm{Lie}\,(G_{i},\,\Om^{1}_{U_{i}/S})$ on 
 $\underline{\Ga}_{i}$  respectively  correspond 
 to the left ({\it resp.} right) composition with 
 $1_{x_{i}}$-pointed automorphisms of $p^{\ast}_{0}x_{i}$ (\re 
  $\epsilon p^{\ast}_{1}x_{i})$. The following lemma asserts that 
  these two actions coincide.
  \begin{lemma}
      The bitorsor 
      structure induced on  $\underline{\Ga}_{i}$ 
      by the $(p^{\ast}_{0}G_{i},\,p^{\ast}_{1}G_{i})$-bitorsor 
      structure on $\Ga_{i}$ is the obvious one determined by the 
      left action on $\UGa_{i}$ of the sheaf of abelian 
      groups
      $\mathrm{Lie}\,(G_{i},\,\Om^{1}_{U_{i}/S})$.
      \end{lemma}
      {\bf Proof:} Let  $\gamma_{i}$  be a local section of $\Ga_{i}$, 
      which we may assume to be pointed. 
      The right $\prb G_{i}$-torsor 
      structure on $\Ga_{i}$ is described in terms of its left $\pra 
      G_{i}$-torsor structure by 
       the $\Delta^{1}_{U_{i}/S}$-group isomorphism 
      \bee
      \label{gcon1}
      m_{i}: \prb G_{i} \la \pra G_{i}\end{equation}
      defined by 
      \bee
      \label{bitorchi}
      \ga_{i} g = m_{i}(g)\ga_{i}\end{equation}
      for all $g \in G_{i}$. Since $\gamma_{i}$ is pointed, 
      $m_{i}$ is in fact a connection on the 
      $U_{i}$-group $G_{i}$.
      The induced bitorsor structure on $\UGa_{i}$ is then determined 
      by the restriction of $m_{i}$ to the abelian subsheaf 
      $\De_{\ast}(\mathrm{Lie}(G_{i},\, \Om^{1}_{U_{i}/S})$ of $\prb 
      G_{i}$. It is a general fact that for any connection $m$ on an 
      $U$-group $G$, the restriction of $m$ to the corresponding 
      subsheaf $\De_{\ast}(\mathrm{Lie}(G_{},\, \Om^{1}_{U/S})$ is 
      trivial since
       the diagram 
      \bee
      \label{mufunct}
      \xymatrix{
 0 \ar[r] & \Delta_{\ast}(\mathrm{Lie}\,(G,\,\Om^{1}_{U/S}))
 \ar@{=}[d]
 \ar[r]& \prb G \ar[d]^{m}\ar[r]&  \Delta_{\ast}(G) 
 \ar@{=}[d]\\
 0 \ar[r]& \Delta_{\ast}(\mathrm{Lie}\,(G,\,\Om^{1}_{U/S})) \ar[r]
 & \pra G \ar[r] &  \Delta_{\ast}(G) 
 }\end{equation}
above $\Delta^{1}_{U_{i}/S}$ commutes  by functoriality of the 
exact sequence  of \cite{cdf} proposition 
       2.2. 
\begin{flushright}
    $\Box$
    \end{flushright}
 
   \bigskip
   
    The previous discussion makes it clear that the 
    $\mathrm{Lie}\,(G_{i},\,\Om^{1}_{U_{i}/S})$-torsor $\underline{\Ga}_{i}$
   carries most of the information embodied in the $(\pra G_{i}, \, 
   \prb G_{i})$-bitorsor $\Ga_{i}$, since the latter's underlying 
   left torsor 
   structure is obtained from $\Delta_{\ast}(\underline{\Ga}_{i})$ by 
   the extension of the structural group associated to the monomorphism
   \[\Delta_{\ast}(\mathrm{Lie}\,(G_{i},\,\Om^{1}_{U_{i}/S})) 
   \hookrightarrow \pra G_{i}\:.\] 
  In the  context of non-abelian gerbes, 
  the right multiplication by $\prb G_{i}$ provides an extra element 
   of structure, which   lives above $\Delta^{1}_{U_{i}/S}$, and 
  is  described, for a chosen section $\gamma_{i}$ of
   $\Gamma_{i}$, by  the connection $m_{i}$ \eqref{gcon1}. 
 This can be expressed more intrinsically by viewing the bitorsor 
 $\Ga_{i}$ as a left $\pra G_{i}$-torsor
 endowed with a pointed morphism
 \[ v_{i}: \Ga_{i} \la \mathrm{Isom}(\prb G_{i},\, \pra G_{i})\]
 above  $\De^1_{U_{i}/S}$,  and which is left  equivariant with respect
 to the inner conjugation morphism 
 \[ \pra G_{i} \la \mathrm{Aut}(\pra 
 G_{i})\,.\]
 Since $\xi_{i}$ \eqref{def:triv} trivializes $\Ga_{i}$ as a bitorsor,
the morphism $v_{i}$ is pointed, and it therefore 
corresponds to a $U_{i}$-morphism
 \bee
 \label{def:vi1}
 \uv_{i}: \UGa_{i} \la \mathrm{Co}(G_{i})\,,\end{equation}
 where $\mathrm{Co}(G_{i})$ is the sheaf of connections on the $U_{i}$-group 
 $G_{i}$. The full bitorsor structure of $\Ga_{i}$ is  now 
 described by the  torsor $\UGa_{i}$ above $U_{i}$, together 
 with a $U_{i}$-morphism $\uv_{i}$ which is equivariant with respect to the 
 homomorphism 
 \begin{equation}
     \label{Liei1}
     \mathrm{Lie}(i,\, \Om^{1}_{U_{i}/S}): 
     \mathrm{Lie}(G, \Om^{1}_{U_{i}/S}) \la \mathrm{Lie(Aut}(G)),\,
 \Om^{1}_{U_{i}/S})
 \end{equation}
 where $i$ is the inner 
 conjugation homomorphism from $G$ to $\mathrm{Aut}(G)$.
 Similarly,  the bitorsor  
${}^{ P_{ij}\ast\,\,}\!\Gamma_{j}$  is 
pointed above $U_{ij}$  and so reduces  to a torsor
${}^{ P_{ij}\ast\,\,}\!\underline{\Gamma}_{j}$ under the abelian group
$\mathrm{Lie}(G_{i},\,\Om^{1}_{U_{ij}/S})$,  defined  above
$U_{ij}$ and endowed with a morphism ${}^{ P_{ij}\ast\,\,}\!\uv_{j}$ 
to $\mathrm{Co}(G_{i})_{\mid U_{ij}}$. Despite the notation, this torsor
does not depend only on $\UGa_{j}$ but on the full bitorsor structure 
of $\Ga_{j}$ ({\it i.e.}  on the pair $(\UGa_{j},\, \uv_{j\, \mid 
U_{ij}}))$. The bitorsor morphisms 
\eqref{defdij1} now corresponds  to  morphisms  of
$\mathrm{Lie}(G_{i},\, \Om^{1}_{U_{ij}/S})$-torsors
\bee
  \label{dired}\xymatrix{
\underline{\Ga}_{i} \ar[r]^(.4){\underline{g}_{ij}}&
{}^{P_{ij}\ast\,\,}\!\underline{\Gamma}_{j}\:.}\end{equation}
above $U_{ij}$, which are compatible with the morphisms $\uv_{i}$ and 
${}^{ P_{ij}\ast\,\,}\!\uv_{j}$, and for which 
for which the diagrams
  \bee
\label{compudij3}
\xymatrix{
\UGa_{i} \ar[rrr]^{\underline{g}_{ij}}\ar[d]_{\underline{g}_{ik}}&&&
{}^{P_{ij}\ast\,}\!\UGa_{j}
\ar[d]^{{}^{P_{ij}\ast\,\,}\! \underline{g}_{jk} }
\\
{}^{P_{ik}\ast\,\,}\! \UGa_{k} 
\ar[rrr]_<<<<<<<<<<<<<<<<<<<{{}^{\Psi_{ijk}\ast\,}\!(1_{\UGa_{k}})}
& &&
{}^{(P_{ij} \wedge P_{jk})\ast\,\,}\!
\UGa_{k}
}\end{equation}
above $U_{ijk}$ commute. 

\bigskip

The bitorsor isomorphisms $\zeta_{i}$ \eqref{equivzi} associated to an 
equivalence between a pair of connections may similarly be described 
by a family of $\mathrm{Lie}(G_{i},\, \Om^{1}_{U_{i/S}})$-torsor isomorphisms
\[\underline{\zeta}_{i}: \UGa_{i}^{\epsilon'} \la 
\UGa_{i}^{\epsilon}\]
compatible with the corresponding morphisms $\uv_{i^{\epsilon}}$ and 
$\uv_{i^{\epsilon'}}$, and for which the diagram
  \bee \label{comuzeta}
     \xymatrix{
\UGa_{i}^{\epsilon'} 
\ar@<-2pt>[r]^(.4){\underline{g}_{ij}^{\epsilon'}} 
\ar[d]_{\underline{\zeta}_{i}}
& {}^{ P_{ij}\ast\,\,}\!\UGa_{\,j}^{\epsilon'}
\ar[d]^{{}^{ P_{ij}\ast\,\,}\!\underline{\zeta}_{\,j}}\\
\UGa_{i}^{\epsilon}
\ar@<-2pt>[r]_(.4){\underline{g}_{ij}^{\epsilon}} 
&{}^{ P_{ij}\ast\,\,}\!\UGa_{\,j}^{\epsilon}
}
\end{equation}
induced by 
\eqref{comzeta} commutes.
 
\begin{remark}
    {\rm   
        For $G= G_{m,\, \mathbb{C}}$ and $\epsilon$ a conncection on 
        the abelian $G_{m,\, \mathbb{C}}$-gerbe $\pc$,  the  expression 
  \eqref{gamunder}  is   the 
 $\Om^{1}_{U_{i}/\sp(\mathbb{C})}$-torsor 
 $Co(x_{i})$ associated by Brylinski (\cite{bry} definition 5.3.1) 
to an object $x_{i}$ in 
 $\pc_{U_{i}}$. The  $\mathrm{Lie}\,(G_{i},\,\Om^{1}_{U_{i}/S})$-torsor 
 $\underline{\Gamma}_{i}$, together with its morphism $\uv_{i}\,$ 
 \eqref{def:vi1}, 
 should therefore be viewed as
 a non-abelian generalization of Brylinski's notion of a
 connective structure.
}  \end{remark}

 \bigskip

\subsection{}
\setcounter{equation}{0}%
 We will  now carry out a 
   parallel discussion for the fake curvature $\kappa$  and 
   the curving  $K$ \eqref{defkK}.
   The fake curvature 
    is described by the family $\pra G_{i}$-bitorsors $\Delta_{i}$ 
    defined above $\Delta^{2}_{U_{i}/S}$ by
    \bee
    \label{dinonu}
    \Delta_{i} := \mathrm{Isom}(\kappa(\pra x_{i}),\, 
    \pra x_{i})\end{equation}
   and  endowed with a canonical trivialization above the degenerate 
   subsimplex  $s\Delta^{2}_{U_{i}/S}
     \hookrightarrow \Delta^{2}_{U_{i}/S}$.
     For any pair of indices $(i,j)$, we define the $\pra 
     G_{i}$-bitorsor  ${}^{ P_{ij}\,}\Delta_{j}$ above $\De^{2}_{U_{ij}/S}$ by 
     \[{}^{ P_{ij}\,}\Delta_{j} := 
   p_{0}^{\ast}P_{ij}  \,\wedge \,\Delta_{j}\wedge
(p_{0}^{\ast}P_{ij})^{0}\/.\]
     For varying $i$, we have
    a family of    isomorphisms
 \bee
     \label{ardij3}
\xymatrix{
\Delta_{i} \ar[r]^(.25){d_{\,ij}}&
p_{0}^{\ast}P_{ij}\, \wedge \,\De_{j}\wedge
(p_{0}^{\ast}P_{ij})^{0}}
\end{equation}
above $\Delta^{2}_{U_{ij}/S}$
which are pointed above $s\De^{2}_{U_{ij}}$
(or equivalently of isomorphisms
 \begin{equation}
\label{ardij6}
\xymatrix{\Delta_{i} \wedge  \pra P_{ij} \ar[rr]^{\tilde{d}_{ij}}&&
\pra P_{ij} \wedge \Delta_{j}}\:),
\end{equation}
 and 
for which the following analogue
     of diagram \eqref{compdij3} commutes:
        \bee
\label{compdij5}
\xymatrix{
\De_{i} \ar[rrr]^{d_{ij}}\ar[d]_{d_{ik}}&&&
{}^{ P_{ij}\,}\!\Delta_{j}
\ar[d]^{ {}^{ P_{ij}\,}\!d_{jk}}
\\
{}^{ P_{ik}\,}\!\Delta_{k}
\ar[rrr]_<<<<<<<<<<<<<<<<<<<{{}^{\Psi_{ijk}\,}\!(1_{\De_{k})}}
& && 
{}^{ P_{ij}\wedge P_{jk}\,}\!\Delta_{k}\:.
}\end{equation}
 It is a diagram of pointed $\pra G_{i}$-bitorsors above $\De^{2}_{U_{ijk}}$.
 To the pointed bitorsor $\De_{i}$ corresponds a  
 left  $\mathrm{Lie}(G_{i},\, 
    \Om^{2}_{U_{i}/S})$-torsor $\underline{\Delta}_{i}$ above 
    $U_{i}$,  such that 
    \bee
    \label{def:di}\underline{\Delta}_{i}:=
     \mathrm{Isom}_{\text{pt}}(\kappa(\pra x_{i}),\, 
    \pra x_{i})\:.\end{equation} It stands in 
    the same relation to the bitorsors $\Delta_{i}$ as does the torsor
    $\underline{\Ga}_{i}$\eqref{gamunder}  to the 
    bitorsor $\Ga_{i}$  \eqref{def:gai}. In particular, the 
    $\mathrm{Lie}(G_{i},\, 
    \Om^{2}_{U_{i}/S})$-bitorsor structure on $\UDe_{i}$ is the 
    obvious one, induced from the left action of the abelian group
     $\mathrm{Lie}(G_{i},\, 
    \Om^{2}_{U_{i}/S})$. The sheaf $\De_{i}$ may be recovered from
    the 
    direct image $\Delta_{\ast}(\underline{\Delta}_{i})$ of 
    $\underline{\Delta}_{i}$ under the embedding $U_{i} 
    \hookrightarrow \De^{2}_{U_{I}/S}$ 
    as the  torsor above $\Delta^{2}_{U_{i}/S}$ obtained 
    by extension of 
    the structural group  from 
    $\mathrm{Lie}(G_{i}, \Om^2_{U_{i}/S})$ to $\pra G_{i}$. The 
    right $\pra G_{i}$-torsor structure of  $\De_{i}$ is described by 
    a
    $\De^{2}_{U_{i}/S}$-morphism
  \bee
  \label{def:uui1}
  u_{i}: \De_{i} \la \mathrm{Aut}(\pra G_{i})
  \end{equation}
 which is pointed above $s\De^{2}_{U_{i}/S}$,
 and  equivariant with respect to the conjugation map for 
 $\pra G_{i}$. 
 This corresponds
to a $U_{i}$-morphism
\begin{equation}
    \uu_{i}:\UDe_{i} \la \mathrm{Lie (Aut}(G) ,\,\Om^{2}_{U_{i}/S})
    \end{equation}
    which is equivariant with respect to the induced homomorphism
   \[ \mathrm{Lie}(i,\, \Om^{2}_{X/S}): \mathrm{Lie}(G,\,
   \Om^{2}_{X/S}) \la \mathrm{Lie\, (Aut}(G),\, \Om^{2}_{X/S})\,.\]
    The  underlying  
    $\mathrm{Lie}(G_{i},\, 
    \Om^{2}_{U_{ij}/S})$-torsor above $U_{ij}$
    associated in the same way to the pointed
    bitorsor ${}^{ P_{ij}\,}\!\De_{j}$
    will be denoted ${}^{ P_{ij}\,}\!\underline{\De}_{j}$.
It once again depends on the full pair $(\UDe_j,\,\uu_{j\,\mid U_{ij}})$.
   Since the morphisms   $d_{ij}$ \eqref{ardij3} are pointed above 
$s\De^{2}_{U_{ij}/S}$, 
they induce
 $\mathrm{Lie}(G_{i}, \Om^2_{U_{ij}/S})$-torsor morphisms
\bee
     \label{ardij4}
\xymatrix{
\underline{\Delta}_{i} \ar[r]^(.4){\underline{d}_{\,ij}}&
{}^{ P_{ij}\,}\!\underline{\Delta}_{j}}\end{equation}
above $U_{ij}$ compatible with the morphisms $\uu_{i}$, and for which 
the diagram of $\mathrm{Lie}(G_i,\,  \Om^{2}_{U_{ijk}/S})$-torsors
\begin{equation}
    \label{compudij5}
\xymatrix{
\UDe_{i} \ar[rrr]^{\underline{d}_{ij}}\ar[d]_{\underline{d}_{ik}}&&&
{}^{ P_{ij}\,}\!\UDe_{j}
\ar[d]^{ {}^{ P_{ij}\,}\!\underline{d}_{jk}}
\\
{}^{ P_{ik}\,}\!\UDe_{k}
\ar[rrr]_<<<<<<<<<<<<<<<<<<<{{}^{\Psi_{ijk}\,}\!(1_{\UDe_{k})}}
& && 
{}^{ P_{ij}\wedge P_{jk}\,}\!\UDe_{k}\:.
}\end{equation}
induced by \eqref{compdij5} commutes. The giving of a fake 
curvature arrow $\kappa$ on the  locally trivialized gerbe $(\pc,\, 
(x_{i})_{i \in I})$ is equivalent to that of the pairs
 $(\UDe_{i},\, \uu_{i})$, and 
of  morphisms $\underline{d}_{ij}$ compatible with $\uu_i$ and
${}^{P_{ij}\,}\!\uu_j$
and  for which 
the diagrams \eqref{compudij5} commute.

\bigskip

In order to describe in similar terms the curving morphism $K$,
we introduce some additional notation. 
For any  sheaf of groups $G$ on an $S$-scheme $U$  and $(\pra G,\, \prb 
G)$-bitorsor $\Gamma$ 
  on $\Delta^{1}_{U/S}$, let us  define the 
induced $(\pra G,\, \pra G)$-bitorsor $\delta^1(\Gamma)$ above 
$\Delta^{2} _{U/S}$  by  the formula
\bee
\label{defd2gai}
\delta^{1}(\Ga) := 
\prab \Ga \wedge \prbc \Ga \wedge
(\prac \Ga)^0\:.\end{equation}
When $\Ga$ is trivialized as a bitorsor above $U$,
so is $\delta^1(\Gamma)$ above 
$s\De^{2}_{U/S}$, and the latter therefore corresponds  to a 
$\mathrm{Lie}(G,\, 
\Om^{2}_{U/S})$-torsor above $U$ which we will denote by 
$\de^{1}(\UGa)$ even though it doesn't  depend only on 
$\UGa$ but rather on the full bitorsor structure of $\Ga$, {\it i.e.} 
on 
the pair $(\UGa,\, \uv)$ with $\uv$ as in \eqref{def:vi1}.
The induced  right 
$\pra G$-torsor structure on $\de^{1}(\Ga)$ is described by a 
$\mathrm{Lie}(G,\, 
\Om^{2}_{U_{i}/S})$-equivariant 
morphism
\[ \de^{1}(\uv): \de^{1}\UGa \la 
\mathrm{Lie}(\mathrm{Aut}(G) ,\,
\Om^{2}_{U/S})\:.\]

\bigskip

A curving 2-arrow $\widetilde{K}$ \eqref{defkK}  is  described by a 
family of pointed bitorsor isomorphisms 
\begin{equation}
\label{def:Ki0}
\xymatrix{\prab \Ga_i\, \wedge\, \prbc \Ga_i
\ar[rr]^{K_i} &&  \Delta_i\,  \wedge
 \, \prac \Ga}\:.\end{equation}
or by the corresponding family of pointed $\pra G_{i}$-bitorsor isomorphisms
 \bee
\label{def:Ki1}
\xymatrix{ \delta^{1}(\Ga_{i}) 
\ar[rr]^{\widetilde{K}_{i}}&& 
\Delta_{i} \:.}\end{equation}
To $\tilde{K}_{i}$ is associated  a $U_{i}$-morphism of 
$\mathrm{Lie}(G_{i},\, \Om^{2}_{U_{i}/S})$-torsors
\bee
\label{def:UKi1}
\xymatrix{ \delta^{1}(\UGa_{i}) 
\ar[rr]^{\underline{\widetilde{K}}_{i}}&& 
\UDe_{i} }\end{equation}
compatible with the morphisms $\de^{1}(\underline{v}_{i})$ and $\uu_{i}$.
The compatibility condition \eqref{cyl} for the various
 morphisms $\widetilde{K}_{i}$ 
 is expressed by the commutativity of the 
following diagram 
(in which the arrow 
 \[\delta^{1}({}^{ P_{ij}\ast\,\,}\!\underline{\Gamma}_{j}) \simeq
{}^{ P_{ij}\,}\!(\delta^{1}\underline{\Gamma}_{j})\]
 is
the canonical isomorphism)
 \bee
   \label{comp:kikj}\xymatrix{
\delta^{1}
(\Ga_{i}) \ar[rrr]^{\widetilde{K}_{i}} 
\ar[d]_{\delta^{1}(g_{ij})} &&&
\Delta_{i}\ar[d]^{d_{ij}}\\
\delta^{1}({}^{ P_{ij}\ast\,\,}\!\Gamma_{j})\ar[r]^(.5){\sim}
&{}^{ P_{ij}\,}\!(\delta^{1}\Gamma_{j})
\ar[rr]_{{}^{P_{ij}\,} \!(\widetilde{K}_{j})}&&
{}^{P_{ij}\,}\!\Delta_{j}\:.} \end{equation}
The left-hand vertical   arrow in this diagram is induced  by the 
morphism
\[ \delta^{1}(g_{ij}) := \prab g_{ij}
\,\wedge \, \prbc g_{ij}\,\wedge \,\prac 
(g_{ij}^{-1})^{0}\:.\]
with
\[\xymatrix{
(\Ga_{i})^{0} \ar[rr]^(.35){(g_{ij}^{-1})^{0}}
&&\prb P_{ij}
\, \wedge \,
(\Ga_{j})^{0} \,\wedge \,(\pra P_{ij})^{0}}\]
 the ``opposite''  \eqref{def:po}
 of the arrow $g_{ij}^{-1}$ \eqref{defdij1}. It is equivalent to 
 require the commutativity of the induced diagram
 above $U_{ij}$ whose vertices are built from the $\UGa_{i}$ and 
 $\UDe_{i}$.

\bigskip

For any $(\pra G,\,\prb 
G)$-bitorsor $\Gamma$ above $\Delta^1_{U/S}$ and any $(\pra 
G,\, \pra 
G)$-bitorsor $\Delta$ above $\Delta^2_{U/S}$, we define a  bitorsor
$ {}^{\Ga\,}\!\De$ above $\dec$ by
\[ {}^{\Ga\,}\!\De:= \Ga_{01} \, \De_{123}\, \Ga_{01}^{-1}\]
where in order to emphasize the analogy with the combinatorial formulas of 
\cite{cdf}
we have preferred to denote by $\Gamma^{-1}$ rather than by
 $\Gamma^0$ the opposite bitorsor of a bitorsor $\Gamma$.
We define  the induced $(\pra G,\, \pra 
G)$-bitorsor $\delta^{2}_{i_{\Gamma}}(\Delta)$ above  $\Delta^3_{U/S}$
 by
\begin{equation}
    \label{def:de2de}
\delta^{2}_{i_{\Ga}}(\Delta) := 
 ({}^{\Ga_{01}\,}\!\De_{123})\, \, \Delta_{013}\, \,
\Delta _{023}^{-1}\,\, \Delta_{012}^{-1}\:.
\end{equation}
We will  abbreviate this expression  to
$\delta^{2}_{\Ga}(\De)$. This construction is functorial in the pair
$(\Ga,\, \Delta)$. 
When the bitorsors $\Gamma$ and $\Delta$ are respectively pointed above 
$U$ and 
$s\Delta^{2}_{U/S}$, then  $\delta^{2}_{i_{\Ga}}(\Delta)$ is pointed 
above $s\Delta^3_{U/S}$ 
and therefore determines a $\mathrm{Lie}(G, \, \Om^{3}_{U/S})$-torsor
above  $U$, which we will denote by $\delta^{2}_{i_{\UGa}}(\UDe)$.
Similarly, for any $(\pra G,\, \pra G)$-bitorsor 
$\Delta$ above $\De^{3}_{U/S}$,
we define the  $(\pra G,\, \pra G)$-bitorsor 
$\delta^{3}_{i_{\Ga}}(\Delta)$ above $\ded$
by the 
analogous  formula
\begin{equation}
    \label{def:de3de}\delta^{3}_{i_{\Ga}}(\Delta):= 
 {}^{\Ga_{01}\,}\De_{1234}\,\, 
\De_{0134}\, \,\De_{0123}\,\,
\De_{0234}^{-1}\,\,
\De_{0124}^{-1}\,.\end{equation}
Once more, when $\Ga$ and $\De$ 
are pointed above their degenerate subsimplexes,
so is the bitorsor
$\delta^{3}_{i_{\Ga}}(\Delta)$ above  $s\De^{4}_{U_{i}/S}$ and it 
therefore determines a $\mathrm{Lie}(G,\, \Om^{4}_{U/S})$-torsor
$\delta^{3}_{i_{\UGa}}(\UDe)$
above $U$. Both $\delta^{2}_{i_{\UGa}}(\UDe)$ and 
$\delta^{3}_{i_{\UGa}}(\UDe)$ are endowed with morphisms 
analogous to the morphism $\uu_{i}$ \eqref{def:uui1}  with values in the 
corresponding sheaves of $\mathrm{Aut}(G)$-valued forms, and which express
the full 
bitorsor structure of  $\delta^{2}_{i_{\Ga}}(\Delta)$ 
and $\delta^{3}_{i_{\Ga}}(\Delta)$.
By the standard  combinatorial  argument of \cite{AK:Bianchi}
theorem 2, carried 
out 
in the present geometric context,
 $\delta^{2}_{i_{\Gamma}}\,(\delta^{1}\Gamma)$ is canonically  
trivialized as a $(\pra G,\, \pra G)$-bitorsor
 by a global  section $t_{\Gamma}$ on $\dec$, whose construction is 
 functorial in $\Ga$. When $\Ga$ is pointed above 
 $s\De^{3}_{U_{i}/S}$, it is compatible with this pointing.
 We denote by
  $t_{\UGa}$ the corresponding
  canonical section  of 
  $\delta^{2}_{i_{\UGa}}\,(\delta^{1}\UGa)$.

\bigskip

By Morita theory, the entire global discussion of section \ref{sec:curv-Bian}
may now expressed  in the language of bitorsors. 
The restriction $\Om_i$ of the 3-curvature $\Om$
 \eqref{defomeg}
 to the open sets $U_i$ above which
$\pc$  trivializes   is expressed as a bitorsor isomorphism\footnote{
We have lifted the indices $i$ in order to make way for  the additional
ones. In the sequel, we will indiscriminately raise or lower indices according
to typographical convenience.}
\[\Delta^i_{012}\,\Delta^i_{023} \la
 {}^{\Gamma^i_{01}} \Delta^i_{123}\,
 \,\Delta^i_{013}\]
or equivalently as  a trivialization $\UOm_{i}$ of the torsor
 $\delta^2_{\underline{\Gamma}_i}(\underline{\Delta}_i)$. The arrow 
$\mu_{01}(\kappa_{123})$ in the
 cube \eqref{cube}   is locally represented by the bitorsor 
$ {}^{\underline{\Gamma}^i_{01}\,}\! \underline{\Delta}^i_{123}$,
so that   the right-hand
2-arrow $M_{01}(\kappa_{123})$ in the cube  corresponds to a bitorsor
 isomorphism induced by the canonical isomorphism $P^0 \wedge P \la T_{G_{i}}$,
  which  can   be neglected here. Since the curving data has been
 locally described
 by the morphism  $\widetilde{K}_{i}$ 
\eqref{def:Ki1}, it   follows  that  $\UOm_i$ is  represented by the image
\begin{equation}
\label{def:omegi}
\UOm_{i} = (\delta_{\UGa_{i}}^{2})(\underline{\widetilde{K}}_{i})(t_{\UGa_i})\end{equation}
of the canonical section 
$t_{\underline{\Gamma}_i}$ 
of  $\delta^{2}_{i_{\UGa_{i}}}\,(\delta^{1}\UGa_{i})$
under the bitorsor isomorphism 
\[\xymatrix{
(\delta_{\underline{\Ga}_{i}}^{2}\delta^1)(\underline{\Ga}_{i}) 
\ar[rr]^{\delta_{\underline{\Ga}_{i}}^{2}(\underline{\widetilde{K}}_{i})}&&
\delta_{\underline{\Ga}_{i}}^{2}(\underline{\Delta}_{i})}\,.\]
Consider the composite morphism
\bee
\label{ijcompat}\delta^{2}_{\Gamma_{i}} \Delta_{i} \stackrel{\sim}{\la} 
\delta^2_{\Ga_{i}}({}^{P_{ij}\,}\!\Delta_{j})
\stackrel{\sim}{\la} 
\delta^2_{({}^{P_{ij}\,\ast\,}\!\Gamma_{j})}({}^{P_{ij}\,}\!\Delta_{j})  
\stackrel{\sim}{\la}{}^{P_{ij}\,}\!(\delta^2_{\Ga_{j}}\Delta_{j})\,  
\end{equation}
above $U_{ij}$, (where the  first and second arrow are respectively
induced by the morphisms $d_{ij}$ and $g_{ij}$, and the third
one is the  canonical isomorphism).
If we   apply the functor
$\delta^2_{\UGa_{i}}$
to diagram \eqref{comp:kikj}, and take into account the functoriality
of the canonical section $t_{\underline{\Ga}_{i}}$, we find that the 
composite map  \eqref{ijcompat} sends the restriction 
to $U_{ij}$
of the section $\UOm_{i}$ of $\delta^2_{\UGa_{i}} \UDe_{i}$
to the section ${}^{P_{ij}\,}\UOm_{j}$ of 
${}^{P_{ij}\,}\!(\delta^2_{\UGa_{j}}\UDe_{j})$.
 This relation between  the sections $\UOm_{i}$ and $\UOm_{j}$ expresses 
the fact that 
the various  $\Om_{i}$  
glue to a  global arrow $\Om$ \eqref{defomeg} 
in the fiber on $X$ of the stack of 
  $\mathrm{(Lie}(\pac))$-valued 3-forms.

\bigskip

The following lemma gives us a very compact description of the image
of a bitorsor
under the composite functor $\de^3_{\UGa^i}\, \de^{2}_{\UGa^i}$:

\vspace{.3cm}
\begin{lemma}
\label{caniso}
There exists a canonical isomorphism
 \begin{equation}
 \label{highBian}
\de^3_{\UGa^i}\, \de^{2}_{\UGa^i}(\UDe_i^{-1})
\stackrel{\sim}{\la}
 {}^{\UDe^i_{012}\, \UGa^i_{02}\,}\!\UDe^i_{234}\, \, 
 {}^{\UGa^i_{01}\, \UGa^i_{12}\,}\!(\UDe^i_{234})^{-1} 
\end{equation}
\end{lemma}

\bigskip

\noindent {\bf Proof:} The following canonical isomorphisms $\alpha$ and $\beta$
are  obtained by the  insertion  for each term $\de^2_{\Ga}(\De)$ of  its value
 according
 to  definition \eqref{def:de2de} and cancellation of  the appropriate factors:
\begin{multline*}
\de^2_{\Ga}(\De^{-1})_{0123}\,\,\, {}^{{}^{\Ga_{01}} \De_{123}}
(\de^2_{\Ga}(\De^{-1}))_{0134}\,\,\,
 {}^{\Ga_{01}}(\de^2_{\Ga}(\De^{-1}))_{1234}\,\,\,{}^{\Ga_{01}\Ga_{12}}
\De_{234}\,\,\,\,\de^2_{\Ga}(\De)_{0124}\,\,\,\,\De_{012}\,\,\,\De_{024}
\stackrel{\alpha}{\la} \\
\la  \De_{012}\,\,\, \De_{023}\,\,
\, \De_{034}
\end{multline*}
\[
\xymatrix{
{}^{\De_{012}}(\de_{\Ga}^2(\De^{-1})_{0234})\,\,\De_{012}\,\,
{}^{\Ga_{02}}\De_{234}\,\, \De_{024}\, \ar[r]^(.65){\beta}
&  \De_{012}\,\,\, \De_{023}\,\,
\, \De_{034}
}
\]
We  then  cancel  the right-hand  factor $\De_{024}$ in both  source
 and target
 of  the composite arrow $\beta^{-1} \, \alpha$. In the context
of the  lemma, one may  then  permute  factors in the source or target of the
 induced isomorphism when they are  indexed by
 a pair 
of common  variables in the set of integers  $[0,4]$. This allows one 
 to permute the factors
${}^{\Ga_{01}\Ga_{12}}
\De_{234}$ and $\de^2_{\Ga}(\De)_{0124}$ in the source, and to neglect the
conjugation action of 
${}^{\Ga_{01}} \De_{123}$ on $\de^2_{\Ga}(\De^{-1}))_{0134}$,
and of  $\De_{012}$ on $\de_{\Ga}^2(\De^{-1})_{0234}$. An elementary
 rearrangement of the factors then yields  the result.
\begin{flushright}
$\Box$
\end{flushright}

The  higher Bianchi identity, 
in its incarnation  \eqref{bianchicat}, now  asserts that
 the diagram
\begin{equation}
\label{biantrian}
\xymatrix{
\de^3_{\UGa^i}\, \de^{2}_{\UGa^i}(\UDe_i^{-1}) \:\: 
\ar[rr]^(.4){\sim}\ar[rd]_{\de^3_{\UGa^i}(\underline{\Om}_i)} &&
\:\:\:\: {}^{\UDe^i_{012}\, \UGa^i_{02}\,}\!\UDe^i_{234}\, \, 
 {}^{\UGa^i_{01}\, \UGa^i_{12}\,}\!(\UDe^i_{234})^{-1} \ar[ld]
\\& T_G&
}
\end{equation}
\medskip
\noindent     is commutative, 
 where  the horizontal map is the isomorphism
\eqref{highBian} and the  the right-hand one is  the composite
\[
\xymatrix{
{}^{\UDe^i_{012}\, \UGa^i_{02}\,}\!\UDe^i_{234}\, \, 
 {}^{\UGa^i_{01}\, \UGa^i_{12}\,}\!(\UDe^i_{234})^{-1}
\ar[r] &{}^{\UDe^i_{012}\, \UGa^i_{02}\,}\!\UDe^i_{234}\, \, 
{}^{\UDe^i_{012}\, \UGa^i_{02}\,}\!(\UDe^i_{234})^{-1} \ar[r]^(.77){\sim} & T_G
}
 \]
induced by the morphism  $K_i$ \eqref{def:Ki0}.

 \bigskip

\subsection{}
\setcounter{equation}{0}%
\label{partly-3}
In section \ref{sec:fulldec}, we will restate the previous discussion in 
purely cocyclic terms, once trivializations have been chosen for 
all the  torsors which occur. Here we 
will begin this process 
by   choosing trivializations
of the torsors $\Ga_{i}$ and 
$\De_{i}$, but without doing so for the bitorsor cocycles $P_{ij}$. In 
this case we will see in proposition \ref{prop:omi}
 that a gerbe $\pc$ with a curving pair 
$(\epsilon,\, K)$ and local objects $x_{i}$
can  be described in  geometric terms, 
involving only differential forms and the remaining non trivialized 
bitorsors $P_{ij}$. The same holds for 
 the associated fake curvature $\kappa$  and the 3-curvature $\Omega$.

\bigskip

We  refine,
as  we did earlier in a related context (\ref{upath}),  the 
given open cover $\mathcal{U}$ of $X$ to one for which  the 
torsors $\UGa_{i}$ and $\UDe_{i}$    have global sections 
$\uga_{i}$
and $\ude_{i}$. By construction,  the corresponding section    $\ga_{i}$
of $\Ga_i$ is a pointed  arrow
\bee
\label{defchii}
\xymatrix{\epsilon (\prb x_{i}) \ar[r]^{\ga_{i}} & \pra 
x_{i}}\end{equation}
in $\pra \pc_{\mid\, U_{i}}$, which trivializes 
$\Ga_{i}$ as a left 
$\pra G_{i}$-torsor. above $\De^{1}_{U_{i}/S}$.
To the quasi-inverse $\epsilon^{-1}$ of $\epsilon$ is associated
the corresponding arrow 
\begin{equation}
\label{defchii1}
\ga'_i : \epsilon^{-1}(\pra x_i)   \la \prb x_i
\end{equation}
inverse to 
\[
\xymatrix{
  \prb x_i \ar[r]^(.4){\sim} & \epsilon^{-1}\, \epsilon\,\, \prb x_i
  \ar[rr]^{\epsilon^{-1}(\ga_i)} && \epsilon^{-1}\, \pra x_i\,.
}
\]
Let us now  pass  from  \eqref{defdij1} to 
the corresponding
$(p^{\ast}_{0}(G_{i}),\,p^{\ast}_{1}(G_{j}))$-bitorsor  isomorphism 
\bee\label{eij0}\xymatrix{
\Ga_{i} \wedge^{p^{\ast}_{1}G_{i}} p^{\ast}_{1}(P_{ij}) 
\ar[rr]^{\tilde{g}_{ij}}&& p_{0}^{\ast}P_{ij} \,
\wedge^{p_{0}^{\ast}G_{j}}\,
\Gamma_{j}
}\end{equation}
above  $\Delta^1_{U_{ij}/S}$. The  pointed  sections 
$\ga_{i}$ and $\ga_{j}$ of $\Ga_{i}$   and  $\Ga_{j}$ give  us another
description of the map $\tilde{g}_{ij}$, in terms of
the $\Delta^{1}_{U_{ij}/S}$-morphism 
\bee
\label{def:eij}
\epsilon_{ij}: p^{\ast}_{1}P_{ij} \la  
p^{\ast}_{0}P_{ij}\end{equation}
defined 
by the equation
\[ \tilde{g}_{ij}(\ga_{i} \wedge \phi) = 
\epsilon_{ij}(\phi)\wedge\ga_{j}\]
for all sections $\phi \in \prb P_{ij}$.
It is readily verified that $\epsilon_{ij}$ is a connection 
as defined in section \ref{sec:contors}    on the 
bitorsor $P_{ij}$ under the $U_{ij}$-groups $G_{i}$ and $G_{j}$
(endowed with their connections $m_{i}$ and $m_{j}$ \eqref{gcon1}).

         \bigskip                     
                     
        The commutativity of diagram \eqref{compudij3}
is now equivalent to that of the diagrams
\bee
\label{hordiag}\xymatrix{\prb P_{ik} \ar[rr]^{\epsilon_{ik}} \ar[d]_{\prb 
\Psi_{ijk}}& &
\pra P_{ik} 
\ar[d]^{\pra \Psi_{ijk}} \\
\prb P_{ij} \wedge \prb P_{jk} \ar[rr]_{\epsilon_{ij}\wedge 
\epsilon_{jk}} && \pra P_{ij} \wedge \pra P_{jk}\:.}\end{equation}
 We may now interpret the bitorsor 
cocycle morphism $\Psi_{ijk}$ \eqref{1bcoc}  as a section of the bitorsor
\bee
\label{chechp}
d^1(P)_{ijk}:= P_{ij} \wedge P_{jk} \wedge 
P_{ik}^{0}
\end{equation}
\label{horpsi}
above $U_{ijk}$,  with its induced connection \[\epsilon_{ijk}:= 
\epsilon_{ij} \wedge \epsilon_{jk}\wedge (\epsilon_{ik}^{-1})^{0}\:. \]
 The commutativity 
of diagram \eqref{hordiag} is equivalent to the equation
\begin{equation}
\label{horpsi-1}
\epsilon_{ijk}(\prb \Psi_{ijk}) = \pra \Psi_{ijk}\,,
\end{equation}
 in other words to 
the assertion that the section   $\Psi_{ijk}$ of $d^1(P)$
is horizontal with respect to
the connnection $\epsilon_{ijk}$.

                    \bigskip
                    
        Similarly, the choice of a  section $\ude_{i}$ of $\UDe_{i}$
        corresponds to that of a pointed arrow
        \bee
        \label{ardei}
        \xymatrix{
        \ka(\pra x_{i}) \ar[r]^{\de_{i}} & \pra x_{i}}
    \end{equation}
together with the corresponding arrow 
 \bee
        \label{ardeiprime}
        \xymatrix{
        \ka^{-1}(\pra x_{i}) \ar[r]^(.55){\de'_{i}} & \pra x_{i}}
    \end{equation}
associated to the quasi-inverse $\kappa^{-1}$ of $\kappa$.
 The  automorphism 
        $\nu_{i}: \pra G_{i} \la \pra G_{i}$ above $\Delta^{2}_{U_{i}/S}$
       which describes in terms of $\de_{i}$ 
       the full bitorsor structure  of $\Delta_{i}$
       is given explicitly by the pointed automorphism $\nu_{i}$  defined by 
        \bee
\label{bitordi}
\delta_{i}\, g = \nu_{i}(g)\,\delta_{i}\end{equation}
 for all $g \in \pra G_{i}$. It corresponds to a 2-form
$\underline{\nu}_{i}$, element of  
$\mathrm{Lie}(\mathrm{Aut}(G_{i}), \Om^{2}_{U_{i}/S})$ and 
bitorsor 
$\De_{i}$ is entirely described in Lie-theoretic terms 
by the 2- form  $\underline{\nu}_{i}$.

\bigskip

The morphism 
$\underline{d}_{ij}$ \eqref{ardij4} corresponds  to a pointed $(\nu_{i},\, 
\nu_{j})$-equivariant  
isomorphism 
\bee\label{def:kij1} \kappa_{ij}:\pra P_{ij} \la \pra P_{ij}\end{equation}
defined in terms of the arrow $\tilde{d}_{ij}$ \eqref{ardij6} by the formula
 by 
 \[\tilde{d}_{ij}(\delta_{i} \wedge  \phi) = \kappa_{ij}( \phi) \wedge 
 \delta_{j}\]
 for all sections  $\phi \in \pra P_{ij}$.
The morphism $\kappa_{ij}$ is the analogue for the fake curvature 
$\kappa$ of the morphism
$\epsilon_{ij}$ \eqref{def:eij} associated to a connnection 
$\epsilon$ on $\pc$, and the commutativity of  
diagram \eqref{compdij5} now asserts that
the diagrams of bitorsors
 \bee
\label{hordiag1}\xymatrix{\pra P_{ik} \ar[rr]^{\ka_{ik}} \ar[d]_{\pra 
\Psi_{ijk}}& &
\pra P_{ik} 
\ar[d]^{\pra \Psi_{ijk}} \\
\pra P_{ij} \wedge \pra P_{jk} \ar[rr]_{\ka_{ij}\wedge 
\ka_{jk}} && \pra P_{ij} \wedge \pra P_{jk}}\end{equation}
commutes.
To the  family of isomorphisms $\kappa_{ij}$ we 
associate the element  $\kappa_{ijk}$ in
$ \mathrm{Lie}((d^{1}P)^{\mathrm{ad}}_{ijk},\, \Om^{2}_{U_{ijk}/S})$
  defined by
\[\kappa_{ijk}:= \kappa_{ij} \wedge \kappa_{jk} \wedge 
(\kappa_{ik}^{-1})^{0}\,,\]
 with the \v{C}ech differential $(d^{1}P)_{ijk}$
of the $(G_{i},\, G_{j})$-bitorsor $P_{ij}$
   defined 
by \eqref{chechp}.
The commutativity of diagram \eqref{hordiag1} is equivalent to the 
equation 
\begin{equation}
    \label{hordiag2}
    [\kappa_{ijk},\, \Psi_{ijk}] = 0
    \end{equation}
where the bracket is the pairing 
\[\mathrm{Lie}((d^{1}P)^{\mathrm{ad}}, \Om^{2}_{U_{ijk}}) \times d^{1}P \la 
\mathrm{Lie}(d^{1}P,\, \Om^{2}_{U_{ijk}})\]
associated  as in \eqref{lie1} to 
the canonical  action 
$(d^{1}P)^{\mathrm{ad}} \times d^{1}P \la d^{1}P$ of the group 
$(d^{1}P)^{\mathrm{ad}}$ on the torsor $d^{1}P$.
 This equation is 
the  analogue for $\kappa$ of the assertion that the section
$\Psi_{ijk}$ of $d^1P$ is horizontal.

\bigskip

Since the source and target of the 
        arrows $\widetilde{K}_{i}$ \eqref{def:Ki1} which locally 
        determines 
        the curving data $K$  are now both trivialized, the map $K_{i}$     
     is described by the form $B_{i} \in \mathrm{Lie}(G_{i},\, 
 \Om^{2}_{U_{i}/S})$
for which 
\bee
\label{defbigki}
\underline{\widetilde{K}}_{i}(\delta^{1}\uga_{i}) =
B_{i}\, \ude_{i}\:. \end{equation}
In order to  express 
the commutativity of the square \eqref{comp:kikj} in Lie-theoretic terms, it is
necessary to observe that while both upper vertices of this square now have 
global sections, and are therefore trivialized as torsors, this is not the
 case for the lower two. Since $\De_j$ is isomorphic to 
$(T_{G_j},\, \nu_j)$ as a $G_j$-bitorsor, we find that the section  $\de_j$
 determines a
bitorsor isomorphism 
{\renewcommand{\arraystretch}{1.5}
\[
\begin{array}{ccl}
{}^{P_{ij}\,}\!\De_j &\simeq &\pra P_{ij} \wedge (T_{G_j},\, \nu_j) \wedge
\pra P_{ij}^0
\\&\simeq &\nu_{j\,\ast}(P_{ij})\wedge P_{ij}^0 
\end{array}
\]
}

For the same reason the left-hand vertex (identified with the $G_i$-bitorsor
${}^{P_{ij\,}}\!(\de^1\Ga_j)$), is now isomorphic to
 $(\de^1\mu_{j})_\ast(P_{ij})\wedge P_{ij}^0$.
The commutativity of \eqref{comp:kikj} is expressed by the  equality
\begin{equation}
\label{relKeps} 
({}^{P_{ij}\,}\!B_{j})  \,\text{curv}(\epsilon_{ij})
=  B_{i} \,{\kappa}_{ij}\:.
\end{equation}
between  pairs of  pointed sections of the $G_i$-torsor
$ \nu_{j\,\ast}(P_{ij})\wedge P_{ij}^0$. 
Here  ${}^{P_{ij}\,}\!B_{j}$ is interpreted as  the pointed section  of 
the $G_i$-torsor
${}^{P_{ij}\,} ((\de^1\mu_{j})_\ast(P_{ij})^0 \wedge   \nu_{j\,\ast}(P_{ij}))
   $ which  describes 
the corresponding lower map ${}^{P_{ij}\,}\!(\widetilde{K}_{j})$ in
\eqref{comp:kikj}. In 
additive notation, this equation becomes
\begin{equation}
    \label{relKeps1}
    {}^{P_{ij}\,}\!B_{j} + \mathrm{curv}(\epsilon_{ij}) = B_{i} + 
\kappa_{ij}\,.\end{equation}

\bigskip

The local 3-curvature 
section $\underline{\Om}_{i}$ 
of $\de^{2}_{\UGa_{i}}(\UDe_{i})$ \eqref{def:omegi} 
is   described by the 3-form $\om_{i} \in \mathrm{Lie}(G_{i},\, 
\Om^{3}_{U_{i}/S})$ defined by the equation
\bee
\label{def:3curfor}
\Om_{i} = \om_{i}\,\de^{2}_{\ga_{i}}(\de_{i})\:. 
\end{equation}
Applying the functor $\de^{2}_{\Ga_{i}}$ to the morphism $K_{i}$
evaluated as in \eqref{defbigki} at $\de^{1}\ga_{i}$ yields the section
\[\Om_{i} =  \delta^{2}_{\ga_{i}}(B_{i}\,\de_{i}) =  
\delta^{2}_{m_{i}}(B_{i})\, \delta^{2}_{\ga_{i}}(\de_{i})\]
of $\de^{2}_{\Ga_{i}}(\De_{i})$, with the last equality following from 
\cite{cdf} lemma 2.8 and remark 2.9. It
follows that
\bee
\label{omidef1}\om_{i} =  \delta^{2}_{m_{i}}(B_{i}) \:.
\end{equation}
Applying $\de^{2}_{m_{i}}$ to  equation 
\eqref{ifi} below,  we find that
\begin{equation}
    \label{ificonj}
i_{\om_{i}} = (\de^2_{m_{i}}\nu_{i})^{-1}
\end{equation}
since $\de^2_{i_{m_{i}}}\ka_{m_{i}} = 1$ by the Bianchi identity 
(\cite{cdf} 
lemma  3.5) for the 1-form $m_i$.
This very simple equation is all that is left in the 
present localized context of the complicated relation between 
$j_{\Om}$ and the various pullbacks of $\kc$ occuring in  \eqref{relKkad}.

\bigskip

We end this paragraph, with the following summary of  the previous 
discussion.
\begin{definition}
    \label{decdef}
    Let  $\pc$ be a gerbe with connection triple $(\epsilon, K, 
    \kappa)$ and  local 
    sections $(x_{i})_{i \in I}$  above some open cover 
    $\mathcal{U} = (U_{i})_{i \in I}$,
    with associated $U_{i}$-groups
 $G_{i}:= \mathrm{Aut}(x_{i})$ and bitorsor cocycle
 data $(P_{ij}, \Psi_{ijk})$. Let $\Gamma_{i}$ be the
 $(\pra G_{i},\, \prb G_{i})$-bitorsor, 
 \eqref{def:gai} which 
 locally describes  the connection $\epsilon$ and
 $\De_i$ the $\pra G_i$-bitorsor\eqref{dinonu} 
 above $\De^2_{U_i/S}$ which describes the
 fake curvature $\kappa$.
 We  say that the  connection triple  $( \epsilon ,\, K,\, \kappa)$
on the gerbe $\pc$  is  
partly 
decomposed when   pointed global sections $\ga_{i}$ and
 $\delta_{i}$ of 
 $\Gamma_{i}$ and $\Delta_{i}$  have been chosen for all $i \in I$.
 \end{definition}
The bitorsor 
 structures on $\Gamma_{i}$ and $\Delta_{i}$ are then respectively described
 by the family of connections  \eqref{gcon1}
 \[m_{i}: \prb G_{i} \la \pra G_{i}\] on the group schemes $G_{i}$,
 and the
 $\mathrm{Aut}(G_{i})$-valued 2-forms  \eqref{bitordi}
 \[\nu_{i}: \pra G_{i}\la \pra 
 G_{i}\:,\]
 and the curving datum morphism $\widetilde{K}_{i}:\de^1(\Ga_{i}) 
\la \De_{i}$ corresponds to the $G_{i}$-valued 2-form $B_{i}$
\eqref{defbigki}. 
\begin{proposition}
        \label{prop:omi}
        Let $\pc$   be gerbe  with  a partly decomposed connection triple
 $(\epsilon,\,K,\,\kappa)$,  and suppose that the associated 
       coefficient groups $\mathrm{Aut}(x_{i})$ are represented 
      by flat  $U_{i}$-group schemes $G_{i}$. To this 
decomposition are associated     the
        isomorphisms
$m_{i}$ \eqref{gcon1}and $\nu_{i}$ \eqref{bitordi} and the pair  $( \epsilon_{ij},
\, B_{i})$, with
$\epsilon_{ij}$  a  connection 
\eqref{def:eij}
on the bitorsor cocycle $(P_{ij},\, \Psi_{ijk})$, and $B_{i}$ a {\rm Lie} 
$G_{i}$-valued 2-form  \eqref{defbigki} on $U_{i}$ satisfying equation
\eqref{relKeps1}. The fake curvature $\kappa$ is described by the
 family of $(\nu_{i},\, \nu_{j})$-equivariant
pointed isomorphisms ${\kappa}_{ij}$ \eqref{def:kij1} 
above $\Delta^{2}_{U_{ij}/S}$, which satisfies equations   \eqref{hordiag2}
and \eqref{relKeps1}. 
The associated local  3-curvature $\Om$ 
is locally described by the $G_{i}$-valued 3-curvature forms $\om_{i}$
\eqref{def:3curfor} (which may also be defined by  \eqref{omidef1}) 
and the restrictions 
to $U_{ij}$ of  $\om_{i}$ and $\om_{j}$
 are compatible under  the composite  map  \eqref{ijcompat}. 
\end{proposition}
\begin{flushright}
    $\Box$
    \end{flushright}

     We refer to \eqref{comoioj} and \eqref{relnufi}  below 
for an explicit description of the compatibility between the forms 
$\om_{i}$ and $\om_{j}$,  and  for a cocyclic description of the
 higher Bianchi  identity.

\begin{remark}
    {\rm We now suppose, as we also will  in \S \ref{special3}, that $\pc$ is 
    an abelian $G$-gerbe on $X$, for some  abelian $S$-group $G$, and 
    that the connection $\epsilon$ is a morphism of abelian gerbes, 
    compatible with the canonical connection on $G_{X}$. 
    In that case, the bitorsor structure of $P_{ij}$ is the obvious one,
    determined by the underlying right torsor structure, and the 
    morphism
     $\epsilon_{ij}$ \eqref{def:eij} is simply a connection on this 
     right torsor. Let us  suppose in addition that the fake curvature
     is trivial. The data attached in proposition \ref{prop:omi}
     to the abelian $G$-gerbe $\pc$  with  partly decomposed
     curving  pair $(\epsilon,\, K)$ reduces to the giving of a $G$-torsor
     $P_{ij}$ endowed with a connection $\epsilon_{ij}$, a
horizontal  torsor 
     isomorphism $\Psi_{ijk}$ \eqref{1bcoc}, and a  family 
     of forms  $B_{i}$
    for which the equation \eqref{horpsi-1} and the simplified version 
    \begin{equation}
    \label{simpl}
    \mathrm{curv}(\epsilon_{ij}) = B_{i} - B_{j}
    \end{equation}
    of \eqref{relKeps1} are satisfied.
    Equation \eqref{omidef1} here becomes simply 
    \[\om_{i} = dB_{i}\]
       and it follows directly, or from equation \eqref{comoioj} 
below, that the $\om_{i}$ glue to a global $G$-valued curvature 3-form on
$X$.  For $G= G_{m,S}$, we recover here the description of an abelian 
gerbe with connective structure given by Hitchin in \cite{hitch} \S 1.3.
    }
    \end{remark}

%% file: dgg-ch6.tex
\section{Cocycles and coboundaries  for  gerbes with curving pairs}
\label{sec:fulldec}
\subsection{}   
\setcounter{equation}{0}%
The aim in this section is to give a local, more explicit, 
description of the global cocycles and
 coboundary relations for a gerbe with connection   $(\pc, \, \epsilon )$
 of section \ref{sec:curv-Bian}.   We will obtain combinatorial relations
involving  $\mathrm{Aut}(G)$-valued forms and their combinatorial
differentials. In translating them into classical terms, we will
implicitly assume that $\mathrm{Aut}(G)$ is a smooth $S$-group
scheme.  This is for example true if the $S$-group scheme $G$ is
reductive. The discussion remains essentially valid when $G$ is simply flat
 over $S$, for the reasons mentioned in footnote \ref{Gflat}.

\bigskip

 We have already seen in section \ref{gerbes} that
 the choice of objects $x_i$ and arrows $\phi_{ij}$ \eqref{path}    determines
a cocycle pair $(\lambda_{ij}, g_{ijk})$ satisfying the cocycle conditions
\eqref{def:lij}, \eqref{defgijk}.
 For 
reasons explained in \cite{lb:2-gerbes} \S 2.7, or as recalled here in
\S \ref{gerbes1} when $X$ is quasi-projective over an affine scheme, 
we may indeed choose these
families of
paths $\phi_{ij}:x_{j} \la x_{i}$ above 
$U_{ij}$, without having to resort to local families of paths 
$\phi_{ij}^{\alpha}$ above members of an open cover 
$(U_{ij}^{\alpha})$ of each $U_{ij}$. In addition, we now suppose 
 as in \S \ref{partly-3} that the connection triple on $\pc$ is partly
decomposed. That is, we now assume that we have chosen
global 
sections  the bitorsors $P_{ij}$ have chosen global sections
of the bitorsors $P_{ij}$, as well as pointed global
sections of  the bitorsors $\Ga_{i}$ and $\De_{i}$. The induced global
 sections  of $\UGa_i$ and  $\UDe_i$ are respectively denoted by 
$\uga_i$ and $\ude_i$.

\bigskip

 Just as the bitorsor structures of the trivial torsors $\Ga_{i }$ 
 and $\De_{i}$ are respectively 
 determined by the isomorphism $m_{i}$ \eqref{gcon1} and $\nu_{i}$ 
 \eqref{bitordi}, that of  $P_{ij}$ corresponds to the isomorphism
 $\lambda_{ij}$ \eqref{def:lij0} above $U_{ij}$.
Since the bitorsor $P_{ij}^{0}$ opposite to $P_{ij}$ 
may be identified with
$\mathrm{Isom}_{\pc_{U_{ij}}}(x_{i},\, x_{j})$ under the map which sends
a section $\phi$ to $\phi^{-1}$, 
the target of the arrow  \eqref{defdij1} is  a $(\pra G_{i},\, \pra 
G_{i})$-bitorsor with chosen global section   $\pra \phi_{ij}  
\wedge \ga_{j} \wedge \prb \phi_{ij}^{-1}$, and which can be 
identified with the restriction to $U_{ij}$ of
the sheaf 
$\mathrm{Isom}_{\pc}(\epsilon \prb x_{i},\, \pra x_{i})$ under the map 
which sends this section to the composite arrow
\[ \xymatrix{
\epsilon \prb x_{i} \ar[rr]^{\epsilon(\prb \phi_{ij}^{-1})}&&
\epsilon \prb x_{j} \ar[r]^(.45){\ga_{j}} & \pra x_{j} \ar[r]^{\pra 
\phi_{ij}}& \pra x_{i}
\:.}\]
The arrow $g_{ij}$ 
\eqref{defdij1} is therefore determined, according 
to the recipe \eqref{def:f}, by the section $\ga_{ij} $ of $\pra G_{i}$
defined by 
\begin{equation}
    \label{gaijdef1}
 g_{ij}(\ga_{i}) =    \ga_{ij}\,(\pra \phi_{ij})\, \ga_{j}\, 
\,\epsilon(\prb \phi_{ij})^{-1}\:.
\end{equation} We have seen that  the morphism $g_{ij}$ is pointed, and 
induces a morphism $\underline{g}_{ij}$  \eqref{dired}, so that $\ga_{ij}$ is 
in fact an element of $\mathrm{Lie}(G_{i},\, 
\Om^1_{U_{ij}/S})$. 
In addition, $g_{ij}$ is a bitorsor morphism so that by 
\eqref{defuprime} the term $\ga_{ij}$ also satisfies the equation 
\begin{equation}
    \label{cocep13}
     i_{\ga_{ij}}\, (\pra \lambda_{ij})\,m_{j} \, (\prb 
     \lambda_{ij})^{-1}=  m_{i}\:.
        \end{equation}
In the notation introduced in \eqref{displ1}, the arrow $\epsilon$ may
 be displayed as 
\begin{equation}
   \label{disple}
\xymatrix{
  (\prb \lambda_{ij},\,\prb  g_{ijk}) \ar[rrr]^{(m_i,\,\,
 \gamma_{ij}^{-1})}&&&
 (\pra \lambda_{ij},\,\pra  g_{ijk})}
\end{equation}
Let us now introduce the notation 
     \[{}^{ \lambda_{ij}\,\ast\,}\!m_{j}:= 
(\pra \lambda_{ij})\,m_{j} \, (\prb 
     \lambda_{ij})^{-1}\]
   for   the twisted conjugate of $m_{j}$. 
Equation \eqref{cocep13} may  be rewritten (in additive notation) as
   \begin{equation}
\label{cocep1}
{}^{\lambda_{ij}\,\ast\,}\!m_{j} = - \,i_{\gamma_{ij}} + 
    m_{i} \:.
    \end{equation}

    \begin{remark}
    {\rm 
    
     When the local groups $G_{i}$ are defined over  the base scheme 
    $S$ (for example when $\pc$ is a $G$-gerbe with $G$ an $S$-group), 
  then  they are endowed with  canonical connections so that $m_{i}$ may 
    be viewed as  an $\mathrm{Aut}(G)$-valued  1-form on $U_{i}$.
    In that case one finds  that in additive notation
    \[
    {}^{ \lambda_{ij}\,\ast\,}\!m_{j}=  {}^{
    \lambda_{ij}\,\,}\!m_{j} + \delta^{0}(\lambda_{ij}^{-1}) 
    \]
    so that equation \eqref{cocep1}
    becomes the 
   compatibility condition
    \begin{equation} \label{cocep4}
 {}^{ \lambda_{ij}\,\,}\!m_{j} + \delta^{0}(\lambda_{ij}^{-1}) 
    = - \,i_{\gamma_{ij}} +  m_{i}\:. \end{equation} 
We now set
    \[
       \tilde{\delta}^{0}(g)(x,\,y)  := g(y)\, \, g(x)^{-1}
       \:= {}^{g\,\,}\!\delta^{0}(g)(x,\,y)\,.
   \]
   with $(x,\, y) \in \dea$.  One then finds that 
   \[ \delta^{0}(-g) = - \,\tilde{\delta}^{0}(g)\:.\]
  so that  equation \eqref{cocep13} may finally 
  be rewritten as the   \v{C}ech-de Rham 
  compatibility
  condition
   \begin{equation}
       \label{cocep4-1}
       {}^{  \lambda_{ij}\,}\!m_{j}  - m_{i} =   - i_{\gamma_{ij}} +  
   \tilde{\delta}^{0}(\lambda_{ij}) \:.\end{equation}
 }  \end{remark}

         \vspace{.75cm}

         Let us now consider once more
         diagram \eqref{compdij3}.
           The sections $\ga_{i}$ of $\Ga_{i}$ and $\phi_{ij}$ of 
           $P_{ij}$ now 
         provide us with global sections for each of its four corner terms.  
The top horizontal and
 left vertical arrows in the diagram are  respectively described 
 according to the prescription \eqref{defmor} by the
 sections  $\ga_{ij}$  and $\ga_{ik}$ of $\pra G_i$, and the right vertical
 map
corresponds to the section $\pra (\lambda_{ij})(\gamma_{jk})$. A more elaborate
computation shows that the lower horizontal arrow is described by the
section
\[\pra g_{ijk}^{-1}\,\,\,(\pra (\lambda_{ij}
\lambda_{jk})\,m_k \,\prb \lambda_{ik}^{-1})
(\prb g_{ijk})\] and by \eqref{coclam} this expression
can be rewritten as 
${}^{\lambda_{ik}\ast\,}\!m_k (\prb g_{ijk})\, \pra g_{ijk}^{-1}$. 
     The commutativity of diagram  \eqref{compdij3}
     is now  expressed by the equation
     \[ \gamma_{ij}\, (\pra \lambda_{ij})(\gamma_{jk})
     = \gamma_{ik} ({}^{ \lambda_{ik}\,\ast\,}\!m_{k})(\prb g_{ijk}) 
     \,\,
     \pra g_{ijk}^{-1}\:.\]
     Applying equation \eqref{cocep13} for the indices $i$ and $k$,
     this may be restated as
     \begin{equation}
\label{cocep5}
\ga_{ij}\, \,(\pra \lambda_{ij})(\gamma_{jk})\,\,(\pra g_{ijk})
\:= \:  m_{i}(\prb g_{ijk})\,\, \ga_{ik}\:.
\end{equation}
 This  equation 
is equivalent to 
\[
\ga_{ij}\, \,(\pra \lambda_{ij})(\gamma_{jk})\,\,(\pra g_{ijk}\,\, 
\ga_{ik}^{-1} \, \pra g_{ijk}^{-1})
\:= \:  m_{i}(\prb g_{ijk})\, \pra g_{ijk}^{-1}\, \:\]
which by \eqref{coclam} may be 
expressed as
\[\gamma_{ij}\,\, \lambda_{ij}(\gamma_{jk})\,\, 
n(\lambda_{ij}\lambda_{jk} \lambda_{ik}^{-1})(\gamma_{ik}^{-1}) = 
{}^{g_{ijk}\,}\!(\delta^{0}_{m_{i}}(g_{ijk})) \:.
\]
  The left-hand term in this equation,
  which is now patterned  on  \eqref{def:d1mu},
   may be thought of as the    differential 
   $d^{1}_{\lambda_{ij}}(\gamma_{ij})$ of the 2-cochain $\gamma_{ij}$
   in the \v{C}ech complex of the 
   open cover $\mathcal{U}$ of $X$ 
   with values in the local coefficient groups  $(G_{i},\, \lambda_{ij})$. 
   Setting
    \begin{equation}
        \label{del1alt}
       \tilde{\delta}^{0}_{\mu}(g)(x,\,y)  := 
       \mu(x,y)(g(y))\, \, g(x)^{-1}
       \:= {}^{g\,\,\,}\!\delta^{0}_{\mu}(g)(x,\,y) \end{equation}
  we see that \eqref{cocep5}
   may finally be rewritten  in \v{C}ech-de Rham 
   terms as
   \begin{equation}
       \label{cocep2}
       d^{1}_{\lambda_{ij}}(\gamma_{ij}) 
       =\tilde{\delta}^{0}_{m_{i}}(g_{ijk})\:.
       \end{equation}

\bigskip

We now carry out a similar discussion for the morphisms $d_{ij}$ 
\eqref{ardij3} and the corresponding commutative diagram 
\eqref{compdij5}.
    The source and target of $\underline{d}_{ij}$
    \eqref{ardij4}  are now  respectively trivialized by the sections
    $\ude_{i}$ and $\pra \phi_{ij}\, \ude_{j}\, 
    \kappa(\pra (\phi_{ij}^{-1})$
    so that this morphism is described by the element $\de_{ij} \in 
    \mathrm{Lie}(G_{i},\, \Om^{2}_{U_{ij}/S})$ defined by 
   \bee
   \label{deijdef}
d_{ij}(\delta_{i})  = \delta_{ij}\,(\pra \phi_{ij})\, 
\delta_{j}\,(\kappa(\pra 
\phi_{ij}))^{-1}
\end{equation}
The compatibility of the morphism $\underline{d}_{ij}$ with the 
bitorsor structure is expressed by the analogue
\bee \label{cockap01}
i_{\de_{ij}}\,  (\pra \lambda_{ij})\,\nu_{j} \:=\:  \nu_{i}\,
        (\pra \lambda_{ij})
        \end{equation}
of equation \eqref{cocep13},
an equation which may be written additively as
\begin{equation}
    \label{cockap1}
    {}^{\lambda_{ij}\,}\!\nu_{j} = \nu_{i} - i_{\delta_{ij}}\:.
    \end{equation}
The upper and left-hand arrows in the 
square \eqref{compdij5} are then 
respectively described  by 
the elements $\de_{ij}$ and $\de_{ik}$, and the right-hand arrow 
corresponds to the expression $\lambda_{ij}(\de_{jk})$.  One verifies 
that  the lower 
arrow is  described in terms of the trivializations  of its source 
and  target by  the  expression $\pra g_{ijk}^{-1}\,\, 
 (\lambda_{ij}\lambda_{jk}\nu_{k}\lambda_{ik}^{-1})(\pra g_{ijk})$,
 so that  equation
 \[\de_{ij} \lambda_{ij}(\de_{jk}) = \de_{ik} \pra  g_{ijk}^{-1}
(\lambda_{ij}\lambda_{jk}\nu_{k}\lambda_{ik}^{-1})(\pra g_{ijk})\]
expresses the commutativity of the diagram \eqref{compdij5}.
 By \eqref{coclam},
this is equivalent to 
\bee
\de_{ij} \lambda_{ij}(\de_{jk}) \, \pra g_{ijk} =
\de_{ik}({}^{\lambda_{ik}\,}\!\nu_{k})(\pra g_{ijk})\:.
\end{equation}
Applying \eqref{cockap01} with indices $i$ and $k$ to this equation,
we find that it
may be expressed  as
\label{cockap02}
\begin{equation}
    \de_{ij}\, \lambda_{ij}(\de_{jk})\,  (\pra g_{ijk})\: = \: 
        \nu_{i}(\pra g_{ijk}) \, \de_{ik}\:,
        \end{equation}
or even, when  equation \eqref{coclam} is taken into account, 
by
\[\delta_{ij} \, \lambda_{ij}(\delta_{jk})\, \lambda_{ij}
 \lambda_{jk}\lambda_{ik}^{-1}(\delta_{ik}^{-1}) 
 = \nu_{i}(\pra g_{ijk})\, \pra g_{ijk}^{-1}\,.\]
The right-hand side of this equation may be restated  as
\[ \nu_{i}(m_{i}^{02}(\prc g_{ijk}))\, (m_{i}^{02}(\prc g_{ijk}))^{-1}\]
so that  the commutativity of diagram
    \eqref{compdij5} may finally be expressed (in the notation of 
    \eqref{lie1})
as 
\begin{equation}
    \label{cockap2}
    d^{1}_{\lambda_{ij}}(\delta_{ij}) = [\nu_{i},\,g_{ijk}]\, .
    \end{equation}

    \bigskip

We have seen that the pointed arrow  $\widetilde{K}_{i}$ \eqref{def:Ki1}.
 is described by the 2-form $B_{i}$  \eqref{defbigki}.
The bitorsor structure  on the source torsor  is described according to 
\eqref{defu} by the automorphism $\delta^1(m_{i})$ of $G_{i}$, in 
other words by the curvature $\kappa_{m_{i}}$ of the connection $m_{i}$
\eqref{gcon1} on the group $G_{i}$.
 Since $\widetilde{K}_{i}$ is a morphism of bitorsors, the rule 
\eqref{defuprime} asserts that the equation 
\[\kappa_{m_{i}} = i_{B_{i}}\,\,\nu_{i}  \]
is satisfied. 
We display this in additive notation as
\begin{equation}
   \label{ifi}
   \kappa_{m_{i}} = i_{B_{i}} + \nu_{i}
    \end{equation}
This formula can either be viewed as the definition of the 
2-form    $\nu_{i} \in \mathrm{Lie(Aut\,}G_{i}, \, \Om^{2}_{U_{i}/S})$ 
or else,  as we have done here, as a   condition satisfied by
 the 2-form $B_{i}$.
 
 \bigskip
 
 We will  now express  in terms of these cocycles the 
commutativity of diagram  \eqref{comp:kikj}. We have seen that the 
upper arrow in that diagram  is described by the 2-form $B_{i}$
 and the right-hand 
vertical one  by $\de_{ij}$. Of the two arrows on the bottom line, 
the unlabelled first one is a canonical isomorphism which contributes 
nothing to our equations, and  the second one is easily seen to be
represented by the $G_{i}$-valued 2-form $\lambda_{ij}(B_{j})$. It can 
also be verified that the left-hand vertical arrow is represented
by the term
 \[\gamma_{ij}^{01}\,\,\, 
({}^{\lambda_{ij}\,\ast\,\,}\!m_{j}^{01})(\gamma_{ij}^{12})\,\,\,
{}^{\lambda_{ij}\,\ast\,}\!(m_{j}^{01}m_{j}^{12})(\gamma_{ij}^{02})^{-1}\]
so that the equation which asserts the commutativity of \eqref{comp:kikj}
is
\[ B_{i}\,\,  \de_{ij} =  \gamma_{ij}^{01}\,\, 
({}^{\lambda_{ij}\,\ast\,\,}\!m_{j}^{01})(\gamma_{ij}^{12})\,\,
{}^{\lambda_{ij}\,\ast\,\,}\!(m_{j}^{01}m_{j}^{12})(\gamma_{ij}^{02})^{-1}
\,\, \lambda_{ij}(B_{j}) \:.\]
This  expression  is better stated (in additive notation) as
\begin{equation}
    \label{relKeps2}
    \de_{ij}+ B_{i} = \lambda_{ij}(B_{j}) + \de^{1}_{({}^{\lambda_{ij}\, 
\ast}m_{j})}(\ga_{ij})\,.
\end{equation}
It is the fully cocyclic form of equation \eqref{relKeps1}.
By \eqref{cocep1}, it is equivalent to  
\[ \de_{ij}+ B_{i} = \lambda_{ij}(B_{j}) + \de^{1}_{(- i_{\gamma_{ij}} + 
m_{i})}(\gamma_{ij})\:.
\]
One verifies that 
\[ \de^{1}_{(- i_{\gamma_{ij}} + 
m_{i})}(\gamma_{ij}) = -\de^1_{m_{i}}(-\ga_{ij})\]
so that  equation \eqref{relKeps2} can be restated as   
\begin{equation}
    \label{compfifj2}
   \de_{ij}+ B_{i} = \lambda_{ij}(B_{j}) - \de^{1}_{m_{i}}(- \ga_{ij})
    \end{equation}
    For any $X$-group $G$ with connection $\mu$, and any 
    $\gamma \in \mathrm{Lie}(G,\, 
    \Om^{1}_{X/S})\:,$ 
   \begin{equation}
       \label{d1rule}
       \delta^{1}_{\mu}(-\gamma) = - \delta^{1} _{\mu}(\gamma) + 
   [\gamma,\, \gamma]_{\mu}  \end{equation}
   so that equation \eqref{compfifj2} may also be expressed  as
    \begin{equation}
       \label{compfifj} 
       \lambda_{ij}(B_{j}) = B_{i} + \delta_{ij}
       - \delta^{1}_{m_{i}}(\gamma_{ij}) + 
   [\gamma_{ij},\, \gamma_{ij}]_{m_{i}}\:.
   \end{equation}

\bigskip

The cocyclic description $\om_{i}$ of  the  3-curvature form $\Om_{i}$  was 
given in \eqref{def:3curfor}, but it remains to interpret cocyclicly 
the compatibility condition between a pair of 
 local forms $\Om_{i}$ and $\Om_{j}$
under the morphism \eqref{ijcompat}. Let us begin by considering the 
left-hand  arrow 
\begin{equation}
    \label{firsteq}
    \de^2_{\Ga_{i}} \De_{i} \la 
\de^2_{\Ga_{i}}({}^{P_{ij}\,}\!\De_{j})\end{equation}
in this composite morphism. We have  seen 
that our
chosen trivializations $ 
\de_{i}$  of $ \De_{i}$ and  $\phi_{ij}$ of $P_{ij}$ determine   
trivializations  of  $\De_{i}$ and ${}^{P_{ij}\,}\!\De_{j}$ in terms 
of which the morphism $d_{ij}$ \eqref{ardij3} was described by the 
2-form $\de_{ij}$ \eqref{deijdef}. With $\Ga_{i}$ trivialized by the 
section $\ga_{i}$ the source of \eqref{firsteq} is trivialized by 
the induced section $\de^{2}_{\ga_{i}}(\de_{i})$, defined by the 
ususal combinatorial formula \eqref{def:d2mu}, with $\ga_{i}$ acting 
on $\de_{i}$ by conjugation. The target of \eqref{firsteq} is 
trivialized by the corresponding expression, but with the term $\de_{i}$
replaced by the trivialization $\pra \phi_{ij}\, \de_{j}\, 
    \pra (\phi_{ij}^{-1})$ of ${}^{P_{ij}\,}\!\De_{j}$.
   The pullback of the 
   morphism $d_{ij}$ to  each of the four factors of the expression
   \eqref{def:de2de}   yields a corresponding pullback above 
   $\Delta^{3}_{U_{ij}/S}$ of the 2-form $\de_{ij}$. In order to get 
   an explicit description of \eqref{firsteq}, we must still regroup 
   these 2-forms. 
  This {\it a priori} involves the bitorsor structure of
   ${}^{P_{ij}\,}\!\De_{j}$, and 
   thus the action on  these sections of the corresponding conjugates of
the  automorphism $\nu_{j}$ 
   which describes this structure. By the analogue  of lemma 2.8 of 
   \cite{cdf},  this action of the $\nu_{j}$ 
   is trivial in each instance.  This gives us the 
   following explicit description of morphism \eqref{firsteq}:
   \begin{lemma}
       For the trivializations of its source and target induced from 
       the chosen trivializations  $\ga_{i}, \,\de_{i},\, \phi_{ij}$ of 
       the torsors
       $\Ga_{i},\,\De_{i},\,P_{ij}$, the arrow $\eqref{firsteq}$ is 
       described by the element 
 $\de^{2}_{m_{i}}(\de_{ij})$ in   \linebreak
 $\mathrm{Lie}(G_i,\, \Om^3_{U_{ij}/S})$.   
       \end{lemma}
       \begin{flushright}
   $\Box$
   \end{flushright}

The second arrow in \eqref{ijcompat}  has the following
 local description, 
as can be seen from its definition: 
  \begin{lemma}
      With the trivializations of its source and target induced  from 
the      chosen trivializations  $\ga_{i}, \,\de_{i},\, \phi_{ij}$, 
the second arrow  in \eqref{ijcompat} is described by the term 
$[\gamma_{ij},\, {}^{\lambda_{ij\,}}\!\nu_{j}]$.  
\end{lemma}
\begin{flushright}
   $\Box$
   \end{flushright}
Since the third morphism in \eqref{ijcompat} is a canonical one, 
and can therefore be ignored in this discussion,  we find that the 
sought-after relation between the local 3-curvature forsm $\om_{i}$ 
and $\om_{j}$ is finally

  \begin{equation}
     \label{comoioj1}
      \lambda_{ij}(\om_{j}) = \omega_{i} + \delta^{2}_{m_{i}}(\delta_{ij})
      +  [\gamma_{ij},\, {}^{\lambda_{ij}\,}\!\nu_{j}] 
     \end{equation}
   Taking into account the relation \eqref{cockap1}, this can also be 
   written as  
   \begin{equation}
     \label{comoioj}
     \lambda_{ij}(\om_{j}) = \om_{i} +
 \delta^{2}_{m_{i}}(\delta_{ij})  + [\gamma_{ij},\,\nu_{i} ]
 - [\gamma_{ij},\, \delta_{ij}]\,.
     \end{equation}

     \bigskip

     The final condition to be interpreted is the higher Bianchi 
     identity.  Its local form can be directly derived from
     \eqref{neomegeq}, or from  the triangle \eqref{biantrian}.
     All three vertices in this triangle are  now trivialized, and the 
     horizontal arrow, being canonical, has no effect  on the 
     cocyclic description of the higher Bianchi identity. 
     It is apparent that the 
     right-hand one, which is induced by applying the morphism $K_{i}$ 
     to objects which conjugate $\De_i$ is represented by the form 
     $[B_{i},\, \nu_{i}]$.
     The higher Bianchi identity now reads 
     \begin{equation}
         \label{relnufi1}
          \delta^{3}_{m_{i}}(\om_{i}) \,+ \,[B_{i},\, \nu_{i}] = 0\,,
          \end{equation}
          or equivalently
     \begin{equation}
         \label{relnufi}
         \delta^{3}_{m_{i}}(\om_{i}) = [\nu_{i},\, B_{i}]\,.
     \end{equation}
     This equation is the local form of 
     \eqref{neomegeq}, from which it could have directly been derived.
     
   \bigskip 

 We  restate the previous discussion as follows:  

\pagebreak
  
\begin{theorem}
        \label{th:omi}
        Consider   a gerbe  $\pc $ with local objects 
        $x_{i}$ 
        defined  above some open cover
      $\mathcal{U}= (U_{i})_{i \in I}$ of $X$, and suppose that the 
      associated coefficient groups $\mathrm{Aut}(x_{i})$ are represented 
      by flat $U_{i}$-group schemes $G_{i}$. 
To local paths $\phi_{ij}$ \eqref{path}  are associated 
 cocycles pair $(\lambda_{ij},\, g_{ijk})$  defined by 
{\rm (\ref{def:lij}), (\ref{defgijk})}
and satisfying the cocycle conditions
{\rm (\ref{coclam}), (\ref{cocg})}. Let $\pc$ be in addition endowed  with 
       a  partly decomposed curving pair $(\epsilon,\, K)$, with given
        paths $\ga_{i}$ \eqref{defchii} and  $\gamma'_{i}$ 
        \eqref{defchii1}
        as well as   paths $\de_{i}$ \eqref{ardei} and $\de'_i$
\eqref{ardeiprime}.
        The induced  isomorphisms  $m_{i}$ \eqref{gcon1} and $\nu_{i}$
        \eqref{bitordi}, and elements  $\ga_{ij}$ \eqref{gaijdef1} and 
        $\de_{ij}$ \eqref{deijdef}  respectively satisfy equations 
        {\rm (\ref{cocep1}), (\ref{cocep2})} and {\rm (\ref{cockap1}),
        (\ref{cockap2})}.
        Let  $B_{i}$  be  given 
        by \eqref{defbigki}. 
        The 3-curvature form $\om_{i}$, defined  by \eqref{omidef1}, satisfies 
 equations  \eqref{compfifj2} and \eqref{ifi}.
      The curvature triple $(\nu_{i},\, \de_{ij},\, \om_{i})$
 satisfies, in addition to 
        \eqref{cockap1} and
        \eqref{cockap2}, the equations
  \eqref{comoioj}, \eqref{ificonj} and  \eqref{relnufi}.
    \end{theorem}
    \begin{flushright}
        $\Box$
        \end{flushright}

We now display the functions and differential 
  forms occuring in the statement of 
  the theorem,  with  the proviso that $\lambda_{ij}$ is actually a section of 
$\mathrm{Isom}(G_{j},\, G_{i})$ and $m_{i}$  a connection on the 
group scheme $G_{i}$. 
   \begin{table}[ht]
     \begin{center}
        \renewcommand{\arraystretch}{1.3}
        \begin{tabular}{|l|c|c|c|c|} \hline
 & functions& 1-forms  &2-forms & 
 3-forms\\\hline
  $G_{i}$-valued& $g_{ijk}$ & $\gamma_{ij}$ & 
  $\delta_{ij},\,\, B_{i}$ & $\om_{i}$\\ 
$\mathrm{Aut}(G_{i})-\mathrm{valued}$&$\lambda_{ij}$& 
$m_{i}$ & $\nu_{i}$&\\\hline
\end{tabular}

\caption{}\label{Ta:1}
\end{center}
\end{table}

The equations which they satisfy are

  \[\renewcommand{\arraystretch}{1.3}
  \begin{array}{rl}
      &
\left\{ \begin{array}{rll}
        \lambda_{ij}(g_{jkl})\,g_{ijl}  &= g_{ijk} \, 
        g_{ikl} \qquad\qquad\qquad \qquad\qquad\qquad\qquad\quad
 & \qquad \hspace{1.5ex}\eqref{cocg}\\
        \lambda_{ij} \: \lambda_{jk} &= i_{g_{ijk}} 
        \lambda_{ik}    \qquad\qquad\qquad \qquad\qquad\qquad\qquad\quad
 & \qquad \hspace{1.5ex}\eqref{coclam}
    \end{array}\right.
  \\&\\&   
\left\{ \begin{array}{rll} 
{}^{\lambda_{ij}\, \ast\,}\!m_{j} 
  &= - \,i_{\gamma_{ij}} +   m_{i} &
 \qquad\qquad \qquad\qquad\qquad\qquad\qquad\qquad\qquad 
 \eqref{cocep1}-\eqref{cocep4-1}\\
  d^{1}_{\lambda_{ij}}(\gamma_{ij}) 
       &=\tilde{\delta}^{0}_{m_{i}}(g_{ijk})
        & \qquad\qquad \qquad\qquad\qquad\qquad\qquad\qquad \qquad 
\eqref{cocep2}
\end{array}\right.
  \\&\\
 & \hspace{.35cm}  \begin{array}{rll}
 \lambda_{ij}(B_{j}) &= B_{i} + \delta_{ij} + 
   \delta^{1}_{m_{i}}(-\,\gamma_{ij})
 & \quad\quad\qquad \qquad \qquad \qquad\qquad\qquad
 \eqref{compfifj2}\end{array}\
 \\&\\&
\left\{  \begin{array}{rll} 
{}^{\lambda_{ij}\,}\!\nu_{j} &= \nu_{i} - i_{\delta_{ij}}
\qquad \qquad \qquad \qquad 
\qquad & \qquad\qquad\qquad\qquad\hspace{2.5ex}\:
        \eqref{cockap1}
        \\
        d^{1}_{\lambda_{ij}}(\delta_{ij}) &= [\nu_{i},\, g_{ijk}]
        \qquad \qquad \qquad \qquad 
\qquad & \qquad\qquad\qquad\qquad\hspace{2.5ex}\: 
        \eqref{cockap2}
        \end{array}\right.
  \\&\\
  &\left\{  \begin{array}{rlr}
   \lambda_{ij}(\om_{j})  &= 
\omega_{i} + \delta^{2}_{m_{i}}(\delta_{ij}) +
[\gamma_{ij},\,\nu_i] - [\gamma_{ij},\,\delta_{ij}] 
    &\qquad \qquad\qquad\qquad\hspace{1.2ex}
\eqref{comoioj}\\
 \delta^{3}_{m_{i}}(\om_{i}) &= [\nu_{i},\, B_{i}]&
\qquad\qquad\qquad\hspace{1.2ex}
  \eqref{relnufi}
\end{array}\right.
  \\&\\&
  \left\{  \begin{array}{rlr}
 \nu_{i} &=  \kappa_{m_{i}} - i_{B_{i}}  \qquad \qquad \qquad \qquad 
 \qquad\qquad \qquad \qquad \qquad  & \qquad\hspace{3ex}\eqref{ifi}\\
 \om_{i}  &=    \delta^{2}_{m_{i}}(B_{i}) & \qquad\hspace{3ex}
 \eqref{omidef1}\\
  i_{\om_{i}}
 &= - \,\de^2_{m_{i}}\nu_{i}
 & \qquad\hspace{3ex}\eqref{ificonj}
 \end{array}\right.
 \end{array}
\]

\bigskip

\pagebreak

\noindent These equations have been regrouped according to their origin, 
as in the following table.

\begin{table}[ht]
    \begin{center}
    \renewcommand{\arraystretch}{1.3}
    \begin{tabular}{|l|c|c|c|} 
 \hline
 global & name & local & equations
 \\\hline
  $\pc$& gerbe  &$(\lambda_{ij},\, g_{ijk})$ &
  (\ref{coclam}, \ref{cocg})
 \\\hline
  $\epsilon$ & connection & $(m_{i},\,\gamma_{ij})$ &  
  (\ref{cocep1}, \ref{cocep2})
 \\ \hline 
$K$&curving datum& 
$B_{i}$ &(\ref{compfifj})\\
\hline
$\kappa$ &fake curvature & $(\nu_{i},\, \de_{ij})$ & (\ref{cockap1},
\ref{cockap2})\\
\hline
$\Omega $& 3-curvature & $\om_{i}$ & (\ref{comoioj}, \ref{relnufi})\\ 
\hline
\end{tabular}

\caption{}\label{Ta:2}
\end{center}
\end{table}

Equations \eqref{ifi},  \eqref{omidef1} and \eqref{compfifj2}
   have been omitted from this
table, since they may be viewed as the definitions  of $\nu_{i}$, 
$\om_{i}$ and  $\delta_{ij}$
in terms of $B_i$, $m_i$ and $\gamma_{ij}$. So has   \eqref{ificonj}, 
since this is a consequence of 
\eqref{ifi} and the Bianchi identity.

\bigskip

\subsection{}  
\setcounter{equation}{0}%
\label{disc:cobound}

We now pass from cocycle to coboundary relations and
 consider a  morphism $(x,h,a)$ \eqref{def:xa}-\eqref{def:a}
between two  curving data pairs $(\epsilon, \,
K)$ and $(\epsilon',\, K')$, on a gerbe $\pc$. In giving 
 a fully combinatorial description of  these coboundary relations, we
 will at times find it expedient  to reason  directly in terms of the 
 composition rules
 for  1- and 2-arrows displayed in paragraph \ref{algar}, instead of
 making use of  commuting diagrams of 
bitorsors as in the last section.

\bigskip

Just as the connection  $\epsilon$ \eqref{gercon}  
was displayed, for 
a chosen pointed arrow  $\gamma_i$ \eqref{defchii}, by the diagram
\eqref{disple}, we   may
 associate to the 1-arrow $h$  \eqref{def:xa}
the corresponding family of $\pra \,G_i$-bitorsors 
\bee
\label{def:H}
H_i:= \mathrm{Isom}_{\pc_{\mid \, U_i}}(h(x_i),\, x_i)
\end{equation}
and the family of bitorsor isomorphisms
\begin{equation}
    \label{defhij1}
    h_{ij}: H_i \la {}^{P_{ij}\,}\!H_j
    \end{equation}
for which the diagrams
 analogous to \eqref{compdij5} commute.
For a chosen section $h_i$ of $H_i$, the arrow $h$
may be therefore
displayed as
\[
\xymatrix{(\pra \lambda_{ij},\,\pra g_{ijk}) \ar[rr]^{(\pi_i,\, \eta_{ij}^{-1})}
&&(\pra \lambda_{ij},\,\pra g_{ijk}) 
 } \]
where  $\pi_{i} \in \mathrm{Lie (Aut}(G_{i}),\,
 \Om^{1}_{U_i/S})    $  is     
induced from $h_i$ by conjugation, and $\eta_{ij} \in
 \mathrm{Lie}(G_{i},  \Om^1_{U_{ij}/S}) $ 
is the 1-form defined by 
\[
h_{ij}(h_i) = \eta_{ij} h_j
\]
in much the same way as $\gamma_{ij}$ was defined in  \eqref{gaijdef1}. 
The following equations, similar
to \eqref{cocep13} and \eqref{cocep5} are satisfied:
{\renewcommand{\arraystretch}{1.5}
 \[
\begin{array}{rcl}
\pi_i& = &i_{\eta_{ij}}\, (\pra \lambda_{ij})\, \pi_j \,
 (\pra \lambda_{ij})^{-1}\\
\eta_{ij}\,\,(\pra \lambda_{ij})(\eta_{jk})\,\, \pra g_{ijk} &=& 
\pi_{i}(\pra g_{ijk})\,\,
\eta_{ik}
    \end{array}
\]
}
The first of these may be displayed additively as
\begin{equation}
    \label{pij1}
\pi_{i }= i_{\eta_{ij}}  + {}^{\lambda_{ij}}\pi_{j}  
\end{equation}
and the second may be rewritten, after some combinatorial
simplifications modelled on the proof of \eqref{cocep2}, 
as
\begin{equation}
    \label{d1eij}
d^1_{\lambda_{ij}}(\eta_{ij}) = [\pi_i,\, g_{ijk}]\:.
\end{equation}

\medskip

To the  diagram \eqref{def:xa}
corresponds a family of bitorsor isomorphisms
 \[ \xymatrix{
H_i \wedge \Ga_i \ar[r]^(.6){\xi_{i}} & \Ga'_i\:.
}\]
These are represented, in terms of the induced trivializations of source and
target, by the 1-forms $  E_i \in \mathrm{Lie}(G_i,\, \Om^1_{U_i/S})$ defined
by 
\[
\xi_i(h_i \wedge \ga_i)  = E_i \, \gamma'_i\,.
\]
Diagram  \eqref{def:xa} may now be displayed as

\begin{center}

\[
\includegraphics{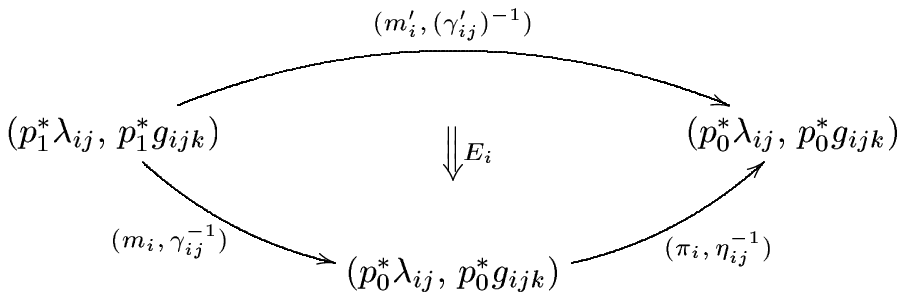}
\]

\end{center}

\bigskip

\noindent Composing the two lower 1-arrows, it follows that the 1-forms $E_{i}$
satisfy the equations
{\renewcommand{\arraystretch}{1.5}
\begin{align*}
i_{E_{i}}\, m'_{i} &= \pi_{i}\, m_{i}\\
E_{i}\, \gamma'_{ij}& = \pi_{i}(\gamma_{ij})\, \eta_{ij}\, 
\lambda_{ij}(E_{j})\,.
\end{align*}}
\noindent In the second equation, the $\mathrm{Aut}(G)$-valued
1-form $\pi_{i}$ acts trivially on the
$G_{i}$-valued form $\gamma_{ij}$, and the various factors 
commute with each other. The two equations may  be 
restated additively as
{\renewcommand{\arraystretch}{1.5}
\begin{align}
    \label{cobe1}
i_{E_{i}} + m'_{i} &= \pi_{i} + m_{i}\\
\label{cobe2}
\gamma'_{ij} - \gamma_{ij}&= \eta_{ij}+ \lambda_{ij}(E_{j}) - E_{i}\,.
\end{align}
}

The 2-arrow $a$ \eqref{def:a} may be described in similar terms. It 
corresponds to a family of bitorsor isomorphisms
\begin{equation}
    \label{bitorsai}
    H^{i}_{01}\, {}^{\Ga^{i}_{01}\,}\!H^{i}_{12}\,\,\De^{i}_{012} 
\stackrel{a_{i}}{\la} {\De'}^{i}_{012}\,\, H^{i}_{02}\,.\end{equation}
For the chosen trivialisations of source and target, the isomorphisms
$a_{i}$ are described by the family of elements $\alpha_{i}\in 
\mathrm{Lie}(G_{i},\,\Om^{2}_{U_{i}/S})$
 defined
by
\begin{equation} 
    \label{def:alphi}
  a_{i}( h^{i}_{01} 
{}^{\ga^{i}_{01}\,}\!h^{i}_{12}\,\, \de^{i}_{012})  = 
\alpha_{i}\,\, ({\de')}_{012}^{i}\, \,h^{i}_{02}\,.\end{equation}   
Assuming the $G_i$ are flat $U_i$-group schemes,
 proposition \ref{psigc}, together with  \cite{cdf} (lemma 2.8, remark 2.9),
 allows  us to 
 rewrite the isomorphisms \eqref{bitorsai}
as 
\bee
\label{permHD}
 (\de^{1}_{\Ga_{i}}H_{i})\, \Delta_{i}  \stackrel{a_i}{\la}    \Delta'_{i}
\end{equation}
and in that form the compatibility condition which the morphisms 
$a_{i}$ must satisfy is given by the commutativity of the bitorsor 
diagram
\begin{equation}
    \label{com:hijk}
 \xymatrix{
 (\de^{1}H_{i})\, \De_{i}  
 \ar[rrr]^{a_i}\ar[d]_{(\de^{1}_{m_{i}}h_{ij})\, d_{ij}} 
 &&& \De'_{i}  \ar[d]^{d'_{ij}} \\
 {}^{P_{ij}\,}\!(\de^{1}H_{j}) \,\, {}^{P_{ij}\,}\!(\De_{j})
 \ar[r]^(.55){\sim} & {}^{P_{ij}\,}\!(\de^{1}H_{j}\, \De_{j}) 
 \ar[rr]^{ {}^{P_{ij}\,}\!a_j} && {}^{P_{ij}\,}\!(\De'_{j})\,.
 }\end{equation}
The   fact that the map  $a_i$  \eqref{bitorsai} is a morphism of bitorsors
and the  commutativity of  diagram \eqref{com:hijk} 
each  determine  an
  equation which the associated 2-forms  $\alpha_i$  must satisfy.
 Instead of  working these equations out directly, we can represent 
the 2-arrow $a_{i}$  \eqref{bitorsai} by a display. Composing the
1-arrows determining the arrow   $\mu_{01}(h_{12})$
\eqref{def:a}  as a conjugate of $h_{12}$
 according the composition rules 
for displays, we  observe  that this  arrow is displayed as 
\[
\xymatrix{(\pra \lambda_{ij},\, \pra g_{ijk})
\ar[rrrrrrr]^{({}^{m_{01}\,}\!\pi_{i}^{12},\, (\gamma_{ij}^{01})^{-1}
m_{i}^{01}(\eta_{ij}^{12})^{-1}\, 
({}^{m_{01}\,}\!\pi_{i}^{12})(\gamma_{ij}^{01})
)}
&&&&&&& (\pra \lambda_{ij},\, \pra g_{ijk})}\:.\]
It follows that the  2-arrow $a_{i}$ 
is  displayed as 
\[
    \xymatrix@C=14pt{&&& (\pra \lambda_{ij},\, \pra g_{ijk}) 
    \ar@/^1pc/[rrrd]^{(\nu'_{i},\, (\de'_{ij})^{-1})} &&&\\
    (\pra \lambda_{ij},\, \pra g_{ijk}) 
    \ar@/_1pc/[rd]_{(\nu_{i},\, (\de_{ij})^{-1})\,\,\,\,\:}
    \ar@/^1pc/[rrru]^{(\pi_{i}^{02},\, (\eta_{ij}^{02})^{-1})}&&&&&&
    (\pra \lambda_{ij},\, \pra g_{ijk})\,.\\
    &  (\pra \lambda_{ij},\, \pra g_{ijk})
    \ar@/_1.3pc/[rrrr]_{({}^{m_i^{01}\,}\!{\pi }^{12}_{i},\,\, 
    (\ga_{ij}^{01})^{-1}\, m_{i}^{01}(\eta_{ij}^{12})^{-1}\,\, 
    ({}^{m_{i}^{01}\,}\!\pi_i^{12})(\gamma_{ij}^{01})
    )} &&
    \ar@{}[uu]_(.35){\,}="1"  
     \ar@{}[uu]_(.65){\,}="2"     &&
    (\pra \lambda_{ij},\, \pra g_{ijk}) 
    \ar@/_1pc/[ru]_{\:\:\:(\pi_{i}^{{01}},\, (\eta_{ij}^{01})^{-1})}&
    \ar@{=>}"2";"1"^{\alpha_{i}}
    }\]

As we saw in \eqref{permHD}, the factors $H$ and $\De$ in \eqref{bitorsai} 
can be permuted, and in the present context this means the the first of the
lower arrows in this display may be permuted past the other two.
It follows that  the 2-forms $\alpha_{i}$ 
 satisfy the two equations
\bee
\label{ai1}
 \nu_{i} \, \de^{1}_{m_{i}}(\pi_{i}) = i_{\alpha_{i}}\, \nu'_{i}
\end{equation}
and
\bee
\label{ai2}
\lambda_{ij}(\alpha_{j}) = 
\de_{ij}^{-1}\,\nu_{i}[(\eta_{j}^{01})^{-1}\,
\pi_{i}^{01}(
(\gamma_{ij}^{01})^{-1}\,m_{i}^{01}((\eta_{ij}^{12})^{-1})\, 
({}^{m_{i}^{01}\,}\!\pi_{i}^{12})(\ga_{ij}^{01})
]\, \alpha_{i}\,\,\nu'_{i}(\eta_{ij}^{02} \, \de'_{ij})
\end{equation}
Equation \eqref{ai1}
 expresses the fact that
\eqref{bitorsai} is a morphism of bitorsors. Since its factors 
commute, it may be written additively as

\begin{equation}
    \label{cob2-i}\nu'_{i} - \nu_{i} + i_{\alpha_{i}} =
    \de^{1}_{m_{i}}(\pi_{i})\:.\end{equation}
It follows from the combinatorial definition \eqref{def:d1mualt} of
    $ \delta^{1}_{m_{i}}$
     that 
  $ \delta^{1}_{m_{i}}(\pi_{i})$ may be described  in classical term, 
 when  $G_i$ is the pullback to  $U_i$ of an $S$-group  and  $m_i$ is
  therefore  represented 
 by a $G_i$-valued 1-form (which we still denote by $m_i$), as
     \begin{alignat}{2}
         \label{def:d1mupi-i}
         \delta^{1}_{m_{i}}(\pi_{i}) &=\delta^{1}(\pi_{i}) +
         [m_{i},\, \pi_{i}]\\
    &=  d \pi_{i} + [\pi_{i}]^{(2)} + [m_{i},\,\pi_{i}]\:.
    \end{alignat}
The corresponding  classical expression for   \eqref{cob2-i} is
         \begin{equation}
       \label{def:5-i}
       \nu'_{i} = \nu_{i} + d\pi_{i} + [\pi_{i}]^{(2)} + [m_{i},\, 
       \pi_{i}] -
       i_{\alpha_{i}}\:.\end{equation}
       Equation  \eqref{ai2}
expresses the commutativity of \eqref{com:hijk}.
The automorphism 
$\nu_{i}$ and $\nu'_{i}$ act trivially here. Further combinatorial 
simplifications imply that equation \eqref{ai2} may be successively
additively rewritten  as
\[ \delta_{ij} - \delta'_{ij} + \lambda_{ij}(\alpha_{j}) - 
\alpha_{i} =
\de^{1}_{\pi_{i}+ m_{i}}(-\eta_{ij}) + [\gamma_{ij},\, 
\eta_{ij}]_{\pi_{i}+ m_{i}} - [\gamma_{i},\, \pi_{i}]_{m_{i}}\]
and as
\[
    \delta_{ij} - \delta'_{ij} + \lambda_{ij}(\alpha_{j}) - 
\alpha_{i} = \de^{1}_{ m_{i}}(-\eta_{ij}) - [\pi_{i},\, 
\eta_{ij}]_{m_{i}}
+ [\gamma_{ij},\,  \eta_{ij}]_{m_{i}} - [\gamma_{i},\, 
\pi_{i}]_{m_{i}}\:.
\]
Taking into account \eqref{d1rule}, we finally find
\begin{equation}
    \label{alpheqij}
     \delta_{ij} - \delta'_{ij} + \lambda_{ij}(\alpha_{j}) - 
\alpha_{i} = - \de^{1}_{ m_{i}}(\eta_{ij}) + [\eta_{ij},\, 
\eta_{ij}]_{m_{i}}
- [\pi_{i},\, 
\eta_{ij}]_{m_{i}}
+ [\gamma_{ij},\,  \eta_{ij}]_{m_{i}} - [\gamma_{i},\, 
\pi_{i}]_{m_{i}}\:.
     \end{equation}
     
     \medskip
     
The  commutativity of diagram \eqref{prism}, 
 may also be stated in  local 
terms. By making a cut in the cube \eqref{prism}
along appropriate edges, this commutativity is equivalent to the 
assertion that the  
composed 2-arrows determined by each of 
 the following two diagrams are equal:
  
 \[   \begin{array}{ccccc}
       \xymatrix{
    &&\ar@/^1pc/[drr]^{\kappa'}
     \ar@{}[d]_(.25){\,}="1"  
      \ar@{}[d]_(.75){\,}="2"  
    &&\\
      \ar@/^1pc/[urr]^{\epsilon'_{02}}\ar[rr]^{\epsilon'_{12}}
       \ar@/_1pc/[dr]_{\epsilon_{12}}&
        \ar@{}[d]_(.25){\,}="3" 
         \ar@{}[d]_(.75){\,}="4" 
       & \ar@/_1pc/[dr]^{\epsilon_{01}}
    \ar[rr]^{\epsilon'_{01}}&
     \ar@{}[d]_(.25){\,}="5" 
     \ar@{}[d]_(.75){\,}="6" 
    &\\&
     \ar@/_1pc/[ur]^{h_{12}} 
     \ar@/_1pc/[dr]_{\epsilon_{01}}&
     &\ar@/_1pc/[ur]_{h_{01}}&\\
     &&\ar@/_1pc/[ur]_{m_{01}(h_{12})}&&
     \ar@{=>}"1";"2"^{K'}
      \ar@{=>}"3";"4"_{x_{12}}
       \ar@{=>}"5";"6"^{x_{01}}
     }&&  &&     
  \xymatrix{
&\ar@{}[d]^(.2){\,}="1"
\ar@{}[d]^(.8){\,}="2"
& \ar[rr]^{\kappa'}&
\ar@{}[dd]^(.3){\,}="3"
\ar@{}[dd]^(.6){\,}="4"
&\\
\ar@/^1pc/[rru]^{\epsilon'_{02}}
\ar[r]_{\epsilon_{02}}
\ar@/_1pc/[rd]_{\epsilon_{12}}
&
\ar@{}[d]^(.2){\,}="5"
\ar@{}[d]^(.8){\,}="6"
\ar@/^1pc/[rd]_{\kappa} \ar@/_1pc/[ru]_{h_{02}}&&&\\
&\ar[r]_{\epsilon_{01}}&
\ar[r]_{m_{01}(h_{12})}&
\ar@/_1pc/[ruu]_{h_{01}}&
 \ar@{=>}"1";"2"^{x_{02}}
 \ar@{=>}"3";"4"^{a}
  \ar@{=>}"5";"6"_{K}
}
\end{array}
\]
These diagrams are displayed as
\[   \begin{array}{ccccc}
       \xymatrix{
    &&\ar@/^1pc/[drr]
     \ar@{}[d]_(.25){\,}="1"  
      \ar@{}[d]_(.75){\,}="2"  
    &&\\
      \ar@/^1pc/[urr]\ar[rr]^{\,}
       \ar@/_1pc/[dr]_{\,}&
        \ar@{}[d]_(.25){\,}="3" 
         \ar@{}[d]_(.75){\,}="4" 
       & \ar@/_1pc/[dr]^{\,}
    \ar[rr]^{{m'_{i}}^{01}}&
     \ar@{}[d]_(.25){\,}="5" 
     \ar@{}[d]_(.75){\,}="6" 
    &\\&
     \ar@/_1pc/[ur]^{\,} 
     \ar@/_1pc/[dr]_{\,}&
     &\ar@/_1pc/[ur]_{\,}&\\
     &&\ar@/_1pc/[ur]_{\,}&&
     \ar@{=>}"1";"2"^{B'_{i}}
      \ar@{=>}"3";"4"_{E_{i}^{12}}
       \ar@{=>}"5";"6"^{E_{i}^{01}}
     }&&  &&     
  \xymatrix{
&\ar@{}[d]^(.2){\,}="1"
\ar@{}[d]^(.8){\,}="2"
& \ar[rrr]^{(\nu'_{i},-)}&
\ar@{}[dd]^(.3){\,}="3"
\ar@{}[dd]^(.6){\,}="4"
&&\\
\ar@/^1pc/[rru]^{\,}
\ar[r]_{\,}
\ar@/_1pc/[rd]_{\,}
&
\ar@{}[d]^(.2){\,}="5"
\ar@{}[d]^(.8){\,}="6"
\ar@/^1pc/[rd]_{\,} \ar@/_1pc/[ru]_{\,}&&&\\
&\ar[r]_{\,}&
\ar[rr]_{({}^{m_i^{01}}\pi_{i}^{12},-)}&&
\ar@/_1pc/[ruu]_{(\pi_{i}^{01},-)}&
 \ar@{=>}"1";"2"^{E_{i}^{02}}
 \ar@{=>}"3";"4"^{\alpha_{i}}
  \ar@{=>}"5";"6"_{B_{i}}
}
\end{array}
\]
so that they are equal if and only if 
\[ E_{i}^{01}\,\, {m'_{i}}^{01}(E_{i}^{12})\,\, B'_{i} = 
\pi_{i}^{01}\,\,\, {}^{m_{i}^{01}}\!\pi_{i}^{12}(B_{i}) \,\,\, \alpha_{i}\,\, 
\, \nu'_{i}(E_{i}^{02})\,.\]
Both the automorphisms $\nu'_{i}$ and
 $\pi_i^{01}\,\,\,{}^{m_{i}^{01}\,}\!\pi_{i}^{12}$
act trivially, so that  after commuting various terms we obtain the
equation
\[ \alpha_{i} = B'_{i} - B_{i} + \delta^1_{m'_{i}}(E_{i})\]
which characterizes $\alpha_{i}$.
With the same assumption on $G_i$  as in \eqref{def:d1mupi-i} we find  
     that in classical terms
      \[  \de_{m'_{i}}^{1}(E_{i}) = dE_{i} + [E_{i}]^{(2)} + [m'_{i},\, 
      E_{i}]\:.\]  
      Substituting for $m'_{i}$ in this equation its value
      determined  by \eqref{cobe1},
      we finally obtain the following relation between the 2-forms $B_{i}$
     and $B'_{i}$:
     \begin{equation}
         \label{cob3-i}
  B'_{i} = B_{i} + \alpha_{i} - dE_{i} + [E_{i}]^{(2)} - [m_{i},\, E_{i}]
  - [\pi_{i},\, E_{i}]\:.
  \end{equation}
\bigskip

  Since the 3-curvatures  $\Omega$ and $\Omega'$ are determined by
the corresponding curving data
pairs $(\epsilon,\, K)$ and $(\epsilon', \, K')$, the associated local
3-forms $\omega_{i}$ and $\omega'_{i}$  are  determined by
the induced  forms $(m_{i},\,\gamma_{ij},\,B_{i})$ and 
$(m'_{i},\,\gamma'_{ij},\,B'_{i})$,  so that
the coboundary
relation between $\omega_{i}$ and $\omega'_{i}$ 
may  be similarly deduced from those
between $(m_{i},\,\gamma_{ij},\,B_{i})$ and 
$(m'_{i},\,\gamma'_{ij},\,B'_{i})$. In fact the terms $\gamma_{ij}$
and $\,\gamma'_{ij}$ do not affect the computation, and one finds:
\begin{proposition}
     \label{prop:cobom}
The coboundary relation between $\omega_{i}$ and $\omega'_{i}$ is given by
\begin{equation}
     \label{cob4-i}
     \om'_{i} = \om_{i} + \de^{2}_{m_{i}}(\alpha_{i}) - [\nu'_{i},\, 
     E_{i}] + [\pi_{i},\, B_{i}] +
   [\pi_{i},\, \alpha_{i}]\:.\end{equation}
   \end{proposition}
  {\bf Proof:} By
   \cite{cdf}  theorem 3.7,
  the  3-form $\de^{2}_{m_{i}}(\alpha_{i})$
  defined as in \cite{cdf} (3.3.1) is classically expressed  as
   \[ \de^{2}_{m_{i}}(\alpha_{i}) = d\alpha_{i} + [m_{i},\, \alpha_{i}]\:.\]
  The 2-form $\nu'_{i}$  is described in terms of $\nu_{i}$ by \eqref{cob2-i}.
  The
  relation to be verified therefore becomes the following one,
  whose right-hand side only involves
  the original quadruple of local forms  $(m_{i}, B_{i}, \nu_{i}, \om_{i})$
  associated to $(\epsilon,\, K)$
and the
coboundary triple $(E_{i}, \,\pi_{i},\, \alpha_{i}) $ defined by the
  morphism triple $(x,\,h,\, a)$ :
  \begin{multline}
       \label{def:6-i}
      \om'_{i} =\\ \om_{i} + d\alpha_{i} + [m_{i},\alpha_{i}] -[\nu_{i},E_{i}]
      - [d\pi_{i} ,E_{i}]
  -[[\pi_{i}]^{(2)},E_{i}]
  -[[m_{i},\pi_{i}],\, E_{i}] + [\alpha_{i},\, E_{i}] + [\pi_{i},B_{i}] + 
  [\pi_{i},\alpha_{i}]\\
  \end{multline}
  A direct computation starting from the equation
  \[\omega'_{i} = dB'_{i} + [m'_{i},\, B'_{i}]\]
  and using equations  \eqref{cobe1},\ \eqref{cob3-i} and
  the graded Jacobi identity \cite{cdf} (2.5.11) yields this result.

  \begin{flushright}
      $\Box$
      \end{flushright}
      
  %    \pagebreak 

 The previous discussion may be summarized by the following table:

{\Small
    \begin{table}[ht]
      \begin{center}
         \renewcommand{\arraystretch}{1.3}
         \begin{tabular}{|l|l|l|c|c|c|c|} \hline
             & source & target & \multicolumn{3}{c |}{
             transformation data} & transformation equations\\
             \hline
             geom. data &$(\epsilon',\, K')$ &$ (\epsilon,\, K)$ &
             $x$&$h$&$a$ & \\ \hline
             diff. forms & $(m'_{i},\, \gamma'_{ij}, \, B'_{i})$ & 
             $(m_{i},\,\gamma_{ij},\, B_{i})$ &
             $E_{i}$&$(\pi_{i},\, \eta_{ij})$ &$\alpha_{i}$&
\eqref{cobe1},\,\eqref{cobe2},\,\eqref{cob3-i}
\\\hline assoc. forms & $(\nu'_{i},\,\delta'_{ij},\, \omega'_{i})$ &
              $(\nu_{i},\,\delta_{ij},\, \omega_{i})$ & $E_{i}$&
              $(\pi_{i},\, \eta_{ij})$ &$\alpha_{i}$&
              \eqref{def:5-i},\,\eqref{alpheqij},\,\eqref{def:6-i} 
              \\ \hline
              degree &&&1 & 1 & 2&
              \\\hline
\end{tabular}

\caption{}\label{Ta:5-i}
\end{center}
\end{table}
}

\bigskip

\subsection{}  
\setcounter{equation}{0}%
\label{disc:cob-cob}

We now express in similar terms the transformation between a pair of 
triples $(x,\,h,\, a)$ and $(x',\, h',\, a')$ determined by a
 2-arrow $r$ \eqref{def:r}
  defined above $\dea$. This 2-arrow may be represented
  by a family of bitorsor isomorphisms
  \[ H'_{i} \stackrel{r_{i}}{\la} H_{i}\,,\]
  compatible with the morphisms \eqref{defhij1} for $H_i$ and
$H'_i$.
It is  described in terms of our chosen trivializations $h'_{i}$ and $h_{i}$
of 
  source and target by the family of 1-forms
  $\rho_{i} \in \mathrm{Lie}(G_{i},\,
  \Om^{1}_{U_{i}/S})$ 
  defined by
  \[ r_{i}(h'_{i}) = \rho_{i}\,h_{i}\]
 We may display  diagram \eqref{def:r}  as 
   \[
   \xymatrix{
  (\pra \lambda_{ij}, \pra g_{ijk}) \ar@/^1pc/[rr]^{(\pi_{i},\, 
  \eta_{ij}^{-1})}_{}="1"
\ar@/_1pc/[rr]_{(\pi'_{i},\, 
  (\eta'_{ij})^{-1})}_{}="2" &&   (\pra \lambda_{ij}, \pra g_{ijk}) 
\ar@{=>}"1";"2"^{r}}
\]
  so that the forms $\rho_{i}$  satisfy (in additive notation) the equations 
  { \renewcommand{\arraystretch}{1.3}
      \begin{align}
  \label{equ:irho-i}
   \pi'_{i} &=\pi_{i} + i_{\rho_{i}} \\
   \lambda_{ij}(\rho_{j}) - \rho_{i} &= \eta_{ij} - \eta'_{ij}
   \label{rhoij}
   \end{align}
   }
    The  compatibility between the composite 2-arrow
       \eqref{cond1:r} and $x'$ may be stated as
       \begin{equation}
           \label{def:rho-i}
           E'_{i} = E_{i} + \rho_{i}\
           \end{equation}
          and the compatibility of the 
           pullbacks of $r$ to $\deb$  {\it via} the three projections    with
           the arrows $a$ and $a'$ yields the equation
       \begin{equation}
           \label{eqrho1-i}
           \alpha'_{i}= \alpha_{i} + \de^{1}_{m_{i}} \rho_{i} + 
           [\pi_{i}, \,\rho_{i}]
           \end{equation}
          with $ \de^{1}_{m_{i}} \rho_i$ defined as in \eqref{def:d1mupi-i},
          in other words the equation
           \begin{equation}
           \label{eqrho2-i}
       \alpha'_{i} =  \alpha_{i} + d\rho_{i} + [\rho_{i}]^{(2)}+ 
       [m_{i},\, \rho_{i}] + [\pi_{i}, \, \rho_{i}]\:.
       \end{equation}
       These rules for the equivalence between transformation triples
       may be summarized by
  the following table:

  {\small  
   \begin{table}[ht]
      \begin{center}
         \renewcommand{\arraystretch}{1.3}
         \begin{tabular}{|l|l|l|c|c|}\hline
             & source & target & equivalence datum & equivalence equations\\
             \hline
             geometric data & $(x,\, h,\, a)$ & $(x',\, h',\, a')$&$r$
             & \\ \hline
             differential forms & $( E_{i},\,\pi_{i},\eta_{ij},\, 
             \alpha_{i})$ &
             $( E'_{i},\,\pi'_{i},\eta'_{ij}\, 
             \alpha'_{i})$ & $\rho_{i}$ &
\eqref{equ:irho-i},\,\eqref{rhoij},\,\eqref{def:rho-i},\,\eqref{eqrho2-i}
            \\ \hline
            degree &&& 1&
             \\ \hline
         \end{tabular}

\caption{}\label{Ta:6}
\end{center}
\end{table}
}

%% file: dgg-ch7.tex
\section{Some special cases}
\label{special}
     \subsection{}
\setcounter{equation}{0}%
 \label{special2}
 We now suppose that $\pc$ has a global section, which 
 identifies it with the
  trivial $G$-gerbe $\mathbb{T}_{G}$ of $G$-torsors on $X$, for some 
$S$-group $G$.  In that case,
  the associated cocycles $(\lambda_{ij},\,
  g_{ijk}) $  are trivial.
  Let us assume in addition
  that the arrow $\epsilon$ \eqref{gercon}  defining
  a connection on $\mathbb{T}_{G}$ is pointed, in
  other
  words that the image by  $\epsilon$ of the trivial torsor  $\prb
  G$-torsor $T_{\prb G}$ on $\dea$ is a trivial $\pra G$-torsor.
   In that situation,
  we can    choose an
      arrow
    \begin{equation}
      \label{gl:gam}
      \epsilon(T_{\prb G}) \stackrel{\ga}{\la} T_{\pra G}
      \end{equation}
   which is   globally defined on $X$. The induced connection $m$
   \eqref{gcon1} is then also  globally defined on the group 
   $G$ and  the  forms  $\gamma_{ij}$
   \eqref{gaijdef1}
   which   embody the gluing data for local versions  $\ga_{i}$ of
   $\gamma$, may  be taken to be trivial. Let us further require that
the bitorsor which describes   the fake curvature morphism $\kappa$ has
 a global section.
 We may then choose  a global arrow
  \begin{equation}
      \label{gl:del}
       \kappa (T_{\pra G}) \stackrel{\delta}{\la} T_{\pra G}
       \end{equation}
and this       determines by conjugation  an $\mathrm{Aut}(G)$-valued 2-form
       $\nu$ on $X$.
The curving datum $K$
 is  then
  described by  a global $G$-valued 2-form $B$. The
  induced 3-curvature $\om$ is now a global $G$-valued 3-form.
Equations \eqref{cocep1}, \eqref{cocep2}, \eqref{cockap1},
\eqref{cockap2}, \eqref{compfifj}, \eqref{comoioj}   are now  vacuous.
The differential forms occuring in this case can be displayed by the
following simplified  version of table \ref{Ta:1}:

  \bigskip

\begin{table}[ht]
      \begin{center}
         \renewcommand{\arraystretch}{1.3}
         \begin{tabular}{|l|c|c|c|c|} \hline
  & functions& 1-forms  &2-forms &
  3-forms\\\hline
   $G$-valued&  &  &
   $ B$ & $\om$\\
$\mathrm{Aut}(G)-\mathrm{valued}$&&
$m$ & $\nu$&\\\hline
\end{tabular}

\caption{}\label{Ta:3}
\end{center}
\end{table}

\bigskip

\noindent  The form  $\nu$ and  $\om$  are now defined in terms of $(m,\, B)$
by the simplified versions
   \begin{align}
     \om &= dB + [m,\, B]  \label{omidef3}\\
    dm + [m]^{(2)} &= i_{B} +\nu   \label{ifi3}
   \end{align}
   of \eqref{omidef1} and  \eqref{ifi},
   and the cocycle conditions \eqref{ificonj} and \eqref{relnufi}
   respectively reduce to the equations
     \begin{align}
     i_{\om} &= - \,(d\nu + [m,\, \nu])
     \label{ificonj3}\\
     d\om + [m,\, \om] &= [\nu,\, B]\:. \label{relnufi3}
   \end{align}

    \bigskip

  \noindent  Table \ref{Ta:2} therefore simplifies, in the present context,
  to

\pagebreak

   \begin{table}[ht]
     \begin{center}
     \renewcommand{\arraystretch}{1.3}
     \begin{tabular}{|l|c|c|}
  \hline
  global & name & differential form
  \\\hline
   $\mathbb{T}_{G}$&trivial  gerbe  &
  \\\hline
   $\epsilon$ & connection & $m$
  \\ \hline
$K$&curving data &
$B$ \\
\hline
$\kappa$ &fake curvature & $\nu$\\
\hline
$\Omega $& 3-curvature & $\om$ \\
\hline
\end{tabular}

\caption{}\label{Ta:4}
\end{center}
\end{table}

\noindent These equations for a connective structure, curving data and
associated 3-curvature on a trivial $G$-gerbe $\pc$ should be
compared with the description of the connection and curvature on a
trivial $G$-torsor reviewed in remark \ref{abcase} {\it i}. If we further
assume  that $m$ is an integrable connection on $G$, and 2 is
invertible on $S$, then
$[\nu,\,B] = - [B, \, B] = 0$ and the equation \eqref{relnufi} reduces to the
higher Bianchi identity of \cite{cdf} proposition 3.11. In particular,
when   $\epsilon$ is the canonical connection on the trivial
$G$-gerbe $\mathbb{T}_{G}$ on $X$ (so that $m = 0$), the first
two previous equations reduce to
  \[     \renewcommand{\arraystretch}{1.3}
   \left\{ \begin{array}{lcl}
     \om &= &dB \\
    \nu &= & - i_{B}\:. \end{array}\right.
   \]
     The first of these is the analogue for $G$-gerbes of the
   relation between a connection 1-form for
   the trivial $G$-bundle $T_{G}$ on $X$  and the associated curvature
   2-form  \eqref{glmc}.

  \bigskip

\subsection{}
\setcounter{equation}{0}%
We now examine the   coboundary relations in the case of a trivial 
gerbe.
Let us  consider a  morphism $(x,h,a)$ \eqref{def:x}-\eqref{def:a}
between two such curving data pairs $(\epsilon, \,
K)$ and $(\epsilon',\, K')$, on the trivial gerbe $\pc:= \mathbb{T}_{G}$,
for which the connections and the fake curvatures are pointed,
 as in \ref{special2}.
We now suppose that the torsor $H$, defined as in \eqref{def:H}, has
a global section. We may therefore choose a global path
\[h(T_{\pra G}) \la   T_{\pra G}\]
and this determines by conjugation a global
$\mathrm{Lie} (\mathrm{Aut} \ G)$-valued 1-form $\pi$ on $X$.
 The 2-arrow
  $x$ \eqref{def:xa}
is  determined by a relative $\mathrm{Lie}(G)$-valued
1-form $E$ on
$X$. The equations \eqref{cobe1} reduce to the single equation
\begin{equation}
    \label{cob1}
    m' + i_{E} = m + \pi
     \end{equation}
     between $\mathrm{Lie}(\mathrm{Aut}(G))$-valued
1-forms on $X$.
     Similarly, since the 2-arrow $a$ \eqref{def:a} is defined above
     $\deb$, it now corresponds  to a globally defined
     $G$-valued relative 2-form $\alpha$. The equations 
     \eqref{alpheqij} vanishes and \eqref{cob2-i} reduces to 
     equation 
     \begin{equation}
         \label{cob2}
         \nu' + i_{\alpha} = \nu + \delta^{1}_{m}(\pi)\:.
         \end{equation}
or equivalently to
         \begin{equation}
       \label{def:5}
       \nu' = \nu + d\pi + [\pi]^{(2)} + [m,\, \pi] -
       i_{\alpha}\:.\end{equation}
       Finally, the equations   \eqref{cob3-i} reduce to 
       the single equation 
    \[ B' = B + \alpha  - \delta^1_{m'}(E)\,\]
   in other words to 
     \begin{equation}
         \label{cob3}
  B' = B + \alpha - dE + [E]^{(2)} - [m,\, E] - [\pi,\, E]\:.
  \end{equation}
  Together with \eqref{cob1}, it expresses the pair $(m',\, B')$
  associated to the curving data pair $(\epsilon',\,K')$ in terms of the
  corresponding pair $(m,\, B)$ associated to $(\epsilon,\,K)$ and of 
  the transformation triple $( E,\,\pi,\, \alpha)$ associated to
  $(x,\, h,\, a)$.

  \bigskip

The family of equations \eqref{cob4-i} or \eqref{def:6-i}
which describe the effect of a transformation triple on the 
3-curvature reduce to 
\begin{equation}
     \label{cob4}
     \om' = \om + \de^{2}_{m}(\alpha) - [\nu',\, E] + [\pi,\, B] +
   [\pi,\, \alpha]\:,\end{equation}
   in other words to 
 \bee   
   \label{def:6}
      \om' = \om + d\alpha + [m,\alpha] -[\nu,E] - [d\pi ,E]
  -[[\pi]^{(2)},E]
  -[[m,\pi],\, E] + [\alpha ,\, E] + [\pi,B] + [\pi,
  \alpha]\:.
\end{equation}

\medskip

       The previous discussion may be summarized by the following table:

    \begin{table}[ht]
      \begin{center}
         \renewcommand{\arraystretch}{1.3}
         \begin{tabular}{|l|l|l|c|c|c|c|} \hline
             & source & target & \multicolumn{3}{c |}{
             transformation data} & transformation equations\\
             \hline
             geometric data &$(\epsilon,\, K)$ &$ (\epsilon',\, K')$ &
             $x$&$h$&$a$ & \\ \hline
             differential forms & $(m, \, B)$ & $(m',\, B')$ &
             $E$&$\pi$ &$\alpha$& \eqref{cob1},\ \eqref{cob3}\\ \hline
             associated forms & $(\nu,\, \omega)$ &
              $(\nu',\, \omega')$ & $E$&$\pi$ &$\alpha$&
              \eqref{def:5}, \ \eqref{def:6}  \\ \hline
              type &&&$G$&$\mathrm{Aut}(G)$ & $G$ & \\ \hline
              degree &&&1 & 1 & 2&
              \\\hline
\end{tabular}

\caption{}\label{Ta:5}
\end{center}
\end{table}

The equations relating the differential
forms associated to a pair of  transformation triples are
hardly changed in this situation, the family of forms $\rho_{i}$ 
being simply replaced here by a single global 1-form $\rho$.
  therefore be  represented by a global relative
  $G$-valued 1-form $\rho$.
 The equations \eqref{equ:irho-i} reduce to the single equation
  \begin{equation}
  \label{equ:irho}
   \pi' =\pi + i_{\rho}\:.
   \end{equation}
     The  compatibility between the composite 2-arrow
       \eqref{cond1:r} and $x'$ is now given by the single equation
       \begin{equation}
           \label{def:rho}
           E' = E + \rho\:.
           \end{equation}
           which replaces the family of equations \eqref{def:rho-i}
           Finally, the compatibility of the three
           pullbacks of $r$ to $\deb$ with
           the arrows $a$ and $a'$ yields the  equation
       \begin{equation}
           \label{eqrho1}
           \alpha'= \alpha + \de^{1}_{m} \rho + [\pi, \,\rho]
           \end{equation}
          or equivalently
           \begin{equation}
           \label{eqrho2}
       \alpha' =  \alpha + d\rho + [\rho]^{(2)}+ [m,\, \rho] +
       [\pi, \, \rho]\:.
       \end{equation}
which respectively replace \eqref{eqrho1-i} and \eqref{eqrho2-i}. 
   This equivalence between transformation triples 
is summarized by
  the table:

  \begin{table}[ht]
      \begin{center}
         \renewcommand{\arraystretch}{1.3}
         \begin{tabular}{|l|l|l|c|c|}\hline
             & source & target & equivalence datum & equivalence equations\\
             \hline
             geometric data & $(x,\, h,\, a)$ & $(x',\, h',\, a')$&$r$
             & \\ \hline
             differential forms & $( E,\,\pi,\, \alpha)$ &
             $(E',\, \pi',\, \alpha')$ & $\rho$ &
            \eqref{equ:irho},\ \eqref{def:rho},\ \eqref{eqrho2} \\ \hline
type &&& $G$&
            \\ \hline
            degree &&& 1&
             \\ \hline
         \end{tabular}

\caption{}\label{Ta:6a}
\end{center}
\end{table}

\subsection{}
\setcounter{equation}{0}%
\label{special3}
The final special case which we will examine  is that
in which  $\pc$
is an abelian gerbe, in the sense of \cite{lb:2-gerbes}
definition 
2.9. The cocycle pair $(\lambda_{ij},g_{ijk})$ associated to
$\pc$ is now of the form $(1, g_{ijk})$, and so defined by a
standard  2-cocycle $g_{ijk}$  with values in the abelian group $G$.
The lien of the gerbe
$\pc$ is the lien $\mathrm{lien}(G)$ associated to the
abelian $X$-group scheme $G$.
We may extend the discussion in \cite{lb:2-gerbes}
by defining for any given $X$-group homomorphism $G \la G'$ a morphism 
\[ u: \pc \la \pc'\]
between the  abelian
$G$-gerbe
and the abelian $G'$-gerbe $\pc'$ 
as a morphism of gerbes  for which the
corresponding diagrams 
  \[\xymatrix{
G \ar[r]  \ar[d]_{\zeta_{x}} & G' \ar[d]^{\zeta_{u(x)}}\\
\mathrm{Aut}(x) \ar[r]_(.45){u}&\mathrm{Aut}(u(x))
}\]
attached to objects $x \in \pc$ commute.

\bigskip

Let  $G$ be an $X$-group  with connection $m$ and let 
$\epsilon$ be  a  connection on the abelian $G$-gerbe $\pc$ compatible 
with $m$
in this sense. 
The connection $\epsilon$ is therefore described by the cocycle pair
  $(1,\, \ga_{ij}^{-1})$.  Condition \eqref{cocep1} is vacuous and
  \eqref{cocep2} simplifies to the relation
  \[ \partial\ga_{ij} =  \delta^{0}_m(g_{ijk})\]
  between $G$-valued 1-forms on $U_{ijk}$.
We now suppose in addition that $G$ is the pullback to $X$ of an
$S$-group, with its canonical connection.
 In that case  $\delta^{0}$ is the
  first differential in the $G$-valued de Rham complex of $X$, as in
  \cite{cdf} 3.1,  and $\partial$ is the \v{C}ech coboundary map.
  In particular, when $G$ is the multiplicative group
  $G_{m,S}$, this equation becomes  
  \[\ga_{jk} - \ga_{ik} + \ga_{ij} =  g_{ijk}^{-1}\,dg_{ijk}\:.\]
  Let us now make the additional assumption
 that the curvature morphism $\kappa$ is
  trivial, so that the expressions $\nu_{i}$ and $\delta_{ij}$  both
  vanish. The curving 2-arrow  $K$ again defines a family of 2-forms
  $B_{i}$. Equation
  \eqref{compfifj} now simplifies to
  \[B_{j} = (\prab \ga_{ij}^{-1})\,
(\prbc \ga_{ij}^{-1})\,
(\prac \ga_{ij})\, B_{i}\,,\]
in other words, by   $loc.\;cit.$  theorem 3.3, to
\[B_{j} - B_{i} = - \, d\gamma_{ij}
\:.\]
  The definition \eqref{omidef1} of the
local 3-curvature form $\om_{i}$  now boils down to the
simplified version
\[\om_{i} = dB_{i}\]
of \eqref{omidef1}, and \eqref{comoioj} asserts that this closed
$G$-valued 3-form is globally defined on $X$. The corresponding 
triple $(g_{ijk},
\,\ga_{ij},\, B_{i})$ is now a 2-cocycle in the standard \v{C}ech-de
Rham complex of the abelian group scheme $G$, with associated
3-curvature defined by the remaining differential in the
corresponding  non-truncated complex. In particular, when  $G:= G_{mS}$
is the multiplicative group scheme, it is simply a 2-cocycle with
values in
the truncated multiplicative de Rham complex, as in \cite{bry} 
theorem 5.3.11.

%% file: dgg-appendixa.tex
{\appendix
\section{Group-valued differential forms}
\label{appendixA}
\subsection{}
%\addtocounter{subsec}{1}%
\setcounter{equation}{0}%
%{\bf \noindent \thesection.\thesubsec} \hspace{1ex}
     In \cite{cdf}, which we will refer to in the appendices  as [cdf],
     we gave a combinatorial development of some aspects
of differential calculus for $X/S$ a relative scheme and $G$ an
$S$-group. In order to interpret our combinatorial
differential forms in classical terms,  we assumed that $G$ was  pushout
reversing ([cdf] 2.1).  The  sheaf
$\Psi^{n}_{X/S}(G)$  of relative $G$-valued differential $n$-forms on
$X$ is the sheaf of $S$-morphisms $\alpha: \De^n_{X/S} \la G$ which
vanish on the degenerate subsimplex $s\De^n_{X/S}$ of $\De^n_{X/S}$.
When $X/S$ is smooth, it  can be identified
with the sheaf  $\mathrm{Lie}(G) \ot_{\mathcal{O}_{X}} \Omega^{n}_{X/S}$. We introduced
combinatorial differentials
\begin{equation}
     \delta^{i}:\Psi^{i}_{X/S}(G) \la \Psi^{i+1}_{X/S}(G)
     \end{equation}
which we also interpreted in classical terms.  Here we extend this
theory to the case where $G$ is an $X$-group endowed with a
connection $\mu$. A prime example of such an $X$-group
 is given by $G:= \pa$ where $P$ is
a torsor under an $S$-group endowed with a connection $\epsilon$ and
$\mu = \epsilon^{\mathrm   ad}$.

\bigskip

If $F$ is a
pointed sheaf on $X$, an $F$-valued combinatorial $n$-form is by
definition an $X$-morphism
\[\Delta^{n}_{X/S} \stackrel{\omega}{\la} F\]
which maps $s\Delta^{n}_{X/S}$ to the point. We assume that $F$ is
pushout reversing so that, by the discussion following
definition 2.4 of [cdf], we may identify $\Psi^{n}_{X/S}(F)$ with
$\mathrm{Lie}(F,\, \Omega^{n}_{X/S})$.
           Observe first of all that the discussion
       in  [cdf] lemma 2.7 extends to the present context,
         since its main ingredient is the
        interpretation of the maps from $\De^{n}_{X/S}$,   to itself which
        permute a pair of adjacent infinitesimally close points.
        This situation is once more governed by the
        $\oxc$-module map which sends an $n$-form $\om \in \Om^{n}_{X}$
        to $-\, \om$. Formula (1.8.2) of
        [cdf]  is now replaced by
        the assertion that for any $G$-valued $m$-form $\om$ on $X$,
        and any permutation $\sigma \in S_{n+1}$,
        the formula
         \begin{equation}
           \label{signsig}
           \mu_{\,0\, \sigma(0)}(\sigma \,\om) = \om^{\,\mathrm{sign}(\sigma)}
       \end{equation}
       is satisfied.
       In particular, the multiplicative version of formula (1.8.2) of
       [cdf] remains valid
       for any permutation $\sigma$
       of $(x_{0}, \ldots, x_{m}) \in \De^{n}_{X/S}$ which leaves fixed the
       first variable $x_{0}$, and   for a transposition
       $\sigma :=(0i)$ involving the first variable one has
       \begin{equation}
           \label{transp}
           \om(x_{0}, \ldots, x_{m})\,
       \, \,\mu(x_{0}, x_{i})\,(\om (x_{i}, x_{1}, \ldots
       x_{i-1}, x_{0}, x_{i+1}, \ldots x_{m}))= 1\:.\end{equation}

\bigskip

    We now define the de Rham differential maps
    for combinatorial forms with values in an  $X$-group  $G$
    with connection $\mu$.
  The  differential
      \[\delta_{\mu}^{0}: G \la \mathrm{Lie}(G,
\Om^{1}_{X/S}) \]
is defined by
  \bee
\label{def:d0mu}
\delta_{\mu}^{0}(g)(x,y) := g(x)^{-1} \, \mu(x,y) (g(y) )
\end{equation}
and a section $g$ of $G$ is said to be horizontal if  it lies
in the kernel of $\de^{0}_{\mu}$.
The same computation as in [cdf] lemma 3.1 shows
that $\delta^{0}_{\mu}$
is a crossed homomorphism with respect to the right adjoint
action of $G$ on $\mathrm{Lie}(G,\, \Om^{1}_{X/S})$.

\begin{remark}
     {\rm
The  twisted differential $\de^{0}_{\mu}$ can  be understood
in more classical terms
by considering, for a fixed
a global  section $g$ of $G$, the
map of sheaves on $\De^{1}_{X/S}$
\bee
\label{gmapO}
\begin{array}{ccc}
     \mathrm{Isom}(\prb G, \, \pra G) & \stackrel{\{-,\,g\}}{ \la} & \pra G\\
     \mu(x,y) & \mapsto & g(x)^{-1}\,\, \mu(x,\, y)(g(y))\:.
     \end{array}
     \end{equation}
Pulling back by the diagonal embedding
  $\Delta: X \hookrightarrow \De^{1}_{X/S}$, this
induces the map of $X$-sheaves
\bee
\label{gmap1}\begin{array}{ccc}
\mathrm{Aut}(G) &\stackrel{[-,\,g]}{\la} &  G\\
u & \mapsto &
  g^{-1}\, u(g)\:.
\end{array}\end{equation}
In
particular, $[u,\, g]= 1$ whenever $u= 1$ so that, by functoriality of the
basic deformation lemma of Lie theory, stated in
 [cdf]  as proposition 2.2, one has the  commutative diagram
of sheaves on $\dea$
\bee
\label{gdiag}
\xymatrix{
1 \ar[r] & i_{\ast}\mathrm{Lie}(\mathrm{Isom}(\prb G,\, \pra G),\,
\Om^{1}_{X/S})\ar[r]
\ar[d]_{i_{\ast}\mathrm{Lie}(\{-,\,g\})}
& \mathrm{Isom}(\prb G,\, \pra G) \ar[r]\ar[d]^{\{-,\,g\}}
& i_{\ast}\mathrm{Aut}(G) \ar[d]^{[-,\,g]}
\\
1 \ar[r]  & i_{\ast}\mathrm{Lie}(G,\, \Om^{1}) \ar[r] & \pra G \ar[r] &
i_{\ast}G \:.
}\end{equation}
The hypotheses of [cdf] are not exactly
satisfied here, since the sheaf \[F:= \mathrm{Isom}(\prb G,\, \pra
G)\] to which we now wish to apply that assertion is
not a sheaf of $H$-sets, and {\it a fortiori } not a sheaf  pointed above
$\De^{1}_{X/S}$. It is however pointed above $X$, and  this suffices
in order for
$\mathrm{Lie}(F,\Om^{1}_{X/S})$ to be defined.
  Since the
middle vertical map
$\{- ,\, g\}$  is compatible with this $X$-pointing,
the left-hand vertical map is well defined.
By exactness of the upper short exact sequence,
the term $\mathrm{Lie}(\mathrm{Isom}(\prb G, \, \pra G),\, \Om^{1}_{X/S})$
is the
space of connections on the $X$-group $G$. The clockwise
image in $\pra G$
of a
connection $\mu$ is $\de^{0}_{\mu}(g)$, so that commutativity  of the
left-hand square  interprets  this expression in classical terms as
the image of the connection $\mu$ under the map
obtained  by applying the Lie functor to
the morphism \eqref{gmapO}.}\end{remark}

\vspace{.75cm}

We now define  the $\mu$-twisted de Rham differentials
\begin{equation}
\label{def:dmu}
  \delta^{i}_{\mu}:\mathrm{Lie}(G, \, \Om^{i}_{X/S})
\la \mathrm{Lie}(G, \,\Om^{i+1}_{X/S})
\end{equation}
for  $i > 0$. We set
         \begin{align}
\label{def:d1mu}
(\delta^{1}_{\mu}\om)(x,y,z) &:= \om(x,y)\,\,[\mu(x,y)(\om(y,z))]\,\,
[\mu(x,y)\mu(y,z)(\om(z,x))]\\ \label{def:d2mu}
(\delta^{2}_{\mu}\om)(x,y,z,u) &:= \mu(x,y)(\om(y,z,u))\,
\om(x,y,u)\,\om(x,u,z)\, \om(x,z,y)\\ \label{def:d3mu}
(\delta^{3}_{\mu}\om)(x,y,z,u,v)
&:=
\mu(x,y)(\om(y,z,u,v))\, \om(x,y,u,v) \,\om(x,y,z,u)
\,\om(x,z,v,u)\,\om(x,y,v,z)
\end{align}
  and define $\delta^{i}_{\mu}$  by
similar formulas for $i > 3$.
By \eqref{transp}, the first of these equations can be written
\[
(\delta^{1}_{\mu}\om)(x,y,z) = \om(x,y)\,\,[\mu(x,y)(\om(y,z))]\,\,
[\mu(x,y)\mu(y,z) \mu(z,x) (\om(x,z)^{-1})]\:,
\]
  and also by the simpler formula
  \begin{equation}
      \label{def:d1mualt}
      (\delta^{1}_{\mu}\om)(x,y,z) = \om(x,y)\,
  \,\mu(x,y)(\om(y,z))\,\,\om(x,z)^{-1} \end{equation}
  since the version of  lemma 2.8 of [cdf] pertaining to the pairing  of
  proposition 2.14 implies that
the curvature $\kappa_{\mu}(x,y,z)$ \eqref{def:kmu}
of the connection $\mu$ acts trivially on $G$-valued 2-form
$\omega (x,\, z)^{-1}$.
Formula \eqref{signsig} ensures that the section $\delta^i_{\mu}\om$ of
$G$ is indeed a
differential form on $X$.
It follows from the definitions (\ref{def:d0mu}), (\ref{def:d1mu}) of
$\delta^{0}_{\mu}$ and
$\delta^{1}_{\mu}$ that
\begin{equation}
     \begin{array}{ccl}
     \delta^{1}_{\mu}\delta^{0}_{\mu} (g)(x,y,z) &=&
g(x)^{-1}\,\kappa_{\mu}(x,\,y,\, z)(g(x))\\
&= &[[\kappa_{\mu},\,g]](x,y,z)
\end{array}\end{equation}
with the bracket  $[[\kappa_{\mu},\,g]]$ defined  in \eqref{doubbra} below.
In particular, if the
curvature $\kappa_{\mu}$ of the connection $\mu$  on the group $G$
is trivial, {\it i.e.} if $\mu$
  is integrable, then
\[\delta^{1}_{\mu}\delta^{0}_{\mu}= 0\:.\]
When $G$ is pushout reversing the combinatorial commutation principle for
$G$-valued forms of [cdf] lemma 2.8 and  remark 2.9 ensures that, for
 $i \geq 2$, $\de^i_{\mu}$ is a group homomomorphism. 

\bigskip

The  adjoint action of a group $G$  on its Lie-valued forms, as
in [cdf]
  \S 2.9, also extends to the case of an $X$-group scheme. In
  addition to the naive adjoint action of $G$ on
  $\mathrm{Lie}(G)$-valued forms, we now have a  twisted (right) adjoint action
  of $G$ on $\mathrm{Lie}(G,\, \Om^{1}_{X/S})$
defined by
\[ \om(x,y)^{\ast_{\mu}\,g}:= g(x)^{-1} \,\om(x,y)\, \,\mu(x,\,y)(g(y))\:.\]
  This expression can also be written as
  \[\om(x,y)^{\ast_{\mu}\,g} = [g(x)^{-1} \,\om(x,y)\, g(x)]\,
  [g(x)^{-1}\mu(x,\,y)(g(y))]\]
  so that this twisted adjoint action is related to the naive adjoint
  action by the formula
  expressed additively as
  \begin{equation}
      \label{astmu}
      \om(x,y)^{\ast_{\mu}\,g} = \om(x,y)^{\,g} + \delta^{0}_{\mu}(g)\:.
      \end{equation}

  \vspace{.75cm}

   We now carry over to the context of an $X$-group $G$  with
  connection $\mu$ the discussion of  [cdf] propositions 2.10 and 2.14
pertaining to the definition of the Lie bracket pairings of $G$- and
$\mathrm{Aut}(G)$-valued forms\footnote{Although the combinatorial
definitions require no hypothesis on $\mathrm{Aut}(G)$, the notation
$\mathrm{Lie}(\mathrm{Aut}(G),\, \Omega^{n}_{X/S})$
and the properties of the pairings on forms require $\mathrm{Aut}(G)$
to be pushout reversing. By [cdf] proposition 2.2, this will hold if
$G/X$ is flat.\label{Gflat}} on $X$.
  We define a pairing
   \[
\begin{array}{ccc}
     \mathrm{Lie}(G,\, \Om^{m}_{X/S}) \times  \mathrm{Lie}(G,\, \Om^{n}_{X/S})
     & \la & \mathrm{Lie}(G,\, \Om^{m+n}_{X/S}) \\
     (f,\, g) & \mapsto & [f,g]_{\mu}
     \end{array}
     \] for $m,\,n > 0$ in terms of the
     corresponding bracket pairing of [cdf] 2.5
     by setting
     \begin{equation}
         \label{bra0mu}
         [f,\, g]_{\mu} := [f,\, \mu_{0m}(g)] \:,\end{equation}
     in other words by
     \begin{equation}
         \label{bra1mu}
     [f,g]_{\mu}(x_{0},\ldots,x_{m+n}) := [f(x_{0},\ldots,x_{m}),\,
     \mu_{0m}(g(x_{m},\ldots,x_{m+n}))]
     \end{equation}
     where the right-hand bracket is the commutator pairing
  $[a,\,     b]:= aba^{-1}b^{-1}$ in the group $\pra G$ above
  $\Delta^{m+n}_{X/S}$.
The bracket (2.8.4) of $\, loc.\ cit$ is extended similarly from $S$-
to $X$-group schemes by
     setting
     \begin{equation}
         \label{bra3mu}
     [u,\, g]_{\mu}:= [u,\, \mu_{0m}(g)]
     \end{equation}
     for   $u \in \mathrm{Lie}\:(\mathrm{Aut}(G),\,\Omega^{m}_{X/S})$
and $g \in \mathrm{Lie}(G, \,\Omega^{m}_{X/S})$ with
$m,\,n > 0$,
  so that
         \begin{equation}
             \label{bra2mu}
         [u,\, g]_{\mu}(x_{0},\ldots,x_{m+n}) =
         u(x_{0},\ldots,x_{m})(\mu_{0m}(g(x_{m},\ldots,x_{m+n})))\,\,
         \mu_{0m}(g(x_{m},\ldots,x_{m+n})^{-1})\:.
         \end{equation}
         The proofs of [cdf] \S 2 extend to the present context, and
         show that the pairings \eqref{bra0mu} and \eqref{bra3mu} are both
         bilinear, and that  the former satisfies the identity
         \begin{equation}
             \label{altcommu}
         [f,\, g]_{\mu} = (-1)^{mn+1}[g,\,f]_{\mu}
         \end{equation}
         of $\,loc.\  cit.$ lemma 2.11. The latter satisfies the corresponding
         identity
\begin{equation}
     \label{altcommu1}
     [u,\, g]_{\mu} = (-1)^{mn+1}[g,\,u]_{\mu} \end{equation}
where the pairing $ G \times \mathrm{Aut}(G) \la \mathrm{Aut}(G)$ on
which the right-hand term is based is  defined, for
any section $a \in G$ and $\phi \in \mathrm{Aut}(G)$, by
\begin{align}
  \notag   [a,\,\phi] &:= [\phi,\, a]^{-1}\\
        \label{altcommu2}       & \:\, =  a \, \phi (a^{-1})  \:.
               \end{align}
               The pairing $[g,\,u]_{\mu}$ is therefore  explicitly defined by
                \begin{equation}
                   \label{bra4mu}
                   [g,u]_{\mu}(x_{0},\ldots,x_{m+n}):= g(x_{0},\ldots,x_{m})
  \:({}^{\mu_{0m}\,}\!u(x_{m},\ldots,x_{m+n}))(g(x_{0},\ldots,x_{m})^{-1})\:.
\end{equation}

\bigskip

         The pairings \eqref{bra1mu} and \eqref{bra2mu} are actually independent
         of the choice of a
         connection $\mu$ on $G$, as we now show for the pairing \eqref{bra1mu}.
         Suppose that $\mu'$ is a second connection
         on the $X$-group scheme $G$, so that $\mu'= \eta\, \mu$, with $\eta$
         an $\mathrm{Aut}(G)$-valued 1-form on $X$. It follows that
                 \begin{multline}
           [f,\, g]_{\mu'}(x_{0},\ldots,x_{m+n})   =
             [f(x_{0},\ldots,x_{m}),\, \mu'_{0m}(g(x_{m},\ldots,x_{m+n})] \\
          = [f(x_{0},\ldots,x_{m}),\,
             (\eta_{0m}\,\mu_{0m})(g(x_{m},\ldots,x_{m+n})]\\
             = [f(x_{0},\ldots,x_{m}),\,
             (\eta_{0m}\,\mu_{0m})(g(x_{m},\ldots,x_{m+n})\,\,
            \mu_{0m}(g(x_{m},\ldots,x_{m+n})^{-1})] \\
            [f(x_{0},\ldots,x_{m}),\,\mu_{0m}
             (g(x_{m},\ldots,x_{m+n}))]
             \end{multline}
         with the last equation following from the standard commutator identity
         $[a,\,bc]= [a,b]\,
         {}^{b\,}\![a,\,c]$. The assertion is now proved, since the first
         factor in this last expression vanishes by [cdf] remark 2.9.
         We will at times simply denote the pairings \eqref{bra1mu} and
         \eqref{bra2mu} by $[f,\, g]$ and by  $[u,\, g]$ just as
         when the group $G$ is defined over $S$.

\bigskip

  The formula \eqref{bra2mu} for $n$ or $m$ equal to zero   
 respectively defines pairings
 \begin{align}
\label{lie1} \mathrm{Lie}(\mathrm{Aut}(G),\, \Om^{m}_{X/S})
  \times G& \la \mathrm{Lie}(G,\, \Om^{m}_{X/S})\\
\mathrm{Aut}(G)  \times  \mathrm{Lie}(G,\, \Om^{n}_{X/S}) &
 \la  \mathrm{Lie}(G,\, \Om^{n}_{X/S})\label{lie2}
\end{align}
which we will both denote once more   by $ [u,\,g]_{\mu}$.
Formula \eqref{altcommu1} or the explicit
formula  \eqref{bra4mu}  similarly determine for $m=0$ or $n=0$
   corresponding  pairings $ [g,\, u]_{\mu}$:
  \begin{align}
\label{lie3}
 G \times \mathrm{Lie}(\mathrm{Aut}(G) ,\, \Om^{n}_{X/S})
&  \la \mathrm{Lie}(G, \,  \Om^{n}_{X/S})\\
\mathrm{Lie}(G,\, \Om^{m}_{X/S}) \times \mathrm{Aut}(G) & \la
 \mathrm{Lie}(G,\, \Om^{m}_{X/S}) \label{lie4}
\end{align}
 These four  pairings are still linear in the positive degree variable,
  but only twisted
  linear in the degree zero one.  For example,  the pairings
 \eqref{lie1}
  and \eqref{lie3} respectively satisfy the equations 
  \begin{alignat*}{2}
      [u,\, gg'] \:&= &[u,g]\:\:{}^{g\,}\![u,\,g']\\
      [gg',\, u] \: &=&{}^{g\,}\![g',\, u]\, \:[g,\,u]
      \end{alignat*}
       analogous to those satisfied by group
      commutators.
We will at times encounter when $n= 0$ instead of \eqref{bra3mu} the
expression $[[u,\,g]]$ in $\mathrm{Lie}(G,\, \Om^{m}_{X/S})$ defined
by
\begin{equation}
     \label{doubbra}
     [[u,\, g]](x_{0},\,\ldots\, x_{m}) =g(x_{0})^{-1}u(x_{0},\,\ldots\,
     x_{m})(g(x_{0}))
     \end{equation}
     The twisted linearity of the bracket pairing, and the relation
     \eqref{altcommu} then imply that the bracket $ [[u,\, g]]$ may
     also be expressed in various other ways  as
     \[ [[u,\, g]] = [g^{-1},\, u] = ([g,\,u]^{g})^{-1} = [u,\, g]^{g}\:. \]
    

%% file: dgg-appendixb.tex
\section{The Lie theory of $gr$-stacks}
\label{sec:conger}
\subsection{}
\setcounter{equation}{0}%
\label{stacks}
In the following paragraph we will discuss the Lie theory of
$gr$-stacks, but will first give a short summary,
for the reader's convenience, of the relevant
portions of the theory of  stacks.

     \bigskip

     By a prestack (of groupoids) $\mathcal{C}$ on
     $X$,
     we mean a
     sheaf of groupoids, {\it i.e.} a split  fibered category on $X$,
     whose objects and arrows form  a sheaf of groupoids on $X$.
     For any open set $U \hookrightarrow X$,
     the set objects in the fiber category  $\cc_{U}$ are sections above
     $U$ of a
     given
      sheaf
     $F_{0}$ on $X$. The arrows in $\cc_{U}$ are sections on $U$ of
     a sheaf
$F_{1}$, endowed with source and target morphisms
$d_{1},\,d_{0}:F_{1}\la F_{0}$ and an identity map $s_{0}: F_{0}\la
F_{1}$.
Once the inverse law, and the composition laws for arrows in $\cc$  have
been specified, the truncated simplicial sheaf
$\xymatrix{F_{1} \ar@<.8ex>[r]
\ar@<-.8ex>[r] & F_{0}\ar[l]}$
extends uniquely to a simplicial sheaf $F_{\ast}$ on  $X$ satisfying the
additional conditions of \cite{illusie} VI propositions 2.2.3 and
2.6.1, which ensure that $F_{\ast}$ is the nerve of a groupoid.
By construction, the arrows in $\mathcal{C}$ glue. So do the objects
of $\cc$ in the restrictive sense  provided by the sheaf property for
$F_{0}$.

\bigskip

When
$\cc$ is endowed with a group-like monoidal structure $\cc \times \cc
\stackrel{\otimes}{\la} \cc$, the
situation is even simpler, since we may then replace the simplicial
sheaf $F_{\ast}$ by one in which the associativity and the unit constraint are
strict, as is the inverse law, so that the tensor law yields 
 a group structure on $F_{0}$. In that case, the simplicial
sheaf $F_{\ast}$ has a very simple description. It is
  determined by the
sheaf of groups of objects $G_{0}:= F_{0}$,
the  sheaf of groups  of arrows $G_{1}:= \mathrm{ker}(d_{1})$
in $\cc$
sourced at the
identity element $I$ of $G_{0}$, and a target
homomorphism
$\delta:G_{1} \la G_{0}$ which is simply the restriction to $G_{1}$ of
the target map
$d_{0}$. A final element of structure is the  left action
\bee
\label{cract}
\begin{array}{ccc}
     G_{0} \times G_{1} & \la & G_{1}\\
     (g,\, f) & \mapsto & {}^{g\,}\!f
     \end{array}
     \end{equation}
   of $G_{0}$ on $G_{1}$,  which sends the arrow
   $f: I \la \delta \,f$  in $\cc$ to
   the composite arrow $gfg^{-1}$:
   \[ I \la g \,I\, g^{-1} \stackrel{1\,f\,1}{\la}
   \:g(\delta \,f) g^{-1}\:.\]
   Once the compatibility of this action with the morphism $\delta$ has
   been fully specified, the  complex of sheaves of non-abelian
   groups
   \[\delta:G_{1} \la G_{0}\]
   possesses the structure of a crossed module (\cite{lb:2-gerbes} 1.2)
   in the category of
   sheaves on $X$ and this structure  completely determines
   the  monoidal prestack
   $\cc$.
   In particular, arrows in $\cc$ correspond to  pairs
   $(f,g) \in G_{1} \times G_{0}$,
   with $g$ the source
   of the arrow, and $(\delta f)\,g $ the target. The tensor law for
   arrows is  defined by the semi-direct
   product  structure $G_{1} \rtimes G_{0}$ on the set $G_{1} \times G_{0}$
   defined by the given action of $G_{0}$ on $G_{1}$.
   A prime example of such a construction
   is the prestack
   \bee
   \label{giag}
   G\: \stackrel{i}{\la} \:\mathrm{Aut}(G)\end{equation}
   associated as in \eqref{inconj} to any  sheaf of $X$-groups
   $G$, together with  the obvious  left action of $\mathrm{Aut}(G)$ on $G$.

   \bigskip

  Let $\tilde{\cc}$ be the stack associated to a given
  prestack  in groupoids $\cc$
  on  $X$. This stack $\tilde{\cc}$ is obtained from $\cc$ by a
  sheafification process, which forces the effectivity of the
  descent data for
  objects in $\tilde{\cc}$. An object $\underline{x}$ in the fiber category
  $\tilde{\cc}_{U}$ is  defined by a family $(x_{i})$ of objects
  in $\cc_{V_{i}}$, for some open cover $\mathcal{V} = (V_{i})_{i \in
  I}$ of $U$, together with arrows
  \[\phi_{ij}:x_{j \,\mid \,V_{ij} }\la x_{i\, \mid \,V_{ij} }\]
  in $\cc_{V_{ij}}$ satisfying the classical descent condition
  \[\phi_{ij}\, \phi_{jk} = \phi_{ik}\]
  above the open set $V_{ijk}$. An arrow $\underline{f}:\underline{x}
  \la \underline{y}$ in $\cc_{U}$ is defined, on a common refinement of
  the covers $\mathcal{V}$ and $\mathcal{V}'$ on  which $\underline{x} $
  and $\underline{y} $ were expressed, by a family of arrows $f_{i}:
  x_{i} \la y_{i}$ whose compatibility with the descent data morphisms
  $\phi_{ij}^{x}$ and $\phi_{ij}^{y}$ for $\underline{x}$ and $\underline{y}$
  is expressed by the commutativity of the diagrams
   \[\xymatrix{
    x_{j} \ar[r]^{\phi^{x}_{ij}} \ar[d]_{f_{j}} & x_{i}
    \ar[d]^{f_{i}}\\
    y_{j}\ar[r]_{\phi^{y}_{ij}} & y_{i}\:.
   }\]
Finally, a monoidal structure $\otimes: \cc \times \cc \la \cc$ on the
prestack $\cc$ induces a monoidal structure on $\tilde{\cc}$, with
tensor law defined by
\[\underline{x} \otimes \underline{y}: = (x_{i} \otimes y_{i})_{i \in
I}\:.\]
The stack on $X$ associated to the crossed module \eqref{giag} is the
stack of $G$-bitorsors on $X$.

   \bigskip

\subsection{}
         % \addtocounter{subsec}{1}%
\setcounter{equation}{0}%
%{\bf \noindent \thesection.\thesubsec} \hspace{1ex}
  The rudiments of a Lie theory for  $gr$-stacks  will now be set forth.
For any open set  $U$  in $X$, we will denote by $U[\epsilon]$ the
$U$-scheme of dual numbers on $U$, with its  canonical
     closed immersion $i:U \hookrightarrow U[\epsilon]$.

\begin{definition}
   Let  $(\gc,\, \ot,\, I)$ be a $gr$-stack on the big \'etale site of $X$.
   The Lie stack
   $\mathrm{Lie}\,
     \gc$ associated to  $\gc$ is the stack on $X$
      whose fiber above an open set $U$
     of $ X$ is the category $(\mathrm{Lie}\,
     \gc)_{U}$ of pairs $(Z, \eta)$, with $Z$
     an object in the fiber  category $\gc_{U[\epsilon]}$
     and
     \[\eta: I \la i^{\ast}Z \]  an arrow in $\gc_{U}$.
     An arrow $ (Z,\, \eta)  \la (Z',\, \eta') $ in
     $(\mathrm{Lie}\, \gc)_{U}$ is an arrow $\phi: Z \la Z'$ in
     $\gc_{U[\epsilon]}$ for which  $ (i^{\ast}\phi )\,\, \eta = \eta'$.
     \end{definition}

The gluing properties in $\gc$ ensure that the fibered category
  $\mathrm{Lie}\,  \gc$  is indeed a stack. Any morphism of $gr$-stacks
  $F: \gc \la \gc'$  on $X$ induces  a morphism of stacks
  \[\mathrm{Lie}(F): \mathrm{Lie}\, (\gc) \la \mathrm{Lie}\, (\gc')\]
  on $X$.
  For any positive integer $n$, the corresponding theory
  for the square-zero immersion
  \bee
  \label{sq0}
  t: s\De^{n}_{X/S} \hookrightarrow
  \De^{n}_{X/S}\end{equation}
  of [cdf] 1.12 defines the  stack of relative  $\gc$-valued
  (one could also say ``Lie
  $\gc$-valued'') $n$-forms on $X$.

  \begin{definition}
      \label{def:b2}
      Let $\gc$ be a $gr$-stack on $X$.
      The stack $\mathrm{Lie}(\gc,\, \Om^{n}_{X/S})$ of
      relative  $\gc$-valued $n$-forms on $X$ is the stack on $X$ whose
      fiber  $\mathrm{Lie}(\gc,\, \Om^{n}_{X/S})_{U}$
      is the category of
      pairs $(Z,\, \eta)$, with $Z$ an object in $\gc_{\De^{n}_{U/S}}$
      and $\eta: I \la t^\ast Z$ an arrow in the category
      $\gc_{s\De^{n}_{U/S}}$ . An arrow $\phi:
     (Z,\, \eta)  \la (Z',\, \eta')$
     in $\mathrm{Lie}(\gc,\, \Om^{n}_{X/S})_{U}$ is an arrow
     $\phi: Z \la Z'$
     in
     $\gc_{\De^{n}_{U/S}}$
      for which    $(t^{\ast}\phi )\,\, \eta = \eta'$.
      An object ({\it resp.} an arrow) in the stack
     $\mathrm{Lie}(\gc,\, \Om^{n}_{X/S})$ will be called a relative
   {\rm  Ob} $\mathrm{Lie}\: \gc$-valued ({\it resp.} an {\rm Ar}
  $\mathrm{Lie}\: \gc$-valued)
     $n$-form on $X$.
      \end{definition}

     \bigskip

   Consider  the $gr$-stack $\tilde{\gc}$  associated to the prestack
   $\gc$ described   by a
   crossed module \bee
   \label{crossF}
   \delta: F_{1} \la F_{0}\end{equation}
   on $X$. Assume that the sheaves $F_{i}$ are universally pushout
   reversing in the sense of [cdf] definition 2.1.
   Applying the Lie functor
   to the map $\delta$,
   we obtain a crossed module
   \bee
   \label{lief}
  \mathrm{Lie}(\delta):
  \mathrm{Lie}(F_{1}) \la \mathrm{Lie}(F_{0})\end{equation}
   in the category of $\oxc$-Lie algebras. The corresponding prestack in
   Lie algebras
   will be denoted $\mathrm{Lie}(\gc)$. Crossed modules
   in the category of Lie algebras have already been
   considered by C. Kassel and  J.-L. Loday (appendix to
   \cite{kl}). We will
   neglect their  Lie structure here. The action of $\mathrm{Lie}(F_{0})$
on $ \mathrm{Lie}(F_{1})$ induced by the action of $F_0$ on $F_1$
  is trivial,so that 
    the complex of
   sheaves of abelian groups \eqref{lief} defines
a group-like symmetric monoidal prestack.  By functoriality of the
   Lie construction, there
   exists a canonical
   morphism of symmetric monoidal prestacks
   \bee
   \label{defstackmap}
   \phi_{\gc}:  \mathrm{Lie}\,(\gc) \la
   \mathrm{Lie}\,(\tilde{\gc})\:,\end{equation}
and  the universal property of the  associated stack implies that
   $\phi_{\gc}$ factors through a  monoidal functor
     \[\Phi_{\gc}: (\mathrm{Lie}\,\gc)^{\sim} \la \mathrm{Lie} \,
  \tilde{\gc}\:.\]
   The same construction, applied to the square-zero
immersion $t$ \eqref{sq0} constructs a monoidal functor
  \[\Phi(\gc, \, \Om^{n}_{X/S}):  (\mathrm{Lie}(\gc,\,
       \Om^{n}_{X/S}))^{\sim} \la \mathrm{Lie} (\tilde{\gc}, \,
       \Omega^{n}_{X/S})\:.\]
       The target of this functor will be called the  stack of
       $\mathrm{Lie}\, (\tilde{\gc})$-valued (or $\tilde{\gc}$-valued)
       relative $n$-forms on $X/S$.
       \begin{proposition}
        \label{psigc}
    Let $\gc$ be the prestack on  $X$ associated to a crossed module
$F_{\ast}$  \eqref{crossF}, and suppose that the components sheaves $F_{i}$
of $F_{\ast}$ universally reverse
   pushouts in the sense of [cdf] definition 2.1.   The functor
    $\Phi_{\gc}$ is an equivalence of monoidal
       stacks, and so are, for all $n > 0$, the
       functors
       $\Phi(\gc, \, \Om^{n}_{X/S})$.
       \end{proposition}
       {\bf Proof:} It follows from the definition of the term
prestack and the construction of its associated stack \cite{lmb} lemme 3.2
    that
       the functor $\Phi_{\gc}$ is fully faithful. All
       that must be proved in order to verify the assertion  for the
       functor $\Phi_{\gc}$ is the essential surjectivity. An
       object of $(\mathrm{Lie} \,\tilde{\gc})_{U}$ is determined
       by a family of
       objects $g_{j}$ in $\tilde{\gc}_{V_{j}[\epsilon]}$, for some
       open cover $\mathcal{V}:= (V_{j})_{j\in J}$ of $U$, together
       with gluing
       data $\phi_{ij}$, and  arrows $\phi_{j}:1 \la  i^{\ast}g_{j}$
       compatible with the induced gluing data
       $i^{\ast}\phi_{ij}$.  We may refine the open cover
       $\mathcal{V}$, and assume  that $g_{j}$ and $\phi_{j}$
       are defined by sections of $F_{0}$ and $F_{1}$ above
       $V_{j}[\epsilon]$.
       Since the structural map $p$ of the $U$-scheme $U[\epsilon]$ is
       a retraction of the immersion $i$, the objects
       \[g'_{j}:= g_{j}\, (i\,p)^{\ast}g_{j}^{-1}\]
       of $(F_{0}{)}_{V_{j}[\epsilon]}$
       are  section  of  $\mathrm{Lie}\, F_{0}$ above $V_{j}$.
       Together with the arrows $\phi_{ij}': g'_{j}\la g'_{i}$ induced by the
       gluing data $\phi_{ij}$, they determine an object
        $\underline{g}'$ in the fiber category of
       $(\mathrm{Lie}\,\gc)^{\sim}$ above $U$, and the sections
       $p^{\ast}\phi_{j}$ then define an arrow
       $p^{\ast}\underline{\phi}: \underline{g}' \la \underline{g}$ in
       $\mathrm{Lie}\,(\tilde{\gc})$.

       \bigskip

       The same reasoning applies to  the square-zero
       immersion $t: U \hookrightarrow
       \De^{1}_{U/S}$ with retraction $p_{0}$, and  proves that
       $\Phi(\gc, \, \Om^{1}_{X/S})$
      is an equivalence. It cannot be applied directly to forms
       of higher degree, since for $n>1$ there is no retraction of
       $ \De^{n}_{U/S}$ onto  $s \De^{n}_{U/S}$.
       Once more, only the essential surjectivity of
       $\Phi(\gc, \, \Om^{n}_{X/S})$ need be verified, and a
       $\tilde{\gc}$-valued
       $n$-form on $U$ can represented, for each member $V$ of an open
       cover $\mathcal{V}$ of $U$,  by a family of $n+1$-tuples
       $(g, \, \phi_{0},\ldots, \phi_{{n-1}})$,
      for sections  $g :\Delta^{n}_{V/S} \la F_{0}$ and
       $\phi_{i}: \De^{n-1}_{V/S} \la F_{1}$ of $F_{0}$  and $F_{1}$
       satisfying the relations

         \begin{alignat}{2}
       g(x_{0}, \ldots, x_{i},x_{i}, \ldots x_{n-1}) &= \delta
       \phi_{i}(x_{0}, \ldots x_{n-1} ) && \qquad \qquad \qquad \qquad
       \quad (A_{i})
       \label{ai}\end{alignat}
      for all $ i \leq n-1$. Such conditions  are  not  in themselves
      sufficient in order to ensure that the $(n+1)$-tuple
      $(g, \, \phi_{0},\ldots, \phi_{{n-1}})$ determines an object   of
      $\mathrm{Lie}(\gc)_{V}$.
     The degenerate subscheme
      $s\Delta^{n}_{V/S}$ of $\Delta^{n}_{V/S}$ may be viewed as
      the colimit
      of the diagram
    \[
      \xymatrix{
      \bigvee_{i,j}\Delta^{n-2}_{V/S} \ar@<.5ex>[r]  \ar@<-.5ex>[r]
      &
      \bigvee_{i} \Delta^{n-1}_{V/S} \ar@{.>}[r] &\Delta^{n}_{V/S}
      }\]
      embodying the relations $s_{i}\,s_{j} = s_{j+1}\, s_{i}$ for
      $i\leq j$ between the degeneracy maps in the simplicial $S$-scheme
      $\Delta^{\ast}_{V/S}$. This follows directly
      from the corresponding assertion
      for  the points of $s\Delta^{n}_{V/S}$ with values in an
      arbitrary $S$-scheme $T$. The morphisms $\phi_{i}$
    therefore  assemble to a  section above
      $s \Delta^{n-1}_{V/S}$  of the pushout reversing sheaf
      $F_{1}$, if and only if
      they satisfy the additional relations
      \begin{alignat}{2}
       \phi_{i}(x_0, \ldots, x_{j}, x_{j}, \ldots x_{n-2})
       &=\phi_{j+1}(x_{0}, \ldots, x_{i}, x_{i}, \ldots x_{n-2})&&
       \qquad \qquad  (B_{i,j})\label{bij}
       \end{alignat}
       for $i \leq j \leq n-2$.
       Finally, we are given  sections
       $\gamma^{VW}: \Delta^{n}_{V\cap W} \la F_{1}$ for any pair of
      sets $V, W$ in $\mathcal{V}$, which compare the
      restrictions of $g^{V}$ and $g^{W}$ to $V \cap W$ {\it via}
      the equation
      \[\delta ( \gamma^{VW}) \, g^{W}  = g^{V}\,.\]
      These must satisfy the usual cocycle conditions
      $\gamma^{TV}\gamma^{VW} = \gamma^{TW}$ on triple
      intersections $T \cap V \cap W$, and  also  the compatibilities
      \[ \delta
      \gamma^{TV}(x_{0}, \ldots, x_{i},x_{i}, \ldots x_{n-1})\,\,
      \phi_{i}^{V}(x_{0}, \ldots , x_{n-1})   =
      \phi_{i}^{T}(x_{0}, \ldots , x_{n-1})\:.\]
     We have now made explicit the full data which determines an object
       in the fiber category   \linebreak    $\mathrm{Lie} (\tilde{\gc}, \,
       \Omega^{n}_{X/S})_{U}$.
       \begin{lemma}
           \label{lemdeg}
           Let $F_{1}\la F_{0}$ be a crossed module in the category of sheaves
           on $U$. Consider for any open set $V$ of $U$  an  $n+1$-tuple
           $(g, \, \phi_{0},\ldots, \phi_{n-1})$,
           satisfying the conditions
           (\ref{ai}, \ref{bij}).
There exists a pair
           $(g'(x_{0},\ldots,x_{n}),\, \chi(x_{0}, \ldots,x_{n-1})$
          with $g':\De^{n}_{V/S} \la
           F_{0}$ and $\chi:\De^{n-1}_{V/S} \la F_{1}$, such that
           \[\delta (\chi)  \, g'  = g\]
           and
           \[ g'(x_{0}, \ldots, x_{i}, x_{i}, \ldots x_{n-1}) = 1\]
           for all $i$.
           \end{lemma}
{\bf Proof of lemma \ref{lemdeg}:}  We will
define, for each fixed integer $k \geq 0$,
sections $g^{k}: \De^{n}_{V/S} \la F_{0}$ and $\phi_{i}^{k}: \De^{n-1}_{V/S}
\la F_{1}$  which satisfy the conditions
       \begin{alignat}{2}
       g^{k}(x_{0}, \ldots, x_{i},x_{i}, \ldots x_{n-1}) &= \delta
       \phi^{k}_{i}(x_{0}, \ldots x_{n-1} ) && \qquad \qquad (A^{k}_{i})
       \label{aik}\\
       \intertext{for all $i \leq n-1$,}
       \phi_{i}^{k}(x_0, \ldots, x_{j}, x_{j}, \ldots x_{n-2})
       &=\phi^{k}_{j+1}(x_{0}, \ldots, x_{i}, x_{i}, \ldots x_{n-2})&&
       \qquad \qquad  (B^{k}_{i,j})\label{bijk}\\
      \intertext{for all $i \leq j \leq n-2$, as well as the additional
      conditions}
      \phi^{k}_{i}(x_{0},\ldots,x_{n-1}) &=1  && \qquad \qquad (C^{k}_{i})
       \label{cik}
       \end{alignat}
       for all $i<k$. It follows that
       \[ g^{k}(x_{0}, \ldots x_{i}, x_{i}, \ldots x_{n-1}) = 1\]
       for all $i <k$.  We set $g^{0} := g$ and $
       \phi^{0}_{i} :=  \phi_{i}$ for all $i$, and will define inductively
       an
       $n+1$-tuple of sections
       $(g^{k+1}, \phi_{0}^{k+1},\ldots, \phi^{k+1}_{n-1})$.
       We set
        \bee
       \label{def:gk}
       g^{k+1}(x_{0}, \ldots, x_{n}) :=   \delta \phi^{k}_{k}(x_{0}, \ldots,
       \widehat{x_{k+1}}, \ldots , x_{n})^{-1} \, g^{k}(x_{0}, \ldots,
      x_{n})\,,\end{equation}
       so that
       \bee
       g^{k+1}(x_{0}, \ldots, x_{n}) =
       g^{k}(x_{0}, \ldots, x_{k}, x_{k}, x_{k+2}, \ldots
       x_{n})^{-1} \,\, g^{k}(x_{0}, \ldots,   x_{n})\:.\end{equation}
       The sections $\phi^{k+1}_{i}$ are  defined for $0 \leq i \leq n-1$ by
     \bee
\label{def:phiijk1}
   \phi^{k+1}_{i}(x_{0}, \ldots, x_{n-1}):=
   \begin{cases}
     1 &  i < k+1\\\phi^{k}_{k}(x_{0},
     \ldots, x_{n-1})^{-1}
      \,\phi^{k}_{k+1}(x_{0}, \ldots, x_{n-1})  &
     i=k+1\\  \phi^{k}_{k}(x_{0},\ldots,
     \widehat{x_{k+1}},\ldots,x_{i},x_{i},\ldots,x_{n-1})^{-1}
     \,\phi^{k}_{i}(x_{0}, \ldots, x_{n-1}) & i > k+1\:.
     \end{cases}
     \end{equation}
  By condition $B^{k}_{k,i-1}$,  one also has
  \bee
  \label{def2phik1}
   \phi^{k+1}_{i}(x_{0}, \ldots, x_{n-1}) =
    \phi^{k}_{i}(x_{0}, \ldots, x_{n-1})\, \phi^{k}_{i}(x_{0}, \ldots,
     x_{k}, x_{k}, x_{k+2}, \ldots, x_{n-1})^{-1}
  \end{equation}
  for $i > k+1$. It is straightforward to verify conditions $A^{k+1}_{i},
  B^{k+1}_{i,j}$ and $C^{k+1}_{i}$,
  though the variety of definitions of $\phi^{k+1}_{i}$ makes the
  proof of the  conditions $B^{k+1}$ somewhat cumbersome to carry out,
  and we will not spell
  it  out here. The proof of the lemma is then complete, since it
  suffices to set $g':= g^{n}$ and
   \[\chi(x_{0}, \ldots, x_{n}):=  \prod_{k=n-1}^{0} \phi^{k}_{k}(x_{0},
  \ldots, \widehat{x_{k+1}}, \ldots x_{n})\:.\]
  \begin{flushright}
      $\Box$
      \end{flushright}

  \bigskip

   Proposition \ref{psigc} now follows, since lemma \ref{lemdeg}
  allows us to replace the local objects $g^{V}(x_{0}, \ldots x_{n})$
  by genuine local sections $(g')^{V}$ of $\mathrm{Lie}\,(F_{0},
  \Om^{n}_{V/S})$ , accompanied by arrows $\chi^{V}: (g')^{V}\la g^{V}$.
  By functoriality, this construction glues correctly on the
  intersections $V \cap W$, and defines the sought-after object
  $\underline{g}'$
  in $\mathrm{Lie}\,(\gc,\,
  \Om^{n}_{U/S})^\sim$, together with the corresponding arrow
  $\underline{\chi}:
  \underline{g}' \la \underline{g}$.
       \begin{flushright}
           $\Box$
       \end{flushright}

\begin{remark}
     {\rm
{\it i})
      Group-like symmetric monoidal stacks
    are often
   referred to as  Picard stacks (\cite{sga7} 1.4).
  By transport of structure by $\Phi_{\gc}$ and $ \Phi (\gc,\, \Om^n_{X/S})$,
   the monoidal stacks $\mathrm{Lie}(\tilde{\gc})$ and
$\mathrm{Lie}\,(\tilde{\gc},\, \Om^n_{X/S})$ are  Picard
stacks on $X$ (and even strict Picard stacks in the sense of $
loc.\;cit.$).

\bigskip

\qquad \qquad \quad {\it ii}) Since the diagram (1.12.1) of [cdf]
is cocartesian, and since its lower map $s$ is
  the closed immersion  of $X$ into the generalized dual number
$X$-scheme determined  by the $\oxc$-module $\Om^n_{X/S}$, it follows
from the corresponding
assertions for the prestack $\gc$ that  an object in
$\mathrm{Lie}\,(\tilde{\gc},\, \Om^n_{X/S})$ corresponds to pairs
$(Z, \,\eta)$ with $Z$ an object in $\tilde{\gc}_{X[\Om^n_{X/S}]}$
and $\eta: I \la s^{\ast}Z)$ an arrow in $\tilde{\gc}_{X}$.
As in [cdf] (2.3.4), an $n$-vector
  field $D \in \Gamma(X,\, \wedge^n T_{X/S})$  contracts the object
  $Z$,
  to a global object of the stack $\mathrm{Lie}\,\tilde{\gc}$.
  This is consistent with the  terminology of `` stack of
  $\mathrm{Lie} \, \tilde{\gc}$-valued $n$-forms'' which we have adopted for
$\mathrm{Lie}\,(\tilde{\gc},\, \Om^n_{X/S})$.

  \bigskip

  \qquad \qquad \quad {\it iii})  We have ignored here the Lie algebra
  structure
  on the components of  $\mathrm{Lie}\, (F_{\ast})$. The Lie brackets
  on $F_{0}$ and $F_{1}$
  determine  bracket pairing on $\mathrm{Lie}\,(\tilde{\gc})$,
  and  pairings of Picard  stacks
  \[ [\:,\:]:\: \mathrm{Lie}\,(\tilde{\gc},\, \Om^m_{X/S}) \times
  \mathrm{Lie}\,(\tilde{\gc},\, \Om^n_{X/S})  \la
  \mathrm{Lie}\,(\tilde{\gc},\, \Om^{m+n}_{X/S}) \]
  for  $m,n
  >0$ which are bilinear (in the categorical sense).
  The Jacobi identity for objects should only be valid in this
  context up to an arrow  which will itself satisfy a higher coherence
  condition, patterned on the corresponding term in the definition of
a homotopy Lie algebras \cite{lada}.
  }
\end{remark}
   }